\documentclass[11pt,leqno]{amsart}
\usepackage{amsthm,amsfonts,amssymb,amsmath,oldgerm,color,bm,multicol}
\usepackage{epsfig}
\usepackage{tikz}
\usetikzlibrary{decorations.markings}

\usepackage{pgfplots}

\usepackage{hyperref}

\usepackage{mathtools}
\mathtoolsset{showonlyrefs}
\numberwithin{equation}{section}

\renewcommand\d{\partial}
\renewcommand\a{\alpha}
\renewcommand\b{\beta}
\renewcommand\o{\omega}
\newcommand\s{\sigma}
\renewcommand\t{\tau}
\newcommand\R{\mathbb R}\newcommand\N{\mathbb N}\newcommand\Z{\mathbb Z}
\newcommand\C{\mathbb C}

\newcommand\dsp{\displaystyle}

\def\g{\gamma}
\def\G{\Gamma}
\def\t{\tau}

\def\th{\theta}
\def\diag{\mbox{diag}\,}

\def\l{\lambda}

\def\e{\varepsilon}

\newcommand\br{\begin{rem}}
\newcommand\er{\end{rem}}
\newcommand\bp{\begin{pmatrix}}
\newcommand\ep{\end{pmatrix}}
\newcommand\be{\begin{equation}}
\newcommand\ee{\end{equation}}
\newcommand\ba{\begin{equation}\begin{aligned}}
\newcommand\ea{\end{aligned}\end{equation}}

\def\fl1{{\mathcal F} L^1}

\pgfplotsset{compat=1.13}

\newcommand{\Id}{{\rm Id }}

\newcommand{\sgn}{\text{\rm sgn}}

\newtheorem{defi}{Definition}[section]
\newtheorem{theo}{Theorem}
\newtheorem{prop}[defi]{Proposition}
\newtheorem{lem}[defi]{Lemma}
\newtheorem{cor}[defi]{Corollary}
\newtheorem{rem}[defi]{Remark}

\newtheorem{assump}[defi]{Assumption}

\newtheorem{nota}[defi]{Note}

\def\op{{\rm op} }

\numberwithin{equation}{section}
\begin{document}

\renewcommand{\refname}{References}

\title[Space-time resonances and Raman instabilities in Euler-Maxwell]{Space-time resonances and high-frequency Raman instabilities in the two-fluid Euler-Maxwell system}
\author{Eric Dumas}
\address{Institut Joseph Fourier UMR CNRS 5582, Universit\'e Grenoble-Alpes}
\email{eric.dumas@univ-grenoble-alpes.fr}
\author{Yong Lu}
\address{Department of Mathematics, Nanjing University, Nanjing, China}
\email{luyong@nju.edu.cn}

\author{Benjamin Texier}
\address{Institut Camille Jordan UMR CNRS 5208, Universit\'e Claude Bernard - Lyon 1} \email{texier@math.univ-lyon1.fr} 

\date{\today}

\begin{abstract} We 
show that space-time resonances induce high-frequency Raman instabilities
in the two-fluid Euler-Maxwell system describing laser-plasma interactions.
 A consequence is that the Zakharov WKB approximation to Euler-Maxwell is unstable for non-zero group velocities. A key step in the proof is the reformulation of the set of resonant frequencies as the locus of weak hyperbolicity for linearized equations around the WKB solution. We analyze those linearized equations with the symbolic flow method. Due to large transverse variations in the WKB profile, the equation satisfied by the symbolic flow around resonant frequencies is a linear {\it partial} differential equation. At space-time resonances corresponding to Raman instabilities, we observe a fast growth of the symbolic flow, which translates into an instability result for the original system. The strongest instability is caused by backscattered Raman waves.
\end{abstract}

\maketitle

\setcounter{tocdepth}{1}

\tableofcontents

\section{Introduction}

Resonances may destabilize WKB approximate solutions, no matter how accurate the approximation is. This was proved and made precise in the monograph \cite{em4}, for semi-linear hyperbolic systems and large-amplitude solutions.

We extend here the results of \cite{em4} in two distinct directions, as we consider {\it quasi-linear} systems, and {\it fast transverse} variations of the WKB profile. Quasi-linear terms model fluid-dynamical convective phenomena; Gaussian laser beams used in current experiments present fast transverse variations.

For quasi-linear systems of partial differential equations,
we describe how space-time resonances may destabilize highly accurate WKB approximations issued from highly-oscillating initial data with fast transverse variations.

The instability is due to incompatible nonlinearities -- that is, {\it non-transparent} nonlinearities in the sense of Joly, M\'etivier and Rauch \cite{JMR-TMB}, or nonlinearities which do not satisfy {\it null forms} in the sense of Klainerman \cite{Kla}.

Here as in the work of Germain, Masmoudi and Shatah \cite{GMS,GM}, time resonances are defined as the set $\{ \Phi = 0\}$ in frequency space, where $\Phi$ is some relevant characteristic phase; space-time resonances are time resonances which in addition belong to $\{ \nabla_\perp \Phi = 0 \},$ where $\nabla_\perp$ is the gradient in the frequency directions corresponding to transverse spatial directions.

Our results apply to (indeed, are stated for) the two-fluid Euler-Maxwell system describing laser-plasma interactions. We show that the Zakharov WKB approximation to Euler-Maxwell is unstable in the case of a non-zero group velocity. The leading instability is caused by the growth of a backscattered Raman component of the electrical field. The growth takes place, in the frequency domain, at a space-time resonance between electromagnetic waves and electronic plasma waves. Thus {\it our analysis rigorously describes the Raman instability (the growth of backscattered Raman waves) in the context of the Euler-Maxwell equations.}

The instability phenomenon in Theorem \ref{th:other} is important in applications, and known as {\it backscattered Raman instability.} According to the physics literature, high-frequency instabilities play an important role in the current failure of large-scale inertial confinement fusion experiments to deliver significant amounts of energy. This is described for instance in \cite{BF}, and in the 2019 {\it White paper on opportunities in plasma physics}, which we quote here: ``Laser-plasma instabilities inhibit the deposition of energy ... [New broad-bandwith lasers could potentially] suppress high-frequency instabilities like [...] stimulated Raman scattering'' (\cite{WP}, pages 2-3).

Forward Raman instabilities are also observed in our analysis. 

A consequence of our results is that the model system famously postulated by Zakharov \cite{Zak+,Zak} does not always accurately describe the nonlinear interactions between the envelope of the electric field of a laser and the mean mode of fluctuation of density in a plasma.

\subsection{The Euler-Maxwell equations} \label{sec:EM}

The non-dimensional form of the equations introduced in \cite{em1}, based on orders of magnitude of physical parameters observed in experiments targeting inertial confinement fusion (see \cite{em1} and \cite{Chen,DB}), is
$$\mbox{(EM}) \left\{ \begin{aligned}  \d_t B & + \nabla \times E  =  0,  \nonumber \\
       \d_t E & -  \nabla \times  B   =  j, \nonumber \\
   \d_t v_e & + \theta_e (\sqrt \e v_e \cdot \nabla) v_e   = - \theta_e  \nabla n_e -  \frac{1}{\e} L_e,  \nonumber \\
  \d_t n_e & + \theta_e \nabla \cdot v_e   +  \theta_e ( \sqrt \e v_e \cdot \nabla) n_e   =  0 , \nonumber \\
\d_t v_i & +  \e ( v_i \cdot \nabla) v_i   =  - \alpha_{ie}^2 \sqrt\e \nabla n_i +  \frac{1}{\theta_e \sqrt \e} L_i,\nonumber \\
  \d_t n_i & +   \sqrt \e \nabla \cdot  v_i   +  \e (v_i \cdot \nabla) n_i   =  0, \nonumber
\end{aligned} \right. $$
where $(B,E) \in \R^6$ is the electromagnetic field. The field $v_e \in \R^3$ is the velocity of the electrons, and $v_i \in \R^3$ is the velocity of the ions. For notational simplicity, the ions are assumed to have charge $+e$ (that is, $+1$ in non-dimensional variables). The fluctuations of density are $n_e \in \R$ and $n_i \in \R.$ The current density $j$ is defined by
$$ j = \frac{1}{\e} e^{\sqrt \e n_e} v_e - \frac{1}{\theta_e \sqrt \e} e^{\sqrt \e n_i} v_i.$$ The Lorentz forces are
$$ L_e = E + \theta_e \sqrt \e v_e \times B, \qquad L_i = E +  \e v_i \times B.$$
The parameters $\e > 0,$ $\theta_e > 0$ and $\a_{ie} > 0,$ are defined as
 $$ \e = \sqrt\frac{m_e}{4 \pi e^2 n_0 t_0^2}, \quad \theta_e = \sqrt\frac{\g_e T_e}{m_e c^2}, \quad \a_{ie} = \frac{T_i}{T_e},$$
 where $n_0$ is the density at equilibrium, $m_e$ the mass of the electrons, and $t_0$ the duration of the incident light pulse. The high-frequency regime corresponds to $\e \to 0.$ The parameter $\g_e$ is the specific heat ratio of the electrons, $T_e$ their temperature, and $c$ the speed of light in vacuum. The temperature of the ions is $T_i.$ We have $\theta_e \ll 1$ and $\a_{ie} \ll 1,$ and both are much smaller than $\e$ in applications. We will consider $\theta_e$ and $\a_{ie}$ to be fixed in the limit $\e \to 0.$

 The divergence equations are
\be \label{div:eq}  \nabla \cdot B  = 0, \qquad \nabla \cdot E = \frac{1}{\e \theta_e} (n_i - n_e).\ee
  We will make sure that the divergence equations are satisfied by our data. The solutions will then satisfy \eqref{div:eq} for all times.

 The form of (EM) as given above derives from two assumptions: the first is a {\it cold ions} hypothesis, in the sense that we posited the ratio $\theta_i$ of the ionic acoustic velocity $\sqrt{T_e/m_i}$ to the speed of light $c$ to be $O(\sqrt \e);$ the second is a choice of amplitude $O(\sqrt \e),$ which can be seen in particular in the convective terms (but in the present form (EM) of the equations, we consider solutions of size $O(1)$ with respect to $\e$).

 For more details on these parameters, and the passage from the physical equations to their non-dimensional forms, we refer to Section 2 in \cite{em1}.

 For fixed $\e > 0,$ local-in-time existence of solutions issued from data in Sobolev spaces $H^s(\R^3),$ with $s > 1 + 3/2,$ derives from the quasi-linear, symmetric hyperbolic nature of (EM). The maximal existence time a priori depends on $\e.$

 For data with amplitudes $O(1),$ the $\e^{-1/2}$ prefactor in front of the nonlinearity suggests an existence time $O(\e^{1/2}).$ We prove here that some small initial perturbations generate solutions which are defined and amplified in time $O(\e^{1/2} |\ln \e|).$ 
 
 In a different scaling, and for non-oscillating data (corresponding to $k = 0;$ by constrast we need $|k| > \sqrt 3$ here: see Theorem \ref{th:other} below), the result of \cite{em3} shows existence and stability of solutions to (EM) up to time $O(1).$
 
 The global existence results of \cite{GM,GIP} concern small initial data, while in our context the data are $O(1)$ with respect to the small parameter $\e.$  

\subsection{The Zakharov approximation to Euler-Maxwell} \label{sec:ZEM}

 In Section 6.2 of \cite{em1}, it is shown how an approximate WKB solution $u_a$ to the above (EM) system can be constructed with an arbitrarily high degree of precision.  
 
\subsubsection{The WKB approximate solution} The approximate solution $u_a = u_a(\e,t,x,y)$ given in \cite{em1} for small enough $\e > 0,$ solves (EM) with a small $O(\e^{K_a})$ error, over a fixed time interval $[0, T_a],$ with $T_a > 0$ independent of $\e.$ We denote $\e^{K_a} R_a(\e,t,x,y)$ the error: it is the extra term in the right-hand side of (EM) when we put $u = u_a.$ The exponent $K_a,$ that is the order of the approximation, can be made arbitrarily large provided the WKB datum has a large enough Sobolev regularity. The remainder $R_a$ is bounded,  pointwise in $(t,x,y)$ and in appropriately weighted Sobolev norms, uniformly in $\e.$ Details are given in Section \ref{sec:wkb}. 

%

\subsubsection{The leading term in the WKB datum} The leading term in the datum for $u_a$ is a Gaussian laser pulse: 
 \be \label{data} u_a(\e,0,x,y) = \cos\left(\frac{k x}{\e}\right) a\left(x, \frac{y}{\sqrt \e}\right) + O(\sqrt \e), \qquad \e > 0, \quad x \in \R, \quad y\in \R^2,\ee
 with a wavenumber $k \in \R \setminus \{ 0 \}.$ 
 
 The initial amplitude $a$ is for instance smooth and compactly supported: $a \in C^\infty_c(\R^3_{x,y}),$ but our analysis holds in the more general setting described in Section \ref{sec:other} below. We assume that $a$ is not identically equal to $0,$ and that is satisfies the {\it polarization} condition $\Pi_1(k) a(x,y) = a(x,y),$ where $\Pi_1$ is an appropriate spectral projector (see Section \ref{sec:wkb}). This polarization condition is a necessary condition for the existence of a WKB solution issued from \eqref{data}.

The datum \eqref{data} defines the characteristic length scales. There are three scales: the fast scale $1/\e$ in the longitudinal $x$ direction, the intermediate scale $1/\sqrt\e$ in the transverse $y$ directions, and the slow $O(1)$ scale of variation of the $\e$-independent amplitude $a.$ In particular, characteristic frequencies in the $x$ direction are $O(1/\e).$ Equivalently, the characteristic scale of variation in the $x$ direction, or wavelength, is $O(\e).$ 

\subsubsection{Corrector terms} 
 The approximate solution $u_a$ is a WKB solution, in the sense that it decomposes as
  $$ u_a = \sum_{0 \leq k \leq K'_a} \e^{k/2} u_k,$$
for some $K'_a \in \N$ depending on $K_a.$ The leading term is $u_0,$ and the first corrector is $u_1.$ 
 Due to the oscillations in \eqref{data}, all the profiles $u_k$ further decompose into oscillating and non-oscillating terms:
 $$ u_k = \sum_{q \in {\mathcal H}_k} e^{i q (k x - \o t)/\e} u_{k,q}(t,x,y/\sqrt \e),$$
 where ${\mathcal H}_k$ is a finite subset of $\Z$ containing $0,$ and $\o = \sqrt{1 + k^2}$ is a characteristic time frequency associated with the initial spatial frequency $k.$ WKB computations (see Section \ref{sec:wkb} and \cite{em1}) show that
 \be \label{wkb:intro} u_0 = \Re e \, \big(  e^{i (k x - \o t)/\e} u_{0,1}(t,x,y/\sqrt \e) \big),\quad u_1 = u_{1,0}(t,x,y/\sqrt \e) + \mbox{fast oscillating terms}\ee
 and that the higher-order corrector terms $u_k,$ for $k \geq 2,$ are functions of the $u_{k'},$ with $k' < k,$ and their spatial derivatives. 
%

\subsubsection{The Zakharov system} In view of (EM), the amplitudes $u_{0,1}$ and $u_{1,0}$ decompose into electromagnetic fields, electronic and ionic velocities, and fluctuations of density. Let $\tilde E$ be the electrical component of $u_{0,1}$ and $\bar n$ be the (electronic or ionic) fluctuation of density of $u_{1,0}.$ Let also ${\bf E}(t,x,Y)$ and ${\bf n}(t,x,Y)$ be such that
 $$ \tilde E(t,x,y) = {\bf E}(t,x,y/\sqrt \e), \quad \bar n(t,x,y) = {\bf n}(t,x,y/\sqrt \e).$$
 Then, $({\bf E}, {\bf n})$ solves a Zakharov system of the form 
 $$ \mbox{(Z)} \left\{\begin{aligned} \big(i(\d_t + c(k) \d_x) - \Delta_Y\big) {\bf E} & = {\bf n} {\bf E}, \\ (\d_t^2 - \Delta_Y) {\bf n} & = \Delta_Y |{\bf E}|^2,\end{aligned}\right.$$
where $c(k)$ is a group velocity, different from 0 if $k \neq 0.$ (For legibility, we did not include physical parameters in (Z) other than the group velocity $c(k).$ The exact (Z) system derived by a WKB expansion from (EM) in \cite{em1} is precisely given in Section \ref{sec:wkb}.)

 The Schr\"odinger equation for the envelope of the electrical field in (Z) is typical of three-scale approximations (see for instance \cite{JMR1, Lt, sw, Lu}, and \cite{Du} for a survey of geometric optics). The wave operator in (Z) is an acoustic (two-scale) approximation of the electronic compressible Euler system. We refer to Section 6.2 of \cite{em1}, and Section \ref{sec:wkb} of the present article, for details on the WKB approximation of (EM) by (Z). A formal derivation of (Z) from (EM) is also found in \cite{SS}.

 There is an extensive  mathematical literature on the Zakharov system. The short-time well-posedness result in Sobolev spaces that we use here is due to Linares, Ponce and Saut \cite{LPS}. It ensures the existence of a solution to (Z) at least for some small (positive) time -- of course, independent of $\e$ --, and then the whole WKB approximate solution can be constructed up to this time, as sketched in Section \ref{sec:wkb}. Notable earlier work on the Cauchy problem for the (Z) system with zero group velocity ($c(k) = 0$) include \cite{SW} and \cite{OT}.

 Since ${\bf E}$ is a component of $u_0,$ the leading term in the WKB expansion, and ${\bf n}$ a component of $u_1,$ a corrector term in the expansion, the (Z) system is an example of a {\it ghost effect}, in which a corrector has a measurable effect on an observable, but is not itself an observable quantity \cite{Sone, TA}. The ghost effect is a consequence of a lack of transparency (in Joly, M\'etivier and Rauch's terminology \cite{JMR-TMB}) of the Euler-Maxwell system, a property that is described in detail in \cite{em1}.

\subsubsection{Background on the Zakharov approximation to Euler-Maxwell}  The stability theory of the Zakharov approximation to Euler-Maxwell was initiated in \cite{em1, em2, em3}. In \cite{em1}, transparency properties of the (EM) system were written out, and a WKB approximate solution was constructed. The article \cite{em2} gave a stability proof of the Zakharov approximation to Klein-Gordon-waves systems, which have the same characteristic variety as Euler-Maxwell (see Section \ref{sec:char} and in particular Figure \ref{fig1}). For the Zakharov approximation to Euler-Maxwell, the {\it stability} result of \cite{em3} addresses specifically the case $k = 0$ (then the group velocity in (Z) vanishes: $c(0) = 0$); by contrast the present {\it instability} Theorem \ref{th:other} holds only for $k$ away from 0. 

\subsection{Main result: instability of the Zakharov approximation to Euler-Maxwell}

The WKB approximate solution $u_a$ to {\rm (EM)}, as described above in Section {\rm \ref{sec:ZEM}}, is unstable in the following sense. Consider an initial datum of the form
\be \label{initial:datum}
u(0,x,y) = u_a(\e,0,x,y) + \e^K \phi_\e(x,y),
\ee
with
\be \label{varphi:initial:perturbation}
\phi_\e(x,y) := \Re e \, \Big(\, \big( e^{i x (\xi_0 + k)/\e} \vec e_0 + e^{i x \xi_0/\e} \vec f_0 \big) \,  \phi\left(x, \frac{y}{\sqrt \e}\right) \, \Big),
\ee
where $\phi$ is scalar, smooth and compactly supported around $(x_0, 0),$ with $x_0 \in \R$ and $\xi_0 \in \R,$ and $\vec e_0$ and $\vec f_0$ are fixed vectors. The perturbation profile $\phi$ depends on $y/\sqrt \e$ just like the WKB datum \eqref{data}. Note that with \eqref{initial:datum}-\eqref{varphi:initial:perturbation}, we have
\be \label{initial:distance} \| (u - u_a)(0, \cdot) \|_{L^\infty(\R^3)} + \| (\e \d_x)^\a (\sqrt \e \d_y)^\b (u - u_a)(0,\cdot) \|_{L^2(\R^3)} \lesssim \e^K,\ee
for any $\a, \b \in \N^3.$ 

\begin{theo} \label{th:other} Given $k > \sqrt 3,$ for $\theta_e$ small enough: given $a \in C^\infty_c(\R^3),$ with $a \not\equiv 0$ satisfying the polarization condition \eqref{pola:a}, given $K > 0$ and $K' > 0:$ for some $\phi_\e \in C^\infty_c(\R^3)$ as in \eqref{varphi:initial:perturbation}, for some $T > 0,$ for any small enough $\e > 0,$ the initial-value problem {\rm (EM)}-\eqref{initial:datum} has a unique solution $u \in \cap_{s > 0} C^0([0, T(\e)], H^{s}(\R^3))$ which satisfies
$$ \| (u - u_a)(T(\e),\cdot) \|_{L^\infty(\R^3)} \geq \e^{K'},$$
with $T(\e) \sim T \sqrt \e |\ln \e|.$ 
%
 \end{theo}

Theorem \ref{th:other} describes the growth of the distance $\|u - u_a\|$ between exact and approximate solutions: initially it is $O(\e^K)$ in pointwise ($L^\infty$) and weighted Sobolev \eqref{weighted:norm} norms, as seen in \eqref{initial:distance}; in short time $O(\sqrt \e |\ln \e|)$ it becomes greater than $\e^{K'},$ pointwise in space. The key is that $K$ can be chosen arbitrarily large, and $K'$ arbitrarily small. We can allow for a weaker condition bearing on the WKB amplitude $a:$ see Theorem \ref{th:main} below.

%
%

\subsection{A quick overview of the proof} \label{sec:quick} The proof of Theorem \ref{th:other} consists in changes of variables, a representation formula, and estimates in weighted Sobolev and $\fl1$ norms. We work with the perturbations equations about the WKB approximate solution. Once an initial instability direction is found (that's $\xi_0$ and $\vec e_0$ and $\vec f_0$ in \eqref{varphi:initial:perturbation}), the main task 
is to control the solution to the Cauchy problem, in order to show that the growth  
does indeed occur. To this aim, the instability direction must be chosen with the optimal amplification rate.

The tools we use to carry this out are {\it reductions to normal forms} and 
the {\it symbolic flow method}. The normal form procedure is a (pseudo-differential) 
change of unknown which eliminates fictitious coupling coefficients 
in the linearized system (i.e., coupling coefficients associated with null forms, or coupling coefficients far from resonances). We use two distinct normal forms. For the systems in normal forms, further changes of variables which exploit a property of {\it separation of resonances} reveal that the fast oscillations in the WKB solution can be factored out to first order. The symbolic flow method then provides an approximate solution 
operator leading to appropriate estimates.

Sections \ref{sec:outline}, \ref{sec:a:little:more} and \ref{sec:intro:spacetime} in this introduction contain much more on the ideas and key concepts of the proof.

\subsection{Plan} The proof (of Theorem \ref{th:main} below, which generalizes Theorem \ref{th:other}) covers Sections \ref{sec:mainproof} to \ref{sec:interaction:coefficients}. In Sections \ref{sec:mainproof} and \ref{sec:char}, 
we set up notation for the perturbation system about the WKB solution and describe the spectral elements (characteristic variety, resonances and space-time resonances) that are involved. 
Normal form reductions are performed in Sections \ref{sec:normal} and \ref{sec:normal:2}. 
Sections \ref{sec:weak:Sobolev:bound} to \ref{sec:lower} are devoted to further changes of variables and estimates, and the 
end of the proof of Theorem~\ref{th:other} is given in Section~\ref{sec:endgame}. 
Sections \ref{sec:wkb} to \ref{sec:interaction:coefficients} contain the computations 
concerning the WKB approximate solution, the abovementioned spectral elements and 
the optimal growth rate. 
Section \ref{app:symb} gathers some notation and classical results pertaining to 
pseudo-differential operators. The symbolic flow theorem (Theorem \ref{th:duh}) is given 
in Section \ref{app:duh}. Sections \ref{sec:notation:index} 
to \ref{sec:lists} provide a notation index, a terminology index, and lists of most parameters used in the text. 

This introduction concludes with three subsections in which we expound on the ideas in the proof of Theorem \ref{th:other}. Before that, we formulate a slight generalization of Theorem \ref{th:other}. 


\subsection{A more general condition bearing on the WKB amplitude} \label{sec:other}

We prove Theorem \ref{th:other} not under the assumption $a \in C^\infty_c(\R^3),$ but under the weaker assumptions that the amplitude $a:$ 
 \begin{itemize}
 \item[(a0)] is not identically equal to zero: $a \not\equiv 0,$
 \item[(a1)]  satisfies the classical polarization condition \eqref{pola:a},
 \item[(a2)] has Sobolev regularity: $a \in H^{s_a}(\R^3),$ 
 \item[(a3)] and satisfies the exponential decay bound \label{page:a3}: for some $R_a > 0,$ some $r_a > 0,$ some $\kappa_a > 0:$ for all $|\a| \leq s_a - d/2,$ for some $C_\a > 0,$ for all $|(x,y)| > R_a:$
 $$  |\d_{x,y}^\a a(x,y)| \leq C_\a \exp\big( - r_a |y|^{\kappa_a}\big).$$
 \end{itemize}
 In condition (a3), the radius $R_a$ can be arbitrarily large, and the radius $r_a$ and the exponent $\kappa_a$ can be arbitrarily small.

 Theorem \ref{th:main} generalizes Theorem \ref{th:other} in the situation that the assumption $a \in C^\infty_c \setminus \{ 0 \}$ is replaced by (a2)-(a3). In Theorem \ref{th:main}, we also formulate precisely the conditions bearing on the initial wavenumber $k \in \R$ and the physical parameter $\theta_e > 0,$ and we describe the final observation time $T \sqrt \e |\ln \e|.$ 

\begin{theo} \label{th:main} Given $k_{max} > k_{min} > \sqrt 3,$ for $\theta_e$ small enough: given $|k| \in [k_{min}, k_{max}],$ given an initial WKB amplitude $a$ satisfying conditions {\rm (a0)-(a3)}, given $K > 0$ and $K' > 0:$ if $s_a > 0$ is large enough, for an appropriate choice of $\phi_\e \in C^\infty_c(\R^3)$ as in \eqref{varphi:initial:perturbation}, for any small enough $\e > 0,$ the initial-value problem {\rm (EM)}-\eqref{initial:datum} has a unique solution $u \in C^0([0, T(\e)], H^{s}(\R^3))$ for some $s > d/2,$ which satisfies 
$$ \| (u - u_a)(T(\e),\cdot) \|_{L^\infty} \geq \e^{K'},$$
for $T(\e) \sim T \sqrt \e |\ln \e|,$ with $T = (K - K')/\g,$ where $\g$ is the optimal linear amplification rate, which depends only on the Euler-Maxwell system and the initial WKB amplitude $a.$ 
 \end{theo}
 
The key unstable frequencies (that is $\xi_0$ in the initial perturbation \eqref{varphi:initial:perturbation}) are pictured on Figure \ref{fig:resonances:unstable}.

\subsection{Main ideas in the proof of Theorem \ref{th:main}: resonances as the locus of weak hyperbolicity, and the symbolic flow method} \label{sec:outline}

There are two key steps in the proof. We believe that each is of interest beyond the scope of this article. 

The first is a series of changes of variables which reveal that {\it resonances are the locus of weak hyperbolicity} for a linearized operator about the WKB approximate solution. That operator is deduced from the linearized (EM) system by admissible changes of variables. By admissible changes of variables, we mean that bounds for the unknown after changes of variables will translate into bounds for the original unknown (the solution of linearized (EM) we started from). 

The second is an application of the {\it symbolic flow method} to our context, in which the symbolic flow equation is a linear partial differential equation (as opposed to an {\it ordinary} differential equation, in previous uses of the symbolic flow method \cite{duh,LNT,em4}). 

\medskip

{\it Resonances as the locus of weak hyperbolicity.} The focus is on the linearized Euler-Maxwell operator about the WKB approximate solution, which, simplifying a bit\footnote{In \eqref{em:discussion}, we overlook in particular the fast dependence in $y.$ We point out the role played by fast transverse variations later on, in Section~\ref{sec:intro:spacetime}; see also the paragraph ``The role of space-time resonances" below.} for the sake of this discussion, takes the form 
\be \label{em:discussion} \d_t + \frac{1}{\e} L_1(\d_x,\d_y) + \frac{1}{\sqrt \e} L_0\big(t,\frac{x}{\e}, x,y\big).\ee
The operator $L_1$ is order-one, hyperbolic and constant-coefficient. The operator $L_0$ is order-zero, has no symmetry and has coefficients which depend on $x/\e.$ This operator comes from the quadratic terms in (EM), and thus represents ``interaction coefficients'' (that is, coupling terms after a decomposition onto the eigenmodes of the linear hyperbolic operator). By hyperbolicity of $L_1,$ and lack of symmetry of $L_0,$ the normal existence time is $O(\sqrt \e),$ in the sense that an $L^2$ estimate for \eqref{em:discussion} a priori shows a growth in time of the form $e^{t \g/\sqrt \e},$ for some $\g > 0.$

We first find {\it normal forms} for this operator, an idea that dates back to Poincar\'e (and for which a good reference is chapter 5 from \cite{Arn}). Somehow there are {\it two distinct relevant normal forms} here, due to the convective terms.

The reductions to normal forms allow to remove the interaction coefficients in $L_0$ away from resonant frequencies, defined as the zeros of the phase functions
\be \label{Phi_jj'} \Phi_{jj'}(\xi,\eta) = \l_j(\xi + k ,\eta) - \l_{j'}(\xi,\eta) - \o,\ee
where $\xi \in \R$ represent longitudinal frequencies, $\eta \in \R^2$ transverse frequencies, the $i \l_j$ are characteristic modes of $L_1$ (pictured on Figure \ref{fig1}, so that $j$ and $j'$ range from $1$ to $5$), and $\o = \sqrt{1 + k^2},$ where $k$ is the initial wavenumber in the WKB oscillations. Examples of resonances are pictured on Figures \ref{fig:resonances:stable} and \ref{fig:resonances:unstable}. The normal form reductions are performed in Sections \ref{sec:normal} and \ref{sec:normal:2}.


The first normal form involves a pseudo-differential operator of order $0.$ It allows for the removal of auto-interaction coefficients, but due to the presence of convection, it creates at the same time a spurious order-one and non-singular-in-$\e$ ``source" term. The other normal form involves a pseudo-differential operator of order $-1.$ No spurious order-one terms are created, but auto-interaction coefficients cannot be removed from $L_0$ in this way. Details are given in Section \ref{sec:normal} and \ref{sec:normal:2}.  

The second normal form reduction is used in order to derive a {\it rough Sobolev bound} for the solution operator to \eqref{em:discussion}. The bound is rough in the sense that it features a rate of exponential growth in time that is far from optimal.

The first normal form reduction is used in most of the rest of the analysis. The key is that the spurious order-one term is not singular in $\e.$ This term implies a loss of derivatives in Sobolev estimates. This issue is circumvented by use of the rough bound. 


 At this point, the non-symmetric ``source" $L_0$ is supported near the resonant set. Thus, by hyperbolicity of $L_1,$ away from the resonant set the solution operator for \eqref{em:discussion} is bounded in time from $L^2$ to itself. This implies that we can restrict the analysis to a neighborhood of the resonant set, which happens to be bounded.
  The operators, $L_1$ and $L_0,$ are thus now compactly supported in frequency, in particular they are order zero.

Next we exploit some ``separation" properties of the resonant set, which essentially state that the zeros of $\Phi_{j_1j_2}$ and $\Phi_{j'_1j'_2}$ are distinct if $(j_1,j_2) \neq (j'_1,j'_2),$ in order to show that an {\it equivalent} operator has the form 
\be \label{em:discussion:2}
 \d_t + \frac{1}{\e} \tilde L_1(\d_x,\d_y) + \frac{1}{\sqrt \e} \tilde L_0\big(t,x,y\big).\ee
The operator \eqref{em:discussion:2} has three key features:
\begin{itemize}
\item it is equivalent to \eqref{em:discussion} in the sense that bounds for the solution operator for \eqref{em:discussion:2} imply bounds for the solution operator for \eqref{em:discussion};
\item the resonant set (the zeros of the $\Phi_{jj'}$ phase functions) is the locus of weak hyperbolicity for $\tilde L_1$ (more about this below). 
\item the operator $\tilde L_0$ does not depend singularly on $x$ (as opposed to $L_0,$ see \eqref{em:discussion}).
\end{itemize}
The equality between resonant set and locus of weak hyperbolicity shows that the operator \eqref{em:discussion:2} is especially relevant in this analysis: at frequencies for which two eigenvalues of $\tilde L_1$ coalesce, the lower-order (in $\e$) operator $\tilde L_0$ comes into play and may possibly cause the spectrum to bifurcate away from the imaginary axis into the rest of the complex plane.

Indeed, simplifying again we consider that \eqref{em:discussion:2} takes the form
\be \label{em:discussion:3}
 \d_t + \frac{1}{\e} \op_\e \left( \begin{array}{cc} i (\l_{j}(\cdot + k) - \o) & \sqrt \e b_{jj'}^+ \\ \sqrt \e b_{j'j}^- & i \l_{j'} \end{array} \right),
\ee
where $\op_\e$ denotes an appropriate anisotropic and semi-classical quantization of operators, and the interaction coefficients $b_{jj'}^+$ and $b_{j'j}^-$ depend on $(x,y,\xi,\eta).$ The spectrum of the symbol in \eqref{em:discussion:3} is 
\be \label{em:discussion:4}
 \frac{1}{2 \e} \left( i \big( \Phi_{jj'} \pm \Big( - \Phi_{jj'}^2 + 4 \e b_{jj'}^+ b_{j'j}^- \Big)^{1/2} \, \right).
\ee  
Away from the zeros of $\Phi_{jj'},$ the operator \eqref{em:discussion:3} is hyperbolic. At frequencies such that $\Phi_{jj'} = 0,$ the spectrum stays purely imaginary if and only if $b_{jj'}^+ b_{j'j}^- < 0.$ If $b_{jj'}^+ b_{j'j}^- > 0$ at a resonant frequency, then the operator \eqref{em:discussion:4} is elliptic and instabilities may arise\footnote{For the related but different problem of the {\it onset of instability} occurring from spectra bifurcating away from $i \R$ as time reaches a critical value, see \cite{Luslow} for a study in geometric optics; for general first-order operators, we refer to \cite{LMX} for scalar equations; the results of \cite{LMX} were extended to systems in Sobolev spaces in \cite{LNT,NT1,NT2}, and to Gevrey spaces in \cite{Mo1}.}. A more refined analysis shows however that the relevant spectral computation is not \eqref{em:discussion:4}: due to the presence of fast transverse variations, it is {\it space-time resonances} that matter. See the discussion on the role of space-time resonances, below.

\medskip

 {\it The symbolic flow.} In this second key step in the proof, we show how to translate the above spectral analysis of the symbol in \eqref{em:discussion:3} or \eqref{em:discussion:2} into bounds for the solution operator for \eqref{em:discussion:2}, by application of the symbolic flow method\footnote{For more on the symbolic flow method, and its connection with G\r{a}rding's inequality, which is the standard way of translating spectral information for the symbol into bounds for the operator, see \cite{duh}.}.

 The symbolic flow method gives us an approximate solution operator for \eqref{em:discussion:2} (in which the operators are order zero), in short time $O(\sqrt \e |\ln \e|),$ which is precisely the relevant time window. 
 
 The approximate solution operator is given by the pseudo-differential operator in anisotropic quantization
  $$ \op_\e(S(0;t)),$$
  where the symbol $S = S(0;t,x,y,\xi,\eta)$ solves the linear partial differential equation in $(t,y):$ 
 \be \label{em:discussion:5} \d_t S +  \frac{1}{\e} \tilde L_1(i \xi , i \eta) S + \frac{1}{\sqrt \e} \tilde L_0(x,y)  S  + \frac{(-i)}{\sqrt \e} \d_\eta \tilde L_1(i \xi, i \eta) \cdot \d_y S = 0,\ee 
  with the datum $S(\t;\t) \equiv {\rm Id}.$

  Thus bounds for $S$ translate into bounds for the linearized (EM) system. We show that {\it at space-time resonances} (and not just resonances) $(\xi,\eta),$ defined as the solutions of
 \be \label{em:discussion:3.5}
  \Phi_{jj'}(\xi,\eta) = 0, \qquad \d_\eta \Phi_{jj'}(\xi,\eta) = 0,
 \ee
  the solution $S$ to the above linear system grows exponentially in time. We also find the optimal rate of growth for $S.$  This translates into an optimal rate of exponential growth for the solution to the linearized (EM) system \eqref{em:discussion}, and concludes the proof. The key space-time resonance is a Raman resonance between fast and slow Klein-Gordon modes, and the growth rate involves the interplay between the current density and Lorentz force terms.

 Above, the pseudo-differential operator $\op_\e(S)$ is defined by its action on the Schwartz class by
 $$ \op_\e(S(0;t)) u(t,x,y) = \int_{\R^3} e^{ i x \xi + i y \cdot \eta} S(0;t,x,y,\e \xi, \sqrt \e \eta) \hat u(\xi,\eta) \, d\xi \, d\eta,$$
 where $\hat u$ is the Fourier transform of $u.$ The above anisotropy reflects the anisotropy in the WKB amplitude \eqref{data}. In the next paragraph, we explain in more detail how fast transverse variations shift the focus from resonances to space-time resonances.

\bigskip

{\it The role of space-time resonances.}
Consider the symbol of the differential operator in the symbolic flow equation \eqref{em:discussion:5}:  
$$ {\mathcal M}(\e,t,x,\xi,\eta;y,\hat y) := \frac{1}{\e} \Big(\tilde L_1 + \sqrt \e \tilde L_0 + \sqrt \e i \hat y \cdot \d_\eta \tilde L_1\Big).$$
 Here $\e,t,x,\xi,\eta$ are parameters, and $\hat y$ is the dual Fourier variable of $y.$ With $\tilde L_1$ and $\tilde L_0$ in the form \eqref{em:discussion:3}, we have 
 $$ {\mathcal M} = \frac{1}{\e} \left( \begin{array}{cc} i (\l_j(\cdot + k) - \o + \sqrt \e \hat y \cdot \d_\eta \l_1(\cdot + k)) & \sqrt \e b^+_{jj'} \\ \sqrt \e b^-_{j'j} & i (\l_{j'} + \sqrt \e \hat y \cdot \d_\eta \l_{j'}) \end{array}\right),$$
 with eigenvalues
 $$ \l^\pm_{\mathcal M} = \frac{1}{2} \mbox{tr}\, {\mathcal M} \pm \frac{1}{2\e} \delta_{\mathcal M}^{1/2},$$
 where
 $$
 \delta_{\mathcal M} := 
  - \Phi_{jj'}^2 
  - \sqrt \e \hat y \cdot \d_\eta(\Phi_{jj'}^2) 
  + \e (4 b_{jj'}^+b_{j'j}^- - (\hat y \cdot \d_\eta\Phi_{jj'})^2),$$
with $\Phi_{jj'}$ the phase function defined in \eqref{Phi_jj'}.
 The trace of ${\mathcal M}$ is purely imaginary. The nature of its spectrum $\l_{\mathcal M}^\pm$ depends on the sign of $\delta_{\mathcal M}:$

\medskip

 $\bullet$ {\it away from resonances}, that is for $(\xi,\eta) \in \R^3$ such that $\Phi_{jj'}(\xi,\eta) \neq 0:$ for small $\e,$ we have $\delta_{\mathcal M} < 0,$ since then the term $- \Phi_{jj'}^2$ dominates in $\delta_{\mathcal M}.$ At such frequencies, the symbol ${\mathcal M}$ is {\it hyperbolic}. 
 
 \medskip
 
 $\bullet$ {\it at a resonant frequency $(\xi,\eta)$ which is not a space-time resonant frequency:} $\Phi_{jj'}(\xi,\eta) = 0$ but $\d_\eta \Phi_{jj'}(\xi,\eta) \neq 0.$ We observe that for such frequences 
 $$ \delta_{\mathcal M} = - \e (\hat y \cdot \d_\eta \Phi_{jj'}(\xi,\eta))^2 + 4 \e b_{jj'}^+b_{j'j}^-.$$
 Relevant spatial coordinates $y$ are small(-ish), since the initial WKB amplitude is assumed to be compactly supported. (In the slightly more general version of our argument, we can invoke the spatial decay of the initial WKB amplitude: assumption (a3) on page \pageref{page:a3}). By the Heisenberg uncertainty principle,  relevant frequencies $\hat y$ are thus large. In particular, away from a small cone in the $\hat y$ plane around the line orthogonal to $\d_\eta \Phi_{jj'}(\xi,\eta),$ the term $(\hat y \cdot \d_\eta \Phi_{jj'}(\xi,\eta))^2$ dominates in $\delta_{\mathcal M},$ hence $\delta_{\mathcal M} < 0,$ corresponding again to a hyperbolic symbol.

 \medskip
 
 $\bullet$ {\it at a space-time resonance}: at $(\xi,\eta)$ such that \eqref{em:discussion:3.5} holds true, then 
 $$ \delta_{\mathcal M} = 4 \e b_{jj'}^+ b_{j'j}^-,$$
 and the symbol is hyperbolic if $b_{jj'}^+ b_{j'j}^- \leq 0$ and elliptic otherwise.

 \medskip

The conclusion of this informal discussion is that instabilities for the solution $S$ to \eqref{em:discussion:5} are susceptible to develop only at space-time resonances, and only in the elliptic case. Our analysis puts these arguments on a rigorous footing (see in particular Section \ref{sec:flow}).

\subsection{A little more on the proof of Theorem \ref{th:main}: on a loss of derivatives in Sobolev estimates, a rough bound, and overcoming the large norm of the Sobolev embedding in a high-frequency setting} \label{sec:a:little:more}

While Section \ref{sec:outline} exposed mostly the algebraic aspects of our proof, we sketch here how bounds are derived. 

\medskip

{\it On a loss of derivatives.} As described in the above Section \ref{sec:outline}, our first normal form reduction brings out an order-one, non-skew-symmetric operator, due to the presence of convective terms in the (EM) system. This results in a loss of derivatives, and, after much of the analysis sketched out in Section \ref{sec:outline}, estimates which in a considerably simplified form look like 
\be \label{d:21:0} \| \dot u \|_{\e,s_1} \lesssim \e^K e^{t \g/\sqrt \e}  + \e^{-1/2 + K'} \int_0^t e^{(t - t') \g/\sqrt \e} \| \dot u (t') \|_{\e,s_1 + 1} \, dt'\ee
and
\be \label{d:21:.5}  
\| D^m \dot u \|_{\fl1} \lesssim \e^K e^{t \g/\sqrt \e}  + \e^{-1/2 + K'} \int_0^t e^{(t - t') \g/\sqrt \e} \| D^{m+1} \dot u (t') \|_{\fl1} \, dt'.\ee

Above,
\begin{itemize}
\item $\dot u$ is the perturbative unknown, which can be thought of as being more or less equal to $u - u_a$ (difference between exact and approximate solution), 
\item the norm $\| \cdot \|_{\e,s_1}$ is a suitably weighted $H^{s_1}$ norm;
\item the index $s_1$ will eventually be chosen large, and smaller than the Sobolev regularity $s$ of the exact solution; 
\item the exponential $e^{t \g /\sqrt \e},$ with $\g > 0,$ is the $\| \cdot \|_{\e,s_1} \to \| \cdot \|_{\e,s_1}$ norm of an appropriate linear solution operator;
\item the prefactor $\e^{-1/2 + K'}$ reflects both the singularity $\e^{-1/2}$ in front of the nonlinear terms in the (EM) system, and the maximal size $O(\e^{K'})$ of the perturbative unknown within the observation window $O(\sqrt \e |\ln \e|);$
\item the norm $\fl1$ is the $L^1$ norm of the Fourier transform; in \eqref{d:21:.5}, the operator $D^m$ denotes up to $m$ suitably weighted derivative. If the perturbative unknown $\dot u$ belongs to $H^s,$ then by Sobolev embedding $D^m \dot u$ belongs to $\fl1$ for $m < s - 3/2.$ 
\end{itemize}
The estimate \eqref{d:21:0} is not in closed form. Fortunately the second normal form reduction allows for an estimate in closed form, but with a much worse rate of exponential growth.

\medskip

{\it On a rough Sobolev bound.} As a consequence of the second normal form reduction, relatively standard estimates yield the bound (see Section \ref{sec:weak:Sobolev:bound} for details) 
\be \label{d:21:1}
\| \dot u \|_{\e,s} \lesssim \e^K e^{t \G/\sqrt \e},
\ee
for $s$ as large as $s_a,$ and some $\G > 0$ which a priori is much larger than the $\g$ that appears in \eqref{d:21:0}. The bound \eqref{d:21:1} is our {\it rough Sobolev bound}. It holds true only {\it so long as} a suitably weighted $W^{1,\infty}$ norm of $\dot u$ is controlled by $\e^{K'}.$

 \medskip
 
{\it A closed and refined Sobolev bound.} We can iterate the loss bound \eqref{d:21:0} a large (equal to $s - s_1 - 1$) number of times, and, taking advantage of the smallness of the final observation time (equal to $O(\sqrt \e |\ln \e|$), obtain something like 
\be \label{d:21:2} \| \dot u \|_{\e,s_1} \lesssim \e^K e^{t \g/\sqrt \e}  + \e^{-1/2 + K'(s-s_1 - 1)/2} |\ln \e|^\star \int_0^t e^{(t - t') \g/\sqrt \e} \| \dot u (t') \|_{\e,s} \, dt'.\ee
Above, $|\ln \e|^\star$ is some power of $|\ln \e|,$ destined to be absorbed into a positive power of $\e.$ 
In the upper bound of \eqref{d:21:2}, we now use the rough bound \eqref{d:21:1}. Under the condition\footnote{This informal discussion only reflects our proof approximately, and condition \eqref{d:21:2.5} is not exactly the one that shows up in our analysis. See the proof of Proposition \ref{cor:the:upper:bound} and condition \eqref{s:s1}.}
\be \label{d:21:2.5}
s \geq 1 +  s_1 + \frac{2 (K - K')}{K'} \cdot \frac{\G - \g}{\g},
\ee
 this leads to the refined Sobolev bound 
\be \label{d:21:3}
 \| \dot u \|_{\e,s_1} \lesssim \e^K e^{t \g/\sqrt \e}.
\ee
(This bound is rigorously derived in Proposition \ref{cor:the:upper:bound}, in one of the final steps of the proof.) We see on \eqref{d:21:2.5} that for $K$ large and $K'$ small, we have the refined bound \eqref{d:21:3} for Sobolev indices much smaller than the ones for which the rough bound \eqref{d:21:1} holds true. 

\medskip

{\it Overcoming the large norm of the Sobolev embedding.} At this point a nagging issue shows up: in weighted norms, the Sobolev embeddings $H^s \hookrightarrow W^{m,\infty}$ for $m < s - 3/2,$ and $H^s \hookrightarrow \fl1,$ for $0 < s - 3/2,$ have large norms. Precisely, given $u \in H^s,$ we only have a priori
$$ \| D^m u \|_{L^\infty} \lesssim \| D^m u \|_{\fl1} \lesssim \e^{-1/2} \| u \|_{\e,s},$$
if $\| \cdot \|_{\e,s}$ involves $(\e \d_x, \d_y),$ with $x \in \R$ and $y \in \R^ 2.$ Thus the use of the Sobolev embedding in \eqref{d:21:.5} results in a loss equal to $\e^{-1/2}.$ That is, iterating \eqref{d:21:.5} as we did above for \eqref{d:21:0}, we obtain
$$ \|D^1 \dot u \|_{\fl1} \lesssim \e^K e^{t \g/\sqrt \e} + \e^{-1/2 + K'(m - 3)/2} |\ln \e|^\star \int_0^t e^{(t - t') \g/\sqrt \e} \| D^{s_1 - 2} \dot u \|_{\fl1} \, dt',$$
and then, using the Sobolev embedding and the refined bound \eqref{d:21:3}:
\be \label{d:21:6}
\| D^1 \dot u \|_{\fl1} \lesssim \e^K e^{t \g/\sqrt \e} + \e^{K -1/2 + K'( s_1 - 3)/2} |\ln \e|^\star e^{t \g/\sqrt \e}.
\ee 
 Thus under a condition of the form (here again, this condition is not exactly the one that appears in the analysis) 
\be \label{cond:s1:intro}
 - \frac{1}{2} + \frac{K'(s_1 - 3)}{2} > 0,
\ee
we find for $\| D^1 \dot u \|_{\fl1}$ the same upper bound as in \eqref{d:21:1}:
\be \label{d:21:7}
 \| D^1 \dot u \|_{\fl1} \lesssim \e^K e^{t \g/\sqrt \e}.
\ee

\medskip

{\it Continuation of the bounds up to the optimal amplification time.} From \eqref{d:21:7}, we see that the refined Sobolev bound \eqref{d:21:1}, which a priori held only conditional to a $W^{1,\infty}$ bound, holds true until the amplification time $T \sqrt \e |\ln \e|,$ with $T = (K - K')/\sqrt \e.$ It is the time at which the $\fl1$ bound \eqref{d:21:7} gives $\| D^1 \dot u \|_{\fl1} \lesssim \e^{K'}.$

\medskip

{\it Optimal rate of growth and lower bound.} A key feature of the upper bound \eqref{d:21:1} is that $\g$ is {\it optimal}, meaning that we are able to show a lower bound of the form
\be \label{d:low} \| \dot u \|_{L^2(B(x_0,y_0,\e^\b))} \geq \e^{K + \b d/2} e^{t \g/\sqrt \e},\ee 
where the exponent $\b$ that is involved in the radius of the ball is not too large, and $(x_0,y_0)$ is a point in space at which the norm of the initial WKB amplitude $a$ is maximal. Thanks to the above upper bounds, the lower bound \eqref{d:low} holds until the amplification time $T \sqrt \e |\ln \e|,$ with $T = (K - K')/\g,$ which concludes the proof.

\subsection{Another viewpoint on fast transverse variations and space-time resonances: points of stationarity for a relevant phase} \label{sec:intro:spacetime}

We put forward another viewpoint here, in which resonances are seen as points of stationarity for a relevant phase. This viewpoint is not the one adopted in the present article, but it might be interesting nonetheless.

 One key difference with \cite{em4} is the consideration of fast transverse variations in the data. We explain here how these fast variations are likely to play a stabilizing role, and why space-time resonances matter here, while the focus in \cite{em4} was on time resonances. Convective phenomena are overlooked in this discussion.

We consider large longitudinal frequencies of size $1/\e$ in $x,$ the direction of propagation. This calls for considering $\e \d_x$ derivatives. Similarly, we consider fast variations in $y,$ with a characteristic length of variation equal to $\sqrt \e.$ This calls for considering $\sqrt \e \d_y$ derivatives.

The focus is on the perturbation variable $\dot u$ defined by $u = u_a + \dot u,$ where $u$ solves Euler-Maxwell and $u_a$ is the approximate solution to Euler-Maxwell described above in Section \ref{sec:ZEM}. The goal is to prove that $\dot u$ grows exponentially in short time.

For the sake of simplicity, we overlook convective terms in this discussion. The perturbation variable $\dot u$ solves a system of the form
\be \label{eq:pert:intro}  \d_t \dot u + \frac{1}{\e} {\mathcal A}(\e \d_x, \sqrt \e \d_y) \dot u = \frac{1}{\sqrt \e} {\mathcal B}(u_a) \dot u,\ee
where ${\mathcal B}(u_a) \dot u$ is linear in $u_a$ and $\dot u.$
The hyperbolic operator ${\mathcal A}$ has a smooth spectral decomposition (see Figure 1 below) into Klein-Gordon modes and null modes:
\be \label{spectral:dec:intro} {\mathcal A}(\e \d_x, \sqrt \e \d_y) = \sum_{1 \leq j \leq 5} i \l_j(D_\e) \Pi_j(D_\e),\ee
where the $\l_j$ are real eigenvalues and the $\Pi_j$ are eigenprojectors. The anisotropic Fourier multipliers are defined by
$$ m(D_\e) f = {\mathcal F}^{-1} \big( m(\e \xi, \sqrt \e \eta) {\mathcal F} f \big),$$
where ${\mathcal F}$ denotes Fourier transform, and $(\xi,\eta)$ are the dual Fourier variables of $(x,y).$

In \eqref{spectral:dec:intro}, the eigenvalues $\l_1$ and $\l_2$ are of Klein-Gordon type:
$$ \l_1(\xi,\eta) = \sqrt{1 + |(\xi,\eta)|^2}, \qquad \l_2(\xi,\eta) = \sqrt{1 + \theta_e^2 |(\xi,\eta)|^2},$$
where $\theta_e > 0$ is a physical parameter in (EM). The mode $\l_3 \equiv 0$ is a degenerate acoustic mode. The other modes are $\l_4 = - \l_2$ and $\l_5 = -\l_1.$ This decomposition of Euler-Maxwell into Klein-Gordon and acoustic modes dates back at least to \cite{em1}.

The constant-coefficient hyperbolic operator \eqref{spectral:dec:intro} has a solution operator
\be \label{solution:op:intro} S(t, D_\e) = \exp\Big( \frac{t}{\e} {\mathcal A}(\e \d_x, \sqrt \e \d_y)\, \Big) = \sum_{1 \leq j \leq 5} e^{i t \l_j(D_\e)/\e} \Pi_j(D_\e).\ee

We may then represent the solution $\dot u$ to the perturbation equations \eqref{eq:pert:intro} in the form
\be \label{dotu:intro} \dot u(t) = S(t,D_\e) \dot u(0) +  \frac{1}{\sqrt \e} \int_0^t S(t - t', D_\e) {\mathcal B}(u_a) \dot u(t') \, dt'.\ee
For very short time $t  \ll \sqrt \e,$ the cumulated response of the large source term is small, and we expect the free solution to approximate $\dot u(t'):$
\be \label{dotu:intro2} \dot u(t) \simeq S(t,D_\e) \dot u(0), \qquad t \ll \sqrt \e.\ee
Plugging \eqref{dotu:intro2} into \eqref{dotu:intro}, we find that
$$ \dot u(t) \simeq S(t, D_\e) \dot u(0) + \frac{1}{\sqrt \e} \int_0^t S(t - t', D_\e) {\mathcal B}(u_a) S(t', D_\e) u(0)\, dt'.$$
The question is whether the perturbation $\dot u$ grows in time. The solution operator $S$ is unitary by hyperbolicity of ${\mathcal A}.$ The initial perturbation $\dot u$ is very small, and so $S(t, D_\e) \dot u(0)$ stays very small in time. Thus the focus is on a bound for the integral
$$ \frac{1}{\sqrt \e} \int_0^t S(t - t', D_\e) {\mathcal B}(u_a) S(t', D_\e) u(0)\, dt'.$$
Using the description of the solution operator in \eqref{solution:op:intro}, and estimating commutators, we find the above integral to be a sum of terms of the form
\be \label{osc:intro} \frac{1}{\sqrt \e} \int_0^t e^{i t' \Phi(D_\e)/\e} f(t') \, dt' = \frac{1}{\sqrt \e} \int_0^t \int_{\R^3} e^{i x \xi + i y \cdot \eta + i t' \Phi(\e \xi, \sqrt \e \eta)/\e} \hat f(t', \xi ,\eta) \, d\xi \, d\eta \,dt'\,,\ee
where the phase $\Phi$ and the interaction coefficient $f$ are defined as, for $p \in \{-1,1\}:$ 
$$ \Phi(\xi,\eta) = \l_j(\xi + p k,\eta) - p \o - \l_{j'}(\xi,\eta), \quad f(t') = \Pi_j(pk + D_\e) {\mathcal B}(u_{0,p}(t')) \Pi_{j'}(D_\e) \dot u(0).$$
We see above that the key phase and interaction coefficients depend on (in the order of apparition):
\begin{itemize}
\item the eigenmodes $\l_j$ of the hyperbolic operator, introduced in \eqref{spectral:dec:intro};
\item the fundamental phase $(\o, k) \in \R \times \R^3$ of the WKB solution $u_a$ (see Section \ref{sec:ZEM} and in particular \eqref{wkb:intro});
\item the eigenprojectors $\Pi_j$ of the hyperbolic operator, introduced in \eqref{spectral:dec:intro};
\item the linearized source ${\mathcal B},$ which comprises convective terms, current densities and Lorentz forces from the physical (EM) equations introduced in Section \ref{sec:EM};
\item the leading WKB profiles $u_{0,1}$ and $u_{0,-1} = u_{0,1}^*$ introduced in \eqref{wkb:intro}.
\end{itemize}

Resonant frequencies are defined as the zeros of $\Phi.$ Far from resonant frequencies, we may integrate by parts in time in \eqref{osc:intro}. To illustrate this point, consider a resonant frequency $(\xi_0, \eta_0) \in \R^3,$ and assume for simplicity that it is the only zero of $\Phi.$ We focus on the restriction of the oscillatory integral in a vicinity of the resonant frequency:
\be \label{osc:intro:disc1} \frac{1}{\sqrt \e} \int_0^t \int_{|(\e \xi, \sqrt \e \eta) - (\xi_0, \eta_0)| \geq 1} e^{i x \xi + i y \cdot \eta + i t' \Phi(\e \xi, \sqrt \e \eta)/\e} \hat f(t', \xi ,\eta) \, d\xi \, d\eta \,dt'\,.\ee
We assumed for simplicity that $\Phi$ is bounded away from zero in the domain in \eqref{osc:intro:disc1}. We may thus integrate by parts, and the integral in \eqref{osc:intro:disc1} appears to be $O(\sqrt \e)$ since $\d_t \hat f$ is uniformly bounded in $\e$ (indeed, the dependence of $f$ in time is through the WKB amplitude, which does not depend on $\e$).

\medskip

Next consider the oscillatory integral \eqref{osc:intro} in the vicinity of a resonant frequency. The question is whether it can be bounded uniformly in $\e,$ in short time. Let $(\xi_0,\eta_0)$ such that $\Phi(\xi_0, \eta_0) = 0.$ That is, we consider
\be \label{osc:intro:2}
\frac{1}{\sqrt \e} \int_0^t \int_{\begin{smallmatrix} |\xi - \xi_0/\e| \leq R \\ |\eta - \eta_0/\sqrt \e| \leq R \end{smallmatrix}} e^{i x \xi + i y \cdot \eta} e^{i t \Phi(\e \xi, \sqrt \e \eta)/\e} \hat f(t',\xi,\eta) \, d\xi \, d \eta \, dt'.
\ee
We expand the phase near the resonant frequency, as we may if the eigenmodes are smooth (they are in Euler-Maxwell, in spite of the coalescing point at the origin): for $(\xi,\eta)$ such that $|\e \xi - \xi_0) \leq \e R$ and $|\sqrt \e \eta - \eta_0| \leq \sqrt \e R:$
$$ \Phi(\e \xi, \sqrt \e \eta) = \Phi(\xi_0, \eta_0) + (\e \xi - \xi_0) \d_\xi \Phi(\xi_0, \eta_0) + (\sqrt \e \eta - \eta_0) \cdot \d_\eta \Phi(\xi_0,\eta_0) + O(\e),$$
that is
$$ \Phi(\e \xi, \sqrt \e \eta) = (\sqrt \e \eta - \eta_0) \cdot \d_\eta \Phi(\xi_0, \eta_0) + O(\e).$$
We plug the above (minus the irrelevant $O(\e)$ remainder) in \eqref{osc:intro:2} to find
$$ \frac{1}{\sqrt \e} \int_0^t e^{-i c t' \cdot \eta_0/\e} \int_{\begin{smallmatrix} |\xi - \xi_0/\e| \leq R \\ |\eta - \eta_0/\sqrt \e| \leq R \end{smallmatrix}} e^{i x \xi + i \eta \cdot( y + ct /\sqrt \e) } \hat f(t',\xi,\eta) \, d\xi \, d\eta \, dt',$$
where $c = \d_\eta \Phi(\xi_0,\eta_0).$
We may ignore the purely oscillatory time factor. We also ignore the frequency domain in this informal discussion, and focus on
\be \label{space-time-res-last} \frac{1}{\sqrt \e} \int_0^t \int_{\R^3} e^{i x \xi + i \eta \cdot(y + c t /\sqrt \e)} \hat f(t',\xi,\eta) \, d\xi\, d\eta \, dt' = \frac{1}{\sqrt\e} \int_0^t f\left(t',x,y + \frac{c t'}{\sqrt \e}\right) \, dt'.\ee
We distinguish two subcases:

\medskip

$\bullet$ {\it The resonant frequency is a space-time resonance:} $c = 0.$ Then the integral is typically $O(t/\sqrt\e),$ and we can expect an exponential amplification in the linearized equations. Note however that if the ratio $\hat f/\Phi$ is bounded near the resonant frequency $(\xi_0,\eta_0)$ then it is possible to integrate by parts in the integral \eqref{osc:intro:2} in spite of $\Phi(\xi_0,\eta_0) = 0,$ and this provides a uniform bound in $\e.$ The condition $\hat f/\Phi$ bounded is a {\it compatibility} condition (called ``null form" in \cite{Kla} and ``transparency" in \cite{JMR-TMB}). Many physical equations exhibit this type of property: in particular Maxwell-Bloch (see \cite{JMR-TMB}), Maxwell-Euler (see \cite{em3}, in a slightly different context than here) and Maxwell-Landau-Lifshitz (see \cite{Lu}). The present result (Theorem \ref{th:main}) reveals a {\it lack of transparency} for Euler-Maxwell: the ratio $\hat f/\Phi$ is not bounded, so that an integration by parts is not possible in \eqref{osc:intro:2}, which leads the rightmost integral in \eqref{space-time-res-last}, and if $c = 0$ we record an instability.

\medskip

$\bullet$ {\it The resonant frequency is a time resonance but not a space-time resonance} or $c \neq 0.$ Then, we may posit $s = t'/\sqrt \e,$ and the above integral becomes
$$ \int_0^{t/\sqrt\e} f( \sqrt \e s, x, y + c s) \, ds.$$
Since $f$ is differentiable in time, and $\d_t f$ is uniformly bounded in $\e,$ in short time the above integral is approximated by
$$ \int_0^{t/\sqrt \e} f(0, x, y + c s) \, ds,$$
and since $c \neq 0,$ this is the integral along a line in the plane $\R^2_y$ of the function $f$ which belongs to $L^1.$ The integral is finite, and we obtain an $\e$-uniform bound.

\section{The perturbation equations}  \label{sec:mainproof}

The proof of Theorem \ref{th:main} starts here.

\medskip

We work with the Euler-Maxwell system (EM) introduced in Section \ref{sec:EM}. In this preparatory step, we set up notations and introduce the perturbation equations \eqref{pert:2}-\eqref{def:mathcaldotF} about the WKB approximate solution.

\subsection{Coordinatization} \label{sec:coord}

We coordinatize solutions $u$ to (EM) as
$$ u = (B, E, v_e, n_e, v_i, n_i) = (B, E |\, v_e, n_e |\, v_i, n_i) \in \R^6 \times \R^4 \times \R^4,$$
with vertical dividers for legibility. Solutions $u$ depend on time $t \geq 0$ and space $(x,y) = (x,y_1,y_2) \in \R \times \R^2,$ and of course on the small parameter $\e > 0.$

\subsection{The $\fl1$ norm} \label{sec:fl1}

 Introduce the norm
 \be \label{def:fl1}
 \| f \|_{{\mathcal F} L^1} := \| \hat f \|_{L^1(\R^3_{\xi,\eta})}.
\ee
We note that the $\fl1$ norm controls the $L^\infty$ norm:
\be \label{Linfty:fl1}
 \| f \|_{L^\infty} \lesssim \| f \|_{\fl1},
\ee
and that a control in $H^{2}(\R^3)$ gives a control in $\fl1:$
\be \label{fl1:sob} 
 \| f \|_{\fl1} \lesssim \| f \|_{H^2(\R^3)}, \qquad \mbox{since $2 > d/2 = 3/2.$} 
 \ee
Also, $\fl1$ is an algebra:
\be \label{algebra}
 \| f g \|_{\fl1} \leq \| f \|_{\fl1} \| g \|_{\fl1}.
\ee
Finally, we will use the elementary interpolation inequality 
\be \label{interpol} 
\| \d_{x,y}^{\a_1} f \d_{x,y}^{\a_2} g \|_{\fl1} \leq \| f \|_{\fl1} \cdot \max_{|\b| = |\a_1| + |\a_2|} \| \d_{x,y}^{\b} g \|_{\fl1} + \| g \|_{\fl1} \cdot \max_{|\b| = |\a_1| + |\a_2|} \| \d_{x,y}^{\b} f \|_{\fl1}.
\ee 
Inequality \eqref{interpol} persists if derivatives are replaced with weighted derivatives. In particular, we will use \eqref{interpol} with $\d_\e$ derivatives; the operator $\d_\e$ is defined below in \eqref{nablae}. 

\subsection{The WKB approximate solution} \label{sec:new:wkb}

 It has the form
 $$ u_a = \sum_{0 \leq k \leq K'_a} \e^{k/2} u_k, \qquad u_k = \sum_{q \in {\mathcal H}_k} e^{i q ( k x - \o t)/\e} u_{k,q}(t,x,y/\sqrt \e),$$
 where the amplitudes $u_{k,q}$ depend on $\e$ only through $y/\sqrt \e.$ Above $\o = \sqrt{1 + k^2}$ is the characteristic frequency associated with the initial wavenumber $k$ via the polarization condition $\Pi_1(k) a = a.$ See Section \ref{sec:char0} for the definition of the eigenprojectors $\Pi_j,$ and Section \ref{sec:wkb} for an explicit version of the polarization condition. 
 
 The leading terms $(u_0, u_1)$ solve the Zakharov (Z) system given explicitly in Section \ref{sec:wkb}, for which an existence $T_a > 0$ is given in \cite{LPS}. In particular, the Sobolev regularity of $u_0$ and $u_1$ is $s_a$ and $s_a -1,$ respectively. 
 
 The higher-order terms $(u_k, u_{k+1})$ solve system (Z) linearized at $(u_0, u_1),$ with additional source terms which depend on spatial derivatives of $u_{k'},$ for $k' \leq k-1.$ 
 In particular, the higher-order terms are all defined over the same existence time $[0, T_a]$ as $(u_0, u_1),$ and the Sobolev regularity of $u_k$ is $s_a - k.$ As a consequence, the full approximate solution $u_a$ is defined over $[0, T_a],$ and has the Sobolev regularity $s_a - K'_a.$ 
 
 In the cascade of WKB equations, the leading term $u_0$ cancels the terms of order $O(\e^{-1}).$ The first correctors cancel the terms of order $O(\e^{-1/2}).$ And so on, so that in order to have a remainder $\e^{K_a} R_a,$ with $R_a$ bounded in $\e,$ we need $K'_a = 2 K_a + 1.$ This means that the Sobolev regularity of $u_a,$ and $R_a$ is $s_a - 2 K_a - 1.$ We have the bound 
 $$ \sup_{\e > 0} \max_{0 \leq t \leq T_a} \| (u_a, R_a)(\e,t,\cdot)\|_{L^\infty(\R^3)} + \| (\e \d_x)^\a (\sqrt \e \d_y)^\b (u_a, R_a)(\e,t,\cdot)\|_{L^2(\R^3)} < \infty,$$
 for any $(\a, \b) \in \N^3$ such that $|\a| + |\b| \leq s_a - 2 K_a - 1,$
 where $T_a > 0$ is an existence time for $u_a.$ 
 
We will find a local-in-time solution $u$ with the same Sobolev regularity as $u_a:$ 
\be \label{def:s}
 s = s_a - 2 K_a - 1,
\ee
The proof goes through if $s_a$ is large enough, depending in particular on $K$ and $K'.$ The constraint bearing on $s_a$ is specified in Remark \ref{rem:Sobolev:indices}.

\subsection{The unique solution to {\rm (EM)} issued from \eqref{initial:datum}} \label{sec:tstar}

The system (EM) is first-order, quasi-linear, symmetric hyperbolic. For fixed $\e > 0,$ the datum \eqref{initial:datum}, which belongs to $H^s$ with $s$ defined in \eqref{def:s}, generates a unique solution $u$ which is defined over some maximal time interval $[0, T_\star(\e)),$ with $T_\star(\e) > 0,$ and has the same Sobolev regularity as the datum:
\be \label{the:solution}
 u \in C^0([0, t], H^s(\R^3)), \qquad \mbox{for all $t < T_\star(\e).$}
 \ee
Moreover, the maximal existence time $T_\star(\e)$ is such that
\be \label{max:time}
\int_0^{T_\star(\e)} \| u(t') \|_{W^{1,\infty}} \, dt' = + \infty.
\ee
Define 
 \be \label{def:Tstar}
 T := \frac{K - K'}{\gamma},
 \ee 
 where $\g$ is the optimal rate of growth defined in \eqref{def:gamma}. Our goal is to show that the solution is defined over $[0, T_\e \sqrt \e |\ln \e|],$ where $T_\e$ is just a little smaller than $T$ (see Section \ref{sec:existence}) and that within that time interval, the solution $u$ diverges from the approximate solution $u_a.$

 At time $t = 0,$ the distance $u - u_a$ is very small. In particular,   
 $$ \| (\e \d_x)^\a (u - u_a)(0) \|_{L^2} +  \| (\sqrt \e \d_y)^\a (u - u_a)(0) \|_{L^2} \lesssim \e^K,$$
 with $|\a| \leq s.$ 

The function $t \to \| (u - u_a)(t) \|_{H^s}$ is continuous over $[0, T_\star(\e)).$ Maps in $\fl1$ belong to $H^s,$ via the embedding \eqref{fl1:sob}, which with weighted derivatives takes the form 
\be \label{embed:fl1:Sob}
 \| f \|_{\fl1} \lesssim \e^{-1} \sum_{|\a| + |\b|  \leq 2} \| (\e \d_x)^\a (\sqrt \e \d_y)^\b f \|_{L^2}. 
 \ee 
 Thus the function $t \to \| ( u - u_a)(t)\|_{\fl1}$ is continuous over $[0, T_\star(\e)),$ and so is the function 
 $$ N_u: \quad t \to \max_{0 \leq t' \leq t} \| ( u - u_a)(t) \|_{\fl1}.$$ Since $N_u$ is continuous and nondecreasing, $[0, \sqrt \e T |\ln \e|] \cap N_u^{-1}([0,\e^{K'}])$ is a closed interval in $[0, T_\star(\e)),$ and since we may assume $K > 1 + K',$ this interval contains $0.$ We denote $t_\star(\e)$ its endpoint, so that $t_\star(\e)$ is the largest time within $[0, \sqrt \e T |\ln \e|] \cap [0, T_\star(\e))$ for which  
 \be \label{a:priori:new}
 \max_{0 \leq t \leq t_\star(\e)} \Big( \| (u - u_a) (t) \|_{\fl1} + \| (\e \d_x, \sqrt \e \d_y) (u - u_a)(t) \|_{\fl1}\Big) \leq \e^{K'}.
 \ee

\subsection{Rescaling the transverse coordinates} \label{sec:y} In view of the form of the initial WKB profile \eqref{data}, we rescale the transverse coordinate, and define
\be \label{def:tildeu}
\tilde u(\e,t,x,y) := u\left(\e, t, x, \sqrt \e y\right).
\ee
 We will use two different weighted gradients: 
\be \label{nablae}
\nabla_\e := (\e \d_x, \sqrt \e \d_y), \quad \mbox{and} \quad \d_\e := (\e \d_x, \d_y).
\ee
 The datum for $\tilde u$ is
$$ \tilde u(0,x,y) = \cos\left(\frac{k x}{\e}\right) a(x,y) + \e^K \Re e \, \Big(\, \big( e^{i x (\xi_0 + k)/\e} \vec e_0 + e^{i x \xi_0/\e} \vec f_0 \big) \,  \phi(x,y) \, \Big).$$
The bound \eqref{a:priori:new} becomes
\be \label{a:priori:tilde:bound}
\sup_{0 \leq t \leq t_\star(\e)} \Big( \| (\tilde u - \tilde u_a)(t) \|_{\fl1} + \|\d_\e (\tilde u - \tilde u_a)(t) \|_{\fl1}\Big) \leq \e^{K'}.
\ee
After the rescaling \eqref{def:tildeu}, the WKB profiles are not singular in $y.$ In particular, 
 $$ \| \d_\e^\a \tilde u_a \|_{L^\infty} \leq C_a,$$ 
 for $|\a|$ less than some function of $s_a,$ where $s_a$ is the Sobolev regularity of the initial WKB amplitude. The constant $C_a$ is uniform in time and $\e,$ for $t \in [0,T_a],$ where $T_a > 0$ is an existence time for the WKB approximate solution. 

\subsection{The linear hyperbolic operator} \label{sec:hypop}

Consider system (EM) from Section \ref{sec:EM}. In the coordinatization introduced in \eqref{sec:coord}, the  linear spatial differential operator in (EM) is
$$ \left(\begin{array}{cc|cc|cc} 0 & \nabla \times & 0 & 0 & 0 & 0 \\ - \nabla \times & 0 & -1/\e & 0 &  1/(\theta_e \sqrt \e) & 0  \\[2pt] \hline &&&&& \vspace{-.3cm} \\ 0 & 1/\e & 0 & \theta_e \nabla & 0 & 0 \\ 0 & 0 & \theta_e \nabla \cdot & 0 & 0 & 0 \\[2pt] \hline &&&&& \vspace{-.3cm} \\ 0 & - 1/(\theta_e \sqrt \e) & 0 & 0 & 0 & \a^2 \sqrt \e \nabla \\ 0 &  0& 0 & 0 & \sqrt \e \nabla \cdot & 0 \end{array}\right).$$
In view of the change of variable of Section \ref{sec:y}, we thus find the linear differential operator in $(t,x,y)$ variables in the equations in $\tilde u$ to be
  $\dsp{\d_t + \frac{1}{\e} A(\nabla_\e)}$
with
$$ A(\nabla_\e) = \left(\begin{array}{cc|cc|cc} 0 & \nabla_\e \times & 0 & 0 & 0 & 0 \\ - \nabla_\e \times & 0 & - 1 & 0 &  \sqrt \e/\theta_e & 0  \\[2pt] \hline &&&&& \vspace{-.3cm} \\ 0 & 1 & 0 & \theta_e \nabla_\e & 0 & 0 \\ 0 & 0 & \theta_e \nabla_\e \cdot & 0  & 0 & 0 \\[2pt] \hline &&&&& \vspace{-.3cm} \\ 0 & - \sqrt \e/\theta_e & 0 & 0 & 0 & \a^2 \sqrt \e \nabla_\e  \\ 0 &  0& 0 & 0 & \sqrt \e \nabla_\e \cdot & 0  \end{array}\right).$$
The characteristic variety of the hyperbolic operator evaluated at $\e = 0$ is pictured on Figure \ref{fig1}. 

\subsection{Convection} \label{sec:conv}

The convection in Euler contributes quasilinear terms to system (EM) introduced in Section \ref{sec:EM} page \pageref{sec:EM}. In view of the change of variable of Section \ref{sec:y}, we denote
$$ {\mathcal J}(\tilde u, \nabla_\e) =  {\mathcal J}_e(\tilde v_e, \nabla_\e) + \sqrt \e {\mathcal J}_i(\tilde v_i, \nabla_\e),$$
with notation
$$ {\mathcal J}_e(\tilde v_e, \nabla_\e) :=  \theta_e (\tilde v_e \cdot \nabla_\e) \underline J_e, \quad \underline J_e :=  \left(\begin{array}{c|c|c} 0 & 0 & 0 \\[2pt] \hline && \vspace{-.3cm}  \\ 0 & {\rm Id}_{\R^4} & 0 \\[2pt] \hline && \vspace{-.3cm} \\ 0 & 0 & 0 \end{array}\right),$$
and
$$ {\mathcal J}_i(\tilde v_i, \nabla_\e) :=  (\tilde v_i \cdot \nabla) \underline J_i, \qquad \underline J_i := \left(\begin{array}{c|c|c} 0 & 0 & 0 \\[2pt] \hline && \vspace{-.3cm} \\ 0 & 0 & 0  \\[2pt] \hline && \vspace{-.3cm} \\  0 & 0 & {\rm Id}_{\R^4} \end{array}\right).$$
in the coordinatization of Section \ref{sec:coord}.

\subsection{The large semilinear source terms} \label{sec:nonlinearterms}

The current density $j$ and the Lorentz forces $L_e$ and $L_i$ contribute semilinear terms to (EM).
The nonlinear current density terms are
 $$ \frac{1}{\e} (e^{\sqrt \e n_e} - 1) v_e - \frac{1}{\sqrt \e \theta_\e} (e^{\e n_i} - 1) v_i = \frac{1}{\sqrt \e} n_e v_e + O(1).$$
The nonlinear electronic Lorentz force term is
$\dsp{-\frac{\theta_e}{\sqrt \e} v_e \times B.}$
We denote, in the coordinatization introduced in Section \ref{sec:coord}:
 \be \label{def:unB} {\mathcal B}(u,u) = \Big(0_{\R^3}, n_e v_e \big|\,  - \theta_e v_e \times B, 0\, \big|\, 0_{\R^3}, 0\,\Big).\ee
 In the equation in $\tilde u,$ the bilinear term ${\mathcal B}$ comes in with a large $1/\sqrt \e$ prefactor.
The other semilinear terms are denoted ${\mathcal F}(\tilde u).$ These incluse the $O(1)$ terms in the current density above, and the  nonlinear ionic Lorentz force term. Precisely,
\be \label{def:F(u)}
{\mathcal F}(\tilde u) := \Big( 0_{\R^3}, \, \sum_{k \geq 2} \frac{\e^{k/2-1}}{k!} \big( \tilde n_e^k \tilde v_e - \frac{\sqrt \e}{\theta_\e} \tilde n_i^k \tilde v_i \big)\,,\,\, 0_{\R^4}, \,\, \frac{\sqrt \e}{\theta_e} \tilde v_i \times \tilde B, \, 0 \, \Big).
\ee

\subsection{The equations in $\tilde u$} \label{sec:eq:tildeu}

Summing up, the initial value problem in $\tilde u$ is
\be \label{ivp:tildeu} \left\{\begin{aligned} \d_t \tilde u + \frac{1}{\e} A(\nabla_\e) \tilde u + \frac{1}{\sqrt \e} {\mathcal J}(\tilde u, \nabla_\e) \tilde u  & = \frac{1}{\sqrt \e} {\mathcal B}(\tilde u, \tilde u) + {\mathcal F}(\tilde u), \\  \tilde u(0,x,y) & = \cos(k x/\e) a(x,y) + \e^K \phi^\e,
\end{aligned}\right.\ee
with notation
\be \label{def:phie}
\phi^\e(x,y) := \Re e \, \Big(\, \big( e^{i x (\xi_0 + k)/\e} \vec e_0 + e^{i x \xi_0/\e} \vec f_0 \big) \,  \phi(x,y) \, \Big).
\ee
We will choose $\xi_0$ and $\phi,$ and vectors $\vec e_0$ and $\vec f_0$ in the course of this proof. In the left-hand side of the evolution equation, the linear symmetric hyperbolic term $A$ a priori generates only fast oscillations in space-time.  Our goal is to show that ${\mathcal B}$  generates short-time instabilities for \eqref{ivp:tildeu}.

\subsection{The perturbation equations} \label{sec:pert} Introduce the perturbation unknown $\dot u:$
\be \label{def:dotu}
 \dot u := \tilde u - \tilde u_a,
\ee
where $\tilde u_a(t,x,y) = u_a(t,x,\sqrt \e y),$ and $u_a$ is the WKB approximate solution solving the Euler-Maxwell system with a remainder $\e^{K_a} R $ for $K_a > K+1/2,$ and $R_a$ bounded in $\| \cdot \|_{\e,s}$ norm, uniformly in $\e$ and in $t$ in the time interval under consideration. Here $s$ is as in the statement of Theorem \ref{th:main}. The construction of $u_a$ is sketched in Section \ref{sec:wkb}. We will work with the perturbation unknown $\dot u$ and subsequent unknowns defined in terms of $\dot u$ in the rest of the proof.

From \eqref{ivp:tildeu} and the definition of $\dot u$ in \eqref{def:dotu}, we deduce the initial-value problem in $\dot u:$  
\be \label{ivp:dotu} \left\{ \begin{aligned}
 \d_t \dot u + \frac{1}{\e} A(\nabla_\e) \dot u + \frac{1}{\sqrt \e} \op_\e(J) \dot u & = \frac{1}{\sqrt \e} \op_\e(\dot {\mathcal B}) \dot u + \dot F, \\
  \dot u(0,x,y) & = \e^K \phi^\e.
 \end{aligned}\right.\ee
 The matrix-valued symbol $\dot {\mathcal B} \in S^0$ (see Section \ref{app:symb} for the definition of the symbol classes $S^m,$ with $m \in \R,$ and the associated operators $\op_\e(\cdot)$ in anisotropic quantization) is defined as
\be \label{def:B0} \begin{aligned}
 \dot {\mathcal B}(\e,t,x,y,\xi,\eta) & :=  {\mathcal B}(\tilde u_a, \cdot) +  {\mathcal B}(\cdot, \tilde u_a) - {\mathcal J}(\cdot, \nabla_\e) \tilde u_a \\ & \qquad - (1 - \chi_{high})(\xi,\eta))  {\mathcal J}(\tilde u_a,i \xi,i \eta) \quad \in \C^{14} \times \C^{14},
\end{aligned} \ee
the dependence in $(\e,t,x,y)$ being through $\tilde u_a.$ In the above definition of $\dot {\mathcal B},$ the smooth frequency cut-off $\chi_{high} \in C^\infty(\R^3_{\xi,\eta};\R)$ is a high-frequency cut-off, such that $0 \leq \chi_{high} \leq 1$ and
 \be \label{def:chiHF}
  \chi_{high} \equiv 0 \quad \mbox{for small frequencies}, \qquad \chi_{high} \equiv 1 \quad \mbox{for very large frequencies.}
  \ee
Above, by ``very" large frequencies, we mean a frequency threshold that depend on parameters $K$ (specifically, which tends to infinity as $K \to \infty$), but not on $\e.$ A precise choice for the support of $\chi_{high}$ is made in Section \ref{sec:prepared}; see in particular Figure \ref{fig:r:chi}. The high-frequency convective symbol $J$ is defined by
\be \label{def:J}
 J(\e,t,x,y, \xi, \eta) = \chi_{high}(\xi,\eta) {\mathcal J}(\tilde u_a(\e,t,x,y), i \xi, i \eta).
\ee
The operator $\op_\e(J)$ is the pseudo-differential operator with symbol $J$ in the anisotropic semiclassical quantization defined in Section \ref{app:symb}. 

Section \ref{sec:remainder:bounds} below gives a description of $\dot F,$ the remainder terms in \eqref{ivp:dotu}.

\subsection{The linear convective terms in the perturbation equations} \label{sec:conv:pert}

By property of the WKB solution (see \eqref{for:interaction:coefficients}), we observe the {\it cancellation}:
 \be \label{cancellation:convection}
  \tilde v_{e0} \cdot (\d_x, 0, 0) = 0,
 \ee
where $\tilde v_{e0}$ is the leading term in $\e$ in the (rescaled in $y$) WKB electronic velocity $\tilde v_{ea}$ (where {\it e} stands for electronic and {\it a} for approximate). This implies 
\be \label{cancellation:convection:2}
 \tilde v_{ea} \cdot \nabla_\e = \tilde v_{e0 y} \cdot \sqrt \e \d_y + \sqrt \e \tilde w_e \cdot \nabla_\e,
 \ee
 where $w_e$ comprises the electronic velocities of all correctors terms in the WKB expansion. 
We can write \eqref{cancellation:convection:2} in the form
\be \label{cancellation:convection:3}
\frac{1}{\sqrt \e} \tilde v_{ea} \cdot \nabla_\e = \tilde v_{e0} \cdot \d_\e + \tilde w_e \cdot \nabla_\e.
\ee
 Besides, the leading ionic velocity term vanishes in the WKB expansion (see \eqref{for:interaction:coefficients} again):
$$ \tilde v_{i0} = 0,$$
so that
$ \tilde v_i = \sqrt \e \tilde w_i,$ where $\tilde w_i$ is bounded in $\|\cdot \|_{\e,s}$ norm. As a consequence, the linear convective terms in \eqref{ivp:dotu} takes the form  
 \be \label{for:convection:this:one} \frac{1}{\sqrt \e} \op_\e(J) =  (\tilde v_{e0} \cdot \d_\e + \tilde w_e \cdot \nabla_\e) \op_\e(\chi_{high} \underline J_e) + \sqrt \e (\tilde w_i \cdot \nabla_\e) \op_\e(\chi_{high}\underline J_i).\ee
 In particular, the linear convective terms are {\it not singular in $\e,$} thanks to the cancellation \eqref{cancellation:convection}. The quasilinear convective terms, however, are singular, as seen in the next Section.   

\subsection{The ``source" terms in the perturbation equations} \label{sec:remainder:bounds}
 The term $\dot F$ in \eqref{ivp:dotu} comprises nonlinear terms and the WKB remainder:
\be \label{def:dotF}
\dot F := {\mathcal F}(\tilde u_a + \dot u) - {\mathcal F}(\tilde u_a) - \e^{-1/2} \big( {\mathcal J}(\dot u, \nabla_\e) \dot u + {\mathcal B}(\dot u, \dot u) \big) - \e^{K_a} \tilde R_a.
 \ee
 In \eqref{def:dotF}, the term $\e^{K_a} R_a$ is the remainder in the WKB approximation, which takes the form $\e^{K_a} \tilde R_a$ after the change of transversal spatial variable of Section \ref{sec:y}.

By definition of $\dot F,$ we have the bound 
$$
\|\dot F \|_{L^2} \lesssim C(\|\tilde u_a + \dot u\|_{L^\infty}) \|\dot u\|_{L^2} + \e^{-1/2} \|\nabla_\e \dot u\|_{L^\infty} \|\dot u\|_{L^2} + \e^{K_a},$$
where $C: \R_+ \to \R_+$ is non-decreasing. 
By the short-time bound \eqref{a:priori:tilde:bound} on $\dot u,$ this translates into
\be \label{bd:dotF:L2}
\|\dot F (t) \|_{L^2} \lesssim \e^{-1/2 + K'} \|\dot u (t) \|_{L^2} + \e^{K_a}, \qquad \mbox{for $t \leq t_\star(\e).$}
\ee
In the upcoming Sections, we will prove weighted Sobolev and $\fl1$ bounds for $\dot F.$ 

\subsection{Weighted Sobolev norms} 

Our weighted Sobolev norms are defined in terms of $L^2$ bounds for $\d_\e^\a$ derivatives, with  
\be \label{def:de}
 \d_\e^\a = ((\e \d_x)^{\a_1}, \d_y^{\a_2}), \qquad \a = (\a_1, \a_2) \in \N \times \N^2,
 \ee 
 in accordance with the definition of $\d_\e$ in \eqref{nablae}. The associated weighted-in-$x$-only Sobolev norm is defined as  
 \be \label{weighted:norm}
 \| f \|_{\e,s}^2 := \sum_{|\a| \leq s} \| \d_\e^\a f \|^2_{L^2(\R^3)}.
 \ee
In particular, in \eqref{weighted:norm} we are considering $\d_y$ derivatives, not $\sqrt \e \d_y$ derivatives in the transverse coordinates $y \in \R^2.$ The definition \eqref{def:de} is in contrast with the operator $\nabla_\e = (\e \d_x, \sqrt \e \d_y)$ that appears in particular in the initial-value problem in $\dot u$ \eqref{ivp:dotu}. Note that when dealing with pseudo-differential operators, we keep, however, the anisotropic $\op_\e(\cdot)$ quantization, in which symbols are evaluated at $(\e \xi, \sqrt \e \eta).$ Both operators, $\nabla_\e$ and $\d_\e,$ are relevant in our analysis. 

Note that the $H^{s_0} \hookrightarrow \fl1$ embedding, for $s_0 > 3/2,$ takes the form
\be \label{embed:loss}
 \| f \|_{\fl1} \lesssim \e^{-1/2} \| f \|_{\e,s_0},
 \ee
 in the weighted Sobolev norm \eqref{weighted:norm}.

\subsection{Sobolev bounds for the source} 
By symmetry of the quasilinear convective terms (see Section \ref{sec:conv}), we have 
$$ \Big| \frac{1}{\sqrt \e} \Re e \, \big(  {\mathcal J}(\dot u, \nabla_\e) \d_\e^\a \dot u, \d_\e^\a \dot u \big)_{L^2} \Big| \lesssim  \e^{-1/2} \| \nabla_\e \dot u \|_{L^\infty} \|  \dot u \|_{\e,s}^2, \qquad |\a| \leq s,$$
where the weighted Sobolev norm $\|\cdot\|_{\e,s}$ is defined in the Section just above. 
Besides, by the Gagliardo-Nirenberg inequality, 
$$  \Big| \frac{1}{\sqrt \e} \Re e \, \big( \big[ \d_\e^\a,  {\mathcal J}(\dot u, \nabla_\e)\big] \dot u), \d_\e^\a \dot u \big)_{L^2} \Big| \lesssim \e^{-1/2} \| \nabla_\e \dot u \|_{L^\infty}  \| \dot u \|^2_{\e,s}, \qquad |\a| \leq s.$$ 
Thus we have, for $\dot F$ defined in \eqref{def:dotF}:
\be \label{bd:dotF:Hs:0} \Re e \, ( \dot F, \dot u)_{\e,s} \lesssim \Big( C(\|\tilde u_a + \dot u\|_{L^\infty}) + \e^{-1/2} \| \dot u \|_{L^\infty} + \e^{-1/2} \| \nabla_\e \dot u \|_{L^\infty} \Big) \| \dot u\|_{\e,s}^2 + \e^{K_a} \|\dot u\|_{\e,s},\ee
where $C: \R_+ \to \R_+$ is non-decreasing. 
By the short-time bound \eqref{a:priori:tilde:bound} on $\dot u,$ and the embedding \eqref{Linfty:fl1}, this translates into
 \be \label{bd:dotF:Hs}
\Re e \, ( \dot F(t), \dot u(t))_{\e,s} \lesssim \e^{-1/2 + K'} \|\dot u(t)\|_{\e,s}^2 + \e^{K_a} \|\dot u(t)\|_{\e,s}, \qquad \mbox{for $t \leq t_\star(\e).$}
\ee
We also have
\be \label{bd:F:sob:new}  \begin{aligned} 
 \| \dot F \|_{\e,s} & \lesssim \e^{-1/2}  \| \dot u \|_{L^\infty} \| \nabla_\e \dot u \|_{\e,s} + \e^{-1/2} (\| \dot u \|_{L^\infty} + \| \nabla_\e \dot u \|_{L^\infty}) \| \dot u \|_{\e,s} ) \\ & \quad + C(\| \tilde u_a + \dot u \|_{L^\infty}) \| \dot u \|_{\e,s} +  \e^{K_a},
\end{aligned}\ee
implying, by the short-time bound \eqref{a:priori:tilde:bound} on $\dot u,$ and the embedding \eqref{Linfty:fl1}: 
\be \label{bd:F:sob}
\| \dot F(t) \|_{\e,s} \lesssim \e^{-1/2 + K'} \| \dot u(t) \|_{\e,s + 1} + \e^{K_a}, \qquad \mbox{for $t \leq t_\star(\e).$}
\ee

\subsection{$\fl1$ bounds for the source} We turn to $\fl1$ bounds for the remainder $\dot F$ defined in \eqref{def:dotF}, where the $\fl1$ norm ($L^1$ norm of the Fourier transform) is defined in \eqref{def:fl1}. 
Since $\fl1$ is an algebra, we have
\be \label{bd:dotF:fl1} \| \dot F \|_{\fl1} \lesssim \big( \, C( \| \dot u \|_{\fl1}) + \e^{-1/2} \| \dot u \|_{\fl1} + \e^{-1/2} \| \d_\e \dot u \|_{\fl1}\, \big) \| \dot u \|_{\fl1} + \e^{K_a},\ee
where $C: \R_+ \to \R_+$ is nondecreasing.  
Using the short-time bound \eqref{a:priori:tilde:bound} on $\dot u,$ this implies
\be \label{bd:dotF:fl1:apriori}
 \| \dot F(t) \|_{\fl1} \lesssim \e^{-1/2 + K'} \| \dot u(t) \|_{\fl1} + \e^{K_a}, \qquad t \leq t_\star(\e).
 \ee 
We bound in a similar way $\d_\e \dot F$ in the $\fl1$ norm:
 \ba \| \d_\e \dot F\|_{\fl1} & \lesssim \big( C(\| \dot u \|_{\fl1} + \e^{-1/2} \| \dot u \|_{\fl1} + \e^{-1/2} \| \d_\e \dot u \|_{\fl1} \, \big) \| \d_\e \dot u (t) \|_{\fl1} \\ & \quad + \e^{-1/2} \| \dot u \|_{\fl1} \| \d_\e^2 \dot u \|_{\fl1} + \e^{K_a},\ea
so that
\be \label{bd:dedotF:fl1:apriori}
\| \d_\e \dot F\|_{\fl1} \lesssim \e^{-1/2 + K'} (\| \d_\e \dot u \|_{\fl1}  + \| \d_\e^2 \dot u \|_{\fl1}) + \e^{K_a}, \qquad t \leq t_\star(\e).\ee
Higher derivatives $\d_\e^\a \dot F$ involve 
$$ \d_\e^\a ({\mathcal F}(\tilde u_a + \dot u) - {\mathcal F}(\tilde u_a)) = \int_0^1 \d_\e^\a ( {\mathcal F}'(\tilde u_a + \t \dot u) \dot u) \, d\t,$$
which takes the form
$$ \sum_{\begin{smallmatrix} \a_1 + \a_2 = \a \\ 1 \leq m \leq |\a_1| \\ \b_1 + \dots + \b_m = \a_1 \end{smallmatrix} } C_{\b_j}\int_0^1 {\mathcal F}^{(m)}(\tilde u_a + \t \dot u) \cdot \prod_{1 \leq j \leq m} \d_\e^{\b_j} \dot u \, \d_\e^{\a_2} \dot u\, d\t,$$
with positive coefficients $C_{\b_j}.$ 
We then use the interpolation inequality \eqref{interpol}: 
$$ \| {\mathcal F}^{(m)}(\tilde u_a + \t \dot u) \cdot \prod_{1 \leq j \leq m} \d_\e^{\b_j} \dot u \d_\e^{\a_2} \dot u\|_{\fl1} \leq C(\| \tilde u_a \|_{\fl1}, \| \dot u \|_{\fl1}) \| \dot u \|^m_{\fl1} \max_{|\b| = |\a|} \| \d_\e^\b \dot u \|_{\fl1},$$
for some nonnegative and nondecreasing function $C.$  
With \eqref{a:priori:tilde:bound}, this gives
$$  \| \d_\e^\a ({\mathcal F}(\tilde u_a + \dot u) - {\mathcal F}(\tilde u_a)) \|_{\fl1} \lesssim \e^{K'} \max_{|\b| = |\a|} \| \d_\e^\b \dot u \|_{\fl1}, \qquad t \leq t_\star(\e).$$
The other terms in $\dot F$ are bounded by the same arguments. In particular, for $t \leq t_\star(\e),$ 
$$ \e^{-1/2} \|  \d_\e^\a {\mathcal B}(\dot u, \dot u) \|_{\fl1} \lesssim \e^{-1/2 + K'} \max_{|\b| = |\a|} \| \d_\e^\b \dot u \|_{\fl1}, $$
and
$$ \e^{-1/2} \|  \d_\e^\a {\mathcal J}(\dot u, \nabla_\e) \dot u ) \|_{\fl1} \lesssim \e^{-1/2 + K'}  \max_{|\b| = |\a|} (\| \d_\e^\b \dot u \|_{\fl1} + \| \d_\e^\b \nabla_\e \dot u \|_{\fl1}).$$ 
Summing up, we have for $|\a| \geq 1$ (and $|\a|$ not greater than some large index depending on high regularity of the WKB solution) the bound 
\be \label{final:F:fl1} \begin{aligned}
  \| \d_\e^\a \dot F \|_{\fl1} & \lesssim \e^{K_a} + \e^{-1/2 + K'}  \max_{|\b| = |\a|} (\| \d_\e^\b \dot u \|_{\fl1} + \| \d_\e^\b \nabla_\e \dot u \|_{\fl1}), \qquad t \leq t_\star(\e).\end{aligned}
  \ee

\subsection{Uniform remainders and equivalent sources} \label{sec:equiv:remainders} 

The solution $\dot u$ to \eqref{ivp:dotu} will go through a series of changes of variables in the course of the proof, and so will the ``source" term $\dot F.$ The changes of variables often involve operators which are bounded in weighted Sobolev norms and in $\fl1$ norm.

\begin{defi} \label{defi:uniform:remainder} Any linear operator $R$ which satisfies the bounds 
 $$ \| R f \|_{\e,s} \lesssim |\ln \e|^\star \| f \|_{\e,s}, \qquad \| \d_\e^\a R f \|_{\fl1} \lesssim |\ln \e|^\star \max_{|\b| \leq |\a|} \|  \d_\e^\b f \|_{\fl1},$$
 for all $f \in H^s$ with $\d_\e^\b f \in \fl1$ for $|\b| \leq |\a|,$
 for all $\a$ not greater than some large index depending on the regularity of the WKB solution, is said to be a \emph{uniform remainder} if the implicit constants do not depend on $\e$ nor on $t,$ for $t \leq t_\star(\e).$
 \end{defi}
 
In Definition \ref{defi:uniform:remainder}, we use notation $|\ln \e|^\star:= |\ln \e|^{N_\star}$ for some $N_\star > 0$ which depends on all parameters in the problem, but neither on $\e$ nor on the space-frequency variables. 

In view of Lemma \ref{lem:tensor}, we see that given $q \in \N,$ given $b_a$ a term in the WKB expansion, given $Q \in L^\infty(\R^3_{\xi,\eta}),$ the operator $\op_\e(e^{i q \theta} b_a(x,y) Q(\xi,\eta))$ is a uniform remainder. Indeed, all terms in the WKB expansion have a high Sobolev regularity, hence belong to $W^{s',\infty}$ for large $s',$ by propagation of regularity for the Zakharov equations. (The leading terms in the expansion solve the nonlinear Zakharov system (Z) given in the introduction, and the correctors terms solve linearized Zakharov equations around the leading terms, with source terms depending on previous terms in the expansion and their derivatives. Details are given in \cite{em1}.)

 \begin{defi} \label{def:equiv:remainder} Consider $u_1$ an avatar of the perturbative unknown $\dot u,$ defined through linear changes of variables in terms of $\dot u,$ and which satisfies a system with a source $F_1.$ (As a starting point: $u_1 = \dot u$ and $F_1 = \dot F.$) Consider a linear change of variable $Q_{1 \to 2}$ which is a uniform remainder, and the associated unknown $u_2 = Q_{1 \to 2} u_1.$ We say that a source $F_2$ in the system satisfied by $u_2$ is \emph{equivalent} to $F_1$ if it has the form 
 $$ F_2 = Q_{1 \to 2} F_1 + R u_1,$$
 where $R$ is a uniform remainder. 
 \end{defi}

  Note that the $(\cdot, \cdot)_{\e,s}$ bound \eqref{bd:dotF:Hs} depends on a symmetry assumption that is a priori lost by changes of variables. In Section \ref{sec:high}, we will verify that our series of changes of variables preserves the symmetry to first order.

\subsection{The leading singular terms in the perturbation equations} \label{sec:Bua}

 The approximate solution $u_a$ has the form $$u_a = u_0 + \sqrt \e (u_1 + \dots),$$ and is such that $(\tilde u_a - \tilde u_0)/\sqrt \e$ is bounded in the weighted Sobolev norm $\| \cdot \|_{\e,s}$ and in the $\fl1$ norm.  The leading term $u_0$ has fast space-time oscillations
 $$ \tilde u_0(t,x,y) = e^{- i \theta/\e} \tilde u_{0,-1}(t,x,y) + e^{i \theta/\e} \tilde u_{0,1}(t,x,y),$$
 in the rescaled spatial frame \eqref{def:tildeu}, where
 \be \label{def:theta}
 \theta = (k x - \o t)/\e.
\ee
The fundamental phase $(k,\o) \in \R \times \R$ (or $(k,\o) \in \R^3 \times \R$ if $k$ is identified with $(k,0,0) \in \R^3$) is pictured on Figure \ref{fig:fund-phase}.

 We further know that the amplitudes $u_{0,p}$ are {\it polarized}, in the sense that there are scalar functions $g_{(p)}$ and fixed spatial directions $\vec e_p$ such that 
 \be \label{def:gp} \tilde u_{0,p} = g_{(p)} \vec e_p.\ee
Details (in particular the definition of vectors $\vec e_p$) are found in Section \ref{sec:wkb}. By symmetry, we have
\be \label{gpm1}
 g_{(-1)} = g_{(1)}^\star, \quad \vec e_{-1} = (\vec e_1)^\star,
\ee
where $z^\star$ denotes the complex conjugate of $z,$ and conjugation for vectors is meant entry-wise.

 Next we evaluate $g_{(p)}$ at $t = 0,$ and introduce a truncation that takes advantage of the assumption that the initial WKB amplitude $a$ decays like $e^{- r_a |y|^{\kappa_a}}$ (this is assumption (a3) below \eqref{data}). Let $\chi_{trans}$ be smooth, with values in $[0,1],$ and such that 
 \be \label{def:chi:infty}
 \chi_{trans}(y) = \left\{\begin{aligned} 1,  & \quad |y| \leq R_{trans} |\ln \e|^{1/\kappa_a}, \\ 0,  & \quad |y| \geq 2 R_{trans} |\ln \e|^{1/\kappa_a}
 \end{aligned}\right. \quad \mbox{with $R_{trans} = \left(\frac{1}{2 r_a}\right)^{1/\kappa_a}.$}
 \ee
 In order to achieve \eqref{def:chi:infty}, we can for instance choose a fixed ``plateau'' function \label{chi0} $\chi_0 \in C^\infty_c(\R^2)$ with values in $[0,1],$ such that $\chi_0(y) = 1$ for $|y| \leq 1$ and $\chi_0(y) = 0$ for $|y| \geq 2,$ and let $\chi_{trans}(y) := \chi_0(R_{trans}^{-1} |\ln \e|^{-1/\kappa_a} y).$ We observe then that $\chi_{trans}$ has somewhat slow variations in $y,$ in the sense that
$$ 
 |\d_y^\a \chi_{trans}|_{L^\infty} \lesssim R_{trans}^{-|\a|} |\ln \e|^{-|\a|/\kappa_a}, \qquad |\a| > 0.
$$ 
  The leading term in $B$ involves only $\chi_{trans}(y) g_{(p)}(0,x,y).$ Indeed, we have, for $|\a| \leq s_a - d/2, t \leq T \sqrt \e |\ln \e|,$ and if $\e$ is small enough, the bound 
\be \label{truncation:WKB}
 \big| \d_{x,y}^\a \big( g_{(p)}(t,x,y) - \chi_{trans}(y) g_{(p)}(0,x,y) \big) \big| \lesssim \e^{1/2} |\ln \e|^\star, 
 \ee
 where notation $|\ln \e|^\star$ is introduced in Definition \ref{defi:uniform:remainder}. 
 
 We now verify \eqref{truncation:WKB}. A starting point is the decomposition 
 \be \label{wkb:bound} \begin{aligned} \d_{x,y}^\a \big( g_{(p)}(t,x,y) & - \chi_{trans}(y) g_{(p)}(0,x,y) \big) \\ & = \d_{x,y}^\a (g_{(p)}(t,x,y) - g_{(p)}(0,x,y)) + \d_{x,y}^\a ((1 - \chi_{trans}(y) ) g_{(p)}(0,x,y)). \end{aligned}\ee 
 The first term in the above right-hand side is controlled pointwise by 
 $$ \left| \int_0^t \d_t g_{(p)}(t', x,y) \, dt' \right| \lesssim t \sup_{0 \leq t' \leq t} \| \d_t g_{(p)} \|_{L^\infty}.$$
 As we see in Section \ref{sec:wkb}, the components of the amplitudes $g_{(p)}$ are expressed in terms of the electric field ${\bf E}$ which solves the (Z) (Zakharov) system of Section \ref{sec:ZEM}. By propagation of Sobolev regularity for the Zakharov system (Z) (proved in \cite{LPS}), this means that $g_{(p)}$ has a high Sobolev regularity, uniformly in time over an $O(1)$ time interval. Via the (Z) equation, we can express $\d_t g_{(p)}$ in terms of spatial derivatives and products of ${\bf E}$ (and ${\bf n},$ which is a corrector term in the WKB ansatz, and has a lower, but still large, regularity). In particular, $\d_t g_{(p)}$ is bounded in $(x,y),$ uniformly in time $O(1).$ 
 
 For the other term in the right-hand side of \eqref{wkb:bound}, by definition of $\chi_{trans}:$  on the support of $1 - \chi_{trans},$ we have $|(x,y)| \geq |y| > R_{trans} |\ln \e|^{1/\kappa_a} > R_a$ if $\e$ is small enough. Thus by assumption (a3) on the WKB amplitude:  
 $$ \| \d_{x,y}^\a ((1 - \chi_{trans}(y) ) g_{(p)}(0,x,y)) \|_{L^\infty} \lesssim e^{- r_a |y|^{\kappa_a}} |\ln \e|^{\star} \lesssim \e^{1/2} |\ln \e|^\star,
 $$ 
 since $|y| > R_{trans} |\ln \e|^{1/\kappa_a}.$ Thus \eqref{truncation:WKB} is proved.

 We denote
 \be \label{def:g}
 g_{p}(x,y) := \chi_{trans}(y) g_{(p)}(0,x,y).
 \ee
By \eqref{truncation:WKB}, the leading term in symbol $\dot {\mathcal B}$  defined in \eqref{def:B0} is
$$ B  := \sum_{p \in \{-1,1\}} e^{i p \theta} g_{p}(x, y) \Big( {\mathcal B}(\vec e_p, \cdot) +  {\mathcal B}(\cdot, \vec e_p) - {\mathcal J}_e(\cdot, i p k) \vec e_p - (1 - \chi_{high,+p}) {\mathcal J}_e(\vec e_p,i \xi,i \eta) \, \Big).$$
 In the above definition of $B,$ we used notation pertaining to frequency shifts of symbols introduced in Remark \ref{rem:shift}. Note that among the convective terms ${\mathcal J},$ only ${\mathcal J}_e,$ which is the contribution of the electrons to the convective terms (defined in Section \ref{sec:conv}), remains in $B.$ The perturbation equations \eqref{ivp:dotu} in $\dot u$ take the form
\be \label{pert:2} \left\{\begin{aligned} \d_t \dot u + \frac{1}{\e} A(\nabla_\e) \dot u  + \frac{1}{\sqrt \e} \op_\e(J) \dot u &  = \frac{1}{\sqrt \e} \op_\e(B) \dot u + \dot G,\\  \dot u(0,x,y) & = \e^K \phi^\e,
 \end{aligned}\right.\ee
where
\be \label{def:mathcaldotF} \dot G := \dot F + \frac{1}{\sqrt \e} \op_\e\big(\dot {\mathcal B} - B\big) \dot u.\ee
Note that the second term in $\dot G$ does {\it not} have a singular prefactor, contrary to what it seems at first sight, since ${\mathcal B} = B + \sqrt \e(...).$ The remainder $\dot G$ is equivalent to $\dot F,$ in the sense of Definition \ref{def:equiv:remainder}.

\section{The characteristic variety and resonant frequencies} \label{sec:char}

The characteristic variety of $A(\nabla_\e) = A(\e \d_x, \sqrt \e \d_y)$ is, by definition, the union of branches of eigenvalues of $A(i \xi, i \eta),$ as $(\xi,\eta)$ ranges in $\R \times \R^2.$ According to the definition of $A$ in Section \ref{sec:hypop}, the symbol is
$$ A(i\xi,i \eta)  = \left(\begin{array}{cc|cc|cc} 0 & i(\xi,\eta) \times & 0 & 0 & 0 & 0 \\ - i(\xi,\eta) \times & 0 & - 1 & 0 &  \sqrt \e/\theta_e & 0  \\[2pt] \hline &&&&& \vspace{-.3cm} \\ 0 & 1 & 0 & i\theta_e (\xi,\eta) & 0 & 0 \\ 0 & 0 & i\theta_e (\xi,\eta) \cdot & 0 & 0 & 0 \\[2pt] \hline &&&&& \vspace{-.3cm} \\ 0 & - \sqrt \e/\theta_e & 0 & 0 & 0 &  \a_{ie}^2 \sqrt \e (\xi,\eta)   \\ 0 &  0 & 0 & 0 & \sqrt \e (\xi,\eta) \cdot & 0 \end{array}\right).$$
The vertical and horizontal dividers follow the convention of Section \ref{sec:coord}.

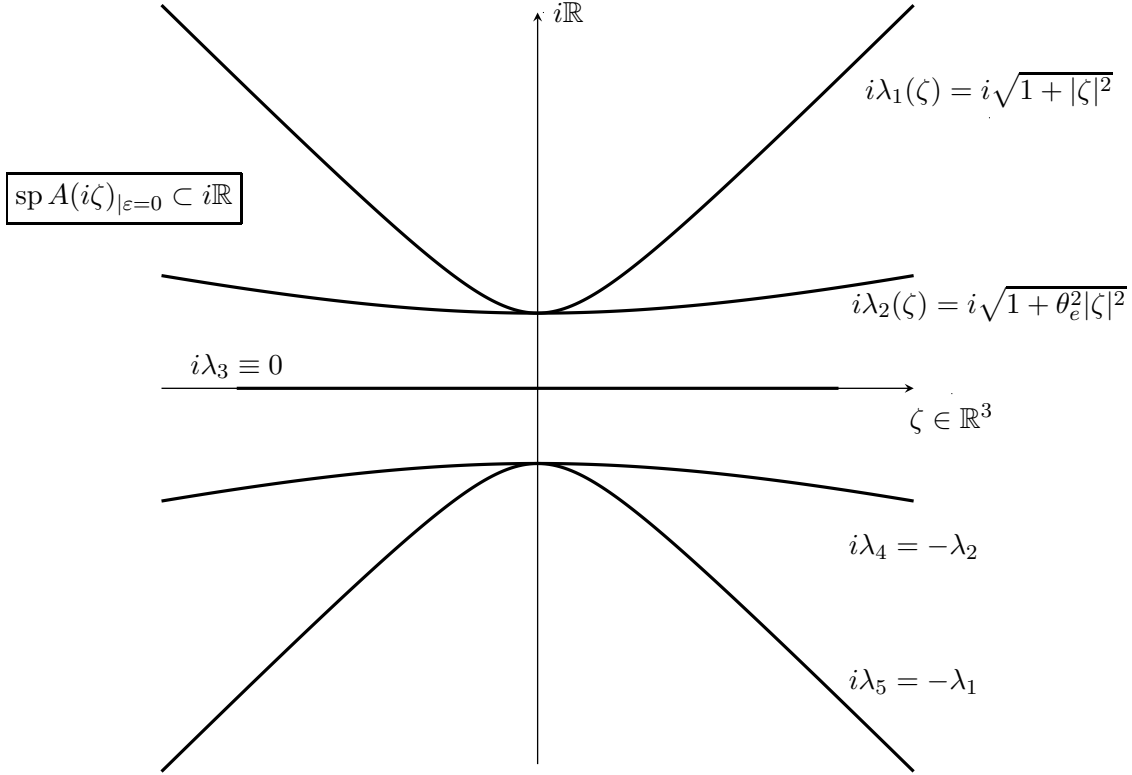
\begin{figure}
\scalebox{1}{
 \begin{tikzpicture}

 \begin{scope}[>=stealth]
  \draw[very thick] (-4,0) -- (4,0);
  \draw[line width=.5pt][->] (-5,0) -- (5,0);

 \draw[line width=.5pt][->] (0,-5) -- (0,5);
 
  \draw (-5.5,3)  node[anchor=north]{$\boxed{\mbox{sp}\,A(i \zeta)_{|\varepsilon = 0} \subset  i {\mathbb R}}$};

  \draw (5.5cm,-2pt) -- (5.5cm,-2pt) node[anchor=north] {$\zeta \in {\mathbb R}^3$};

  \draw (2pt,5cm) -- (2pt,5cm) node[anchor=west] {$i {\mathbb R}$} ;

 \draw (6,3.6) -- (6,3.6) node[anchor=south] {$i \lambda_1(\zeta) = i \sqrt{1 + |\zeta|^2}$} ;

\draw (6,1.5) node[anchor=north] {$i \lambda_2(\zeta) = i \sqrt{1 + \theta_e^2 |\zeta|^2}$} ;

 \draw (-4,0) node[anchor=south] { $i \lambda_{3} \equiv 0$} ;

\draw (5, -2.4) node[anchor=south] {$ i\lambda_4 = -\lambda_2$} ;

\draw (5,-4.2) node[anchor=south] {$i \lambda_5 = -\lambda_1$} ;

 \draw[domain=0:5,samples=100, very thick] plot( {\x}, {(1 + \x^2)^(1/2)} );

 \draw[domain=-5:0,samples=100, very thick] plot( {\x}, {(1 + (-\x)^2)^(1/2)} );
 
  \draw[domain=0:5,samples=100, very thick] plot( {\x}, {(1 + .05*\x^2)^(1/2)} );
 
 \draw[domain=-5:0,samples=100, very thick] plot( {\x}, {(1 + .05*(-\x)^2)^(1/2)} );

  \draw[domain=0:5,samples=100, very thick] plot( {\x}, {-(1 + \x^2)^(1/2)} );
 
 \draw[domain=-5:0,samples=100, very thick] plot( {\x}, {-(1 + (-\x)^2)^(1/2)} );
 
  \draw[domain=0:5,samples=100, very thick] plot( {\x}, {-(1 + .05*\x^2)^(1/2)} );
 
 \draw[domain=-5:0,samples=100, very thick] plot( {\x}, {-(1 + .05*(-\x)^2)^(1/2)} );

  \end{scope}

 \end{tikzpicture}
}
\caption{The linear electronic characteristic variety for the Euler-Maxwell equations. This is the union of branches $i \lambda_j(\zeta)$ of eigenvalues of $A(i \zeta)|_{\epsilon = 0},$ for $\zeta \in \R^3.$ The variety is rotation invariant in $\zeta.$ In the above figure for the sake of representation we let $\theta_e =
\sqrt 5 \cdot 10^{-1}.$ In physical applications, $\theta_e$ is quite smaller. In the physics literature, the fast Klein-Gordon modes are called {\it electromagnetic waves}, while the slow Klein-Gordon modes are called {\it electronic plasma waves}.}
\label{fig1}

\end{figure}

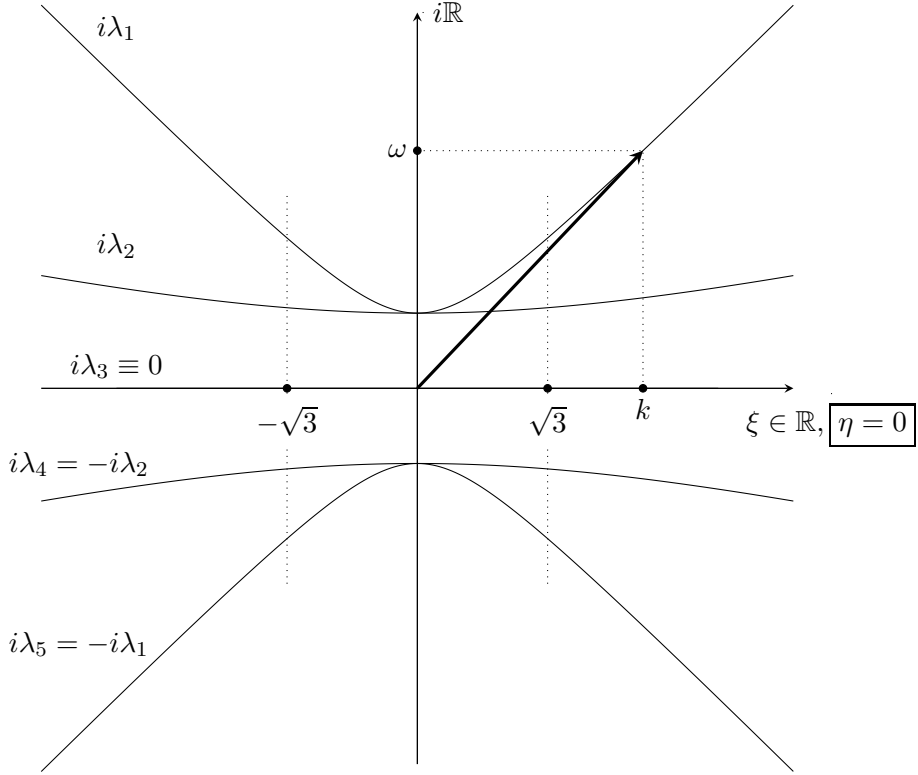
\begin{figure}
\scalebox{1}{
\label{fig:fund-phase}
 \begin{tikzpicture}

 \begin{scope}[>=stealth]
 \draw[line width=.5pt][->] (-5,0) -- (5,0);

 \draw[line width=.5pt][->] (0,-5) -- (0,5);

  \draw (5.5cm,-2pt) -- (5.5cm,-2pt) node[anchor=north] {$\xi \in {\mathbb R}, \boxed{\eta = 0}$};

  \draw (2pt,5cm) -- (2pt,5cm) node[anchor=west] {$i {\mathbb R}$} ;

 \draw (-4,0) -- (4,0);

 \draw[domain=0:5,samples=100] plot( {\x}, {(1 + \x^2)^(1/2)} );

 \draw[dotted] (1.732,-2.6) -- (1.732,2.6) ; 
 
 \draw[dotted] (-1.732,-2.6) -- (-1.732,2.6) ; 
 
 \draw (3,0) node[anchor=north] {$k$} ; 
 
 \draw (0,3.162) node[anchor=east] {$\omega$} ; 
 
 \filldraw (0,3.162) circle (0.05cm) ; 
 
 \filldraw (3,0) circle (0.05cm) ;

 \draw (-4,2.2) node[anchor=north] {$i \lambda_2$} ;
 \draw (-4,0) node[anchor=south] { $i \lambda_{3} \equiv 0$} ;

  \draw (-4.5, -1.3) node[anchor=south] {$ i\lambda_4 = - i \lambda_2$} ;

 \draw (-4.5,-3.7) node[anchor=south] {$i \lambda_5 = - i \lambda_1$} ;

\draw (-4,4.5) node[anchor=south] {$i \lambda_1$} ;

 \draw[domain=-5:0,samples=100] plot( {\x}, {(1 + (-\x)^2)^(1/2)} );
 
  \draw[domain=0:5,samples=100] plot( {\x}, {(1 + .05*\x^2)^(1/2)} );
 
 \draw[domain=-5:0,samples=100] plot( {\x}, {(1 + .05*(-\x)^2)^(1/2)} );

  \draw[domain=0:5,samples=100] plot( {\x}, {-(1 + \x^2)^(1/2)} );
 
 \draw[domain=-5:0,samples=100] plot( {\x}, {-(1 + (-\x)^2)^(1/2)} );
 
  \draw[domain=0:5,samples=100] plot( {\x}, {-(1 + .05*\x^2)^(1/2)} );
 
 \draw[domain=-5:0,samples=100] plot( {\x}, {-(1 + .05*(-\x)^2)^(1/2)} );

 \filldraw (1.732,0) circle (0.05cm) ; 
 
  \filldraw (-1.732,0) circle (0.05cm) ; 
 
 \draw (1.732,-.1) node[anchor=north,fill=white] {$\sqrt 3$} ; 
  
  \draw (-1.732,-.1) node[anchor=north,fill=white] {$-\sqrt 3$} ;

 \draw[dotted] (3,0) -- (3,3.162) ; 
 \draw[dotted] (0,3.162) -- (3,3.162) ; 
 
 \begin{scope}[>=stealth] 
 \draw[domain=0:3,->,very thick] plot( {\x}, { sqrt(10)/3*\x } ); 
 \end{scope}
 
  \end{scope}

 \end{tikzpicture}
}
 \caption{The fundamental phase in the WKB oscillations is $(k,\omega),$ with $|k| > \sqrt 3,$ and $\omega = \lambda_1(k,0) = \sqrt{1 + k^2}.$ It describes the space-time oscillations of the initial (or incident) {\it electromagnetic wave} which interacts nonlinearly with the plasma. We often identify $k \in {\mathbb R}$ with the vector $(k,0,0) \in {\mathbb R}^3.$ Condition $|k| > \sqrt 3$ guarantees the existence of the key unstable $(1,4)$ and $(2,5)$ resonances, see Proposition \ref{prop:res} and Figure \ref{fig:resonances:unstable}.}

\end{figure}

\begin{figure}
\scalebox{1}{
 \begin{tikzpicture}

 \begin{scope}[>=stealth]

 \draw[line width=.5pt][->] (0,-5) -- (0,5);

  \draw (5.5cm,0) -- (5.5cm,0) node[anchor=south] {$\xi \in {\mathbb R}, \boxed{\eta = 0}$};

  \draw (2pt,5cm) -- (2pt,5cm) node[anchor=west] {$i {\mathbb R}$} ;

\draw[domain=0:5,samples=100] plot( {\x}, {(1 + \x^2)^(1/2)} );
 
 \draw[domain=-5:0,samples=100] plot( {\x}, {(1 + (-\x)^2)^(1/2)} );
 
  \draw[domain=0:6,samples=100] plot( {\x}, {(1 + .05*\x^2)^(1/2)} );
 
 \draw[domain=-6:0,samples=100] plot( {\x}, {(1 + .05*(-\x)^2)^(1/2)} );

  \draw[domain=0:5,samples=100] plot( {\x}, {-(1 + \x^2)^(1/2)} );
 
 \draw[domain=-5:0,samples=100] plot( {\x}, {-(1 + (-\x)^2)^(1/2)} );
 
  \draw[domain=0:6,samples=100] plot( {\x}, {-(1 + .05*\x^2)^(1/2)} );
 
 \draw[domain=-6:0,samples=100] plot( {\x}, {-(1 + .05*(-\x)^2)^(1/2)} );

 \begin{scope}[>=stealth]

  \draw[domain=0:3,->,very thick] plot( {\x}, { sqrt(10)/3*\x } ); 
  \draw (2.8, 3.162) node[anchor=east] {$(1,3)$} ; 
   \filldraw (0,0) circle (.05cm) ;  
   \draw (.3,0) node[anchor=north] {\footnotesize $\xi_{13}^+$}; 

  \draw[->,very thick] (-6,0) -- (-3, 3.162) ;
   \draw (-3.2, 3.162) node[anchor=east] {$(1,3)$} ;
   \filldraw (-6,0) circle (.05cm) ;  
   \draw (-6,0) node[anchor=north] {\footnotesize $\xi_{13}^-$}; 
 \draw[domain=0:3,->,very thick] plot( {1.1 + \x}, {sqrt(1 + .05*1.1^2)  + sqrt(10)/3*\x} ) ; 
 \draw (4.1, 4.19) node[anchor=west] {$(1,2)$} ; 
   \draw (1.1,0) node[anchor=north] {\footnotesize $\xi_{12}^+$}; 
 \draw[dotted] (1.1,0) -- (1.1,1.03) ; 
\filldraw (1.1,0) circle (.05cm) ;  
   \draw[domain=0:3,->,very thick] plot( {-4.04 + \x}, {-sqrt(1 + 4.04^2)  + sqrt(10)/3*\x} ) ;
   \draw (-1.04, -.7) node[anchor=east] {$(4,5)$} ; 
\draw[dotted] (-4.04,0) -- (-4.04,-4.17) ; 
   \draw (-4.04,0) node[anchor=north,fill=white] {\footnotesize $\xi_{45}^-$}; 
\filldraw (-4.04,0) circle (.05cm) ;  
 \draw[domain=0:3,->,very thick] plot( {-3 + \x}, {-sqrt(1 + 3^2)  + sqrt(10)/3*\x} ) ;
   \draw (-.3, 0.3) node[anchor=east] {$(3,5)$} ;    
\draw[dotted] (-3,0) -- (-3,-3.162) ; 
   \draw (-3,0) node[anchor=north,fill=white] {\footnotesize $\xi_{35}^-$}; 
  \filldraw (-3,0) circle (.05cm) ;

  \draw[domain=0:3,->,very thick] plot( {3 + \x}, {-sqrt(1 + 3^2)  + sqrt(10)/3*\x} ) ;
   \draw (6, -.5) node[anchor=north] {$(3,5)$} ;  
\draw[dotted] (3,0) -- (3,-3.162) ; 
   \draw (3,0) node[anchor=north,fill=white] {\footnotesize $\xi_{35}^+$}; 
  \filldraw (3,0) circle (.05cm) ;

 \end{scope}
 
  \end{scope}

 \draw[line width=.5pt][->] (-6,0) -- (6.5,0);

 \end{tikzpicture}
}
\caption{Some {\it stable} space-time resonances are pictured here. By Proposition \ref{prop:space-time}, space-time resonances are resonances with $\eta = 0.$ We see above the $(1,3)$ resonances, which in the small $\theta_e$ limit occur at $\xi_{13}^- = - 2 k$ and $\xi_{13}^+ = 0;$ a $(1,2)$ resonance at $\xi_{12}^+ = -k + (\o^2 + 2 \o)^{1/2}$ (the other $(1,2)$ resonance, occurring at $\xi_{12}^- = - k - (\o^2 + 2 \o)^{1/2},$ is outside the scope of the picture), and a $(4,5)$ resonance occurring at $\xi_{45}^- = - (\o^2 + 2 \o)^{1/2}.$ All those frequencies are understood in the small $\theta_e$ limit. The analysis will reveal that $(1,2),$ $(1,3),$ $(3,5)$ and $(4,5)$ resonances are stable. }
\label{fig:resonances:stable}

\end{figure}
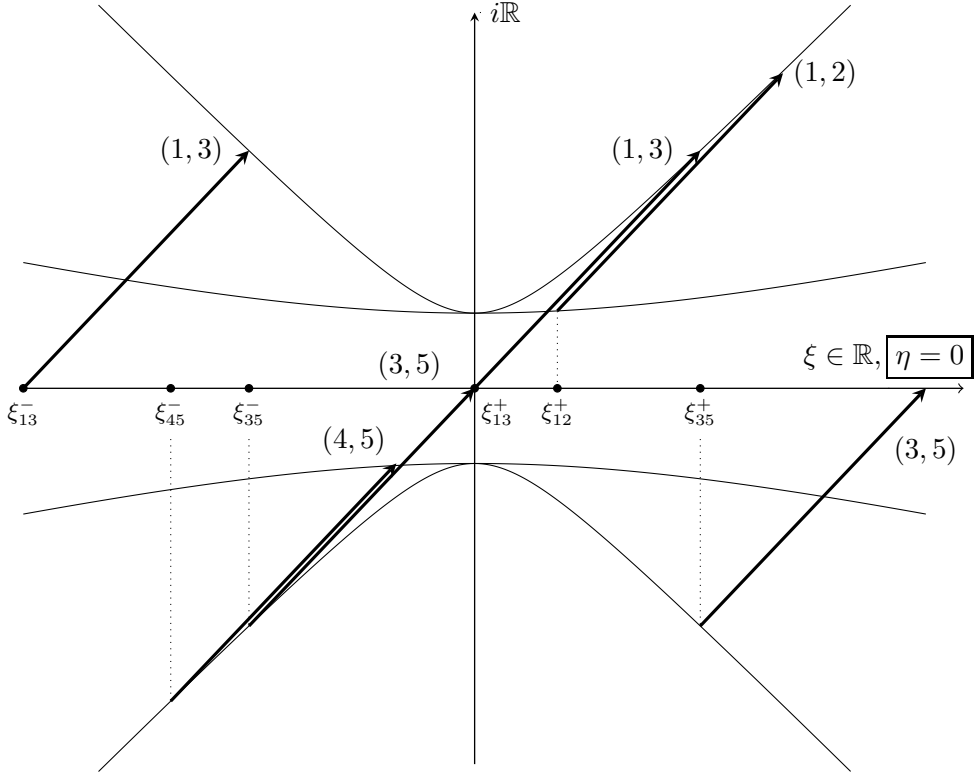

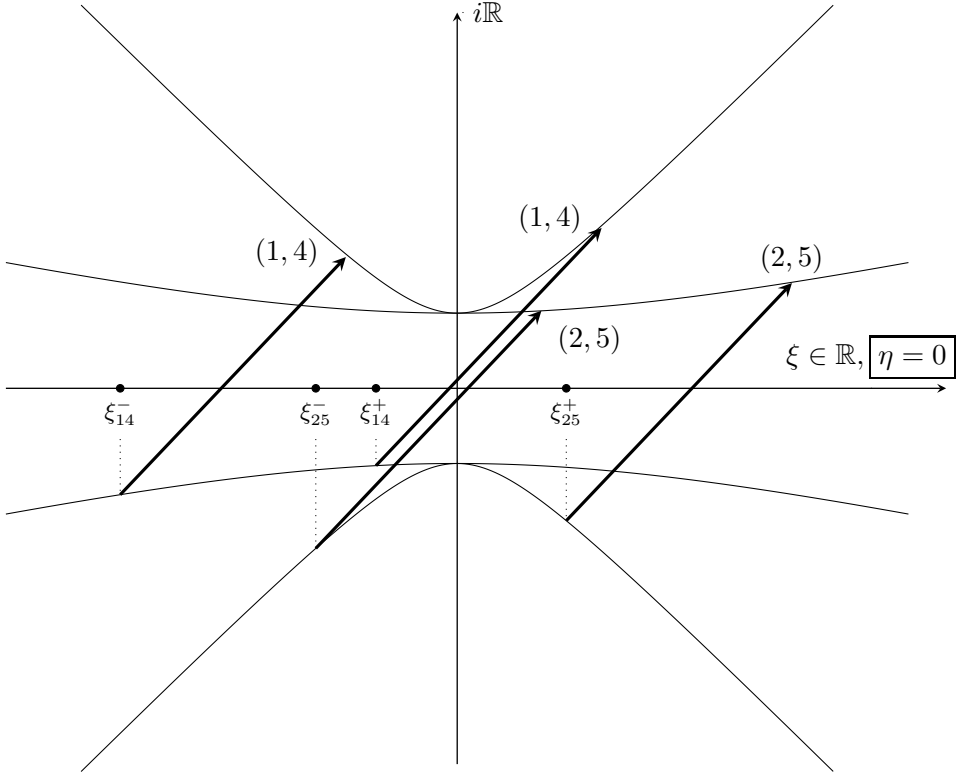
\begin{figure}
\scalebox{1}{
 \begin{tikzpicture}

 \begin{scope}[>=stealth]

 \draw[line width=.5pt][->] (0,-5) -- (0,5);

  \draw (5.5cm,0) -- (5.5cm,0) node[anchor=south] {$\xi \in {\mathbb R}, \boxed{\eta = 0}$};

  \draw (2pt,5cm) -- (2pt,5cm) node[anchor=west] {$i {\mathbb R}$} ;

\draw[domain=0:5,samples=100] plot( {\x}, {(1 + \x^2)^(1/2)} );
 
 \draw[domain=-5:0,samples=100] plot( {\x}, {(1 + (-\x)^2)^(1/2)} );
 
  \draw[domain=0:6,samples=100] plot( {\x}, {(1 + .05*\x^2)^(1/2)} );
 
 \draw[domain=-6:0,samples=100] plot( {\x}, {(1 + .05*(-\x)^2)^(1/2)} );

  \draw[domain=0:5,samples=100] plot( {\x}, {-(1 + \x^2)^(1/2)} );
 
 \draw[domain=-5:0,samples=100] plot( {\x}, {-(1 + (-\x)^2)^(1/2)} );
 
  \draw[domain=0:6,samples=100] plot( {\x}, {-(1 + .05*\x^2)^(1/2)} );
 
 \draw[domain=-6:0,samples=100] plot( {\x}, {-(1 + .05*(-\x)^2)^(1/2)} );

 \begin{scope}[>=stealth]

  (1,4) + ***
 \draw[domain=0:3,->,very thick] plot( {-1.08 + \x}, {-sqrt(1 + .05*1.08^2)  + sqrt(10)/3*\x} ) ; 
 \draw[dotted] (-1.08,0) -- (-1.08, -1.03) ;
 \draw (-1.08,0) node[anchor=north,fill=white] {\footnotesize $\xi_{14}^+$} ; 
 \filldraw (-1.08,0) circle (.05cm) ; 
 \draw (1.8,2.25) node[anchor=east] {$(1,4)$};

 \draw[domain=0:3,->,very thick] plot( {-4.48 + \x}, {-sqrt(1 + .05*4.48^2)  + sqrt(10)/3*\x} ) ; 
  \draw[dotted] (-4.48,0) -- (-4.48, -1.415);
  \draw (-4.48,0) node[anchor=north,fill=white] {\footnotesize $\xi_{14}^-$} ; 
   \filldraw (-4.48,0) circle (.05cm) ; 
  \draw(-1.7,1.8) node[anchor=east] {$(1,4)$} ;

  \draw[domain=0:3,->,very thick,fill=white] plot( {-1.88 + \x}, {-sqrt(1 + 1.88^2)  + sqrt(10)/3*\x} ) ;
  \draw[dotted] (-1.88,0) -- (-1.88, -2.13);
   \draw (-1.88,0) node[anchor=north,fill=white] {\footnotesize $\xi_{25}^-$} ; 
   \filldraw (-1.88,0) circle (.05cm) ; 
 \draw (1.2, .65) node[anchor=west] {$(2,5)$} ; 

  \draw[domain=0:3,->,very thick] plot( {1.45 + \x}, {-sqrt(1 + 1.45^2)  + sqrt(10)/3*\x} ) ;
  \draw[dotted] (1.45,0) -- (1.45, -1.76);
   \draw (1.45,0) node[anchor=north,fill=white] {\footnotesize $\xi_{25}^+$} ; 
   \filldraw (1.45,0) circle (.05cm) ; 
 \draw(4.45, 1.4) node[anchor=south] {$(2,5$)};

    \draw[line width=.5pt][->] (-6,0) -- (6.5,0);

 \end{scope}
 
  \end{scope}

 \end{tikzpicture}
}
\caption{The unstable space-time resonances associated with pairs $(1,4)$ and $(2,5)$ are shown on this figure. In the small $\theta_e$ limit, we have $\xi_{14}^\pm = -k \pm (\o^2 - 2 \o)^{1/2},$ and $\xi_{25}^\pm = \pm (\o^2 - 2 \o)^{1/2}.$ The analysis will show that these space-time resonances, between electromagnetic waves (fast Klein-Gordon) and electronic plasma waves (slow Klein-Gordon), are responsible for the Raman amplification. } 
\label{fig:resonances:unstable}

\end{figure}

The variety for $A|_{\e = 0 }$ is pictured on Figure \ref{fig1}. It shows two Klein-Gordon modes and acoustic modes. The fast Klein-Gordon modes are called {\it electromagnetic waves}, and the slow Klein-Gordon modes are called {\it electronic plasma waves} in the physics literature. In a different context, the decomposition of (EM) into Klein-Gordon with different speeds and acoustic subsystems is central in \cite{GIP}.

\subsection{Eigenmodes and projectors} \label{sec:char0}

Details of the computations for the eigenmodes and eigenvectors of the symbol $A$ of the hyperbolic operator are given in Section \ref{sec:spectral:dec}. We sum up some of the results of Section \ref{sec:spectral:dec} here. 
The matrices $A(i \xi, i\eta)$ have eigenvalues
$i \l^\e_j(\xi, \eta),$ with $j \in \{1,2,3^-,3,3^+,4,5\}.$ By symmetry, $\l_j^\e \in \R.$ There are two fast {\it transverse} Klein-Gordon branches
$$ \l^\e_1 = -\l_5^\e = \big( 1 + \e \theta_e^{-2} + \xi^2 + |\eta|^2 \big)^{1/2},$$
with multiplicity two. There are two slow {\it longitudinal} Klein-Gordon branches 
$$ \l^\e_2 = - \l_4^\e = \big( 1 + \theta_e^2  |\zeta|^2 \big)^{1/2} + O(\e),$$
uniformly in frequency. The slow Klein-Gordon modes have multiplicity one. The eigenvalues $\l_2^\e$ and $\l_4^\e$ have multiplicity one. The other modes are acoustic:
$$ \l^\e_{3^\pm} = O(\sqrt \e |\zeta|), \qquad \l^\e_{3} \equiv 0,$$
uniformly in frequency.
The eigenvalues $\l^\e_{3^\pm}$ are simple, and $\l^\e_3$ has multiplicity 6.
The branches $\l^\e_1$ and $\l^\e_2$ (similarly, $\l^\e_4$ and $\l^\e_5$) intersect at $(0,0) \in \R^3,$ with no loss of smoothness. The acoustic modes also coalesce at the zero frequency, with no loss of smoothness. Associated with $\l^\e_j,$ denote $\Pi_j^\e,$ with $j \in \{1,2,3^-,3,3^+,4,5\},$ the eigenprojector of $A$ onto the eigenspace associated with $\l^\e_j,$ and parallel to the direct sum of the other eigenspaces.

From the above, we deduce the approximate spectral decomposition
\be \label{spectral:dec:symbols}
 A(i \zeta) = \sum_{1 \leq j \leq 5} i \l_j(\zeta) \Pi_j(\zeta) + \sqrt\e R_A(\zeta),\ee 
and the decomposition of identity
\be \label{id:proj}
 \qquad {\rm Id} = \sum_{1 \leq j \leq 5} \Pi_j, \qquad \Pi_j \Pi_{j'} = \delta_{jj'} \Pi_j,
 \ee
where
\begin{itemize}
\item the $\l_j$ and $\Pi_j$ are the eigenvalues and eigenprojectors of $A_{|\e = 0},$ so that
$$ \l_1 = -\l_5 = (1 + \xi^2 + |\eta|^2)^{1/2}, \quad \l_3 \equiv 0, \quad \l_2 = - \l_4 = (1 + \theta_e^2(\xi^2 + |\zeta|^2))^{1/2},$$
and, in particular
$$ \Pi_3 = \big(\Pi^\e_{3^-} + \Pi^\e_{3} + \Pi^\e_{3^+}\big)_{|\e = 0}.$$
The approximate eigenvalues $\l_j$ depend smoothly on $\zeta.$ By inspection (see Section \ref{sec:spectral:dec}), the projectors $\Pi_j,$ for $j \neq 3,$ also depend smoothly on $\zeta.$ By \eqref{id:proj}, this implies that $\Pi_3$ depends smoothly on $\zeta.$ 
\item The remainder $R_A$ is order one, since the acoustic longitudinal modes $\l^\e_3$ are order one in frequency. It is bounded in $\e$ (in particular by the uniform bound \eqref{parallel:modes:3}), and smooth in $\zeta$ by \eqref{spectral:dec:symbols} (since all other symbols in \eqref{spectral:dec:symbols} are smooth).
\end{itemize}
From \eqref{spectral:dec:symbols}, we deduce 
\be \label{approximate:spectral:decomposition}
 A(\nabla_\e) = \sum_{1 \leq j \leq 5} \op_\e(i \l_j) \op_\e(\Pi_j) + \sqrt \e \op_\e( R_A),
\ee
since all the symbols in \eqref{spectral:dec:symbols} are Fourier multipliers, that is independent of the spatial variables. Note that the remainder $\op_\e(R_A)$ is {\it not} a uniform remainder in the sense of Definition \ref{defi:uniform:remainder}, since it is order one. 

\subsection{Resonances} \label{sec:resonances}

The resonant frequencies first appear in our analysis as small divisors in homological equations (see \eqref{homo} in the proof of Proposition \ref{prop:normal:form}).

\begin{defi} \label{def:resonances} The resonant frequencies are defined as the solutions $(\xi,\eta) \in\R^3$ to
$$
 \o = \l_j(\xi + k,\eta) - \l_{j'}(\xi,\eta),
$$
 for $1 \leq j,j' \leq 5,$ where the $\l_j$ are the eigenmodes described in Section {\rm \ref{sec:char0}} evaluated at $\e = 0.$ In \eqref{def:resonances}, $k$ is the fundamental wavenumber that appears in the datum \eqref{data}, and $\o = (1 + k^2)^{1/2}$ is the associated characteristic time frequency.
 If there exists $(\xi,\eta)$ such that \eqref{def:resonances} holds, we say that a $(j,j')$ resonance occurs at $(\xi,\eta),$ and that $(j,j')$ is a resonant pair. We denote ${\mathcal R}_{jj'}$ the set of $(j,j')$ resonant frequencies, and ${\mathcal R}$ the union of all ${\mathcal R}_{jj'}.$
 \end{defi}

\subsection{Structure of the resonant set} \label{sec:resonant:set} Based on the explicit description of the eigenvalues given in Section \ref{sec:char} (see also Figure 1), we make a number of important observations for the resonant set ${\mathcal R}$ introduced in Definition \ref{def:resonances}:

\begin{prop}[Resonant set] ${}^{}$\label{prop:res} 
\begin{enumerate}
\item The phase functions
\be \label{def:phase:function} \Phi_{jj'}: \quad (\xi,\eta) \to  \l_{j}(\xi + k, \eta) - \l_{j'}(\xi,\eta) - \o\ee
are all bounded away from zero in a neighborhood of infinity in $\R^3.$ In particular, the resonant set ${\mathcal R}$ is bounded. Moreover, unless $j = j',$ we have for some $c_{jj'} > 0$ and some $R_{jj'} > 0,$ the lower bound $|\Phi_{jj'}(\xi,\eta)| \geq c_{jj'} |(\xi,\eta)|,$ for all $(\xi,\eta)$ with $|(\xi,\eta)| \geq R_{jj'}.$ Finally, if the phase function $\Phi_{jj'}$ does not have any zero, then it is bounded away from 0 uniformly in $(\xi,\eta).$ 
\item There may occur a $(j,j')$ resonance only if $j < j'.$ 
\item There are no $(1,5)$ resonances. 
\item Given $|k| > k_c = \sqrt 3 + O(\theta_e^2),$ where $k_c$ is defined in \eqref{def:new:kc} in the proof below, then all possible resonances do occur, that is we observe resonances of the type $(1,2),$ $(1,3),$ $(1,4),$ $(2,3),$ $(2,4),$ $(2,5),$ $(3,4),$ $(3,5)$ and $(4,5).$\end{enumerate}
\end{prop}

\begin{proof} (1) We may verify the first two claims by straightforward computations based on the explicit expressions of the eigenvalues $\l_j$ of $A_{|\e = 0},$ given in Section \ref{sec:char0}. See also Figure \ref{fig1}. The third claim, about non-existence of near-resonances, follows from the proof of (2) and (3) below.  

\medskip

(2) Consider first the possibility of $(j,j)$ resonances, sometimes called auto-resonances. The $(j,j)$ phase function is 
$$ \Phi_{jj}(\xi,\eta) = \l_j(\xi + k, \eta) - \l_j(\xi,\eta) - \o,$$
so that
$$ |\Phi_{jj}(\xi,\eta)| \geq |\o| - |\nabla \l_j(\xi,\eta)| |k|.$$
Based on the explicit expressions of the $\l_j,$ we observe that $|\nabla \l_j(\xi,\eta)| < 1$ (euclidian norm), for all $j$ and all $(\xi,\eta).$ Thus
$$ |\Phi_{jj}(\xi,\eta)| \geq |\o| - |k| = (1 + |k|^2)^{1/2} - |k|$$
is uniformly bounded away from zero. This proves that auto-resonances do not occur.

The $(2,1)$ phase function is %
$$ \Phi_{21}(\xi,\eta) = \l_2(\xi + k, \eta) - \l_1(\xi,\eta) - \o = \l_2(\xi,\eta) - \l_1(\xi,\eta) + (\l_2(\xi + k, \eta) - \l_2(\xi,\eta) - \o).$$
In the above right-hand side, the first term $\l_2 - \l_1$ is at most equal to 0. The second term is strictly smaller than zero, by 
$$ |\l_2(\xi + k,\eta) - \l_2(\xi,\eta)| \leq |\nabla \l_2| |k| < |k| < \o.$$
Thus there are no $(2,1)$ resonances. By symmetry, the same argument shows that there are no $(5,4)$ resonances either. 

By a simple consideration of sign, we can rule out the other resonances of the type $(j,j')$ with $j > j',$ that is $(3,1),$ $(3,2),$ $(4,1),$ $(4,2),$ $(4,3)$ and $(5,1),$ $(5,2).$

\medskip

(3) We prove the inequality
 \be \label{15} (1 + (\xi + k)^2)^{1/2} + (1 + \xi^2)^{1/2} > \o, \quad \mbox{for all $\xi \in \R.$} \ee
 It is then easy to put in $\eta$ and deduce from \eqref{15} that there are no $(1,5)$ resonances. For a proof of \eqref{15}, first note that if $|k| \leq |\xi|,$ then the inequality is obvious, since $\o = (1 + k^2)^{1/2}.$ Also, if $k$ and $\xi$ are the same sign, the inequality is obvious. Hence we focus on the case $|\xi| < |k|$ and $\xi k < 0.$ Then, inequality \eqref{15} is implied by 
$$ %
 1 + (\xi + k)^2  > \o^2 + (1 + \xi^2) - 2 \o (1 + \xi^2)^{1/2},
$$ %
 which is equivalent to 
 \be \label{for:15}
 \xi k > \frac{1}{2} - (1 + k^2)^{1/2} (1 + \xi^2)^{1/2}.
\ee 
Recall, $\xi k < 0.$ Thus \eqref{for:15} is equivalent to 
$$ (1 + k^2)^{1/2} (1 + \xi^2)^{1/2} > \frac{1}{2} + |\xi| |k|,$$
which holds true.

\medskip

(4) There are $10$ pairs of the form $(j,j')$ with $1 \leq j < j' \leq 5.$ The $(1,5)$ resonance is ruled out. Thus we have to verify that all $9$ possible resonances do occur under the stated condition on $|k|.$ 

For $(1,2),$ we simply observe that $\Phi_{12}(0,0) = - 1 < 0,$ while $\Phi_{12} \to +\infty$ in the large-frequency limit. By continuity, this proves the existence of $(1,2)$ resonances.

For $(1,3),$ the zero frequency is resonant: $\Phi_{13}(0,0) = 0.$

We turn to $(1,4).$ The associated phase function is 
$$ \Phi_{14}(\xi,\eta) = \l_1(\xi + k,\eta) - \o + \l_2(\xi,\eta).$$
We observe that $\Phi_{14}(0,0) = 1,$ and $\Phi_{14}(-k,0) = 1 - \o + (1 + \theta_e^2 k^2)^{1/2}.$ Thus by continuity, the condition
\be \label{for:14}
 1 + (1 + \theta_e^2 k^2)^{1/2} < (1 + k^2)^{1/2}
\ee
guarantees the existence of $(1,4)$ resonances. Elementary computations show that inequality \eqref{for:14} is equivalent to
 $$ k^4 (1 - \theta_e^2)^2 -2 k^2(1 + \theta_e^2) - 3 > 0.$$  
For $|k| > k_c$ defined by
\be \label{def:new:kc}
k_c :=  \frac{1}{\sqrt 2 \cdot (1 - \theta_e^2)} \Big( 2( 1 + \theta_e^2) + 4 \big( 1 - \theta_e^2 + \theta_e^4 \, \big)^{1/2} \, \Big)^{1/2} = \sqrt 3 + O(\theta_e^2),
\ee the above inequality holds true. 

We have $\Phi_{23}(0,0) = (1 + \theta_e^2 k^2)^{1/2} - \o < 0,$ since $\theta_e < 1.$ Besides, $\Phi_{23}$ is large in the large-frequency limit. Hence the existence of $(2,3)$ resonances by continuity. 

Consider next the $(2,4)$ phase function. We have $\Phi_{24}(0,0) = (1 + \theta_e k^2)^{1/2} - \o + 1,$ which under the condition $|k| > k_c$ is strictly negative. In the large-frequency limit, $\Phi_{24}$ tends to infinity. This proves the existence of $(2,4)$ resonances by continuity. 

The same argument works for the $(2,5)$ resonance as well, since $\Phi_{25}(0,0) = \Phi_{24}(0,0),$ and $\Phi_{25}$ tends to infinity for large $|(\xi,\eta)|.$ 

Resonances of type $(3,4)$ are deduced from resonances of type $(2,3)$ by translation. 

Similarly, resonances of type $(3,5)$ are deduced from resonances of type $(1,3)$ by translation.

Finally, we observe that $\Phi_{45}(0,0) = - (1 + \theta_e^2 k^2)^{1/2} - \o + 1 < 0,$ and $\Phi_{45} \to +\infty$ in the large-frequency limit, implying the existence of $(4,5)$ resonances by continuity.  \end{proof}

\begin{rem} \label{rem:k<3/4} If $|k| < 3/4,$ then $(1,4),$ $(2,4)$ and $(2,5)$ resonances do not occur. Indeed, we have
$$ \begin{aligned} \Phi_{14}(\xi,\eta) & = \l_1(\xi,\eta) + \l_2(\xi,\eta) + (\l_1(\xi + k, \eta) - \l_1(\xi,\eta)) - \o \\ & > 2 - |k| - \o,
\end{aligned}$$
since $|\nabla \l_j| < 1,$ and condition $2 - |k| - \o > 0$ translates into $|k| < 3/4.$ The same holds for $(2,4)$ and $(2,5).$ 
\end{rem}

We now state and prove a proposition which describes {\it separation} properties of the resonant sets. These will be used much further along in the analysis, in Section \ref{sec:6}. There Proposition \ref{prop:separation} will imply the existence of coordinates for the {\it in} subsystem that are associated with a {\it non-oscillating interaction matrix}. This removal of the fast oscillations in the leading source term will turn out to be a crucial step in the analysis.

\begin{prop}[Resonance separation] \label{prop:separation} We observe the following separations properties for the resonant set, if $|k| > 1$ and $\theta_e$ is small enough:
\be \label{for:newsep1}
 {\mathcal R}_{12} \cap \Big( \, \big( {\mathcal R}_{23} + k \big) \cup \big({\mathcal R}_{24} + k \big) \, \Big) = \emptyset,
 \ee 
\be \label{for:newsep2}
 {\mathcal R}_{23} \cap {\mathcal R}_{24} = \emptyset,
 \ee
 \be \label{for:newsep3}
 {\mathcal R}_{25} \cap \big( {\mathcal R}_{23} \cup {\mathcal R}_{24}\big) = \emptyset,
 \ee
 \be \label{for:newsep4}
  \big ({\mathcal R}_{14} \cup {\mathcal R}_{24}\big) \cap {\mathcal R}_{34} = \emptyset,
  \ee
  \be \label{for:newsep5}
  {\mathcal R}_{34} \cap {\mathcal R}_{45} + k = \emptyset.
  \ee
\end{prop}

The proof below gives more details on the conditions bearing on $|k|$ and $\theta_e.$ 

\begin{proof} {\it First separation in \eqref{for:newsep1}:} Let $\zeta = (\xi,\eta)$ belong to ${\mathcal R}_{12},$ and also to ${\mathcal R}_{23} + k.$ We denote $\zeta + k$ for $(\xi + k,\eta).$ This means that $\zeta$ solves
 $$
  (1 + |\zeta + k |^2)^{1/2}  = \o + (1 + \theta_e^2  |\zeta|^2)^{1/2}, \quad \mbox{and} \quad (1 + \theta_e^2 |\zeta|^2)^{1/2}  = \o.
 $$
 Thus, denoting $S(X,R)$ the sphere with center $X$ and radius $R$ in $\R^3,$ we find that $\zeta$ belongs to the intersection
 $$ S\big(-k , \, (3 + 4k^2)^{1/2} \big) \cap S\big( 0, \, |k|/\theta_e\big)$$
 (above $k$ denotes the vector $(k,0) \in \R \times \R^2;$ we often use this slight abuse of notation) which is empty if $|k| > k_{min},$ with $k_{min}$ defined by 
 \be \label{def:kmin}
 k_{min} := \sqrt 3 \cdot \big(\theta_e^{-1} - 1)^2 - 4\big)^{-1/2} \simeq \sqrt 3 \cdot \theta_e, \quad \mbox{for small $\theta_e.$}
\ee
 Indeed, an element in $S(-k, (3 + 4 k^2)^{1/2})$ has norm at most $|k| + (3 + 4 k^2)^{1/2},$ and an element in $S(0,|k|/\theta_e)$ has norm exactly $|k|/\theta_e;$ under condition $|k| > k_{min},$ we have $|k| /\theta_e > |k| + (3 + 4 k^2)^{1/2}.$
 
 \medskip
 
 {\it Second separation in \eqref{for:newsep1}:} 
If $\zeta$ belongs to both ${\mathcal R}_{12}$ and ${\mathcal R}_{24} + k,$ then
$$ \left\{\begin{aligned}
 (1 + (\zeta + k)^2)^{1/2} & = 2 \o - (1 + \theta_e^2 (\zeta - k)^2)^{1/2}, \\
 (1 + \theta_e^2 \zeta^2)^{1/2} & = \o - (1 + \theta_e^2(\zeta - k)^2)^{1/2}.
 \end{aligned}\right.$$
From the first equation we derive
$$ |\zeta + k|^2 < 4 (\o^2 - \o), \qquad \mbox{implying} \quad |\zeta| < |k| + 2 \sqrt{\o^2 - \o}.$$
From the second equation we derive
$$ \o \leq 2 + \theta_e (2 |\zeta| + |k|), \quad \mbox{implying} \quad \frac{\o -   2 - 2 \theta_e|k|}{\theta_e} \leq |\zeta|.$$
Thus the intersection is empty as soon as 
\be \label{new:kmin}
(1 + k^2)^{1/2} > 2 + \theta_e \big( 3 |k| + 2 \sqrt{\o^2 - \o}\big).
\ee
which is a condition bearing on $\theta_e$ and $|k|,$ that can be reformulated as
\be \label{new:cond:k:theta} \sqrt{1 + k^2} \geq 2 + \theta_e(3 |k| + \sqrt{3 + 4 k^2}).\ee 
In the small $\theta_e$ limit, it suffices to have $|k| > 1$ in order to \eqref{new:cond:k:theta} to be satisfied. 

\medskip

{\it Verification of \eqref{for:newsep2} and \eqref{for:newsep3}:} those are straightforward. We have ${\mathcal R}_{23} \cap {\mathcal R}_{2j} = \emptyset,$ with $j = 4$ or $j = 5,$ simply because $|\l_4| \geq 1$ and $|\l_5| \geq 1.$ Besides, if $\zeta \in {\mathcal R}_{25} \cap {\mathcal R}_{24},$ then $\zeta = 0$ (otherwise $\l_4 \neq \l_5$), but then $\o \neq (1 + \theta_e^2 k^2)^{1/2}$ for $\theta_e < 1.$

\medskip

{\it Verification of \eqref{for:newsep4}:} this is straightforward, like the above.

\medskip

 {\it Verification of \eqref{for:newsep5}:} it goes just like the first one in this proof. We find that if $\zeta$ belongs to $\big( {\mathcal R}_{45}  + k \big) \cap {\mathcal R}_{34},$ then it must belong to $S\big((k,0), \, (3 + 4k^2)^{1/2} \big) \cap S\big( 0, \, |k|/\theta_e\big),$ which is empty if $|k| > k_{min}.$

\end{proof}

\subsection{Space-time resonances} \label{sec:space-time} 

We let ${\mathcal S}_{jj'}$ be the frequency set of {\it $(j,j')$-space-time resonances}, defined as 
 \be \label{def:space-time} {\mathcal S}_{jj'} = \big\{ (\xi, \eta) \in \R^3, \quad (\xi,\eta) \in {\mathcal R}_{jj'} \quad \mbox{and} \quad \d_\eta \l_j(\xi + k,\eta) = \d_\eta \l_{j'}(\xi + k,\eta) \big\}.\ee

\begin{prop}[Localization of space-time resonances] \label{prop:space-time} Under condition
\be \label{another:cond:k}
 (1 + k^2)^{1/2} \leq (1 - \theta_e^2)/\theta_e^2,
 \ee
 which amounts roughly to $|k| \leq \theta_e^{-2},$ we have for all resonant pair $(j,j')$ the equality 
 $$ {\mathcal S}_{jj'} = {\mathcal R}_{jj'} \cap \{ \eta = 0 \}.$$
\end{prop}

\begin{proof} Given $\zeta = (\xi,\eta) \in {\mathcal S}_{12},$ we have 
$$ \left\{\begin{aligned} \l_1(\xi + k,\eta) & = \o + \l_2(\xi,\eta), \\ \d_\eta \l_1(\xi + k,\eta) = \frac{\eta}{\l_1(\xi + k, \eta)} & = \d_\eta \l_2(\xi,\eta) = \frac{\theta_e^2 \eta}{\l_2(\xi,\eta)},
\end{aligned}\right.$$
or, equivalently,
$$ \left\{\begin{aligned} \l_1(\xi + k,\eta) & = \o + \l_2(\xi,\eta), \\ \eta\big( (1 - \theta_e^2) \l_2(\xi,\eta) - \theta_e^2 \o \big) & = 0 .
\end{aligned}\right.$$
If $\eta \neq 0,$ then the second equation above implies 
$$ (1 + \theta_e^2 |\zeta|^2)^{1/2} = \frac{\theta_e^2}{1 - \theta_e^2} \o,$$
which is impossible under condition \eqref{another:cond:k}. Thus under that condition we have ${\mathcal S}_{12} = {\mathcal R}_{12} \cap \{ \eta = 0 \}.$  

By symmetry, the same arguments imply that ${\mathcal S}_{45} = {\mathcal R}_{45} \cap \{ \eta = 0 \}.$ 

For the space-time resonances ${\mathcal S}_{3j}$ or ${\mathcal S}_{j3},$ involving $\l_3 \equiv 0$ and a Klein-Gordon mode $j,$ this is even more straightforward: the equality of $\eta$ derivatives implies that $\d_\eta \l_j = 0$ at the space-time resonance. But the Klein-Gordon modes have non-zero $\eta$ derivatives unless $\eta = 0.$

It only remains to prove the statement for $(j,j') \in \{ (1,4), (2,4), (2,5)\}.$ Here a sign argument directly implies the result. For instance, in the case of $(j,j') = (1,4),$ we have  
$$ \d_\eta \l_1(\zeta + k) = \frac{\eta}{\l_1(\zeta + k)} = \d_\eta \l_{4}(\zeta) = \frac{\theta_e^2 \eta}{- \l_2(\zeta)},$$
implying $\eta = 0$ since $\l_1 > 0$ and $\l_2 > 0.$ The same is true for $(2,4)$ and $(2,5).$ 

\end{proof}

\subsection{Further separation properties for the resonances and space-time resonances} \label{sec:further:separation}

Next we need a further separation result for (space-time) resonances. This is best understood in the context of Section \ref{sec:6}. Remark \ref{rem:about:matrix:B} and Figure \ref{fig:tree} explain how the following separation properties will be useful later on in the analysis. 
 
\begin{prop}[Further separation properties] \label{prop:space-time:sep} For any $k_{max} > \sqrt 3,$ if $\theta_e$ is small enough, for any given $k \in (\sqrt 3, k_{max}),$ we have 
\be \label{st:sep:1} 
 {\mathcal R}_{12} \cap {\mathcal R}_{13} = \emptyset, \quad {\mathcal R}_{12} \cap {\mathcal R}_{14} = \emptyset, \quad {\mathcal R}_{13} \cap {\mathcal R}_{14} = \emptyset,
\ee
\be \label{st:sep:2}
 {\mathcal R}_{12} \cap {\mathcal R}_{25} + k = \emptyset,
\ee
\be \label{st:sep:3}
 {\mathcal R}_{23} \cap {\mathcal R}_{24} = \emptyset,
\ee
\be \label{st:sep:4} \left\{ 
 \begin{aligned} {\mathcal R}_{13} \cap {\mathcal R}_{23} = \emptyset, & \quad {\mathcal R}_{13} \cap {\mathcal R}_{34} + k = \emptyset, \\
  {\mathcal R}_{13} \cap {\mathcal R}_{35} + k & = {\mathcal S}_{13} \cap {\mathcal S}_{35} + k = \{ 0 \}, \\ 
 {\mathcal R}_{23} \cap {\mathcal R}_{34} + k \neq \emptyset, & \quad {\mathcal S}_{23} \cap {\mathcal S}_{34} + k  = \emptyset, \\ {\mathcal R}_{23} \cap {\mathcal R}_{35} + k = \emptyset, & \quad {\mathcal R}_{34} \cap {\mathcal R}_{35} = \emptyset,
\end{aligned}\right.\ee
\be \label{st:sep:6}
 {\mathcal R}_{14} \cap {\mathcal R}_{24} = \emptyset, \quad {\mathcal R}_{14} \cap {\mathcal R}_{45} + k = \emptyset, \quad {\mathcal R}_{24} \cap {\mathcal R}_{45} + k = \emptyset,
 \ee
 \be \label{st:sep:7} 
{\mathcal R}_{25} \cap {\mathcal R}_{35} = \emptyset, \quad {\mathcal R}_{25} \cap {\mathcal R}_{45} = \emptyset, \quad {\mathcal R}_{35} \cap {\mathcal R}_{45} = \emptyset.
\ee

\end{prop}

\begin{proof} {\it Verification of \eqref{st:sep:1}:} this is straightforward, since $\l_1 > 0$ and $\l_2 > 0.$ 

\medskip

{\it Verification of \eqref{st:sep:2}:} Given $(\xi,\eta) \in {\mathcal R}_{12},$ we have
$$ (1 + (\xi + k)^2 + |\eta|^2)^{1/2}   = \o + (1 + \theta_e^2 (\xi^2 + |\eta|^2))^{1/2}.$$ 
We verify first that $\xi$ and $|\eta|$ are bounded in the small $\theta_e$ limit. 
We have
$$ 1 + (\xi + k)^2 + |\eta|^2 \leq 2\big( \o^2 + 1 + \theta_e^2 (\xi^2 + |\eta|^2)\big),$$
so that
$$ (1 - \theta_e^2) (\xi^2 + |\eta|^2) + 2 \xi k \leq \o^2 + 2.$$
With $- 2 |\xi k| \geq - 2 \e_0 \xi^2 - \frac{1}{2 \e_0} k^2,$
this gives
$$ (1 - \theta_e^2 - 2 \e_0) \xi^2 + (1 - \theta_e^2) |\eta|^2 \leq \o^2 +2 + \frac{1}{2 \e_0} k^2,$$
and we can choose $\e_0 = (1 - \theta_e^2)/4.$ This proves that if $\zeta = (\xi,\eta)$ belongs to ${\mathcal R}_{12},$ then $|\zeta|$ is bounded by a function of $k,$ uniformly in $\theta_e$ for $\theta_e$ smaller than a universal constant. 

Back to $\zeta = (\xi,\eta)$ in ${\mathcal R}_{12},$ we can now write 
$$ 1 + (\xi + k)^2 + |\eta|^2 = (\o + 1 + O(\theta_e^2))^2,$$
which gives
$$ (\xi + k)^2 + |\eta|^2 = \o^2 + 2 \o + O(\theta_e^2),$$
so that 
$(\xi,\eta) \in S(-k, (\o^2 + 2 \o + O(\theta_e^2))^{1/2}).$ If in addition $(\xi,\eta)$ belongs to ${\mathcal R}_{25} + k,$ then 
$$ (1 + \theta_e^2 (\xi^2 + |\eta|^2))^{1/2}  = \o - (1 + (\xi - k)^2 + |\eta|^2)^{1/2},$$
so that
$$ (\xi - k)^2 + |\eta|^2 = \o^2 - 2 \o + O(\theta_e^2))^2,$$
and $(\xi,\eta)$ belongs to the sphere $S(k, (\o^2 - 2 \o + O(\theta_e^2))^{1/2}).$  
By inspection (see Figure \ref{fig:blueblack}), for any $k_{max} > \sqrt 3,$ 
\be \label{1225:1} \inf_{k \in [\sqrt 3, k_{max}]} (k - R_-(k) - (-k + R_+(k)) > 0,\ee
and
\be \label{1225:2} \inf_{k \in [-k_{max}, -\sqrt 3]} (-k - R_+(k) - (k + R_-(k)) > 0,\ee 
with notation
\be \label{def:R+-}
 R_\pm(k) = (\o^2 \pm 2 \o)^{1/2}, \qquad \o = (1 + k^2)^{1/2}.
\ee
and this implies
$$ S\big(-k, (\o^2 + 2 \o + O(\theta_e^2))^{1/2}\big) \cap S\big(k, (\o^2 - 2 \o + O(\theta_e^2))^{1/2}\big) = \emptyset,$$
no matter the sign of $k,$ for $|k| \in [\sqrt 3, k_{max}]$ and $\theta_e$ small enough, depending only on $k_{max}.$ Hence \eqref{st:sep:2} for $\theta_e$ small enough. 

\begin{figure}
\begin{tikzpicture}[scale=1.3]
\draw (3,0) node[anchor=north] {$k$} ; 
\draw[->] (0,-2) -- (0,2) ;  
\node (1) at (-1,-.5) {${}$} ; 
\node (2) at (1,.5) {${}$} ; 
\node (3) at (-1,.5) {${}$} ; 
\node (4) at (1,-.5) {${}$} ; 
\draw[dotted,shift=(1)] (1.74,-1) -- (1.74,2) ; 
\draw[shift=(1)] (1.74,.5) node[anchor=north,fill=white] {${\footnotesize \sqrt 3}$} ;   
\filldraw[shift=(1)] (1.74,.5) circle (.05cm) ;
\draw[dotted,shift=(1)] (2.8,-1) -- (2.8,2) ; 
\draw[shift=(1)] (2.8,.5) node[anchor=north,fill=white] {${\footnotesize k_{max}}$} ;   
\filldraw[shift=(1)] (2.8,.5) circle (.05cm) ;

\draw[dotted,shift=(2)] (-1.74,-2) -- (-1.74,1) ; 
\draw[shift=(2)] (-1.74,-.5) node[anchor=north,fill=white] {${\footnotesize -\sqrt 3}$} ;   
\filldraw[shift=(2)] (-1.74,-.5) circle (.05cm) ;
\draw[dotted,shift=(2)] (-2.8,-2) -- (-2.8,1) ; 
\draw[shift=(2)] (-2.8,-.5) node[anchor=north,fill=white] {${\footnotesize -k_{max}}$} ;   
\filldraw[shift=(2)] (-2.8,-.5) circle (.05cm) ;

\draw[domain=1.74:3, samples=50,shift=(1)]
plot( \x, {-\x + ( 1 + \x^2 + 2*(1 + \x^2)^(1/2) )^(1/2) } ) ;  
\draw[domain=1.74:3, samples=50,blue,shift=(1)]
plot( \x, { \x - ( 1 + \x^2 - 2*(1 + \x^2)^(1/2) )^(1/2) } ) ;
 
\draw[domain=-1.74:-3, samples=50,shift=(2)]
plot( \x, {\x + ( 1 + (-\x)^2 - 2*(1 + (-\x)^2)^(1/2) )^(1/2) } ) ; 
\draw[domain=-1.74:-3, samples=50,blue,shift=(2)]
plot( \x, {-\x - ( 1 + (-\x)^2 + 2*(1 + (-\x)^2)^(1/2) )^(1/2) } ) ;

\draw[domain=1.74:3, samples=50,red,shift=(3)]
plot( \x, {\x - ( 1 + (\x)^2 + 2*(1 + (\x)^2)^(1/2) )^(1/2) } ) ;
\draw[domain=1.74:3, samples=50,yellow,shift=(3)]
plot( \x, {-\x + ( 1 + (\x)^2 - 2*(1 + (\x)^2)^(1/2) )^(1/2) } ) ;

 \draw[domain=-1.74:-3, samples=50,red,shift=(4)]
plot( \x, {-\x - ( 1 + (-\x)^2 - 2*(1 + (-\x)^2)^(1/2) )^(1/2) } ) ;
\draw[domain=-1.74:-3, samples=50,yellow,shift=(4)]
plot( \x, {\x + ( 1 + (-\x)^2 + 2*(1 + (-\x)^2)^(1/2) )^(1/2) } ) ;

\draw[->] (-3,0) -- (3,0) ;

\end{tikzpicture}

\caption{The blue curve is above the black curve: for $k > \sqrt 3,$ this illustrates \eqref{1225:1}, and for $k < - \sqrt 3$ this illustrates \eqref{1225:2}. Taken together, these inequalities imply the separation property ${\mathcal R}_{12} \cap {\mathcal R}_{25} + k = \emptyset.$ The red curve is above the yellow curve: for $k > \sqrt 3,$ this illustrates \eqref{1445:1}, and for $k < - \sqrt 3$ this illustrates \eqref{1445:2}. These inequalities imply ${\mathcal R}_{14} \cap {\mathcal R}_{24} + k = \emptyset.$} 

\label{fig:blueblack}
\end{figure}
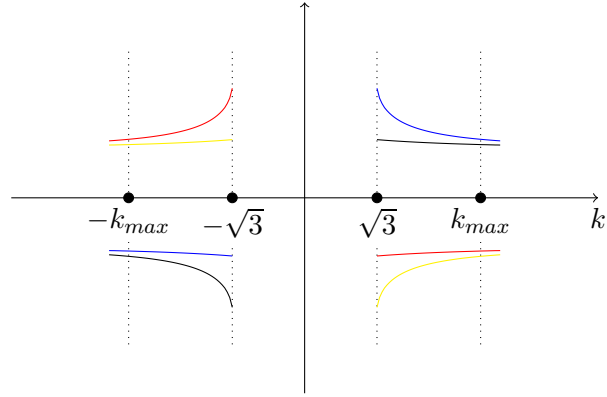

\medskip

{\it Verification of \eqref{st:sep:3}:} it is straightforward, just like the verification of \eqref{st:sep:1}. 

\medskip

{\it First and fourth lines in \eqref{st:sep:4}:} we observe that ${\mathcal R}_{13} \cap {\mathcal R}_{23} = \emptyset$ since $\theta_e < 1.$ Next we observe that if $\zeta$ belongs to the intersection ${\mathcal R}_{13} \cap {\mathcal R}_{34} + k,$ then it belongs to $S( -k , |k|) \cap S( k, |k|/\theta_e),$ which is empty as soon as $\theta_e < 1/3.$ By symmetry, the verification of ${\mathcal R}_{23} \cap {\mathcal R}_{35} + k = \emptyset$ goes exactly as the verification of ${\mathcal R}_{13} \cap {\mathcal R}_{34} + k = \emptyset,$ and the verification of ${\mathcal R}_{34} \cap {\mathcal R}_{35} = \emptyset$ goes exactly as the one for ${\mathcal R}_{13} \cap {\mathcal R}_{23}.$  

\medskip

{\it Second line in \eqref{st:sep:4}:} The intersection ${\mathcal R}_{13} \cap {\mathcal R}_{35} + k$ is equal to $S(-k, |k|) \cap S(k, |k|) = \{ 0 \}.$ By Proposition \ref{prop:space-time}, the intersection of the associated space-time resonant sets is also reduced to $\{ 0 \}.$ 

\medskip

{\it Third line in \eqref{st:sep:4}:} 
The set ${\mathcal R}_{23} \cap {\mathcal R}_{34} + k$ lies at the intersection of the sphere centered at $(-k,0) \in \R^3_{\xi,\eta}$ with radius $|k|/\theta_e$ and the sphere centered at $(k,0) \in \R^3_{\xi,\eta}$ with radius $|k|/\theta_e.$ This intersection is the disc centered at $\eta = 0$ with radius $(k^2/\theta_e^2 - k^2 - 1)^{1/2}$ in the $\{ \xi = 0\}$ plane. In particular, this disc does not intersect $\{ \eta = 0 \}.$ Since space-time resonances are confined to $\{ \eta = 0 \},$ this implies that ${\mathcal S}_{23} \cap {\mathcal S}_{34} + k = \emptyset.$

\medskip

{\it Verification of \eqref{st:sep:6}:} the first intersection is empty simply because $\theta_e < 1.$ Next if $(\xi,\eta) \in {\mathcal R}_{14} \cap {\mathcal R}_{45} + k,$ then we have 
$$ \left\{\begin{aligned} (1 + (\xi + k)^2 + |\eta|^2)^{1/2} & = \o - (1 + \theta_e^2 (\xi^2 + |\eta|^2))^{1/2}, \\ (1 + (\xi - k)^2 + |\eta|^2)^{1/2} & = \o + (1 + \theta_e^2 (\xi^2 + |\eta|^2))^{1/2}.
\end{aligned}\right.$$
Similarly to the analysis for ${\mathcal R}_{12} \cap {\mathcal R}_{25} + k,$ this implies 
$$ (\xi + k)^2 + |\eta|^2 = \o^2 - 2 \o + O(\theta_e^2), \quad (\xi - k)^2 + |\eta|^2 = \o^2 + 2 \o + O(\theta_e^2),$$
so that
$$ (\xi,\eta) \in S\big(-k, (\o^2 - 2 \o + O(\theta_e)^2)^{1/2}\big) \cap S\big(k, (\o^2 + 2 \o + O(\theta_e^2))^{1/2}\big).$$
The fact that these two spheres do not intersect is implied, for a given $k_{max} > \sqrt 3$ and for $\theta_e$ small enough, depending only on $k_{max},$ by the inequalities 
\be \label{1445:1}
 \inf_{k \in [\sqrt 3, k_{max}]} k - R_+(k) - (-k + R_-(k)) > 0,
\ee
and
\be \label{1445:2}
\inf_{k \in [-k_{max}, -\sqrt 3]} -k - R_-(k) - (k + R_+(k)) > 0,
\ee
with the radii $R_\pm$ defined in \eqref{def:R+-}.
Inequalities \eqref{1445:1} and \eqref{1445:2} are illustrated on Figure \ref{fig:blueblack}.

Finally, if $(\xi,\eta) \in {\mathcal R}_{45} + k,$ we find that $\xi$ and $|\eta|$ are bounded in the small $\theta_e$ limit, as we did for the $(1,2)$ resonance. If in addition $(\xi ,\eta) \in {\mathcal R}_{24},$ then from 
$$ (1 + \theta_e^2( (\xi + k)^2 + |\eta|^2))^{1/2} + (1 + \theta_e^2( \xi^2 + |\eta|^2))^{1/2} = \o,$$
we deduce
$ 2 + O(\theta_e^2) = \o^2,$
which is not true for $k > \sqrt 3$ and $\theta_e$ small enough.

\medskip

{\it Verification of \eqref{st:sep:7}:} this is straightforward, just like the verification of \eqref{st:sep:1}, by symmetry. 
\end{proof}

 \begin{figure}
\scalebox{.8}{
 \begin{tikzpicture}

 \begin{scope}[>=stealth]
 \draw[line width=.5pt][->] (-5,0) -- (5,0);

 \draw[line width=.5pt][->] (0,-3) -- (0,3);
 
 \draw (5,0) node[anchor=north] {$\R_\xi$} ;

 \draw (0,3) node[anchor=east] {$\R^2_\eta$} ; 
 
 \draw (1,0) node[anchor=north] {$(k,0)$} ; 
 
  \draw (-1,0) node[anchor=north] {$(-k,0)$} ; 
 
  \filldraw (1,0) circle (0.05cm) ; 
  
  \draw (1,0) circle (2.5cm) ;
  
  \filldraw (-1,0) circle (0.05cm) ; 
  
  \draw (-1,0) circle (2.5cm) ;
  
  \filldraw (0,2.27) circle (0.1cm) ;

    \filldraw (0,- 2.27) circle (0.1cm) ; 
    
   \draw[dotted] (-1,0) circle (1cm) ; 
   
   \draw[dotted] (1,0) circle (1cm) ; 
   
   \filldraw (0,0) circle (0.1cm) ; 
   
   \draw (-3.2,2) node[anchor=north,fill=white] {${\mathcal R}_{23}$} ;
   
   \draw (3.5,2) node[anchor=north,fill=white] {${\mathcal R}_{34} + k$} ;

   \draw (-2,1) node[anchor=north,fill=white] {${\mathcal R}_{13}$} ; 
    
   \draw (2.3,1) node[anchor=north,fill=white] {${\mathcal R}_{35} + k$} ; 
  
  \end{scope}
 
 \end{tikzpicture}

} 
 \caption{We see here the non-empty intersection ${\mathcal R}_{23} \cap {\mathcal R}_{34} + k:$ it is a circle in the $\{ \xi =  0 \}$ plane, which does not pass by the origin. As a consequence, ${\mathcal S}_{23} \cap {\mathcal S}_{34} + k = \emptyset.$ By contrast, ${\mathcal R}_{13}$ and ${\mathcal R}_{35} + k$ intersect precisely that the zero frequency, implying ${\mathcal S}_{13} \cap {\mathcal S}_{35} + k = \{ 0 \}.$}
 
 \label{fig:2334}
 
 \end{figure}
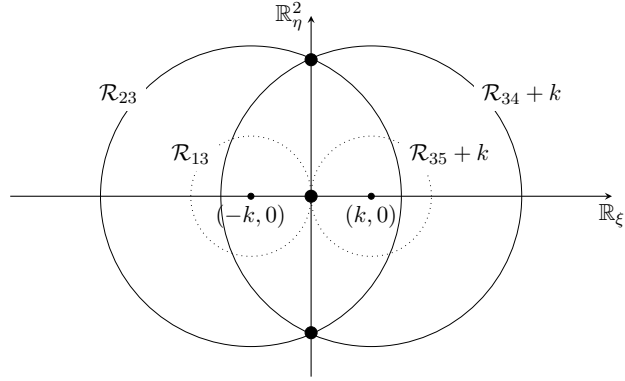

\begin{rem} \label{rem:about:matrix:B} Why do we focus on the separation properties of Proposition {\rm \ref{prop:space-time:sep}}? Here it's useful to jump ahead to the end of Section {\rm \ref{sec:6}.} In Section {\rm \ref{sec:step4},} we introduce matrix ${\bf B}$ which will play a key role in the equation for the symbolic flow. Proposition {\rm \ref{prop:space-time:sep}} is here precisely to describe how the rows and columns of matrix ${\bf B}$ are uncoupled, or coupled. This is illustrated in Figure {\rm \ref{fig:tree}.}  
\end{rem}

  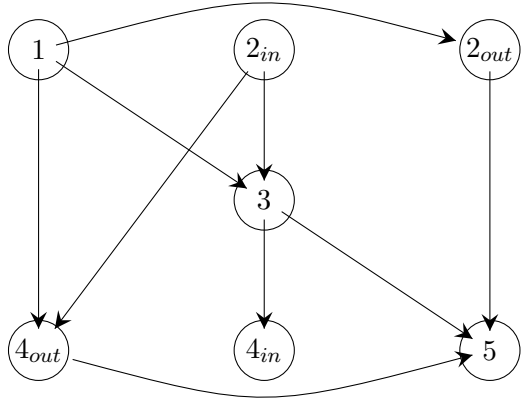
\begin{figure}
\scalebox{1}{
 \begin{tikzpicture}
\tikzset{myptr/.style={decoration={markings,mark=at position 1 with %
    {\arrow[scale=2,>=stealth]{>}}},postaction={decorate}}}
 \begin{scope}[>=latex]
 \node (1) at (-3, 2) {$1$};
\node (2in) at (0, 2) {$2_{in}$};
\draw (-3,2) circle (.4cm) ; 
\draw (0,2) circle (.4cm) ; 
\node (2out) at (3, 2) {$2_{out}$}; 
\draw[shift=(2out)] circle (.4cm) ; 
\node (3) at (0,0) {$3$};
\draw[shift=(3)] circle (.4cm) ; 
\node (4out) at (-3,-2) {$4_{out}$};
\draw[shift=(4out)] circle (.4cm) ; 
\node (4in) at (0,-2) {$4_{in}$}; 
\draw[shift=(4in)] circle (.4cm) ; 
\node (5) at (3,-2) {$5$}; 
\draw[shift=(5)] circle (.4cm) ; 
\end{scope}
\draw[myptr] (1) .. controls (0,2.8) .. (2out) ;
\draw[myptr] (3) .. controls (1.5,-1) .. (5) ; 
\draw[myptr] (2in) .. controls (0,1) .. (3) ;
\draw[myptr] (3) .. controls (0,-1) .. (4in) ; 
\draw[myptr] (1) .. controls (-3,0) .. (4out) ;
\draw[myptr] (2out) .. controls (3,0) .. (5) ;
\draw[myptr] (4out) .. controls (0,-2.8) .. (5) ; 
\draw[myptr] (2in) .. controls (-1.5,0) .. (4out) ; 
\draw (-.8,0.2) node[anchor=south,fill=white] {$\,$} ; 
\draw[myptr] (1) .. controls (-1.5,1) .. (3) ;
 \end{tikzpicture}
 }
 \caption{This is a picture of the ways the rows and columns of matrix ${\bf B}$ from Section \ref{sec:6} are coupled. Seventeen intersections of resonant sets need to be considered. There are 3 arrows going out of the $1$ node, corresponding to 3 intersections: those are described in \eqref{st:sep:1}. There are only 2 arrows connected to the $2_{out}$ node, corresponding to one possible intersection: this one is described in \eqref{st:sep:2}. Etc.}

 \label{fig:tree} 
 \end{figure}

\subsection{Eigenspace decomposition for the linearized singular source.} \label{sec:B}
We may decompose the symbol $B$ onto the eigenspaces of the symbol of $A.$ According to \eqref{id:proj}, we have, pointwise in $(t,x,y,\xi,\eta):$
 $$ B = \sum_{1 \leq j,j' \leq 5} \Pi_{j} B \Pi_{j'},$$
  so that
$$  B(t,x,y,\xi,\eta) = \sum_{p,j,j'} e^{i p \theta} g_p(t,x,y) B_{pjj'}(\xi,\eta),$$
where $g_p$ is defined in \eqref{def:g}, and $B_{pjj'}$ are the $(t,x,y)$-independent matrix-valued symbols
\be \label{def:interaction:coefficient}
 \begin{aligned} B_{pjj'}(\xi,\eta) & := \Pi_{j,+p}(\xi,\eta) {\mathcal B}(\vec e_p, \Pi_{j'}(\xi,\eta) \, \cdot) + \Pi_{j,+p}(\xi,\eta) {\mathcal B}(\Pi_{j'}(\xi,\eta)\, \cdot, \vec e_p) \\ & \qquad - \Pi_{j,+p}(\xi,\eta) {\mathcal J}_e(\Pi_{j'}(\xi,\eta) \, \cdot , i p k) \vec e_p \\ & \qquad - (1 - \chi_{high}) \Pi_{j,+p}(\xi,\eta) {\mathcal J}_e(\vec e_p,i \xi,i \eta) \Pi_{j'}(\xi,\eta).\end{aligned}
 \ee
 Here $\Pi_{j,+p}(\xi,\eta): = \Pi_{j}(\xi + pk , \eta )$.  The subscript $+p$ refers to a frequency shift induced by the fast WKB oscillations (see Remark \ref{rem:shift}).
  
Above in the decomposition of $B$ into a sum over $(p,j,j'),$ we used the shorthand
\be \label{short:sum}
\sum_{p,j,j'} \dots = \sum_{(p,j,j') \in \{-1,1\} \times \{1,2,3,4,5\}^2} \dots
\ee

\subsection{Frequency cut-offs} \label{sec:cutoffs}

 Given any smooth function $\chi \in C^\infty_c$ such that $0\leq \chi \leq 1,$ we use notation (borrowed from \cite{LNT})
 \be \label{sharp:tilde}
 \chi^\flat \prec \chi \prec \chi^\sharp
 \ee
to denote a couple $(\chi^\flat,\chi^\sharp)$ of smooth and compactly supported functions, each taking values in $[0,1],$ which satisfy the relations
$$ (1 - \chi) \chi^\flat \equiv 0, \qquad (1 - \chi^\sharp) \chi \equiv 0.$$
In other words: $\chi \equiv 1$ on the support of $\chi^\flat,$ $\chi^\sharp \equiv 1$ on the support of $\chi.$ Somehow, $\chi$ extends (or ``is above", hence the musical notation) $\chi^\flat,$ and $\chi^\sharp$ extends $\chi.$ When the analysis requires a further extension, we use a tilde:
\be \label{sharp:tilde2} \chi^\sharp \prec \tilde \chi.\ee

Let $\chi_0: \R \to [0,1]$ be a standard smooth ``plateau'' cut-off, that is such that $\chi_0(x) = 1$ if $|x| \leq 1,$ and $\chi_0(x) = 0$ if $|x| \geq 2,$ with $0 \leq \chi_0 \leq 1,$ and $\chi_0$ monotonic and smooth. (In the definition of $\chi_{trans}$ on page \pageref{chi0}, we used a two-dimensional version of this type of cut-off and also denoted it $\chi_0.$)

\begin{defi} \label{def:cut-offs} Given a resonant pair $(j,j')$ and $\delta_{res} > 0,$  we define $\chi_{jj'} \in C^\infty_c(\R^3)$ to be the cut-off around ${\mathcal R}_{jj'}$ defined by 
 \be \label{def:chijj'} \chi_{jj'}(\xi,\eta) := \chi_0\left( \frac{C \sqrt{\Phi_{jj'}(\xi,\eta)^2 + (\delta_{res}/6)^2}}{\delta_{res}}\right), \quad C = 2.\ee
 If $(j,j')$ is not a resonant pair, we let $\chi_{jj'} = 0.$ For $\chi_{jj'}^\flat,$ $\chi_{jj'}^\sharp$ and $\tilde \chi_{jj'},$ we use the right-hand side of \eqref{def:chijj'} where we change $C$ into $C^\flat = 4,$ $C^\sharp = 1,$ $\tilde C = 1/2.$ This guarantees $\chi_{jj'}^\flat \prec \chi_{jj'} \prec \chi_{jj'}^\sharp \prec \tilde \chi_{jj'}.$ %
\end{defi}

With the above definition, we note that $\chi_{jj'}$ is smooth (the $\delta_{res}/6$ term is here to guarantee smoothness). 

We note that $\chi_{jj'}(\xi,\eta) = 1$ is equivalent to $(\Phi_{jj'}(\xi,\eta)^2 + ( \delta_{res}/6)^2)^{1/2} \leq \delta_{res}/2,$ that is $|\Phi_{jj'}(\xi,\eta)| \leq \delta_{res}(1/4-1/36)^{1/2}.$ In particular, the resonant set ${\mathcal R}_{jj'}$ is included in $\chi_{jj'}^{-1}(\{ 1 \}).$ We also note that if $(\xi,\eta)$ belongs to the support of $\chi_{jj'},$ then $|\Phi_{jj'}(\xi,\eta)| \leq \delta_{res}(1 - 1/36)^{1/2}.$ If $\chi_{jj'}(\xi,\eta) < 1,$ then $2 (\Phi_{jj'}^2 + (\delta_{res}/6)^2)^{1/2} > \delta_{res},$ implying 
 $|\Phi_{jj'}(\xi,\eta)| > \delta_{res}(1/4 - 1/36)^{1/2}.$  In particular:
 \be \label{transp}
 \mbox{The phase function $\Phi_{jj'}$ is bounded away from zero on the support of $1 - \chi_{jj'}.$}
 \ee

\begin{lem} \label{lem:cut-offs} Under the assumptions of Proposition {\rm \ref{prop:separation}}, we can choose $\delta_{res}$ to be small enough so that all the separation properties of  Proposition {\rm \ref{prop:separation}} extend to all the cut-offs associated with resonant pairs, as introduced in Definition {\rm \ref{def:cut-offs}}. For instance, we have 
 \be \label{sep:cut-offs}  \tilde \chi_{23}(\xi - k, \eta) \tilde \chi_{12}(\xi,\eta) \equiv 0.\ee
\end{lem}

\begin{proof} By Proposition \ref{prop:separation}, if $|k| > k_{min},$ then the sets ${\mathcal R}_{12}$ and ${\mathcal R}_{23} + k$ do not intersect. By continuity of the phase functions \eqref{def:phase:function}, the resonant sets are closed. Since they are disjoint, there are open neighborhoods $U_{12}$ and $U_{23}$ of ${\mathcal R}_{12}$ and ${\mathcal R}_{23} + k$ which do not intersect. The phase functions associated with resonances are proper (see Proposition \ref{prop:res}(1)), and this implies that we can find $0 < \delta < 1$ such that $\{ |\Phi_{12}| < \delta\} \subset U_{12}$ and $\{|\Phi_{23}(\cdot - k, \cdot) | < \delta\} \subset U_{23}.$ Thus
$$ \{ |\Phi_{12}| < \delta \} \cap \{ |\Phi_{23}(\cdot - k,\cdot)| < \delta \} = \emptyset.$$
It only remains to show that if $\delta_{res}$ is small enough, then the support of $\tilde \chi_{12}$ is included in $\{|\Phi_{12}| < \delta\}$ (and similarly for the $(2,3)$ resonance). By definition of $\tilde \chi_{jj'}$ (see Definition \ref{def:cut-offs}), the support of $\tilde \chi_{12}$ is included in $\{ |\Phi_{12}| \leq \delta_{res}(16 - 1/36)^{1/2}\}.$ Thus the condition is $\delta_{res}(16 - 1/36)^{1/2} < \delta.$ Under this condition, the support of $\tilde \chi_{12}$ is included in  $\{|\Phi_{12}| < \delta\},$ and \eqref{sep:cut-offs} holds true. For each of the other separation properties of Proposition \ref{prop:separation}, we find a similar condition bearing on $\delta_{res}.$ Note that $\delta > 0$ a priori depends on the resonant pair. Since there are a finite number of separation conditions, we can find $\delta_{res}$ small enough which satisfies them all. 
\end{proof}

A further smallness constraint bearing on $\delta_{res}$ will be formulated in Corollary \ref{cor:234}.

\section{Reduction to normal form} \label{sec:normal}

 This is a central, if now classical, piece in the analysis. In the context of ordinary differential equations, the idea dates back to Poincar\'e; in this high-frequency context to Joly, M\'etivier and Rauch \cite{JMR-TMB} (see also \cite{em3, CG, Lu, em4}). With a linear change of variable to normal form, we eliminate here all interaction coefficients away from resonant frequencies.  

\begin{prop}[Normal form reduction] \label{prop:normal:form} For some $Q \in S^{0},$ for $\check u$ defined by
  \be \label{def:checku:good} \check u := (\Id + \sqrt \e \op_\e(Q))^{-1} \dot u,\ee
 with $\dot u$ the solution to \eqref{pert:2}, we have
 \be \label{eq:checku0}
 \d_t \check u + \frac{1}{\e} A(\nabla_\e) \check u + \frac{1}{\sqrt \e} \op_\e(J) \check u = \frac{1}{\sqrt \e} \op_\e(\check B) \check u + \op_\e(\check J) \check u + \check F,
\ee
where
\be \label{def:checkJ}
\op_\e(\check J) := - (\Id + \sqrt \e \op_\e(Q))^{-1} [\op_\e(J + R_A), \op_\e(Q)], \qquad \mbox{so that $\check J  \in S^1,$}
\ee
and
\be \label{def:checkB}
 \begin{aligned} \op_\e(\check B) \check u & = g_1(t,x,y) \sum_{(j,j') \in {\mathcal R}} e^{i \theta} \op_\e(\chi_{jj'} B_{1jj'}) \check u \\ & \qquad + g_{-1}(t,x,y) \sum_{(j,j') \in {\mathcal R}} \op_\e(\chi_{jj'} (B_{-1j'j})_{+1}) (e^{- i \theta} \check u).\end{aligned}
\ee
The scalar amplitudes $g_{\pm 1}$ are introduced in \eqref{def:gp}. In \eqref{eq:checku0}, the source $\check F$ is equivalent to $\dot F,$ in the sense of Definition {\rm \ref{def:equiv:remainder}.}
\end{prop}

In \eqref{def:checkB} below, the sum $\dsp{\sum_{(j,j') \in {\mathcal R}} \dots}$  means the sum over all resonant pairs, as described in Proposition \ref{prop:res}.

In Proposition \ref{prop:normal:form}, the subscript $+1$ in $(B_{0,-1j'j})_{+1}$ refers to a frequency shift induced by the fast WKB oscillations (see notation introduced in Remark \ref{rem:shift}).  The interaction coefficients $B_{pjj'}$ are defined in \eqref{def:interaction:coefficient}. The cut-offs $\chi_{jj'}$ are introduced in Definition \ref{def:cut-offs}.

 Compared to \cite{em4}, the key new term is $\check J,$ which captures the interplay between the convective terms and the normal form, and the error $R_A$ in the approximation of $A$ by $A_{\e = 0}$ (see \eqref{spectral:dec:symbols}) and the normal form. 

\begin{proof} We use the description of $B$ given in Section \ref{sec:B}.

\medskip

{\it Step 1: taking out a single interaction coefficient from $B$ with $p = 1.$}  We single out a pair of indices $(j,j')$ and the oscillation corresponding to $p = 1.$ We work in a neighborhood of the set ${\mathcal R}_{jj'}$ comprising $(j,j')$-resonant frequencies.

 In this first step, we look for a symbol $Q \in S^0,$ supported in a neighborhood of ${\mathcal R}_{jj'},$ such that, letting \eqref{def:checku:good},
the unknown $\check u$ satisfies the reduced perturbation equation
\be \label{eq:checku}
\d_t \check u + \frac{1}{\e} A(\nabla_\e) \check u + \frac{1}{\sqrt \e} \op_\e(J) \check u = \frac{1}{\sqrt \e} \op_\e(\check B_0) \check u  + \op_\e(\check J) \check u + \check F,
\ee
where $\check F$ is equivalent to $\dot F,$ the symbol $\check B_0$ is the reduced source
\be \label{def:checkB:inproof}
\check B_0 = e^{i \theta} g_1 \chi_{jj'} B_{1jj'} + \sum_{ (q,\ell,\ell') \neq (1,j,j')} e^{i q \theta} g_q B_{q\ell\ell'},
\ee
and $\check J$ is an order-one symbol that is not singular in $\e.$ 
Once \eqref{eq:checku}-\eqref{def:checkB:inproof} is achieved, we may sum up Step 1 as follows: $$ \mbox{{\it The term $e^{i  \theta} g_1 (1 - \chi_{jj'}) B_{1jj'}$ is removed from the source, via the change of variable \eqref{def:checku:good}}}.$$ 
The price to pay for this removal is the creation of a possibly non-symmetric, order-one symbol $\check J,$ which is not singular in $\e.$ 

Many pair of indices $(j,j')$ do not correspond to resonant pairs. The associated ``cut-off" is then $\chi_{jj'} \equiv 0$ (see Definition \ref{def:cut-offs}). In the change of variable that we now describe, for those pairs the whole interaction coefficient $B_{1jj'}$ is removed from the source.

\medskip

Thus we posit \eqref{def:checku:good}, and look for $Q$ in the form
\be \label{ansatz:Q} Q(\e,t,x,y,\xi,\eta) = e^{i \theta} g_1(t,x,y) \Pi_{j,+1}(\xi,\eta)  Q_{1jj'}(\xi,\eta) \Pi_{j'}(\xi,\eta), \qquad Q_{1jj'} \in S^{0}.\ee
Then, from the perturbation equations \eqref{pert:2} in $\dot u,$ we deduce that $\check u$ solves
\be \label{for:homo} \begin{aligned}  \d_t \check u  + \frac{1}{\e} A(\nabla_\e) \check u  & = \frac{1}{\sqrt \e} \op_\e(B - J) \check u \\ & + ({\rm Id} + \sqrt \e \op_\e(Q))^{-1} \Big[ - \frac{1}{\sqrt \e} A(\nabla_\e)  + \op_\e(B - J), \op_\e(Q)\Big] \check u \\ & - \sqrt \e \op_\e(Q))^{-1} \op_\e(\d_t Q) \check u + ({\rm Id} + \sqrt \e \op_\e(Q))^{-1} \dot G.
\end{aligned}\ee
Consider first the leading, $O(1/\sqrt\e)$ terms in the above right-hand side. There are three of these:
\begin{itemize}
\item the original ``source" term $(1/\sqrt \e) B;$
\item the commutator $[A(\nabla_\e), \op_\e(Q)];$
\item and also the term involving $\d_t Q,$ since $Q$ contains $O(1/\e)$ time oscillations.
\end{itemize}
With \eqref{approximate:spectral:decomposition} and the a priori form of $Q$ \eqref{ansatz:Q}, we have
$$ [A(\nabla_\e), \op_\e(Q)] = \sum_{\a} \big[\op_\e(i \l_{\a} \Pi_{\a}), \op_\e(e^{i \theta} g_1 \Pi_{j,+1}  Q_{1jj'} \Pi_{j'})\big] + \sqrt \e [\op_\e(R_A), \op_\e(Q)],$$
where the sum in $\a$ ranges over $\{1,2,3,4,5\}.$
The fast spatial oscillations in $Q$ induce a shift, as described in Remark \ref{rem:shift}:
$$ \begin{aligned} {} [A(\nabla_\e), \op_\e(Q)] & = \sum_{\a} e^{i \theta} \op_\e(i \l_{\a,+1} \Pi_{\a, +1}) \op_\e(g_1 \Pi_{j,+1} Q_{1jj'} \Pi_{j'}) \\ & \qquad - \sum_{\a} e^{i \theta} \op_\e(g_1 \Pi_{j,+1} Q_{1jj'} \Pi_{j'}) \op_\e(i \l_{\a} \Pi_{\a}) \\ & \qquad + \sqrt \e [\op_\e(R_A), \op_\e(Q)].\end{aligned}$$
By orthogonality of the eigenprojectors of $A$ and composition of pseudo-differential operators (see Appendix \ref{app:symb}), since the above symbols are not oscillating fast in space (note that the $e^{i\theta}$ factors were pulled out), the leading term above is
 $$ e^{i \theta} g_1(t,x,y) \op_\e\Big(\,  i (\l_{j,+1} - \l_{j'}) \Pi_{j,+1} \tilde Q_{1jj'} \Pi_{j'}\Big).$$
 At this point what is left out, in the leading term in $\e$ in $[A(\nabla_\e), \op_\e(Q)]$ commutator, is a term that is both order 0 as a pseudo-differential operator and non singular in $\e.$ Such a term belongs to the remainder term $\check F.$ Since $Q$ is a priori order zero, the commutator $\sqrt \e [\op_\e(R_A), \op_\e(Q)]$ is order one and cannot be seen as a remainder. We put it in the $\check J$ error term. 

 Taking into account the contribution of the original source term $B$ and of $\d_t Q$ we arrive at a leading term in the form
 $$
\begin{aligned} \frac{1}{\sqrt \e} & \op_\e\Big(e^{i \theta} g_1 \Pi_{j,+1} \Big(B_{1jj'} -  \big(\l_{j,+1} - \l_{j'} - \o \big) Q_{1jj'}  \Big) \Pi_{j'}\Big) \\ & \qquad + \frac{1}{\sqrt \e} \op_\e\Big( - J + \sum_{(q,\ell,\ell') \neq (p,j,j')} e^{i q \theta} B_{q\ell\ell'}\Big).
\end{aligned}$$
We let
\be \label{homo}
 Q_{1jj'} := (1 - \chi_{jj'}(\xi,\eta)) \big(\l_{j,+1} - \l_{j'} - \o\big)^{-1} B_{pjj'}.
\ee
This is an admissible definition: if $(j,j')$ is a resonant pair, then the phase function $\l_{j,+1} - \l_{j'} - \o$ is bounded away from zero on the support of $1 - \chi_{jj'},$ as observed in \eqref{transp}, and if $(j,j')$ is not a resonant pair, then the phase function is uniformly bounded away from zero (see Proposition \ref{prop:res}(1)). Thus the symbol $Q$ is well defined and belongs to $S^{0}.$ With this choice, we have
 $$ B_{1jj'} -  \big(\l^\e_{j,+1} - \l^\e_{j'} - \o \big) Q_{1jj'}  = O(\e),$$
 uniformly in frequency.
The equation in $\check u$ takes the form \eqref{eq:checku}, with
$\op_\e(\check J)$ given by \eqref{def:checkJ}.  This concludes Step 1.
\medskip

{\it Step 2: same as Step 1 but with $p = -1.$} Consider now oscillations associated with $p = -1.$ By symmetry with Step 1, we will consider indices $(j',j).$ Thus the symbol in the source under consideration is $e^{- i \theta} B_{-1j'j}.$ Working as in Step 1, we find a phase
 \be \label{other:resonance} \l_{j',-1} + \o - \l_{j} = - (\l_{j,+1} - \o - \l_{j'})_{-1},\ee
 using notation for shifted symbols introduced in Remark \ref{rem:shift}. Thus, what is left of $B_{-1j'j}$ after the reduction to normal form is
 $$ e^{- i \theta} g_{-1}(t,x,y) \op_\e(\chi_{jj',-1} B_{-1j'j}).$$
 According to Remark \ref{rem:shift}, this operator can also be written
 $$ g_{-1}(t,x,y) \op_\e(\chi_{jj'} B_{-1j'j,+1}) \big( e^{- i \theta} \, \cdot \big).$$

{\it Step 3.} Next we observe that the above procedure of normal form reduction is an additive process, in the following sense: given two distinct triples $(p,j,j')$ and $(p_1,j_1,j'_1),$ Step 1 and 2 provide associated symbols $Q$ and $Q_{1}\in S^{0}.$ If we posit
\be \label{simultaneously} \check u = \Big(\Id + \sqrt \e \big( e^{i p_1 \theta} \op_\e(Q_1) + e^{i p \theta} \op_\e(Q) \Big)^{-1} \dot u,\ee
we obtain a reduced source with truncated $(p,j,j')$ term {\it and} truncated $(p_1,j_1,j'_1)$ term. The same holds if we posit
\be \label{consecutively} \check u = (\Id + \sqrt \e  e^{i p_1 \theta} \op_\e(Q_1) )^{-1} (\Id + \sqrt\e e^{i p \theta} \op_\e(Q) \big)^{-1} \dot u.\ee
In \eqref{simultaneously}, we take out the different components of the source away from their respective resonances simultaneously. In \eqref{consecutively}, we take these components out consecutively. Up to leading order in the source, the result is the same.

\medskip

{\it Conclusion.} We remove all the source terms away from resonant frequencies and arrive at \eqref{eq:checku0}-\eqref{def:checkJ}-\eqref{def:checkB}. By inspection, the source $\check F$ is equivalent to $\dot F,$ in the sense of Definition \ref{def:equiv:remainder}.

\end{proof}

\begin{rem} \label{rem:Q-1} With $Q$ defined as in the proof of Proposition {\rm \ref{prop:normal:form}}, we have $Q_{1jj'} \in S^{-1}$ unless $j = j'.$ Indeed, if $j \neq j',$ the phase function $\Phi_{jj'}$ belongs to $S^1$ (see Proposition {\rm \ref{prop:res}(1)}. Similarly, $Q_{-1j'j}$ belongs to $S^{-1}$ unless $j = j'.$ We will exploit this remark in the upcoming Section.  
\end{rem}

\section{Another reduction to normal form} \label{sec:normal:2}

 We perform here a variation on the normal form reduction of Section \ref{sec:normal}. What was achieved in Proposition \ref{prop:normal:form} was the elimination of terms in the source $B$ away from the frequency set ${\mathcal R}.$ Due to the presence of convective terms, and the high-frequency behavior of the longitudinal acoustic modes\footnote{This high-frequency behavior is responsible for the error term $R_A$ in \eqref{spectral:dec:symbols} being order one; see Sections \ref{sec:char0} and \ref{sec:spectral:dec})}, the price to pay for that elimination was the creation of the non-singular but non-symmetric and order-one term $\op_\e(\check J).$ We exploit here Remark \ref{rem:Q-1} as we show that we may perform a reduction to normal form which does not create order-one terms. The downside is that the reduction is incomplete: auto-interaction coefficients persist.

\begin{prop}[Second reduction to normal form] \label{prop:normal:form:2} For some $Q_{(-1)} \in S^{-1},$ for $\check v$ defined by
  \be \label{def:checkv} \check v := (\Id + \sqrt \e \op_\e(Q_{(-1)}))^{-1} \dot u,\ee
 with $\dot u$ the solution to \eqref{pert:2}, we have
 \be \label{eq:checkv0}
 \d_t \check v + \frac{1}{\e} A(\nabla_\e) \check v + \frac{1}{\sqrt \e} \op_\e(J) \check v = \frac{1}{\sqrt \e} \op_\e(\check B_{(-1)}) \check v + \check F_{(-1)},
\ee
where, informally,
\be \label{def:checkB-1:informal}
 \check B_{(-1)} = \check B + \Big(\mbox{auto-interaction coefficients in restriction to the support of $\chi_{high}^\flat$}\Big),
\ee
that is, precisely, 
\be \label{def:checkB-1}
 \begin{aligned} \op_\e(\check B_{(-1)}) \check u & = \op_\e(\check B)\check v \\ & + g_1(t,x,y) \sum_{1 \leq j \leq 5} e^{i \theta} \op_\e( \chi_{high}^\flat B_{1jj}) \check v \\ & +  g_{-1}(t,x,y) \sum_{1 \leq j \leq 5} \op_\e((\chi_{high}^\flat B_{-1jj})_{+1}) (e^{- i \theta} \check v).\end{aligned}
\ee
The scalar amplitudes $g_{\pm 1}$ are introduced in \eqref{def:gp}. In \eqref{eq:checkv0}, the remainder $\check F_{(-1)}$ 
is equivalent to $\dot F$ in the sense of Definition {\rm \ref{def:equiv:remainder}.}
\end{prop}

\begin{proof} It suffices to follow the proof of Proposition \ref{prop:normal:form} and leave the auto-interaction coefficients untouched on the support of $\chi_{high}^\flat.$  Then, the associated $Q_{(-1)}$ is equal to $Q$ minus the $\chi_{high}^\flat Q_{1jj}$ and $\chi_{high}^\flat Q_{-1jj}$ terms. In particular, as noticed in Remark \ref{rem:Q-1}, we have $Q_{(-1)} \in S^{-1}.$ Indeed, the symbols $Q_{pjj'}$ with $j \neq j'$ belong to $S^{-1}.$ The symbols $(1 - \chi_{high}^\flat) Q_{pjj}$ also belong to $S^{-1}$ because $1 - \chi_{high}^\flat$ is a low-frequency truncation. The associated commutator $[J + R_A, Q_{(-1)}]$ belongs to $S^0,$ hence what was $\op_\e(\check J)$ in Proposition \ref{prop:normal:form} is now an order-zero, non-singular-in-$\e$ pseudo-differential operator which belongs in the remainder term $\check F_{(-1)}.$ 
\end{proof} 

At this point, we have two distinct reductions to normal form for the system in $\dot u:$
\begin{itemize}
\item a system in $\check v$ \eqref{eq:checkv0}, in which some spurious source terms are present in the right-hand side; spurious since they do not correspond to resonances. The system in $\check v$ will be used, in Section \ref{sec:weak:Sobolev:bound}, in order to produce a ``rough" Sobolev bound on the solution.  
\item a system in $\check u$ \eqref{eq:checku0}, in which all source terms have been eliminated away from resonant frequencies, but there appears an order-one non-symmetric term $\op_\e(\check J).$ 
\end{itemize}

\section{A rough Sobolev bound} \label{sec:weak:Sobolev:bound}

Here we prove a rough a priori Sobolev bound for the solution $\dot u$ of the initial-value problem \eqref{pert:2}, based on the system \eqref{eq:checkv0} in $\check v.$ The bound is rough in the sense that it features a rate of exponential growth ($\Gamma$ below) that is far from optimal.
Much of the subsequent analysis after this Section will be devoted to deriving the optimal rate of growth.

 We use the weighted-in-$x$-only Sobolev norm $\| \cdot \|_{\e,s}$ defined in \eqref{weighted:norm}. This norm is naturally associated with the weighted-in-$x$-only differential operator $\d_\e$ introduced in \eqref{nablae}.  
 
We first observe that Sobolev bounds in $\check u$ or $\check v$ imply Sobolev bounds in $\dot u,$ and vice versa:

\begin{lem} \label{lem:dot-checkv} For all $t \leq t_\star(\e),$ we have
  \be \label{equiv:checkv:dotu} \| \check u(t)\|_{\e,s} \leq (1 + \sqrt \e C_1) \| \dot u(t) \|_{\e,s}, \qquad \| \dot u(t) \|_{\e,s} \leq (1 + \sqrt\e C_2) \|  \check u(t) \|_{\e,s},
  \ee
 where $\check u$ is our solution after the first normal form reduction, defined in \eqref{def:checku:good} In \eqref{equiv:checkv:dotu}, the positive constants $C_1$ and $C_2$ are independent of $\e,t.$ The same bounds, with different constants $C_1$ and $C_2,$ hold also with $\check v,$ our solution after the second normal form reduction, defined in \eqref{def:checkv}.
 
 Identical bounds hold in $\fl1$ norms: we have, for $t \leq t_\star(\e),$ and for $|\a|$ constrained only by the regularity of $u_a:$ 
\be  \| \d_\e^\a \check u(t)\|_{\fl1} \leq (1 + \sqrt \e C_1) \max_{|\b| \leq |\a|} \| \d_\e^\b \dot u(t) \|_{\fl1}, \quad \| \d_\e^\a \dot u(t) \|_{\fl1} \leq (1 + \sqrt\e C_2)  \max_{|\b| \leq |\a|} \| \d_\e^\b \check u(t) \|_{\fl1},
 \ee
 and same with $\check v$ instead of $\check u.$ 
  \end{lem}

\begin{proof} By definition of $\check u$ \eqref{def:checku:good}, 
$$ \d_\e^\a \dot u = \d_\e^\a \check u + \e^{1/2} \d_\e^\a \op_\e(Q) \check u.$$
We see in the proof of Proposition \ref{prop:normal:form} that every entry of $Q$ is a tensor product of the form $e^{i p \theta} g_p(x,y) \underline Q(\xi,\eta),$ with $\underline Q \in S^{0}.$ Thus by Lemma \ref{lem:tensor}, we have the bounds for $\dot u$ in terms of $\check u.$ 

Conversely, we expand in Neumann series: 
$$ \check u = \sum_{n \geq 0} (-\e)^{n/2} \op_\e(Q)^n \dot u.$$
For fixed $n,$ we may now use Lemma \ref{lem:tensor} repeatedly, and this proves for the bounds for $\check u$ in terms of $\dot u.$ 

The same arguments hold for $\check v,$ with $\underline Q \in S^{-1}$ instead of $\underline Q \in S^0.$ 
\end{proof}

\begin{prop} \label{prop:weak:bound} For $0 \leq t \leq t_\star(\e),$ so that bound \eqref{a:priori:tilde:bound} holds, the solution $\dot u$ to \eqref{pert:2} satisfies the bound
 \be \label{weak:Sobolev:bound}
  \| \dot u(t) \|_{\e,s} \lesssim \e^K e^{t \Gamma/\sqrt \e}, 
 \ee
 for some $\Gamma > 0$ that is made semi-explicit in the proof below. The weighted-in-$x$-only Sobolev norm $\| \cdot \|_{\e,s}$ is defined in \eqref{weighted:norm}. 
\end{prop}

\begin{proof} For the most part, this is a standard Sobolev estimate. We use the system \eqref{eq:checkv0} in $\check v,$ derived in the second norm form reduction (Section \ref{sec:normal:2}). We carefully follow the dependence in $\e.$ Multiplicative constants do matter since they translate into rates of exponential growth.  
 
 For $|\a| \leq s,$ we let $v_\a := \d_\e^\a \check v,$ where $\d_\e$ is defined in \eqref{def:de},  and consider the real part of the $L^2$ scalar product with $v_\a.$ Since $A(\nabla_\e)$ is a hyperbolic Fourier multiplier, we have
$$ \Re e \, \big( \d_\e^\a A(\nabla_\e) \check v, v_\a\big)_{L^2} = 0.$$
We turn to convective terms. As discussed in Section \ref{sec:conv:pert}, these take the form \eqref{for:convection:this:one}, which we partially reproduce here: 
\be \label{convection} \frac{1}{\sqrt \e}\op_\e(J) = (\tilde v_{e0} \cdot \d_\e  +  \tilde w_e \cdot \nabla_\e) \op_\e(\chi_{high} \underline J_e)  + \mbox{smaller ionic terms.} \ee %
 The ionic convective terms are symmetric in the same sense as the electronic convective terms (see Section \ref{sec:conv}), and smaller than the electronic convective terms. Thus in the convection we focus on the role played by the electrons. %

 We have, using \eqref{cancellation:convection} and symmetry of $\underline J_e,$ 
$$ \begin{aligned} \Big| \Re e \, \big( (\tilde v_{e0} \cdot \d_\e) \op_\e(\chi_{high} \underline J_e) v_\a, v_\a \big)_{L^2} \Big| & = \Big| - \frac{1}{2} \big( (\d_y \cdot \tilde v_{e0y}) \op_\e(\chi_{high}\underline J_e) v_\a, v_\a \big)_{L^2}\Big| \\ & \leq \frac12 \| \d_y \cdot \tilde v_{e0y} \|_{L^\infty} \| v_\a \|_{L^2}^2.\end{aligned}$$
The commutator term is bounded as follows
$$ \Big| \Re e \, \big( \sum_{\begin{smallmatrix} \a_1 + \a_2 = \a \\ |\a_1| > 0 \end{smallmatrix}} ( \d_\e^{\a_1} \tilde v_{e_0}) \cdot \d_\e \op_\e(\chi_{high} \underline J_e) v_{\a_2}, v_\a \big)_{L^2} \Big| \lesssim \max_{|\b| \leq |\a|} \| \d_\e^\b \tilde v_{e0} \|_{L^\infty} \| \check v \|_{\e,s}^2.$$
The $\tilde w_e$ term contributes a similar quantity. The key is that $\d_\e^\b \tilde v_{e0}$ and $\d_\e^\b \tilde w_e$ are bounded, uniformly in $(\e,t,x,y),$ for $t \leq t_\star(\e),$ since we are considering $\e$-derivatives in $x.$ Thus we obtain
$$ \Big| \frac{1}{\sqrt \e} \Re e \, \big( \d_\e^\a \op_\e(J)\check v, \d_\e^\a \check v \big)_{L^2} \Big| \lesssim \| \check v\|_{\e,s}^2,$$
the implicit constant depending on the WKB solution and $s.$

Next onto the singular term in $\check B_{(-1)}$ in the right-hand side of \eqref{eq:checkv0}. Every entry of the symbol $\check B_{(-1)}$ has the form $e^{ip\theta} g_p \underline B,$ where $p \in \{-1,1\}$ and $\underline B = \underline B(\xi,\eta)$ is independent of $(t,x,y)$ and bounded. Thus, by Lemma \ref{lem:tensor},
$$ |\Re e \, \big( \op_\e(\check B_{(-1)}) \check v, \check v)|_{L^2} \leq \Gamma \| \check v \|_{L^2}^2,$$
with \be \label{def:Gamma} \Gamma := |g_\pm |_{L^\infty} \| \check B_{(-1)} \|_{L^\infty}.\ee

We now consider the remainder term $\check F_{(-1)},$ introduced in Proposition \ref{prop:normal:form:2}, with the goal of proving a bound of the form \eqref{bd:dotF:Hs}. That Sobolev bound for $\dot F$ involved a symmetry argument for the order-one convective term. The key is that multiplication by $(\Id + \sqrt \e \op_\e(Q_{(-1)}))^{-1}$ to the left does not harm the symmetry argument, since $Q_{(-1)}$ is order $-1.$ 

Precisely, consider the contribution of $\dot F$ in $\check F_{(-1)},$ that is $$\Re e \, \big( \dot F, \check v\big)_{\e,s} = \Re e \, \big( \dot F, (\Id + \sqrt \e \op_\e(Q_{(-1)}))^{-1} \dot u \big)_{\e,s}.$$
For the leading term, we can use \eqref{bd:dotF:Hs} and Lemma \ref{lem:dot-checkv} and obtain
$$ \big| \Re e \, \big( \dot F, \dot u \big)_{\e,s} \big| \lesssim \e^{-1/2 + K'} \| \check v \|_{\e,s}^2 + \e^{K_a} \| \check v \|_{\e,s}.$$  
The other term is 
 $$ \Re e \, \big( \dot F, \big( (\Id + \sqrt \e \op_\e(Q_{(-1)}))^{-1} - \Id \big) \dot u \big)_{\e,s}.$$ 
 We can expand the operator $(\Id + \sqrt \e \op_\e(Q_{(-1)}))^{-1} - \Id$ in Neumann series, as in the proof of Lemma \ref{lem:dot-checkv}, to find the bound 
\be \label{neumann:0}
  \| \big( \, (\Id + \sqrt \e \op_\e(Q_{(-1)}))^{-1} - \Id\, \big) f \|_{\e,s} \lesssim \sqrt \e \| f \|_{\e,s-1}, \quad \mbox{for all $f \in H^{s-1}.$}
 \ee
 Above in \eqref{neumann:0} we used the fact that $Q_{(-1)}$ is order $-1.$ The $L^2$ adjoint of the operator $(\Id + \sqrt \e \op_\e(Q_{(-1)}))^{-1} - \Id$ satisfies the same bound. Thus %
 \ba \big| \Re e \, \big( \dot F, \big( (\Id + \sqrt \e \op_\e(Q_{(-1)}))^{-1} - \Id \big) \dot u \big)_{\e,s} \big| & \lesssim \| \dot F \|_{\e,s-1} \| \dot u \|_{\e,s} \\ & \lesssim \e^{K'} \| \check v \|_{\e,s}^2 + \e^{K_a + 1/2}  \| \check v \|_{\e,s},\ea
 where we used \eqref{bd:F:sob} and  Lemma \ref{lem:dot-checkv}.

 It remains to handle the contribution of $(\Id + \sqrt \e \op_\e(Q_{(-1)}))^{-1} - \Id \big) \dot F$ to the scalar product involving $\check F_{(-1)}.$ For that term, we use \eqref{neumann:0} and \eqref{bd:F:sob} again. Summing up, we find 
$$ \big| \Re e \, \big( \check F_{(-1)}, \check v\big)_{\e,s}\big| \lesssim \e^{-1/2 + K'} \| \check v\|_{\e,s}^2 + \e^{K_a} \| \check v \|_{\e,s},$$
for $t \leq t_\star(\e).$

At this point we have bounds in $s$ norms for all terms except the leading singular term $\check B_{(-1)},$ for which we have an $L^2$ bound. This gives the $L^2$ estimate 
$$  \frac{1}{2} \d_t \| \check v(t) \|_{L^2}^2 \leq \e^{-1/2} (\Gamma + C \e^{K'}) \| \check v \|_{L^2}^2 + C \e^{K_a} \| \check v \|_{L^2},$$
where $C > 0$ is independent of $\e,$ hence, by Gronwall's lemma and Lemma \ref{lem:dot-checkv}, the rough bound \eqref{weak:Sobolev:bound} in the case $s = 0.$  

Assuming now that the bound \eqref{weak:Sobolev:bound} holds true for some integer $s,$ we endeavour to bound $\| \check v \|_{\e,s+1}.$ The focus is on the commutator between $\d_\e^\a,$ with $|\a| = s + 1,$ and $\op_\e(\check B_{(-1)}),$ since
$$ \| \d_\e^{\a} \big( \op_\e(\check B_{(-1)}) \check v \big) \|_{L^2} \leq \Gamma \| \d_\e^{\a} \check v\|_{L^2} + \big\| [\d_\e^\a, \op_\e(\check B_{(-1)}) ] \check v \|_{L^2}.$$   
We have
$$ \big\| [\d_\e^\a, \op_\e(\check B_{(-1)}) ] \check v \|_{L^2} \leq C_s \sum_{\begin{smallmatrix} \a_1 + \a_2 = \a \\ |\a_1| > 0 \end{smallmatrix} } \| \op_\e(\d_\e^{\a_1} \check B_{(-1)}) \|_{L^2 \to L^2} \| \d_\e^{\a_2} \check v\|_{L^2},$$
 and we may use the induction hypothesis, since $|\a_2| \leq s.$ Besides, every entry of $\d_\e^{\a_1} \check B_{(-1})$ has the form $\d_\e^{\a_1} (e^{i p \theta} g_p) \underline B,$ so that the $L^2 \to L^2$ norm of $\op_\e(\d_\e^{\a_1} \check B_{(-1)})$ is bounded independently of $\e.$ Summing up, we found 
 $$ \| \d_\e^{\a} \big( \op_\e(\check B_{(-1)}) \check v \big) \|_{L^2} \leq \Gamma \| \d_\e^{\a} \check v\|_{L^2} + C \e^K e^{t \Gamma/\sqrt \e},$$
 where $C > 0$ depends on $s,$ not on $\e.$ Together with the Sobolev bounds for the other terms, above, we find 
 $$  \frac{1}{2} \d_t \| \check v(t) \|_{\e,s}^2 \leq \e^{-1/2} (\Gamma + C \e^{K'}) \| \check v \|_{\e,s}^2 + C \e^{K_a} \| \check v \|_{\e,s} + C \e^K e^{t \G/\sqrt \e} \| \check v \|_{\e,s},$$
  and another application of Gronwall's lemma completes the induction and the proof.
 \end{proof}

\begin{rem} It's in the proof of high-frequency estimates (see Section {\rm \ref{sec:high}}) that the weak bound of Proposition {\rm \ref{prop:weak:bound}} will be needed. The need for a weak bound will then come from the necessity to control the order-one term $\op_\e(\check J)$ which enters the system after the first normal form reduction in Section {\rm \ref{sec:normal}}.
\end{rem}

\section{The prepared system} \label{sec:prepared}

 Our starting point is the perturbation equations {\it in normal form}, which we derived in \eqref{eq:checku0} and \eqref{eq:checkv0}.
We perform more changes of unknowns in this Section, as we introduce a partition of the phase space $(x,y,\xi,\eta) \in \R^3 \times \R^3,$ and associated unknowns $u_{in}, u_{out}$ and $u_{high},$ defined in \eqref{def:U:in:out:high} below in term of $\check u$ and $\check v.$ 

 Here ``out" stands for {\it out in space}: the unknown $u_{out}$ is supported in the spatial far-field domain, in particular where the effect of the linearized source ${\mathcal B}(u_a)$ is not much felt, since the amplitude in $u_a$ decays to 0 at spatial infinity (due to its high Sobolev regularity).

 In $u_{high},$ ``high" stands for {\it high-frequency:} the unknown $u_{high}$ is supported far from the (bounded) resonant set.

 In particular, it's around the support of $u_{in}$ that the amplification takes place (see Sections \ref{sec:Duh} and \ref{sec:endgame}).

\subsection{The {\it in/out/high} variables} \label{sec:in:out:high}

Let $\chi_{long}$ be a smooth spatial cut-off in the $x$ variable, with
$$0 \leq \chi_{long}(x) \leq 1, \qquad \chi_{long}(x) \equiv 1 \quad \mbox{for $|x| \leq R_{long}$},$$
where $R_{long} > 0$ will ultimately chosen to be very large, in terms of $K$ and $K',$ and depending on the spatial decay of the initial WKB profile $a$ (which is assumed to have a high Sobolev regularity, hence decays to 0 at spatial infinity). Precisely, we will need $R_{long}$ to be large enough so that the associated rate of growth $\g_{out}$ derived in Section \ref{sec:out} is smaller than the optimal rate of growth derived in Sections \ref{sec:flow} and \ref{sec:interaction:coefficients}.

Let also $\chi_{low}$ be a smooth frequency cut-off, such that
\be \label{def:chilow}
0 \leq \chi_{low} \leq 1, \qquad \chi_{low}(\xi,\eta) \equiv 1 \quad \mbox{on a large enough neighborhood of ${\mathcal R}.$}
\ee
It is specified in the proof of Proposition \ref{prop:prepared} below how large we need the support of $\chi_{low}$ to be.

We introduce the unknowns %
\be \label{def:U:in:out:high} \begin{aligned} u_{in} & := \op_\e(\chi_{low}) (\chi_{long} \check u), \\ u_{out} & := \op_\e(\chi_{low}) ((1 - \chi_{long}) \check v), \\ u_{high} & := (1 - \op_\e(\chi_{low})) \check u,
\end{aligned}\ee
where $\check u$ is defined in Proposition \ref{prop:normal:form} and $\check v$ is defined in Proposition \ref{prop:normal:form:2}.

\begin{lem} \label{lem:reconstruct} Estimates for $u_{in},$ $u_{out}$ and $u_{high}$ imply estimates for $\check u:$ 
 $$ \begin{aligned} \| \check u \|_{\e,s} + \| \nabla_\e^\a \check u \|_{\e,s} & \leq (1 + \sqrt \e C_1) \big( \| u_{in} \|_{\e,s} + \| u_{out} \|_{\e,s} + \| \nabla_\e^\a u_{high} \|_{\e,s} \big) , \\
  \| \check u \|_{\fl1} + \| \nabla_\e^\a \check u \|_{\fl1} & \leq (1 + \sqrt \e C_2) \big( \| u_{in} \|_{\fl1} + \| u_{out} \|_{\fl1} + \| \nabla_\e^\a u_{high} \|_{\fl1}\big), \end{aligned}$$
 with $\a =0$ and $|\a| = 1,$ and constants $C_1 > 0$ and $C_2 > 0$ which are independent of $\e,t.$ 
\end{lem} 

\begin{proof} By \eqref{def:U:in:out:high}, 
we have 
\ba \label{dec:ucheck}  \check u & = u_{high} + u_{in} + u_{out} + \op_\e(\chi_{low}) \big( (1 - \chi_{long}) (\check u - \check v) \, \big),
\ea
so that
$$ (1 +  r_{out}) \check u =  u_{high} + u_{in} + u_{out},$$
with notation
$$ \begin{aligned} r_{out}  =  \op_\e(\chi_{low}) \circ (1 - \chi_{long}) \circ \big( ({\rm Id} + \sqrt \e \op_\e(Q_{(-1)}))^{-1} ({\rm Id} + \sqrt \e \op_\e(Q)) - \Id \,\big). \end{aligned}$$
Here we used
$$
\check v = (\Id + \sqrt \e \op_\e(Q_{(-1)}))^{-1} \dot u =  (\Id + \sqrt \e \op_\e(Q_{(-1)}))^{-1}   (\Id + \sqrt \e \op_\e(Q)) \check u. 
$$

We note that $\nabla_\e^\a \op_\e(\chi_{low})$ is bounded, from $H^s$ to itself, and from $\fl1$ to itself, independently of $\e:$
\be \label{bd:chi:low} \| \nabla_\e^\a \op_\e(\chi_{low}) f \|_{\e,s} \lesssim \| f \|_{\e,s}, \quad \mbox{for all $f \in H^s,$}\ee
and similarly in $\fl1.$  The operator $\e^{-1/2} ({\rm Id} + \sqrt \e \op_\e(Q_{(-1)}))^{-1} ({\rm Id} + \sqrt \e \op_\e(Q)) - \Id \,\big)$ is similarly $\| \cdot \|_{\e,s} \to \| \cdot \|_{\e,s}$ and $\| \cdot \|_{\fl1} \to \| \cdot \|_{\fl1}$ bounded, uniformly in $\e.$ Thus
$$ \| \nabla_\e^\a r_{low} \check u \|_{\e,s} \lesssim \e^{1/2} \| \check u \|_{\e,s}, \quad \| \nabla_\e^\a r_{low} \check u \|_{\fl1} \lesssim \e^{1/2} \| \check u \|_{\fl1}.$$
It now suffices to use \eqref{dec:ucheck} and \eqref{bd:chi:low}, the latter implying $\| \nabla_\e^\a u_{in} \|_{\e,s} \lesssim \| u_{in} \|_{\e,s},$ and the same for $u_{out},$ also in $\fl1$ norms. 
\end{proof}

Choosing the support of the initial amplitude $\phi$ in the initial perturbation \eqref{initial:datum} to satisfy
\be \label{cond:phi:datum}
\mbox{supp}\, \phi \subset \{ \chi_{long} \equiv 1 \},
\ee
we find that the datum for $u_{in}$ satisfies
\be \label{datum:uin}
\| u_{in}(0,x,y) - \e^K \op_\e(\chi_{low}) \phi^\e \|_{\e,s} \lesssim \e^{K + 1/2}.
\ee

\subsection{The prepared system} \label{sec:prepared:prepared}

\begin{prop} \label{prop:prepared} For an appropriate choice of $\chi_{low}$ and $\chi_{high},$ the unknowns $(u_{in}, u_{out}, u_{high})$ defined in \eqref{def:U:in:out:high} solve
$$
 \left\{\begin{aligned}
 \d_t u_{in} + \frac{1}{\e} A(\nabla_\e) u_{in} & = \frac{1}{\sqrt \e} \op_\e(\chi_{long}^\sharp\check B) u_{in} + F_{in}, \\
 \d_t u_{out} + \frac{1}{\e} A(\nabla_\e) u_{out}  & = \frac{1}{\sqrt \e} \op_\e((1 - \chi_{long}^\flat) \check B_{(-1)}) u_{out} + F_{out},\\
 \d_t u_{high} + \frac{1}{\e} A(\nabla_\e) u_{high} + \frac{1}{\sqrt \e} \op_\e(J) u_{high} & =  \op_\e(\check J) u_{high} + F_{high},
\end{aligned}\right.
$$
where we used convention \eqref{sharp:tilde} for the truncation functions. The remainders $F_{in}$ and $F_{high}$ are equivalent to $\check F,$ and $F_{out}$ is equivalent to $\check F_{(-1)},$ in the sense of Definition {\rm \ref{def:equiv:remainder}}.
\end{prop}

\begin{proof} In a first step, we focus on the $A(\nabla_\e)$ terms in the systems in $u_{in},$ $u_{out}$ and $u_{high}.$ We repeatedly use
$$ \big[ \op_\e(\chi_\star),  A(\nabla_\e) \big] = 0,$$
with $\star \in \{ low, high\},$ which holds true since $\op_\e(\chi_\star)$ are scalar Fourier multipliers, and $A(\nabla_\e)$ is a constant-coefficient differential operator. Thus
$$ \frac{1}{\e} \op_\e(\chi_{low}) \big(\chi_{long} A(\nabla_\e) \check u \, \big) = \frac{1}{\e} A(\nabla_\e) u_{in} -  \op_\e(\chi_{low}) (A(\d_x,0) \chi_{long}\big) \, \check u,$$
since $\chi_{long}$ depends only on $x.$ The second term in the above right-hand side satisfies
the remainder bounds
$$ \begin{aligned} \| \op_\e(\chi_{low}) (A(\d_x,0) \chi_{long}\big) \, \check u  \|_{\e,s} & \lesssim \| \check u \|_{\e,s}, \\ \| \d_\e^\a \op_\e(\chi_{low}) (A(\d_x,0) \chi_{long}\big) \, \check u  \|_{\fl1} & \lesssim \max_{|\b| \leq |\a|} \| \check u \|_{\fl1},\end{aligned}$$
so that it is a {\it uniform remainder} in the sense of Definition \ref{defi:uniform:remainder}. 
Similarly,
$$ \frac{1}{\e} \op_\e(\chi_{low}) \big( \, (1 - \chi_{long}) A(\nabla_\e) \check v \, \big) = \frac{1}{\e} A(\nabla_\e) u_{out} + R \check v,$$
where $R$ is a uniform remainder.

 Now onto the contribution of $\check B$ to the prepared system. Since $\chi_{long}$ is independent of the frequency variables, we have
 $$ \chi_{long} \op_\e(\check B) \equiv \op_\e(\chi_{long} \check B),$$
 and, with convention \eqref{sharp:tilde} bearing on the cut-offs,
 $$ \op_\e(\chi_{long} \check B) \equiv \op_\e(\chi_{long} \chi_{long}^\sharp \check B),$$
 and
 $$ \op_\e(\chi_{long} \chi_{long}^\sharp \check B) = \op_\e(\chi_{long}^\sharp \check B) \big( \chi_{long} \cdot \, \big) + \e R,$$
 where $R$ is a uniform remainder. Next, as seen in \eqref{def:checkB}, the term $\chi_{long}^\sharp \check B$  is a sum of terms of the form $e^{i \theta} g_1(t,x,y) \chi_{jj'}  B_+$ and $e^{- i \theta} g_{-1}(t,x,y) \chi_{jj',-1}  B_-,$ where the Fourier multipliers $B_{\pm 1}$ are order zero. We compute
 $$ \op_\e(\chi_{low}) \big( \, e^{i \theta} \op_\e(g_1 \chi_{jj'}  B_+) \, \big) = e^{i \theta} \op_\e(\chi_{low,+1} g_1 \chi_{jj'}  B_+) + \sqrt \e R,$$
 where $R$ is a uniform remainder. Note that we have $\sqrt \e R$ here, as opposed to $\e R$ above, because $g_1$ does depend on $y.$ We may choose the support of $\chi_{low}$ to be a large enough neighborhood of the resonant set ${\mathcal R}$ so that
 $$ \chi_{low,+1} \chi_{jj'} \equiv \chi_{jj'},$$
 for all $j,j'.$ For the terms corresponding to $p = -1,$ the condition is $\chi_{low} \chi_{jj'} \equiv \chi_{jj'},$ for all $j,j'.$
Under these conditions bearing on the support of $\chi_{low},$ we have
$$ \op_\e(\chi_{low}) \op_\e(\chi_{long}^\sharp \check B)  = \op_\e(\chi_{long}^\sharp \check B) + \sqrt \e R = \op_\e(\chi_{long}^\sharp \check B) \op_\e(\chi_{low}) + \sqrt \e R,$$
where $R$ is a uniform remainder. 

There is one more thing to verify for the equation in $u_{in},$ namely that the $J$ and $\check J$ terms do not feature in it. This can be ensured by an adequate choice of the truncation $\chi_{high},$ introduced in \eqref{def:chiHF}. Observe indeed that $J$ is a sum of terms of the form $e^{i p_a \theta} \chi_{high} J_1,$ where $p_a$ is a WKB harmonics and $J_1 \in S^1.$ Associated with $K$ and $K',$ we have a WKB index $K_a$ which describes the order of the WKB approximation. For a given order of the WKB expansion, we have a finite number of harmonics. So for given $K$ and $K',$ the WKB harmonics $p_a$ range in a finite subset of $\Z.$ We choose the support of $\chi_{high}$
to be so high-frequency that
\be \label{low:high:cancellation} \chi_{low, +p_a} \chi_{high} \equiv 0,\ee
for all WKB harmonics $p_a.$ Then,
\be \label{low:J} \op_\e(\chi_{low}) \circ \chi_{long} \circ \op_\e(J) = \sqrt \e R,
\ee
where $R$ is a uniform remainder.
We have similarly
\be \label{low:checkJ}
\op_\e(\chi_{low}) \circ \chi_{long} \circ \op_\e(\check J) = \sqrt \e R,
\ee
where $R$ is a uniform remainder, if the supports of $\chi_{low}$ and $\chi_{high}$ are well separated. Indeed, by definition of $\check J$ in Proposition \ref{prop:normal:form}, 
$$ \op_\e(\chi_{low}) \op_\e(\check J) = - ({\rm Id} + \sqrt \e \op_\e(Q))^{-1} \op_\e(\chi_{low}) [\op_\e(J), \op_\e(Q)] + \sqrt \e R.$$
In view of \eqref{low:J}, we only have to consider $\op_\e(\chi_{low}) \op_\e(Q) \op_\e(J).$ By definition of $Q$ in the proof of Proposition \ref{prop:normal:form}, the composition $\op_\e(\chi_{low}) \op_\e(Q)$ is a sum of terms of the form $e^{i p \theta} \op_\e(\chi_{low,+p} Q_p),$ with $Q_p \in S^0$ and $|p| = 1.$ 
Thus if $\chi_{low,+p+p_a} \chi_{high} = 0,$ with $|p| = 1$ and $p_a$ any WKB harmonics, then 
$$ \op_\e(\chi_{low}) \op_\e(Q) \op_\e(J) = \sqrt \e R,$$
which implies 
$$ \op_\e(\chi_{low}) \op_\e(\check J) = \sqrt \e R.$$
We then handle the extra truncation $\chi_{long}$ for instance by writing 
$$\op_\e(\chi_{low}) \circ \chi_{long} = \op_\e(\chi_{long} \chi_{low}) + \e R,$$ where $R$ is a uniform remainder. This implies \eqref{low:checkJ}.

At this point we verified that the equation in $u_{in}$ is as stated in the Proposition.

For the equation in $u_{out},$ the same arguments give the source term $\op_\e((1 - \chi_{long})^\sharp \check B_{(-1)}) u_{out},$ and we observe that we may write $(1 - \chi_{long})^\sharp = 1 - \chi_{long}^\flat.$

\medskip

For the equation in $u_{high},$ we observe that, with the above choice for $\chi_{low},$ we have
$\op_\e(1 - \chi_{low}) \op_\e(\check B) = \sqrt \e R,$
where $R$ is a uniform remainder. Besides, by the arguments used to prove \eqref{low:J} and \eqref{low:checkJ}, we find that $\op_\e(1 - \chi_{low}) \op_\e(J) = \op_\e(J) + \sqrt \e R,$ where $R$ is a uniform remainder, and an analogous equality with $\check J.$ 
\end{proof}

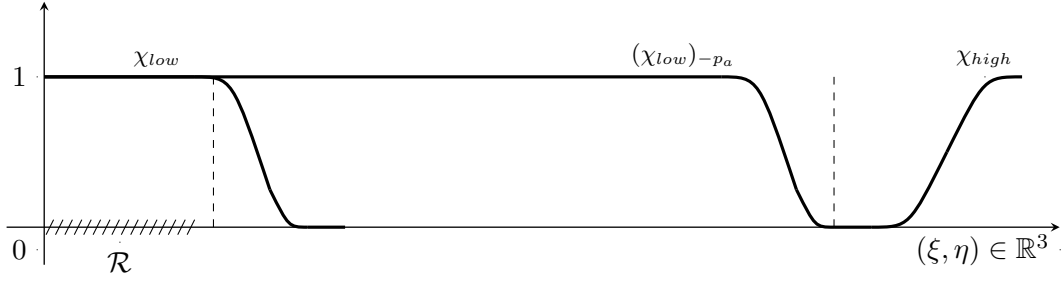
\begin{figure}

 \begin{tikzpicture}

 \begin{scope}[>=stealth]

  \draw[line width=.5pt][->] (-7,-2) -- (7,-2);
  \draw[line width=.5pt][->] (-6.5,-2.5) -- (-6.5,1);
  \draw (-6.6,0) -- (-6.6,0) node[anchor=east] {{$1$}} ;
  \draw (-6.6,-2.3) -- (-6.6,-2.3) node[anchor=east] {$0$} ;
  \draw (7,-2.3) -- (7,-2.3) node[anchor=east] {$(\xi,\eta) \in {\mathbb R}^3$} ;

  \draw (-5,0) -- (-5,0) node[anchor=south] {{\footnotesize $\chi_{low}$}} ;
  \draw[very thick] (-6.5,0) -- (-4.5,0) ;
  \draw[very thick] (-4.5,0) .. controls (-4,0) .. (-3.5,-1.5) ;
  \draw[very thick] (-3.5,-1.5) .. controls (-3.25,-2) .. (-3,-2) ;
  \draw[very thick] (-3,-2) -- (-2.5, -2) ;
  \draw[dashed] (-4.25,-2) -- (-4.25,0) ;

   \draw (2,0) -- (2,0) node[anchor=south] {{\footnotesize $(\chi_{low})_{-p_a}$}} ;
     \draw[very thick] (-6.5,0) -- (2.5,0) ;
     \draw[very thick] (2.5,0) .. controls (3,0) .. (3.5,-1.5) ;
     \draw[very thick] (3.5,-1.5) .. controls (3.75,-2) .. (4,-2) ;
     \draw[very thick] (4,-2) -- (4.5, -2) ;
  \draw[dashed] (4,-2) -- (4,0) ;

   \draw (6,0) -- (6,0) node[anchor=south] {{\footnotesize $\chi_{high}$}} ;
    \draw[very thick] (5.5,-1) .. controls (6,0) .. (6.5,0) ;
        \draw[very thick] (4.5,-2) .. controls (5,-2) .. (5.5,-1) ;

   \foreach \i in {1,...,16}
{
\draw (-6.6 + \i/8,-2.1)-- (-6.5 + \i/8, -1.9);
}

 \draw (-5.5,-2.2) -- (-5.5,-2.2) node[anchor=north] {${\footnotesize {\mathcal R}}$} ;

  \end{scope}

 \end{tikzpicture}

 \caption{The resonant set ${\mathcal R}$ is bounded (Proposition \ref{prop:res}). The frequency cut-off $\chi_{low},$ defined in \eqref{def:chilow} is identically equal to 1 on a large enough neighborhood of ${\mathcal R}.$ How large this neighborhood must be is specified in the proof of Proposition \ref{prop:prepared}. In order to allow for $K$ large and $K'$ small, we need the WKB approximation to be very precise, hence the largest harmonics $p_a > 0$ in the WKB expansion to be very large. Thus the support of the shifted cut-off $(\chi_{low})_{-p_a}$ extends far beyond the support of $\chi_{low}.$ The high-frequency cut-off $\chi_{high}$ is chosen so that $(\chi_{low})_{-p_a} \chi_{high} \equiv 0.$}
 \label{fig:r:chi}

 \end{figure}

\begin{rem} \label{trunc:y} We see in the proof of Proposition {\rm \ref{prop:prepared}} that if the spatial truncation $\chi_{long}$ depended on $y,$ then the prepared system would contain extra singular terms, of the form
 $\dsp{\frac{1}{\sqrt \e} \op_\e(\chi_{low}) \big( A(0,\d_y) \chi_{long} \big) \check u.}$
 That is why we do not truncate in $y$ in \eqref{def:U:in:out:high}.
 \end{rem}

\section{High-frequency estimates} \label{sec:high}

The system in $u_{high}$ defined in \eqref{def:U:in:out:high} is, according to Proposition \ref{prop:prepared}, 
\be \label{eq:u:high} \d_t u_{high} + \frac{1}{\e} A(\nabla_\e) u_{high} + \frac{1}{\sqrt \e} \op_\e(J) u_{high}  =  \op_\e(\check J) u_{high} + F_{high}.
\ee

We obtain here an upper bound for $u_{high},$ in spite of the first-order, non-symmetric term $\op_\e(\check J)$ in the above right-hand side. The term $\check J$ comes from the interplay between the normal form symbol $Q$ and the convection $J.$ It is order one and not symmetric. Hence it causes a loss of one derivative in the Sobolev estimates. By exploiting the smallness of the time window under consideration, the smallness of the datum for $u_{high},$ and the rough Sobolev bound derived in Section \ref{sec:weak:Sobolev:bound}, we are able to derive a closed estimate for Sobolev norms of $u_{high},$ in terms of $u_{in}$ and $u_{out}.$ 

\begin{lem} \label{lem:datum:high} For an appropriate choice of $\chi_{low}$ we have for small enough $\e$ the bound, for all $\a:$ 
 $$ \| u_{high}(0) \|_{\e,s} + \| \d_\e^\a u_{high}(0) \|_{\fl1}  \lesssim \e^{K_a}.
 $$
\end{lem}

\begin{proof} The datum for $u_{high}$ is
 $$ u_{high}(0) = \e^K \op_\e(1 - \chi_{low}) (1 + \sqrt \e \op_\e(Q))^{-1} \phi^\e,$$
 with $\phi^\e$ defined in \eqref{def:phie}.
 We may expand $(1 + \sqrt \e \op_\e(Q))^{-1}$ in Neumann series. Thus we consider terms of the form
$\e^{K + n/2} \op_\e(1 - \chi_{low}) \op_\e(Q)^n \big( e^{i x \xi_0/\e} \phi^\e).$ Now recall that every entry of $Q$ has the form $e^{i p \theta} g_p Q_p,$ where $Q_p$ is an order-zero Fourier multiplier, and $|p|  =1.$

Thus the composition of $\op_\e(1 - \chi_{low})$ with $\op_\e(Q)^n$ results in shifts to $1 - \chi_{low}$ by multiples of $p k/\e.$ If we choose the support of $\chi_{low}$ to be large enough, then the first $n$ terms in the Neumann expansion of $u_{high}(0)$ have the form $\op_\e(\tilde \chi) \phi,$ where $\tilde \chi$ is supported far from the zero frequency. Since $\phi$ is compactly supported, its Fourier transform $\hat \phi$ belongs to the Schwartz class, hence $\hat \phi$ evaluated at frequencies $|\zeta| \geq c/\sqrt \e,$ for some $c > 0,$ is of order $\e^\infty,$ in particular smaller than $\e^{K_a}.$
 \end{proof}

\subsection{High-frequency Sobolev estimates}  \label{sec:high:sob}

In the proposition below, we use the Sobolev index $s$ defined in \eqref{def:s}. In particular, as noted in Remark \ref{rem:Sobolev:indices}, the index $s$ is very large in the large $K$ limit, and in the small $K'$ limit. 

\begin{prop} \label{prop:high}  For $t \leq t_\star(\e) \leq T \sqrt \e |\ln \e|,$ for all $s_1$ and $s_2$ with $s_1 + s_2 \leq s,$ we have for $u_{high}$,  $u_{in}$ and  $u_{out}$ introduced in \eqref{def:U:in:out:high} the bound 
$$ \begin{aligned} \| u_{high}(t)\|^2_{\e,s_1} &  \lesssim \e^{2 K_a} + \e^{2 K + s_2} |\ln \e|^\star e^{t (2 \Gamma + \g_{high})/\sqrt \e} \\ & +   \e^{K'} |\ln \e|^\star e^{t \g_{high}/\sqrt \e} \max_{ 0 \leq t' \leq t} (\| u_{in}(t') \|^2_{\e,s_1} + \| u_{out}(t') \|^2_{\e,s_1}),\end{aligned}
$$%
where $\g_{high} > 0$ is arbitrarily small. %
 The implicit constant is independent of $\e,t,$ but depends on $\g_{high}.$ The rate $\Gamma > 0$ is the rough growth rate from Proposition {\rm \ref{prop:weak:bound}.}
\end{prop}

\begin{proof} In a first step, we prove 
\be \label{loss:estimate} \begin{aligned} 
 \| u_{high}(t)  \|_{\e,s_1}^2 & \lesssim \e^{2 K_a} e^{t \g_{high}/\sqrt\e} \\ & \quad + \e^{-1/2 + K'} \int_0^t e^{(t - t') \g_{high}/\sqrt \e} \max_{0 \leq t'' \leq t'} (\| u_{in}(t'') \|_{\e,s_1}^2 + \| u_{out}(t'') \|_{\e,s_1}^2 \big) \, dt' \\ & \quad  +  \e^{1/2}  \int_0^t  e^{(t - t') \g_{high}/\sqrt\e} \max_{0 \leq t'' \leq t'} \| u_{high}(t'') \|_{\e,s_1 + 1}^2 \, dt', \end{aligned}
 \ee 
where $\g_{high} > 0$ is arbitrarily small (and the implicit constant depends on $\g_{high}$). 

 The terms in the left-hand side of the system in \eqref{eq:u:high} are identical to the left-hand side in the system \eqref{eq:checkv0} in $\check v.$ Thus, following the proof of Proposition \ref{prop:weak:bound}, we have
$$\begin{aligned} \frac{1}{2} \d_t \| u_{high} \|_{\e,s_1}^2& \lesssim  \| u_{high} \|_{\e,s_1}^2 + \Re e \, \big( \op_\e(\check J) u_{high} + F_{high}, u_{high} \big)_{\e,s_1}.\end{aligned}$$
We first examine the contribution of the $\check J$ term, which was defined in Proposition \ref{prop:normal:form}:
$$\op_\e(\check J) = - (\Id + \sqrt \e \op_\e(Q))^{-1} [\op_\e(J + R_A), \op_\e(Q)].$$
We saw in the proof of Lemma \ref{lem:dot-checkv} that $(\Id + \sqrt \e \op_\e(Q))^{-1}$ is bounded from $\| \cdot \|_{\e,s_1}$ to itself, uniformly in $\e.$ Thus we focus on the commutator $[\op_\e(J + R_A), \op_\e(Q)].$ It follows from the definition of the convective term $J$ in \eqref{def:J} and the remainder $R_A$ in Section \ref{spectral:dec:symbols} that 
$$ \| \op_\e(J + R_A) f \|_{\e,s_1} \lesssim \| f \|_{\e,s_1+1}.$$
Besides, $\op_\e(Q)$ is a uniform remainder, in the sense of Definition \ref{def:equiv:remainder}. Thus 
$$ \| [\op_\e(J + R_A), \op_\e(Q)] f \|_{\e,s_1} \lesssim \| f \|_{\e,s_1 + 1},$$
implying 
 \be \label{est:checkJ} \Re e \, \big( \op_\e(\check J) u_{high}, u_{high} \big)_{\e,s_1} \leq C_{\check J} \| u_{high} \|_{\e,s_1+1} \| u_{high} \|_{\e,s_1},\ee
 for some $C_{\check J} > 0.$

We handle the source term $F_{high}$ very much like we took care of $\check F_{(-1)}$ in the proof of Proposition \ref{prop:weak:bound}. We note indeed that $F_{high}$ is equivalent to $\check F,$ as indicated in Proposition \ref{prop:prepared} ({\it equivalent} in the sense of Definition \ref{def:equiv:remainder}). Precisely, in the scalar product $\Re e \, (F_{high}, u_{high})_{\e,s},$ we find the term $\Re e \, (\op_\e(1-\chi_{low}) \dot F, u_{high})_{\e,s},$ in which the leading term in $\e$ is $\Re e \, (\op_\e(1 - \chi_{low}) \dot F, \op_\e(1 - \chi_{low}) \dot u)_{\e,s}.$ For this term, we can use the estimate \eqref{bd:dotF:Hs}, since the scalar operators $\op_\e(1 - \chi_{low})$ do not destroy the symmetry in the nonlinear convection. All the other terms in the contribution of $F_{high}$ are $O(1)$ in $\e.$ For those terms, the symbol $Q \in S^0$ does disrupt the symmetry, and we have to use \eqref{bd:F:sob:new} instead of \eqref{bd:dotF:Hs}. This implies a loss of derivative, but with Lemma \ref{lem:reconstruct}, this loss occurs only for the {\it high} component of the solution. We find 
\be \label{Fhigh:est} |\Re e \, \big( F_{high}, u_{high} \big)_{\e,s_1}| \leq C_{high} \big( \e^{-1/2 + K'} \| \dot u\|_{\e,s_1}^2 + \e^{K_a} \| \dot u \|_{\e,s_1} + \e^{K'} \| u_{high} \|_{\e,s_1 + 1} \| u_{high} \|_{\e,s_1}\big),\ee
for some $C_{high} > 0.$   
Next we bound 
 $$ \begin{aligned} (C_{\check J} + \e^{K'} C_{high}) \| u_{high} \|_{\e,s_1+1} \| u_{high} \|_{\e,s_1} & \leq C'_{high} \e^{1/2} \| u_{high} \|_{\e,s_1 + 1}^2 \\ &  +  \e^{-1/2} \g_{high} \| u_{high} \|_{\e,s_1}^2,
 \end{aligned}$$
 for any $\g_{high} > 0,$ and some $C'_{high} > 0$ which depends on $C_{\check J},$ $C_{high}$ and $\g_{high}.$ 
Summing up, we obtained  
$$ \begin{aligned} \frac{1}{2} \d_t \| u_{high} \|^2_{\e,s_1} & \leq \e^{2 K_a} + C'_{high} \e^{1/2} \| u_{high} \|_{\e,s_1 + 1}^2 + \e^{-1/2 + K'} C''_{high} \| \dot u \|^2_{\e,s_1} \\ & + \e^{-1/2} \g_{high} \| u_{high} \|_{\e,s_1}^2, 
\end{aligned}$$
where $\g_{high} > 0$ is arbitrarily small, and $C'_{high} > 0$ and $C''_{high} > 0$ depend on $\g_{high}.$ With Lemma \ref{lem:reconstruct}, this implies
$$ \begin{aligned} \frac{1}{2} \d_t \| u_{high} \|^2_{\e,s_1} & \leq \e^{2 K_a} + C'_{high} \e^{1/2} \| u_{high} \|_{\e,s_1 + 1}^2 + \e^{-1/2 + K'} \big( \| u_{in} \|^2_{\e,s_1} + \| u_{out}\|^2_{\e,s_1}\big) \\ & + \e^{-1/2} \g_{high} \| u_{high} \|_{\e,s_1}^2, 
\end{aligned}$$
with a different $\g_{high} > 0$ that is still arbitrarily small.
 The bound \eqref{loss:estimate} then follows from Gronwall's lemma and Lemma \ref{lem:datum:high}.

 In a second step, we use the loss estimate \eqref{loss:estimate} iteratively, and the smallness of the time interval: at this point we use $t_\star(\e) \leq T \sqrt \e |\ln \e|.$ This gives %
$$ %
\begin{aligned} 
 \| u_{high}(t)  \|_{\e,s_1}^2  \lesssim \e^{2 K_a} & + \e^{K'} e^{t \g_{high}/\sqrt \e} \e |\ln\e| \max_{0 \leq t' \leq t}  (\| u_{in}(t') \|_{\e,s_1 + 1}^2 + \| u_{out}(t') \|_{\e,s_1 + 1}^2) \\ & \quad  +  \e^{3/2} |\ln \e| \int_0^t  e^{(t - t') \g_{high}/\sqrt\e} \max_{0 \leq t'' \leq t'} \| u_{high}(t'') \|_{\e,s_1 + 2}^2 \, dt'. \end{aligned}
 $$
 We observe that
 $$  \e \| u_{in}(t') \|_{\e,s_1 + 1}^2  \lesssim \e (\| u_{in}(t') \|_{\e,s_1}^2 + \| \d_\e u_{in}(t') \|_{\e,s_1}^2),$$
 and that
 $$ \e \| \d_\e u_{in}(t') \|_{\e,s_1}^2 \lesssim \| \nabla_\e u_{in}(t') \|_{\e,s_1}^2 \lesssim \| u_{in}(t') \|_{\e,s_1}^2,$$
 by definition of $\d_\e,$ $\nabla_\e$ \eqref{nablae} and $u_{in}$ \eqref{def:U:in:out:high}.
 Thus the estimate \eqref{loss:estimate} implies no loss for the {\it in} (and, similarly, for the {\it out}) components of the solution: we obtain 
 $$ %
\begin{aligned} 
 \| u_{high}(t)  \|_{\e,s_1}^2  \lesssim \e^{2 K_a} & + \e^{K'} e^{t \g_{high}/\sqrt \e} |\ln\e| \max_{0 \leq t' \leq t}  (\| u_{in}(t') \|_{\e,s_1}^2 + \| u_{out}(t') \|_{\e,s_1}^2) \\ & \quad  +  \e^{3/2} |\ln \e| \int_0^t  e^{(t - t') \g_{high}/\sqrt\e} \max_{0 \leq t'' \leq t'} \| u_{high}(t'') \|_{\e,s_1 + 2}^2 \, dt'. \end{aligned}
 $$
 Iterating, we obtain
  $$ %
\begin{aligned} 
 \| u_{high}(t)  \|_{\e,s_1}^2  \lesssim \e^{2 K_a} & + \e^{K'} e^{t \g_{high}/\sqrt \e} |\ln\e|^\star \max_{0 \leq t' \leq t}  (\| u_{in}(t') \|_{\e,s_1}^2 + \| u_{out}(t') \|_{\e,s_1}^2) \\ & \quad  +  \e^{s_2 - 1/2} |\ln \e|^\star \int_0^t  e^{(t - t') \g_{high}/\sqrt\e} \max_{0 \leq t'' \leq t'} \| u_{high}(t'') \|_{\e,s_1 + s_2}^2 \, dt'. \end{aligned}
 $$
 With the rough bound of Proposition \ref{prop:weak:bound}, we finally obtain the result.
\end{proof} 

\subsection{Estimates in the $\fl1$ norm} \label{sec:high:fl1}

Next we turn to estimates for $u_{high}$ in the ${\mathcal F}L^1$ norm. This norm is defined in \eqref{def:fl1}. We denote 
\be \label{fl1:m} \| \d_\e^m f \|_{\fl1} := \max_{0 \leq |\a| \leq m} \| \d_\e^\a f \|_{\fl1}, \qquad \mbox{for all $f \in H^{2  +m}(\R^3).$}\ee

\begin{prop} \label{prop:high:fl1} For $t \leq t_\star(\e) \leq T \sqrt \e |\ln \e|,$ we have
for $u_{high}$,  $u_{in}$ and  $u_{out}$ from \eqref{def:U:in:out:high} the bound
 \be \label{fl1:1.5}
\begin{aligned} \| u_{high} (t) \|_{\fl1} &  \lesssim \e^{K_a} + \e^{-1/2 + K'} \int_0^t (\| u_{in}(t')\|_{\fl1} + \| u_{out}(t') \|_{\fl1} ) \, dt' \\ & \quad  + \int_0^t \| \d_\e u_{high}(t') \|_{\fl1} \, dt',
\end{aligned}
\ee 
and, for $1 \leq m \leq s-2:$ 
 \be \label{fl1:2} \begin{aligned}
 \| \d_\e^m u_{high}(t) \|_{\fl1} & \lesssim \e^{K_a} + \e^{-1/2 + K'} \int_0^t (\| \d_\e^m u_{in}(t') \|_{\fl1} + \| \d_\e^m u_{out}(t') \|_{\fl1} ) \, dt' \\ & + \e^{-1/2 + K'} \int_0^t \| \d_\e^{m+1} u_{high} (t')\|_{\fl1} \, dt'.
 \end{aligned}\ee
\end{prop}

\begin{proof} The equation satisfied by $u_{high}$ is \eqref{eq:u:high}. The solution operator associated with the hyperbolic differential operator is defined as 
 $$ e^{t A(\nabla_\e)/\e}: \qquad f \to {\mathcal F}^{-1} \big(e^{t A(i \e \xi, i \sqrt \e \eta)/\e)} \hat f \big), \qquad t \in \R.$$ 
 It is bounded in the $\fl1$ norm, by strong hyperbolicity of $A:$ 
 \be \label{strong:hyperbolicity} \| e^{t A(\nabla_\e)/\e} f \|_{\fl1} \lesssim \| f \|_{\fl1}, \quad \mbox{for all $f \in \fl1.$}\ee
 The above bound holds since the eigenvalues of $A(i \zeta)$ are real and semi-simple, and the eigenprojectors are bounded (see Section \ref{sec:spectral:dec}).
  
 \medskip
 
{\it First step: proof of \eqref{fl1:1.5}.} We start with the description \eqref{for:convection:this:one} of the convective terms, which we reproduce here:
\be \label{linear:conv:cancellation} \frac{1}{\sqrt \e} \op_\e(J) =  (\tilde v_{e0} \cdot \d_\e + \tilde w_e \cdot \nabla_\e) \op_\e(\chi_{high} \underline J_e) + \sqrt \e (\tilde w_i \cdot \nabla_\e) \op_\e(\chi_{high}\underline J_i).\ee
The WKB maps $\tilde v_{e0}$ and $\tilde w_j$ and their $\d_\e^\a$ derivatives are bounded in $\fl1$ norm, uniformly in $\e.$ Indeed, we have 
 $$ \| e^{i p \theta} g_p\|_{\fl1} = \| \hat g_p(x - p k/\e,y)\|_{L^1} = \| g_p \|_{\fl1} \lesssim \| g_p\|_{H^{s_0}},$$
 for $s_0 > 3/2.$ Now, by virtue of the WKB expansion, the maps $\tilde v_{e0}$ and $\tilde w_j$ have the form of polynomials in $e^{i \theta}$ times profiles which are defined in terms of the solution to the Zakharov system. In particular, by propagation of Sobolev regularity by the Zakharov system \cite{LPS}, those profiles and their derivatives are bounded in $\fl1$ norm.
 
Since $\chi_{high} \underline J_e$ is a bounded Fourier multiplier, we have
 $$ \| \op_\e(\chi_{high} \underline J_e) f \|_{\fl1} \lesssim \| f \|_{\fl1}, \quad \mbox{for all $f \in \fl1,$}$$
Thus by the product law \eqref{algebra}, we deduce
 $$ \| (v_{e 0} \cdot \d_\e + w_e \cdot \nabla_\e) \op_\e(\chi_{high} \underline J_e) u_{high} \|_{\fl1} \lesssim \big(\| v_{e0} \|_{\fl1} + \| w_e\|_{\fl1}\big) \| \d_\e u_{high}\|_{\fl1},$$
 The same holds true for the $\tilde w_i$ and $\check J$ terms, so that 
  $$ \frac{1}{\sqrt \e} \| \op_\e(J) u_{high}\|_{\fl1} + \| \op_\e(\check J) u_{high} \|_{\fl1} \lesssim \| \d_\e u_{high} \|_{\fl1}.$$ 
  Note that here, crucially, the upper bound does not have a singular prefactor in $\e,$ thanks to the {\it linear convective transparency} (cancellation) observed in \eqref{linear:conv:cancellation}. 
  
  Thus, by \eqref{strong:hyperbolicity}, Lemma \ref{lem:datum:high} and \eqref{embed:loss}, we obtain
$$ \| u_{high}(t) \|_{\fl1} \leq \e^{K_a} + \int_0^t \| \d_\e u_{high} \|_{\fl1} + \| F_{high}(t') \|_{\fl1} \, dt'.$$
It remains to bound, in $\fl1$ norm, the source (given by Proposition \ref{prop:prepared}) $$F_{high} = (1 - \op_\e(\chi_{low})) \big( (1 + \sqrt \e \op_\e(Q))^{-1} \dot F + R \check u \big),$$
where $R$ is a uniform remainder in the sense of Definition \ref{defi:uniform:remainder}. Since $0 \leq 1 - \chi_{low}(\xi,\eta) \leq 1,$ we have  
$$\| \op_\e(1 - \chi_{low}) f\|_{\fl1} \leq \| f\|_{\fl1}, \qquad \mbox{for all $f \in \fl1.$}$$ We also have 
\be \label{Q:fl1} \| (\Id + \sqrt \e \op_\e(Q))^{-1} \|_{\fl1 \to \fl1} \lesssim 1,\ee
 by the arguments of the proof of Lemma \ref{lem:dot-checkv}. Indeed, we may write $(\Id + \sqrt \e \op_\e(Q))^{-1} f$ as a sum of terms $\e^{n/2} \op_\e(Q)^n f,$ and we use the fact that $Q$ is a tensor product of a function in $\fl1$ by a bounded Fourier multiplier. Thus
$$ \| \op_\e(Q) f \|_{\fl1} \lesssim \| f\|_{\fl1}, \qquad f \in \fl1,$$
 and for $\e$ small enough the Neumann series converges, proving \eqref{Q:fl1}. Summing up, we have 
 $$\| F_{high} \|_{\fl1} \lesssim \| \dot F \|_{\fl1} + \| \check u \|_{\fl1}.$$ Thus we proved %
 \be \label{fl1:1} 
 \| u_{high} (t) \|_{\fl1}  \lesssim \e^{K_a - 1/2} + \int_0^t (\| \d_\e u_{high}(t') \|_{\fl1} + \| \dot F(t') \|_{\fl1}  + \| \check u(t') \|_{\fl1} ) \, dt'. 
 \ee 
  We now use the bound \eqref{bd:dotF:fl1:apriori} for $\dot F,$ which we reproduce here: for $t \leq t_\star(\e):$ 
$$  \| \dot F(t) \|_{\fl1}  \lesssim \e^{-1/2 + K'} \| \dot u(t) \|_{\fl1} + \e^{K_a}.$$
We plug the above in \eqref{fl1:1} and, with Lemmas \ref{lem:dot-checkv} and \ref{lem:reconstruct} and an application of Gronwall's lemma, arrive at \eqref{fl1:1.5}, which ends the first step of this proof.

\medskip

{\it Second step: proof of \eqref{fl1:2}.} We apply $\d_\e^\a$ to the equation in $u_{high},$  and find
$$ \begin{aligned} \| \d_\e^\a & u_{high}(t) \|_{\fl1} \lesssim \e^{K_a - 1/2} \\ & +  \int_0^t (\e^{-1/2} \| \d_\e^\a \op_\e(J) u_{high} (t') \|_{\fl1} + \|\d_\e^\a \op_\e(\check J) u_{high}(t') \|_{\fl1}   + \| \d_\e^\a F_{high}(t') \|_{\fl1})  \, dt'.\end{aligned}$$
Using the same arguments as in the first step for the convective terms, and notation \eqref{fl1:m}, we find 
$$ \e^{-1/2}  \| \d_\e^\a  \op_\e(J) u_{high} ) \|_{\fl1} +  \|\d_\e^\a \op_\e(\check J) u_{high} \|_{\fl1} \lesssim \| \d^{|\a| + 1}_\e u_{high} \|_{\fl1}.$$
Besides, 
$$ \| \d_\e^\a F_{high} \|_{\fl1} \lesssim \| \d^\a_\e \check F \|_{\fl1} + \| \d_\e^{|\a|} \check u \|_{\fl1} \lesssim \| \d_\e^{|\a|} \dot F \|_{\fl1} + \| \d_\e^{|\a|} \check u \|_{\fl1}.$$
Now we appeal to \eqref{final:F:fl1}. That bound involves $\| \d_\e^{|\a|} \nabla_\e \dot u \|_{\fl1},$ which we want to relate to $\fl1$ norms of $u_{in},$ $u_{out}$ and $u_{high}.$ This norm involved a mixed differential operator is not exactly covered by Lemmas \ref{lem:dot-checkv} and \ref{lem:reconstruct}, but a quick look at the proof of Lemma \ref{lem:dot-checkv} shows that
 $$ \| \d_\e^{|\a|} \nabla_\e \dot u \|_{\fl1} \lesssim \| \d_\e^{|\a|} \nabla_\e \check u \|_{\fl1},$$ 
and the proof of Lemma \ref{lem:reconstruct} shows that
$$ \| \d_\e^{|\a|} \nabla_\e \check u \|_{\fl1} \lesssim \| \d_\e^{|\a| + 1} u_{high} \|_{\fl1} + \| \d_\e^{|\a|} u_{in} \|_{\fl1} + \| \d_\e^{|\a|} u_{out} \|_{\fl1}.$$
Thus 
$$ \| \d_\e^\a F_{high} \|_{\fl1} \lesssim \e^{K_a} + \e^{-1/2 + K'} \big( \, \| \d_\e^{|\a| + 1} u_{high} \|_{\fl1} + \| \d_\e^{|\a|} u_{in} \|_{\fl1} + \| \d_\e^{|\a|} u_{out} \|_{\fl1} \, \big),$$
and \eqref{fl1:2} follows.
 \end{proof} 

\section{Estimates in the {\it out}} \label{sec:out}

We prove here upper bounds for the $u_{out}$ component in weighted Sobolev and $\fl1$ norms. The system is (see Section \ref{sec:prepared:prepared}):
$$ \d_t u_{out} + \frac{1}{\e} A(\nabla_\e) u_{out}  = \frac{1}{\sqrt \e} \op_\e((1 - \chi_{long}^\flat) \check B_{(-1)}) u_{out} + F_{out}.$$

We start with a Lemma describing the datum for $u_{out}:$

\begin{lem} \label{lem:datum:out} If the support of the initial perturbation $\phi$ is included in the interior of $\{ \chi_{long} \equiv 1\},$ then $u_{out}$ introduced in \eqref{def:U:in:out:high} satisfies the bounds, for $|\a| \leq s_1 - 2:$ 
 $$ \| u_{out} \|_{\e,s} + \| \d_\e^\a u_{out} \|_{\fl1} \lesssim \e^{K_a}.$$
\end{lem}

\begin{proof} By definition of $\check v$ in \eqref{def:checkv} and the description of the datum for $\dot u$ in \eqref{pert:2},  
$$ u_{out}(0) = \e^K \op_\e(\chi_{low}) \big( (1 - \chi_{long}) ({\rm Id}  + \sqrt \e \op_\e(Q_{(-1)}))^{-1} \phi^\e \big),$$
where $\phi^\e$ is defined in \eqref{def:phie}. Consider a term from the Neumann series:
 $$ \e^{n/2} (1 - \chi_{long}) \op_\e(Q)^n \phi^\e = \e^{n/2} \op_\e(Q)^n ((1 - \chi_{long}) \phi^\e) - \e^{n/2} [\chi_{long}, \op_\e(Q)^n] \phi^\e.$$
 If we choose the support of $\phi$ to be included in $\{ \chi_{long} \equiv 1 \},$ then
 $$ (1 - \chi_{long}) \phi^\e \equiv 0.$$
 Besides, by \eqref{compo:e}, Proposition \ref{prop:composition}, and \eqref{composition:sobolev}, 
 \be \label{13:10:20} [\chi_{long}, \op_\e(Q)^n] \phi^\e = \op_\e(\dot{\chi}_{long} \tilde Q_{n,m}) + \e^{m} R_m,\ee 
where $R_m$ is a uniform remainder, and $\dot \chi_{long}$ is some smooth and compactly supported map  which is supported away from $\{ \chi_{long} \equiv 1\},$ implying $\dot \chi_{long} \phi \equiv 0.$ The above holds for any $m \in \N,$ for some $\tilde Q_{n,m}$ depending on $Q,$ $n,$ and $m,$ and which is a sum of fast oscillations times functions of the WKB amplitudes times Fourier multiplier of order $-1.$ Next, using again the composition result from Section \ref{app:symb}, 
 $$ \op_\e(\dot{\chi}_{sp} Q_{n,m}) = \op_\e(Q_{n,m}) \circ \dot{\chi}_{long} + \cdots + \e^{m'} R,$$
 where $m' \in \N$ is arbitrary, $R$ is a uniform remainder, and the dots represent a finite sum of terms involving spatial truncations supported within the support of $\dot{\chi}_{long}.$ In particular those spatial truncations times $\phi$ vanish identically. This implies that the commutator in \eqref{13:10:20} is arbitrarily small, and 
 $$ \| [\chi_{long}, \op_\e(Q)^n] \phi^\e \|_{\e,s} + \| \d_\e^\a [\chi_{long}, \op_\e(Q)^n] \|_{\fl1} \lesssim \e^m,$$
 with $m$ arbitrarily large, and the result follows.
\end{proof}

Next we have a Sobolev bound {\it with no loss of derivative} for the {\it out} component. This is due to the fact that we used $\check v,$ resulting from the second normal form reduction, in order to define $u_{out}:$  

\begin{prop} \label{prop:out} For any $\g_{out} > 0,$ we can choose the support of the spatial truncation $\chi_{long}$ to be large enough so that, for $t \leq t_\star(\e)$ and $s_1 \leq s,$ we have for $u_{out}$,  $u_{in}$ and  $u_{high}$ from \eqref{def:U:in:out:high}:
$$ \begin{aligned} \| u_{out}(t) \|_{\e,s_1}^2 & \lesssim e^{t \g_{out}/\sqrt \e} \Big(\e^{2 K_a} + \e^{K'} |\ln \e|^\star  \max_{0 \leq t' \leq t} \big( \| u_{in}(t') \|^2_{\e,s_1} + \| u_{high}(t') \|^2_{\e,s_1} \,\big) \,  \Big).
\end{aligned}$$  
\end{prop}

\begin{proof} 
By hyperbolicity of $A,$ we have
$$ \frac{1}{2} \d_t \| u_{out}\|_{\e,s}^2 \lesssim \frac{1}{\sqrt \e} \Re e \, \big( \op_\e((1 - \chi_{long}^\flat) \check B_{(-1)}) u_{out}, u_{out} \big)_{\e,s} + \Re e \, \big( F_{out}, u_{out} \big)_{\e,s}.$$
Recall that by definition of $B$ (see Section \ref{sec:Bua}) and $\check B_{(-1)}$ in terms of $B$ (in Section \ref{sec:normal:2}), every entry of the matrix-valued symbol $(1 - \chi_{long}^\flat) \check B_{(-1)}$ is a sum of terms of the form $e^{ip \theta} (1 - \chi_{long}^\flat) u_{0,p} \check B_p, \ p = \pm 1$, where $\check B_p \in S^0$ is a Fourier multiplier. The $L^2 \to L^2$ norm of the corresponding operators is controlled by $|(1 - \chi_{long}^\flat) u_{0,p} |_{L^\infty_{x,y}} |\check B_p|_{L^\infty_{\xi,\eta}}.$ By Sobolev regularity, the amplitudes $u_{0,p}$ decay to zero at spatial infinity. Hence for any $\g_{out} > 0,$ we can find a spatial truncation $\chi_{long}$ with large enough support so that $|(1 - \chi_{long}^\flat) u_{0,\pm 1}|_{L^\infty} \leq \g_{out}.$ The same holds for a finite number of spatial derivatives $\d_\e^\a (1 - \chi_{long}^\flat) u_{0,p}.$ 
This gives
$$ \frac{1}{2} \d_t \| u_{out}(t)\|_{\e,s}^2 \lesssim \frac{\g_{out}}{\sqrt \e} \| u_{out}\|_{\e,s}^2 + \Re e \, \big( F_{out}, u_{out} \big)_{\e,s}.$$
We handle the $F_{out}$ term similarly to the $F_{high}$ term in the proof of high-frequency estimates in Section \ref{sec:high:sob}. That is, by definition of $F_{out}$ in Proposition \ref{prop:prepared} and estimates \eqref{neumann:0} and \eqref{bd:F:sob}, we may focus on 
\be \label{for:Fout}
 \Re e \, \big( \,\op_\e(\chi_{low}) ((1 - \chi_{long}) \dot F),  (\op_\e(\chi_{low}) ((1 - \chi_{long}) \dot u) \big)_{\e,s}, 
\ee
meaning that the difference between $\Re e \, (F_{out}, u_{out})_{\e,s}$ and \eqref{for:Fout} is controlled by 
\be \label{for:out:bd} (\e^{-1/2 + K'} \| \dot u \|_{\e,s} + \e^{K_a}) \| u_{out} \|_{\e,s}.\ee
Now for \eqref{for:Fout}, we can use the symmetry and Gagliardo-Nirenberg arguments that led to \eqref{bd:dotF:Hs}, since the scalar operators $\op_\e(\chi_{low})$ and $1 - \chi_{long}$ do not disrupt the symmetry. This gives a control of \eqref{for:Fout} also by \eqref{for:out:bd}.

Summing up, we obtained
$$  \frac{1}{2} \d_t \| u_{out}(t)\|_{\e,s}^2 \lesssim \frac{\g_{out}}{\sqrt \e} \| u_{out}\|_{\e,s}^2 +  (\e^{-1/2 + K'} \| \dot u \|_{\e,s} + \e^{K_a}) \| u_{out} \|_{\e,s},$$
and the result then follows from Gronwall's lemma.
     \end{proof}

\begin{prop} \label{prop:out:fl1} For $t \leq t_\star(\e),$ we have, with the same $\g_{out}$ as in Proposition {\rm \ref{prop:out}}:
\be \label{fl1:out} \begin{aligned} \| u_{out}(t) \|_{\fl1} & \lesssim \e^{K_a} e^{t \g_{out}/\sqrt \e} \\ & \quad + \e^{-1/2 + K'} \int_0^t e^{(t - t') \g_{out}/\sqrt\e} \max_{0 \leq t'' \leq t'} \big( \| u_{in}(t'') \|_{\fl1} + \| u_{high}(t'') \|_{\fl1} \big) \, dt',
\end{aligned} \ee
and, for $1 \leq m \leq s-2:$ 
\be \label{fl1:out:der} \begin{aligned} & \| \d_\e^{m} u_{out} \|_{\fl1}   \lesssim \e^{K_a} e^{t \g_{out}/\sqrt\e} \\  & \quad + \e^{-1/2 + K'} \int_0^t e^{(t- t') \g_{out}/\sqrt\e} \big( \max_{ \begin{smallmatrix} 0 \leq t'' \leq t' \end{smallmatrix} }  \| \d_\e^{m} u_{in}(t'') \|_{\fl1} +  \max_{\begin{smallmatrix} 0 \leq t'' \leq t' \end{smallmatrix} } \| \d_\e^{m  +1} u_{high}(t'') \|_{\fl1}\,\big) \, dt',\end{aligned}\ee 
where we used notation \eqref{fl1:m}.
\end{prop}

There is no loss of derivatives for $u_{high}$ in \eqref{fl1:out}. This is due to the estimate for $\dot F$ \eqref{bd:dotF:fl1:apriori}, where no derivatives are lost because of the short-time bound \eqref{a:priori:tilde:bound}.

\begin{proof} We give the details for the proof of \eqref{fl1:out:der}, since the proof of \eqref{fl1:out} is similar and simpler. Just like in the proof of Proposition \ref{prop:high:fl1}, by strong hyperbolicity of $A$ we find the bound 
$$ \| \d_\e^\a u_{out}(t) \|_{\fl1} \lesssim \int_0^t (\e^{-1/2}\| \d_\e^\a \op_\e((1 - \chi_{long}^\flat) \check B_{(-1)}) u_{out} \|_{\fl1} + \| \d_\e^\a F_{out} \|_{\fl1}) \, dt'.$$ We use the fact that $\check B_{(-1)}$ has a tensor structure: $\check B_{(-1)} = e^{i p \theta} g_p \underline B_{(-1)},$ for some Fourier multiplier $\underline B_{(-1)} \in S^0.$ Thus, by Lemma \ref{lem:tensor}, pointwise decay of $g_p,$ and choice of the support of $\chi_{long}$ (as in the proof of Proposition \ref{prop:out}), we find 
$$ \e^{-1/2}\| \d_\e^\a \op_\e((1 - \chi_{long}^\flat) \check B_{(-1)}) u_{out} \|_{\fl1} \leq \frac{\g_{out}}{\sqrt \e} \max_{\b \leq \a} \| \d_\e^\b u_{out} \|_{\fl1}.$$
By definition of $F_{out}$ in Proposition \ref{prop:prepared}, and property of the linear operator that defines $u_{out}$ in terms of $\check v,$ we have 
 $$ \| \d_\e^\a F_{out} \|_{\fl1} \lesssim \max_{\b \leq \a} (\| \d_\e^\b \check F \|_{\fl1} + \| \d_\e^\b \check v \|_{\fl1}).$$
 Next by definition of $\check F$ in Proposition \ref{prop:normal:form:2}, and property of the linear operator that defines $\check v$ in terms of $\dot u,$ we have 
$$ \| \d_\e^\b \check F \|_{\fl1} + \| \d_\e^\b \check v \|_{\fl1} \lesssim \max_{\b' \leq \b} (\| \d_\e^{\b'} \dot F \|_{\fl1} + \| \d_\e^{\b'} \dot u \|_{\fl1}).$$ 
With \eqref{final:F:fl1} and notation \eqref{fl1:m}, this gives  
$$ \| \d_\e^{\a} F_{out} \|_{\fl1} \lesssim \e^{K_a} + \e^{-1/2 + K'} (\| \d_\e^{|\a|} \dot u \|_{\fl1} + \| \d_\e^{|\a|} \nabla_\e \dot u \|_{\fl1}).$$
Just like in the proof of Proposition \ref{prop:high:fl1}, we find that the loss of derivative affects only the high-frequency component $u_{high}$ of the solution $\dot u,$ so that 
$$ \| \d_\e^{\a} F_{out} \|_{\fl1} \lesssim  \e^{K_a} + \e^{-1/2 + K'} \big( (\| \d_\e^{|\a|} u_{in} \|_{\fl1} + \| \d_\e^{|\a|} u_{out} \|_{\fl1}) + \| \d_\e^{|\a| + 1} u_{high} \|_{\fl1}\,\big).$$
Thus we find, using Lemma \ref{lem:datum:out},
$$ \begin{aligned} & \| \d_\e^{m} u_{out} \|_{\fl1} \lesssim \e^{K_a} + \int_0^t \frac{\g_{out}}{\sqrt \e} \| \d_\e^{m} u_{out}(t') \|_{\fl1} \, dt' \\ & \quad + \e^{-1/2 + K'} \int_0^t \big( \| \d_\e^{m} u_{in} \|_{\fl1} + \| \d_\e^{m} u_{out} \|_{\fl1} + \| \d_\e^{m  +1} u_{high} \|_{\fl1}\,\big) \, dt',\end{aligned}$$
whence the result by Gronwall's lemma. 
\end{proof}

\section{Coordinatization for the {\it in} subsystem} \label{sec:6}

The {\it in} subsystem is the equation in $u_{in}$ derived in Proposition \ref{prop:prepared}. We reproduce this equation here:
 \be \label{eq:in} \d_t u_{in} + \frac{1}{\e} \op_\e\big( \chi_{low}^\sharp A  \big) u_{in}  = \frac{1}{\sqrt \e} \op_\e(\chi_{long}^\sharp\check B) u_{in} + F_{in}.
 \ee
 The goal of the changes of variables in this Section is twofold:
 \begin{itemize}
 \item we will {\it eliminate the fast oscillations} in the source $\chi_{long}^\sharp \check B;$
 \item and the resonant set will be revealed to be the {\it locus of weak hyperbolicity} for the leading hyperbolic operator.
 \end{itemize}
These changes of variables
 depend heavily on the structure of the set ${\mathcal R}$ of resonant frequencies (described in Propositions \ref{prop:res} and \ref{prop:separation}).

\subsection{First step: diagonalization} \label{sec:diago}

Using identity \eqref{id:proj} for the eigenprojectors of $A,$ we decompose any function $f$ of $(t,x,y)$ with values in $\C^{14},$ in particular the solution $u_{in}$ to the prepared system \eqref{eq:in}, as
 \be \label{coordinatize}
 f = \sum_{1 \leq j \leq 5} f_j, \qquad f_j  := \op_\e(\Pi_j) f.
 \ee
 We frequently identify $f$ with $(f_1, \dots, f_5)$ in the following, that is, we view the decomposition \eqref{coordinatize} as a coordinatization.

 We denote $U_j$ the coordinates of $u_{in}$ in decomposition \eqref{coordinatize}:
\be \label{dec:uin}
 u_{in} = \sum_{1 \leq j \leq 5} U_j, \qquad \mbox{$U_j := \op_\e(\Pi_j) u_{in}.$}
\ee

\begin{lem}[Approximate diagonalization] \label{lem:diago} In the coordinatization \eqref{coordinatize}, the first-order operator in \eqref{eq:in} takes the form
\be \label{A:coord}
\op_\e(\chi_{low}^\sharp A)  f = \op_\e(\chi_{low}^\sharp) \big( i \op_\e(\l_1) f_1, \,\, i \op_\e(\l_2) f_2, \,\, 0, \,\, i \op_\e(\l_4) f_4, \,\, i \op_\e(\l_5) f_5\big) + \sqrt \e R f,
\ee
where $R$ is a uniform remainder in the sense of Definition {\rm \ref{defi:uniform:remainder},}
 and the source $\check B$ defined in \eqref{def:checkB} takes the form
  \be \label{def:checkB:new}
\check B =  \left(\begin{array}{ccccc} 0 & e^{i \theta} b_{12}^+ & e^{i \theta}  b_{13}^+ & e^{i \theta} b_{14}^+ & 0 \\  b_{21}^- e^{-i\theta} & 0 & e^{i \theta}  b_{23}^+ & e^{i \theta}  b_{24}^+ & e^{i \theta}  b_{25}^+ \\ b_{31}^- e^{-i\theta} &  b_{32}^- e^{-i\theta} & 0 & e^{i \theta}  b_{34}^+ &  e^{i \theta} b_{35}^+ \\ b_{41}^- e^{- i \theta} &  b_{42}^- e^{- i \theta} & b_{43}^- e^{-i\theta} & 0 & e^{i\theta}  b_{45}^+  \\ 0 &  b_{52}^- e^{- i \theta} &  b_{53}^- e^{-i\theta} &  b_{54}^- e^{-i\theta} & 0 \end{array}\right) + \sqrt \e R, 
\ee
where $R$ is another uniform remainder, with notation 
\be \label{def:bjj'}  \left.\begin{aligned} b_{jj'}^+ & :=  g_{1}(x,y) \chi_{jj'}(\xi,\eta) B_{1jj'}(\xi,\eta) \\ b_{j'j}^- & := g_{-1}(x,y) \chi_{jj'}(\xi,\eta) B_{-1j'j}(\xi + k, \eta)\end{aligned}\right\}\qquad \mbox{if $(j,j')$ is a resonant pair,}
\ee
where $g_\pm$ is the truncated WKB amplitude introduced in \eqref{def:g}.
\end{lem}

The diagonal in $\check B$ is equal to 0 since there are no auto-resonances. The $(1,5)$ and $(5,1)$ entries in $\check B$ are equal to 0 since there are no $(1,5)$ resonances. 

\begin{nota} \label{note} Above in \eqref{def:checkB:new}, when we write $\op_\e(\sigma e^{-i\theta})$ we really mean the linear operator $u \to \op_\e(\sigma)(e^{-i\theta} u).$ That is, the fast oscillations $e^{-i \theta}$ should be understood as a multiplication operator to the right. This is clearly an abuse of notation, but a convenient one, which we will use only for composition to the right by $e^{-i\theta}.$
\end{nota}

\begin{proof} The decomposition \eqref{A:coord} directly follows from \eqref{approximate:spectral:decomposition}. Indeed, the symbol $\chi_{low}^\sharp R_A$ is order zero, hence a uniform remainder. 

For the source term, we start from \eqref{def:checkB}. Consider the terms associated with $p  = 1$ in $\check B.$ Those have the form 
 $$ e^{i \theta} g_1 \chi_{jj'} B_{1jj'}, \qquad \mbox{for any resonant pair $(j,j').$}$$ %
In decomposition \eqref{dec:uin}, the $(j_0, j'_0)$ coordinate of the above term in $\check B$ is
$$ \op_\e(\Pi_{j_0}) \big( e^{i \theta} g_1 \op_\e(\chi_{jj'} B_{1jj'}) \big) \op_\e(\Pi_{j'_0}) = e^{i \theta} g_1 \op_\e\big(\chi_{jj'} \Pi_{j_0, +1} B_{1jj'} \Pi_{j'_0} \big) + \sqrt \e R,$$
where $R$ is a uniform remainder. The interaction coefficient $B_{1jj'}$ satisfies
 $$ B_{1jj'} \equiv \Pi_{j,+1} B_{1jj'} \Pi_{j'}.$$
 Thus, by orthogonality of the eigenprojectors \eqref{id:proj}, the $(1,j,j')$ term in $\check B$ takes the form $e^{i \theta} g_1 b_{jj'}^+,$ as in \eqref{def:checkB:new}.

 Consider next the terms with $p = -1$ in $\check B.$ Those have the form
 $$g_{-1}(t,x,y) \op_\e(\chi_{jj'} (B_{-1j'j})_{+1}) (e^{- i \theta} \cdot), \qquad (B_{-1j'j})_{+1} \equiv \Pi_{j} (B_{-1j'j})_{+1} \Pi_{j',+1}.$$ Applying $\op_\e(\Pi_{j'_0})$ to the left and $\op_\e(\Pi_{j_0})$ to the right, we find
 $$ \begin{aligned} \op_\e(\Pi_{j'_0}) g_{-1}(t,x,y) & \op_\e(\chi_{jj'} (B_{-1j'j})_{+1}) (e^{- i \theta} \op_\e(\Pi_{j_0}) \, \cdot) \\ & = \op_\e(\Pi_{j'_0}) g_{-1}(t,x,y) \op_\e(\chi_{jj'} (B_{-1j'j})_{+1}) \op_\e(\Pi_{j_0,+1}) \big( (e^{- i \theta} \, \cdot\big) \\ & = \delta_{jj_0} \delta_{j'j'_0} \op_\e(g_{-1} (B_{-1j'j})_{+1})) (e^{- i \theta} \cdot) + \sqrt \e R, \end{aligned}$$
 hence \eqref{def:checkB:new}.
\end{proof}

\subsection{Second step: shifting the extremal components $U_1$ and $U_5$}

We further simplify equation \eqref{eq:in}, now viewed in the coordinatization \eqref{coordinatize}.

\begin{lem} \label{lem:step2} Let
\be \label{b1}
\tilde U := (\tilde U_1, \dots, \tilde U_5) := \big( e^{-i \theta} U_1, \, U_2, \, U_3, \, U_4, \, e^{i \theta} U_5 \big), \qquad \theta = (k x - \o t)/\e,
\ee
where the coordinates $U_j$ of $u_{in}$ are defined in \eqref{dec:uin}. Then, $\tilde U$ solves
 \be \label{eq:tildeU}
\d_t \tilde U + \frac{i}{\e} \op_\e(\tilde A) \tilde U = \frac{1}{\sqrt \e} \op_\e(\tilde B) \tilde U + \tilde F,
\ee
with
\be \label{tildeA}
 \tilde A := \chi_{low}^\sharp {\rm diag}\, \big( \l_{1,+1} -  \o, \,  \l_2, \, \l_3, \, \l_4, \, \l_{5,-1} + \o \big),
\ee
and, using notation $b_{jj'}^\pm$ introduced in Lemma {\rm \ref{lem:diago}},
\be \label{b2}
\tilde B = \chi_{long}^\sharp \left(\begin{array}{ccccc} 0 & b_{12}^+ &  b_{13}^+ &  b_{14}^+ & 0 \\  b_{21}^- & 0 & e^{i \theta}  b_{23}^+ & e^{i \theta}  b_{24}^+ & (b_{25}^+)_{-1} \\  b_{31}^-  &  b_{32}^- e^{-i\theta} & 0 & e^{i \theta}  b_{34}^+ &  (b_{35}^+)_{-1} \\ b_{41}^- &  b_{42}^- e^{- i \theta} &  b_{43}^- e^{-i\theta} & 0 & (b_{45}^+)_{-1}  \\ 0 & (b_{52}^-)_{-1} & (b_{53}^-)_{-1}& (b_{54}^-)_{-1} & 0 \end{array}\right).
\ee
In \eqref{eq:tildeU}, the remainder $\tilde F$ is equivalent to $\dot F$ in the sense of Definition {\rm \ref{def:equiv:remainder}.} 
\end{lem}

In the definition of $\tilde B$ above, we use notation $\s_{-1}$ to denote a shift of the argument of symbol $\s$ by $-k$ in the $\xi$ direction (see Remark \ref{rem:shift}). The convention introduced in Note \ref{note} applies to $\tilde B.$

\begin{proof} Given a symbol $\mu \in S^0$ and $f \in L^2,$ consider the differential equation
\be \label{edo} \d_t v + \frac{i}{\e} \op_\e(\mu) v = f.\ee
Let $\tilde v := e^{i p \theta} v,$ with $p \in \{-1,1\}.$ Then, using notation introduced in Remark \ref{rem:shift}, the shifted unknown $\tilde v$ solves
$$ \d_t \tilde v + \frac{i}{\e} \op_\e(\mu_{-p} + p \o) \tilde v = e^{i p \theta} f.$$
The equation in $U_1,$ as given by \eqref{eq:in} and Lemma \ref{lem:diago}, has the form \eqref{edo}, with a source term $f = e^{i \theta} g.$ We may thus apply the above, with $p = -1,$ to the equation in $U_1,$ and this gives the equation in $\tilde U_1$ as announced.

The equation in $U_5$ (again, \eqref{eq:in} and Lemma \ref{lem:diago}) has the form \eqref{edo} with a source $f = g (e^{- i \theta} \cdot).$ We apply the above with $p = 1,$ and this gives the equation in $\tilde U_5,$ since $e^{i \theta} g e^{-i \theta} = g_{-1}.$ (Here as before, oscillations $e^{i p \theta}$ are understood as multiplication operators and $g_{-1}$ denotes a shift of the symbol $g$ by $-k$ in the $\xi$ direction.)

Besides, the first column of $\check B$ \eqref{def:checkB:new} involves source terms in $U_1.$ Consider for instance the source $\op_\e(b_{21}^-) e^{- i \theta} U_1$ in the equation in $U_2,$ as given by Lemma \ref{lem:diago}. This term is obviously equal to $\op_\e(b_{21}^-) \tilde U_1.$ Similarly, the last column of $\check B$ involves terms in $U_5.$ Consider for instance the term $e^{i \theta} \op_\e(b_{35}^+) U_5$ in the equation in $U_3.$ This term is equal to $\op_\e((b_{35}^+)_{-1}) \tilde U_5$ (here we use again the notation for translates of symbol introduced in Remark \ref{rem:shift}) and we arrive at $\tilde B$ \eqref{b2}.
\end{proof}

\subsection{Third step: decomposition of $U_2$ into ``in" and ``out" variables} \label{sec:coord:3} In view of the form of the source $\tilde B$ in Lemma \ref{lem:step2}, it would be tempting to define $V_2 := e^{- i \theta} \tilde U_2$ in order to get rid of the fast oscillations in the entries of $\tilde B$ involving the $(2,3)$ resonance, just like defining $\tilde U_1 = e^{-i \theta} U_1$ allowed us, in Lemma \ref{lem:step2}, to get rid of the fast oscillations associated with the $(1,2)$ resonance. But then clearly this would introduce another oscillation in the $\tilde U_1/U_2$ coupling term. Thus we may take out fast oscillations only away from the support of $\chi_{12}.$

This is done in Lemma \ref{lem:step3} below, in which we use notation pertaining to frequency cut-off functions introduced in Section \ref{sec:resonances}.

\begin{lem} \label{lem:step3} Let
$$ U^\sharp := (\tilde U_1, U_2^{out}, U_2^{in}, U_3,  U_4, \tilde U_5),$$
with
$$ U_2^{out} := \op_\e(1 - \chi_{23,-1}^\sharp - \chi_{24,-1}^\sharp)  U_2, \qquad U_2^{in} : = \op_\e(\chi_{23}^\sharp + \chi_{24}^\sharp) \big( e^{- i \theta} U_2\big),$$
where the $U_j$ are the coordinates of $v$ \eqref{dec:uin} and $\tilde U_1$ and $\tilde U_5$ are defined in Lemma {\rm \ref{lem:step2}.}
Then,  $U^\sharp$ solves
\be \label{eq:Usharp}
\d_t U^\sharp + \frac{i}{\e} \op_\e(A^\sharp)  U^\sharp = \frac{1}{\sqrt \e} \op_\e( B^\sharp)  U^\sharp +  F^\sharp,
\ee
with notation (using notation for cut-offs introduced in Section {\rm \ref{sec:cutoffs}})
$$ A^\sharp = \chi_{low}^\sharp {\rm diag}\, \big( \l_{1,+1} - \o, \,\, (1 - \chi_{23,-1} - \chi_{24,-1}) \l_2, \,\, (\tilde \chi_{23} + \tilde \chi_{24})(\l_{2,+1} - \o) , \,\, \l_3, \,\, \l_4, \,\, \l_{5,-1} + \o \big),
$$
and
\be \label{def:Bsharp}
 B^\sharp := \chi_{long}^\sharp \left(\begin{array}{cccccc} 0 &  b_{12}^+ & 0 &  b_{13}^+ & b_{14}^+ & 0 \\[1pt]   b_{21}^- & 0& 0  & 0 & 0 & (b_{25}^+)_{-1}\\[1pt] 0 & 0 & 0 & b_{23}^+ & b_{24}^+ & 0 \\[1pt]  b_{31}^-  & 0 & b_{32}^-  & 0 & e^{i \theta}  b_{34}^+ &  (b_{35}^+)_{-1} \\[1pt] b_{41}^- & 0 &  b_{42}^- & b_{43}^- e^{-i\theta} & 0 & (b_{45}^+)_{-1}  \\[1pt] 0 & (b_{52}^-)_{-1} & 0 & (b_{53}^-)_{-1}& (b_{54}^-)_{-1} & 0 \end{array}\right).
\ee
In \eqref{eq:Usharp}, the remainder $F^\sharp$ is equivalent to $\dot F$ in the sense of Definition {\rm \ref{def:equiv:remainder}.} 
\end{lem}

\begin{proof} The proof relies on Proposition \ref{prop:separation}. We start with the equation in $U_2^{in}.$ First we consider the oscillating term. Just like in the equation in $\tilde U_1,$ the fast oscillation in the definition of $U_2^{in}$ gives
 \be \label{for:u2in} \big(\d_t + \frac{i}{\e} \op_\e(\l_{2,+1} - \o) \big) U_2^{in} = \dots\ee
 Next we remark that
 \be \label{u2in:identity}
 U_2^{in} = \op_\e(\tilde \chi_{23} + \tilde \chi_{24}) U_2^{in}, \qquad \chi_{23}^\sharp \prec \tilde \chi_{23}, \quad \chi_{24}^\sharp \prec \tilde \chi_{24},
 \ee
 where we use notation introduced in \eqref{sharp:tilde}. Thus the operator in \eqref{for:u2in} may be changed into $\op_\e( (\tilde \chi_{23} + \tilde \chi_{24}) (\l_{2,+1} - \o)).$ The source term in the equation in $U_2^{in}$ is, omitting the $1/\sqrt \e$ prefactor: 
 \be \label{source:U2in} \begin{aligned} \op_\e(\chi_{23}^\sharp & + \chi_{24}^\sharp) \\ & \cdot \big( e^{ - i \theta} \big( \op_\e(\chi_{12} b_{21}^-) \tilde U_1 + e^{i \theta} \op_\e(\chi_{23} b_{23}^+) U_3 + e^{i \theta} \op_\e(\chi_{24} b_{24}^+) U_4 + \op_\e(\chi_{25} b_{25}^-)_{-1} \tilde U_5\big)\big).
 \end{aligned}\ee 
 Consider the term in $\tilde U_1$ above. Up to a uniform remainder of the form $\sqrt \e R U^\sharp,$ it is equal to 
 $$ e^{- i \theta} \op_\e( (\chi_{23,-1}^\sharp + \chi_{24,-1}^\sharp) b_{21}^-) \tilde U_1.$$ 
  By Lemma \ref{lem:cut-offs}, we have
\be \label{newsep:1} (\chi_{23,-1}^\sharp + \chi_{24,-1}^\sharp) \chi_{12}  \equiv 0,\ee
hence the term in $\tilde U_1$ in the equation in $U^{in}$ vanishes to first order. 

Consider now the terms in $U_3$ and $U_4$ in \eqref{source:U2in}. Up to a uniform remainder of the form $\sqrt \e R U^\sharp,$ they are equal to 
\be \label{for:source:U2in} \op_\e( (1 + \chi_{24}^\sharp)  b_{23}^+) U_3 + \op_\e( (1 + \chi_{23}^\sharp)  b_{24}^+) U_4.\ee 
By Lemma \ref{lem:cut-offs}, we have
\be \label{newsep:2} \chi_{23}^\sharp \chi_{24} \equiv 0, \quad \chi_{23} \chi_{24}^\sharp \equiv 0,\ee
hence \eqref{for:source:U2in} takes the form
$$ \op_\e( b_{23}^+) U_3 + \op_\e(b_{24}^+) U_4.$$
Finally consider the term in $\tilde U_5$ in \eqref{source:U2in}.  Up to a uniform remainder of the form $\sqrt \e R U^\sharp,$ it is equal to
 $$ e^{i \theta} \op_\e\Big(  \, \big((\chi_{23}^\sharp + \chi_{24}^\sharp\big)_{-1} b_{25}^-)_{-1} \, \Big) \tilde U_5.$$
 We appeal again to Lemma \ref{lem:cut-offs}, which asserts that
\be \label{newsep:3} (\chi_{23}^\sharp + \chi_{24}^\sharp) \chi_{25} \equiv 0,\ee
 so that the term in $\tilde U_5$ in \eqref{source:U2in} vanishes to first order.
 We obtained the third row in $B^\sharp.$

 We turn to the equation in $U_2^{out}.$ The oscillating term is handled as above, by $$U_2^{out} \equiv \op_\e\big(1 - \chi_{23,-1} - \chi_{24,-1}\big) U_2^{out}.$$ The source term involves $\tilde U_1,$ $U_3,$ $U_4$ and $\tilde U_5$ just like in \eqref{source:U2in} but with a different operator to the left. We use again the cancellations \eqref{newsep:1}-\eqref{newsep:2}-\eqref{newsep:3}. From 
 $$ (1 - \chi_{23,-1}^\sharp - \chi_{24,-1}^\sharp) \chi_{12} \equiv \chi_{12}$$ 
 we deduce that the source term in $\tilde U_1$ is $\op_\e(b_{21}^-) \tilde U_1,$ to first order. From  
 $$ (1 - \chi_{23}^\sharp - \chi_{24}^\sharp) \chi_{2j} \equiv 0, \qquad j \in \{3,4\},$$
 we verify that the source term in $U_3$ and $U_4$ vanish to first order, and finally
 $$ (1 - \chi_{23}^\sharp - \chi_{24}^\sharp) \chi_{25} = \chi_{25}$$
 implies that the source term in $\tilde U_5$ is $\op_\e( (b_{25}^+)_{-1}) \tilde U_5.$

We now have a complete description of the second and third rows of $B^\sharp.$ It remains to consider the source terms involving $U_2$ in $\tilde B,$ that is the second column of $\tilde B$ \eqref{b2}. The first of these terms is $\op_\e(b_{12}^+) U_2$ in the equation in $\tilde U_1.$ We use
\be \label{reconstruct:U_2}
U_2 = e^{i \theta} U_2^{in} + U_2^{out} = e^{i \theta} \op_\e(\tilde \chi_{23} + \tilde \chi_{24}) U_2^{in} + \op_\e\big((1 - \chi_{23,-1}^\star -  \chi_{24,-1}^\star)\big) U_2^{out},
\ee
where $\chi_{2j} \prec \chi_{2j}^\star \prec \chi_{2j}^\sharp,$ for $j = 3$ and $j = 4.$ 
Using \eqref{newsep:1} once more, and \eqref{reconstruct:U_2}, we find
$$ \op_\e(\chi_{12} b_{12}^+) U_2 \equiv \op_\e(\chi_{12} b_{12}^+) U_2^{out}.$$
Using \eqref{newsep:2} once more, and \eqref{reconstruct:U_2}, we find
$$ \op_\e(b_{j2}^-) e^{i \theta} U_2 \equiv \op_\e(b_{j2}^-) U_2^{in}, \quad j \in \{3,4\}.$$
Finally, for the source in $U_2$ from the equation in $\tilde U_5,$ we use \eqref{newsep:3} and find 
$$ \op_\e( (b_{52}^-)_{-1}) U_2 \equiv  \op_\e( (b_{52}^-)_{-1}) U_2^{out}.$$
This completes the justification of the third column in $B^\sharp.$
\end{proof}

\subsection{Fourth step: decomposition of $U_4$ into ``in" and ``out" variables} \label{sec:step4}
In a last factorization step, we apply to $U_4$ operations that are (almost) symmetrical to the ones applied to $U_2$ in the previous step:

\begin{lem} \label{lem:step3bis} Let
\be \label{b5}
 v_{in} := (\tilde U_1, U_2^{out}, U_2^{in}, U_3, U_4^{in}, U_4^{out}, \tilde U_5),
 \ee
with notation
$$ U_4^{in} := \op_\e(\chi_{34,-1}^\sharp) (e^{i \theta} U_4), \quad U_4^{out} := \op_\e(1- \chi_{34}^\sharp) U_4,$$
where the $U_j$ are the coordinates of $u_{in}$ \eqref{dec:uin} and $\tilde U_1$ and $\tilde U_5$ are defined in Lemma {\rm \ref{lem:step2}.}
Then, the unknown $v_{in}$ solves
\be \label{eq:Uflat}
\d_t v_{in} + \frac{i}{\e} \op_\e({\bf A}) v_{in} = \frac{1}{\sqrt \e} \op_\e({\bf B}) v_{in} +  {\bf F},
\ee
where
$$ \begin{aligned} {\bf A} = \chi_{low}^\sharp {\rm diag}\, \Big( \l_{1,+1} - \o, & \,\, (1 - \chi_{23,-1} - \chi_{24,-1}) \l_2,  \,\,\,\, (\tilde\chi_{23} + \tilde\chi_{24}) (\l_{2,+1} - \o),   \,\,\, \l_3, \\ & \,\,\,\,\tilde \chi_{34,-1}( \l_{4,-1} + \o), \,\,\, (1 - \chi_{34})\l_4, \,\,\, \l_{5,-1} + \o \Big),\end{aligned}$$
and
$$
{\bf B} :=  \chi_{long}^\sharp \left(\begin{array}{ccccccc} 0 &  b_{12}^+ & 0 &  b_{13}^+ & 0 &  b_{14}^+ & 0 \\[1pt]  b_{21}^- & 0& 0 & 0  & 0 & 0 & (b_{25}^+)_{-1} \\[1pt] 0 & 0 & 0 &  b_{23}^+ & 0 &  b_{24}^+ & 0 \\[1pt]  b_{31}^-  & 0 & b_{32}^-  & 0 & (b_{34}^+)_{-1} & 0 &  (b_{35}^+)_{-1} \\[1pt] 0 & 0 & 0  & (b_{43}^{-})_{-1} & 0 & 0 & 0  \\[1pt] b_{41}^- & 0 & b_{42}^- & 0 & 0 & 0 & (b_{45}^+)_{-1}\\[1pt] 0 & (b_{52}^-)_{-1} & 0 & (b_{53}^-)_{-1}& 0 & ( b_{54}^-)_{-1} & 0 \end{array}\right).
$$
In \eqref{eq:Uflat}, the remainder ${\bf F}$ is equivalent to $\dot F$ in the sense of Definition {\rm \ref{def:equiv:remainder}.}  \end{lem}

\begin{proof} We proceed as in the proof of Lemma \ref{lem:step3} and Proposition \ref{prop:separation}, and its extension to cut-offs in Lemma \ref{lem:cut-offs}.
\end{proof} 

In \eqref{eq:Uflat}, we achieved our goal for this Section, which was to go from equation \eqref{eq:in} with a source containing fast oscillations in space-time to equation \eqref{eq:Uflat}, an equivalent system with a source containing no fast oscillations. The ``equivalence" here means only that bounds for $v_{in}$ translate into bounds for $u_{in}.$

Another key point of the formulation \eqref{eq:Uflat} of the {\it in} equation is that {\it the resonant set now corresponds to the locus of weak hyperbolicity}. Thus at the resonances, the term ${\bf B}$ is susceptible to generate instabilities (as sketched in Section \ref{sec:outline}).

By definition of $v_{in}$ in \eqref{b5}, the first two components of the datum $v_{in}(0)$ satisfy
\be \label{datum:vin}
 v_{in}(0) = \Big( e^{- i k x/\e} \op_\e(\Pi_1) u_{in}(0), \,\, \op_\e(1 - \chi_{23,-1}^\sharp) \op_\e(\Pi_2) u_{in}(0), \dots \Big),
\ee
where $u_{in}(0)$ is described in \eqref{datum:uin}.

\section{The symbolic flow}  \label{sec:rep:flow}

This is the central piece of the analysis, in which we solve \eqref{eq:Uflat} by the symbolic flow method introduced in \cite{duh} and used in similar contexts in \cite{em4,LNT}.

\medskip

In view of proving estimates for \eqref{eq:Uflat}, we introduce the auxiliary linear partial differential equation
\begin{equation} \label{resolvent0-in-proof} \d_t S(\t;t) + \frac{1}{\e} (i {\bf A} - \sqrt \e {\bf B}) S(\t;t) + \frac{1}{\sqrt \e} \d_\eta {\bf A}  \d_y S(\t;t)  = 0, \qquad S(\t;\t) = {\rm Id},\end{equation}
where ${\bf A}$ and ${\bf B}$ are introduced in Lemma \ref{lem:step3bis} and ${\rm Id}$ is the identity matrix with the same size as ${\bf A}$ and ${\bf B}.$ 
 The solution $S$ to \eqref{resolvent0-in-proof} (the {\it symbolic flow}) depends on
$$ (\t;t,x,y,\xi,\eta), \quad \mbox{with $0 \leq \t \leq t \leq T \sqrt \e |\ln \e|, \quad (x,y) \in \R^3, \quad (\xi,\eta) \in \R^3,$}$$
 and is matrix-valued (same dimensions as ${\bf A}$ and ${\bf B}$). In \eqref{resolvent0-in-proof}, by $\d_\eta {\bf A}  \d_y S$ we mean $\d_{\eta_1} {\bf A} \d_{y_1} S + \d_{\eta_2} {\bf A} \d_{y_2} S.$

Theorem \ref{th:duh}, stated in Section \ref{app:duh} below, sums up the symbolic flow method applied to our context. It says roughly that the solution $S$ to \eqref{resolvent0-in-proof} provides a good approximation to the symbol of the solution operator to the initial-value problem for \eqref{eq:Uflat} in the sense that
\be \label{approx:flow} v_{in}(t) \simeq \op_\e(S(0;t)) v_{in}(0).\ee
 The approximation \eqref{approx:flow} is made precise in Section \ref{sec:Duh}. In the remainder of this Section, we show that $S$ is uniquely determined by the initial-value problem \eqref{resolvent0-in-proof}, and describe explicitly its far-field behavior.

 In the course of the analysis, we will have to consider the symbolic equation with a source $f$ and datum $g:$
\be \label{symbolic:flow:source}
\d_t S_{f,g}(\t;t) + \frac{1}{\e} (i {\bf A} - \sqrt \e {\bf B}) S_{f,g}(\t;t) + \frac{1}{\sqrt \e} \d_\eta {\bf A} \d_y S_{f,g}(\t;t)  = f, \qquad S_{f,g}(\t;\t) = g,
\ee
where $f$ and $g$ depend on $(t,x,y,\xi,\eta),$ are bounded, and their $(x,y,\xi,\eta)$ derivatives up to a large degree are also bounded.

\subsection{Separation and block decomposition} \label{sec:flow:blocks}

We use here the separation properties of the resonant set, as given in Proposition \ref{prop:separation}, and the block structure of ${\bf B}$ (defined in Lemma \ref{lem:step3bis}), in order to show that the system \eqref{resolvent0-in-proof} in $S$ decouples into smaller subsystems.

 First, by Proposition \ref{prop:space-time:sep} and the structure of the interaction matrix ${\bf B},$ we identify in \eqref{resolvent0-in-proof} the subsystems
\be \label{2blocks:tilde} 
\d_t S + \frac{1}{\e} \left(\begin{array}{cc} i \tilde \mu_{j} & \sqrt \e \chi_{long}^\sharp b_{jj'}^+ \\ \sqrt \e \chi_{long}^\sharp b_{j'j}^- & i \tilde \mu_{j'} \end{array}\right) S + \frac{1}{\sqrt \e} \left(\begin{array}{cc} \d_\eta \tilde \mu_{j} & 0 \\ 0 & \d_{\eta} \tilde \mu_{j'} \end{array}\right) \cdot \d_y S = 0,
\ee
for the pairs of indices 
\be \label{pairs} (j,j') \in \{ (1,2), (1,4), (2,4), (2,5), (4,5) \},\ee
with notation for the eigenvalues: 
\begin{itemize}
\item for $(1,2),$ we have 
\be \label{def:tildemu:12} \tilde \mu_1 := \chi_{low}^\sharp ( \l_{1,+1} - \o) , \quad \tilde \mu_2 = \chi_{low}^\sharp (1 - \chi_{23,-1} - \chi_{24,-1}) \l_2;\ee 
\item for $(1,4),$ we have 
\be \label{def:tildemu:14} \tilde \mu_1 := \chi_{low}^\sharp ( \l_{1,+1} - \o) , \quad \tilde \mu_4 = \chi_{low}^\sharp (1 - \chi_{34}) \l_4;\ee 
\item for $(2,4):$ 
\be \label{def:tildemu:24} \tilde \mu_2 := \chi_{low}^\sharp (\tilde \chi_{23} + \tilde \chi_{24})(\l_{2,+1} -\o), \quad \tilde \mu_4 = \chi_{low}^\sharp (1 - \chi_{34}) \l_4,\ee
\item for $(2,5):$ 
\be \label{def:tildemu:25} \tilde \mu_2 := \chi_{low}^\sharp (1 - \chi_{23}  -\chi_{24}) \l_2, \quad \tilde \mu_5 := \chi_{low}^\sharp(\l_5 + \o),\ee
\item and for $(4,5)$: 
\be \label{def:tildemu:45} \tilde \mu_4 := \chi_{low}^\sharp (1 - \chi_{34,+1}) \l_{4,+1}, \quad \tilde \mu_5 = \chi_{low}^\sharp(\l_5 + \o).\ee 
\end{itemize}
Thus notation $\tilde \mu_j$ is flexible, hence potentially a little ambiguous. It depends on the {\it resonant pair} that is considered. Precisely, for modes $2$ and $4,$ the definition of $\tilde \mu_2$ and $\tilde \mu_4$ depends on whether the pair that is considered has the form $(2,j')$ or $(j,2),$ and similarly $(4,j')$ or $(j,4).$ This explains the difference between \eqref{def:tildemu:12} and \eqref{def:tildemu:25} for $\tilde \mu_2,$ and between \eqref{def:tildemu:24} and \eqref{def:tildemu:45} for $\tilde \mu_4.$

Next, looking at the fourth row and fourth column of ${\bf B},$ and taking into account Proposition \ref{prop:space-time:sep}, we find within the system \eqref{resolvent0-in-proof} the subsystems
 \be \label{3blocks:tilde} \begin{aligned} 
 \d_t S & + \frac{1}{\e} \left(\begin{array}{ccc} i \tilde \mu_{j_1} + \sqrt \e \d_\eta \tilde \mu_{j_1} \cdot \d_y & \sqrt \e \chi_{long}^\sharp b_{j_1j_2}^+ & 0 \\ \sqrt \e \chi_{long}^\sharp b_{j_2j_1}^- & i \tilde \mu_{j_2} + \sqrt \e \d_\eta \tilde \mu_{j_2} \cdot \d_y & \sqrt\e \chi_{long}^\sharp (b_{j_2j_3}^+)_{-1} \\ 0 & \sqrt \e \chi_{long}^\sharp (b_{j_3j_2}^-)_{-1} & i \tilde \mu_{j_3} + \sqrt \e \d_\eta \tilde \mu_{j_3} \cdot \d_y \end{array}\right) S  = 0, \end{aligned} 
 \ee
with 
\be \label{triplets}
(j_1,j_2,j_3) = (1,3,5) \quad \mbox{or} \quad (j_1,j_2,j_3) = (2,3,4),
\ee
and
\begin{itemize}
\item for the triplet $(1,3,5),$ which has to be considered due to a lack of separation between $(1,3)$ and $(3,5),$ as stated in Proposition \ref{prop:space-time:sep}: 
\be \label{def:tildemu:135} 
\tilde \mu_1 = \chi_{low}^\sharp ( \l_{1,+1} - \o),  \quad \tilde \mu_3 = \chi_{low}^\sharp \l_3,  \quad \tilde \mu_5 = \chi_{low}^\sharp(\l_{5,-1} + \o);
\ee 
\item and for the triplet $(2,3,4),$ which accounts for the fact that the resonances $(2,3)$ and $(3,4)$ are not fully separated (see Proposition \ref{prop:space-time:sep}):  
\be \label{def:tildemu:234}
\tilde \mu_2 = \chi_{low}^\sharp (\tilde \chi_{23} + \tilde \chi_{24}) (\l_{2,+1} - \o),\quad \tilde \mu_3 = \chi_{low}^\sharp\l_3, \quad \tilde \mu_4 = \chi_{low}^\sharp \tilde \chi_{34,-1}( \l_{4,-1} + \o).
\ee 
\end{itemize}
 Note that in the triplet $(2,3,4),$ the mode $2$ is in first position, hence a definition of $\tilde \mu_2$ identical to the one in $(2,4)$ above. Symmetrically, the mode $4$ in second position, hence a definition of $\tilde \mu_4$ identical to the one in $(1,4)$ above. 

 We will consider systems \eqref{2blocks:tilde} and \eqref{3blocks:tilde} separately, hence the fact that the definition of the $\tilde \mu_j$ is context-dependent should be not a factor of confusion.

In \eqref{3blocks:tilde}, the notation $(\cdot)_{-1}$ refers to the translation $\xi \to \xi - k,$ see Remark \eqref{rem:shift}.

 Note finally that by Proposition \ref{prop:space-time:sep} and structure of the matrix ${\bf B},$ the collection of subsystems \eqref{2blocks:tilde} with $(j,j')$ as in \eqref{pairs} and subsystems \eqref{3blocks:tilde} with $(j_1,j_2,j_3)$ as in \eqref{triplets} account for the whole of system \eqref{resolvent0-in-proof}.

\begin{nota} \label{note:dimensions} A note on dimensions: the symbolic flow $S$ is defined as the solution to \eqref{resolvent0-in-proof}: it is a matrix with dimensions equal to those of ${\bf A}$ and ${\bf B}.$ The solution to \eqref{2blocks:tilde} for a given pair $(j_1,j_2),$ is only a submatrix of $S.$ We denote it $S$ rather than $S_{j_1j_2}.$ Thus the $S$ from \eqref{2blocks:tilde} and \eqref{3blocks:tilde} is not the same as the one from \eqref{resolvent0-in-proof}. Furthermore, it will be convenient to handle $S$ (whether the full matrix or one of its distinguished submatrices) as a vector with length equal to the dimensions of ${\bf A}$ and ${\bf B}.$ Precisely, given a fixed vector $\vec u,$ we coordinatize $S \vec u =: (S_1, S_2).$ The computations are linear in $\vec u$ and the upper bounds are linear in $|\vec u|,$ so that we omit $\vec u$ altogether. Sometimes, in particular when we use the far-field description of the upcoming Lemma {\rm \ref{lem:away}}, we go even further and identify $S$ with one of its entries. 
\end{nota}

\subsection{A second look at the subsystems for $S$} \label{sec:second:look} 

We observe here that the subsystems \eqref{2blocks:tilde} and \eqref{3blocks:tilde}, the frequency truncations that appear in the definitions of $\tilde \mu_j$ can be taken out with no loss of information. 

Consider first \eqref{2blocks:tilde} with $(j,j') = (1,2).$ Given $(\xi,\eta)$ outside of the support of $\chi_{12},$ the partial differential system in $(t,y)$ \eqref{2blocks:tilde} becomes constant-coefficient, since the $(x,y)$-dependent terms $b_{12}^+$ and $b_{21}^-$ contain $\chi_{12}$ (see \eqref{def:bjj'}). Thus for $(\xi,\eta)$ outside of the support of $\chi_{12},$ the system \eqref{2blocks:tilde} with $(j,j') = (1,2)$ decouples into subsystems
\be \label{constant:coeff:subsystem}
 \d_t S + \frac{1}{\e} \tilde \mu_j S + \frac{1}{\sqrt \e} \d_\eta \tilde \mu_j \d_y S = 0.
\ee
Besides, by Lemma \ref{lem:cut-offs} (a consequence of Proposition \ref{prop:separation}(1)), we have $\chi_{12} (\chi_{23,-1} + \chi_{24,-1}) \equiv 0.$ Thus given $(\xi,\eta)$ such that $\chi_{12}(\xi,\eta) \neq 0,$ for $\tilde \mu_2$ defined in \eqref{def:tildemu:12}, we have $\tilde \mu_2(\xi,\eta) = \l_2(\xi,\eta).$ 

 By symmetry, and separation of resonances, the same holds for all other resonant pairs in \eqref{pairs}. 
 
 This shows that in the analysis of \eqref{2blocks:tilde}, we will be able to focus on the constant-coefficient equations \eqref{constant:coeff:subsystem}, and on the slightly simpler subsystems 
\be \label{2blocks}
\d_t S + \frac{1}{\e} \left(\begin{array}{cc} i \mu_{j} + \sqrt \e \d_\eta \mu_j \cdot \d_y & \sqrt \e \chi_{long}^\sharp b_{jj'}^+ \\ \sqrt \e \chi_{long}^\sharp b_{j'j}^- & i \mu_{j'} + \sqrt \e \d_\eta \mu_{j'} \cdot \d_y \end{array}\right) S = 0, \qquad (\xi,\eta) \in \mbox{supp}\, \chi_{jj'},
\ee
for $(j,j')$ as in \eqref{pairs},
and 
\be \label{def:mu} \begin{aligned} 
 & \mbox{$\mu_j$ is defined just like $\tilde \mu_j$ in Section \ref{sec:flow:blocks}} \\ & \mbox{but without any truncation, that is $\chi_{low}^\sharp$ changed into 1 and $\chi_{jj',\star}$ changed into 1.} \end{aligned} 
\ee 

 Next we turn to the subsystems \eqref{3blocks:tilde}. First we observe that for
$$ (\xi,\eta) \in \big( \mbox{supp}\,\chi_{j_1j_2} \setminus \mbox{supp}\, \chi_{j_2j_3,-1} \big) \bigcup \big( \mbox{supp}\, \chi_{j_2j_3,-1} \setminus \mbox{supp}\,\chi_{j_1j_2} \big),$$ 
then the system \eqref{3blocks:tilde} decouples into a system of the form \eqref{2blocks:tilde} and an equation of the form \eqref{constant:coeff:subsystem}. This means in particular that we have to consider systems \eqref{2blocks} also for $(j_1,j_2) \in \{ (1,3), (3,5), (2,3), (3,4)\}.$ Now if $(\xi,\eta)$ belongs to the intersection of the supports of $\chi_{j_1j_2}$ and $\chi_{j_2j_3,-1},$ then we observe that the truncations in the $\tilde \mu_j$ \eqref{def:tildemu:135} and \eqref{def:tildemu:234} are identically equal to one. 
 Thus for triplets the systems are 
  \be \label{3blocks} \begin{aligned} 
 \d_t S & + \frac{1}{\e} \left(\begin{array}{ccc} i \mu_{j_1} + \sqrt \e \d_\eta \mu_{j_1} \cdot \d_y & \sqrt \e \chi_{long}^\sharp b_{j_1j_2}^+ & 0 \\ \sqrt \e \chi_{long}^\sharp b_{j_2j_1}^- & i  \mu_{j_2} + \sqrt \e \d_\eta \mu_{j_2} \cdot \d_y & \sqrt\e \chi_{long}^\sharp (b_{j_2j_3}^+)_{-1} \\ 0 & \sqrt \e \chi_{long}^\sharp (b_{j_3j_2}^-)_{-1} & i \mu_{j_3} + \sqrt \e \d_\eta \mu_{j_3} \cdot \d_y \end{array}\right) S  = 0, \end{aligned} 
 \ee
with $\mu_j$ defined in terms of $\tilde \mu_j$ from \eqref{def:tildemu:135} and \eqref{def:tildemu:234}, with the frequency truncations replace by 1, and $(\xi,\eta) \in \mbox{supp}\,\chi_{j_1j_2} \cap \mbox{supp}\, \chi_{j_2j_3,-1}.$ 
  
 {\it Summing up:} besides constant-coefficient equations of the form \eqref{constant:coeff:subsystem}, the symbolic flow equation \eqref{resolvent0-in-proof} reduces to :
 \begin{itemize}
 \item subsystems of the form \eqref{2blocks} for any resonant pair $(j_1,j_2),$ where the definition of $\mu_{j_1}$ and $\mu_{j_2}$ depends on the pair $(j_1,j_2):$ we have $\mu_j = \l_j$ up to translation $\xi \to \xi + k$ and plus or minus $\o,$ in such a way that the frequency set $\{ \mu_{j_1} - \mu_{j_2}  = 0 \}$ is the resonant set ${\mathcal R}_{j_1j_2}.$ 
 \item subsystems of the form \eqref{3blocks} for ``triplets" $(1,3,5)$ and $(2,3,4).$ Those triplets are really pairs of non-separated resonant pairs $(j_1,j_2)-(j_2,j_3)$ as described in Proposition \ref{prop:space-time:sep}.
 \end{itemize}

\subsection{Far-field behavior} \label{sec:far-field}

Given a resonant pair $(j,j'),$ introduce the truncation function 
\be \label{def:chi:Omega:pair} \chi_{\Omega}(x,y,\xi,\eta) := \tilde \chi_{long}(x) {\bm \chi}_{trans}(y) \chi_{jj'}^\sharp(\xi,\eta),\ee
where, in accordance with \eqref{sharp:tilde} and \eqref{sharp:tilde2}, we have $\chi_{long}^\sharp \prec \tilde \chi_{long}$ and $\chi_{jj'} \prec \chi_{jj'}^\sharp,$ and the transverse truncation ${\bm \chi}_{trans}$ is defined by 
\be \label{def:bm:chi:trans}
 {\bm \chi}_{trans}(y) = \left\{\begin{aligned} 1, & \quad |y| \leq  \big( R_{trans} |\ln \e|^{1/\kappa_a} + T |\ln \e|\big), \\ 0, & \quad |y| \geq 2\big(R_{trans} |\ln \e|^{1/\kappa_a} + T |\ln \e|\big),
 \end{aligned}\right.
\ee
where $R_{trans}$ is the radius that appears in the definition of $\chi_{trans}$ in \eqref{def:chi:infty}. The partial definition \eqref{def:bm:chi:trans} is completed by having ${\bm \chi}_{trans}$ smoothly and monotonically transition from $1$ to $0$ as $|y|$ increases, just like for instance in the definition of $\chi_{trans}$ \eqref{def:chi:infty}. We define $\Omega$ to be the support of $\chi_\Omega:$ 
\be \label{def:Omega}
 \Omega := \mbox{supp}\,\chi_\Omega.
 \ee
Note that $\chi_\Omega$ and $\Omega$ depend on the particular resonant pair that we consider, but we denote $\Omega$ and not $\Omega_{jj'}.$ By definition of $\chi_{\Omega},$ we see that $\Omega$ is a product:
\be \label{dec:Omega} \begin{aligned} 
\Omega & =: \Omega_{sp} \times \Omega_{freq} \\ & := \bar B_{\R_x}(0,\tilde R) \times \bar B_{\R^2_y}\big(0, 2\big(R_{trans} |\ln \e|^{1/\kappa_a} + T |\ln \e|\big)\big) \times \overline{\{ (\xi,\eta) \in \R^3, \quad \chi_{jj'}^\sharp(\xi,\eta) > 0 \}}. \end{aligned}
\ee
The projection $\Omega_{freq}$ of $\Omega$ over $\R^3_{\xi,\eta}$ is the support of $\chi_{jj'}^\sharp.$ Above, $\tilde R > 0$ is a radius in $\R_x,$ such that $\tilde R > R,$ where $R$ is the radius associated with the longitudinal cut-off $\chi_{long}$ introduced in Section \ref{sec:in:out:high}.  

\begin{lem}[Far-field behavior for \eqref{2blocks}] \label{lem:away} For the solution $S_{f,g}$ to \eqref{2blocks} with a source $f$ and datum $g,$ we have the explicit description outside of $\Omega:$ for $0 \leq \t \leq t \leq T \sqrt \e |\ln \e|$ and  $(x,y,\xi,\eta) \notin \Omega$,
\be \label{explicit:S} 
\begin{aligned} S^j_{f,g} & = e^{-i (t - \t) \mu_j/\e} g_j\left(x, y - \frac{\d_\eta \mu_j}{\sqrt \e} (t - \t), \xi,\eta\right) \\ & \quad + \int_\t^t e^{- i (t - t') \mu_j/\e} f_j\left(t', y - \frac{\d_\eta \mu_j}{\sqrt \e} (t' - \t), \xi , \eta\right) \, dt', \qquad (x,y,\xi,\eta) \notin \Omega,\end{aligned}
\ee
where $j \in \{1,2\}$ denotes the coordinate in the size-two-by-block system \eqref{2blocks}.
\end{lem}

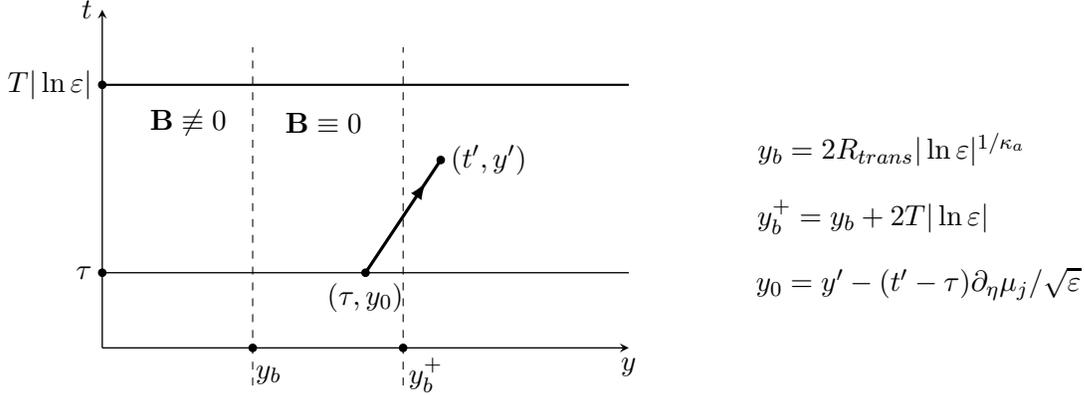
\begin{figure}
\scalebox{1}{
 \begin{tikzpicture}

 \begin{scope}[>=stealth]
 \draw[line width=.5pt] (0,0) -- (7,0);

 \draw[line width=.5pt][->] (0,-1) -- (0,3.5);
 
  \draw[line width=.5pt][->] (0,-1) -- (7,-1);
  
 \end{scope} 
 
 \draw[dashed] (2,-1.5) -- (2,3) ; 
 
 \draw[dashed] (4,-1.5) -- (4,3) ; 
 
 \draw[thick] (0,2.5) -- (7,2.5) ; 
 
 \draw (0,0) node[anchor=east] {$\t$} ; 
 
 \draw (0,2.5) node[anchor=east] {$T |\ln \e|$} ; 
 \draw (7,-1) node[anchor=north] {$y$}; 
 \draw (0,3.5) node[anchor=east] {$t$}; 
 \draw (2.21,-1.1) node[anchor=north]  {$y_b$};
 \draw (4.3,-1) node[anchor=north]  {$y_{b}^+$};

\filldraw[black] (2,-1) circle (0.05cm) ;

\filldraw[black] (0,0) circle (0.05cm) ;
\filldraw[black] (0,2.5) circle (0.05cm) ;

\filldraw[black] (4,-1) circle (0.05cm) ;

\draw (.5,2) node[anchor=west] {${\bf B } \not\equiv 0$};

\draw (2.3,2) node[anchor=west] {${\bf B } \equiv 0$};

\draw (10.5,2) node[anchor=north] {$y_b = 2 R_{trans} |\ln \e|^{1/\kappa_a}$};
\draw (10.25,1.1) node[anchor=north] {$y_b^+ = y_b + 2 T |\ln \e|$};
\draw (10.85,.2) node[anchor=north] {$y_0 = y' - (t' - \t) \d_\eta \mu_j/\sqrt \e$} ; 

\draw (4.5, 1.5) node[anchor=west] {$(t',y')$} ; 
\filldraw[black] (4.5, 1.5) circle (0.05cm) ; 

\draw (3.5,0) -- (4.5,1.5) ; 

\filldraw[black] (3.5, 0) circle (0.05cm) ; 

\draw (3.5,0) node[anchor=north] {$(\t,y_0)$} ; 

\draw[very thick, decoration = {markings, mark = at position .8 with {\arrow{latex}} }, postaction = {decorate} ] (3.5,0) -- (4.5,1.5) ;

 \end{tikzpicture}
}
 \caption{In the $(t,y)$ domain, we find the ${\bf B} \not\equiv 0$ region to the left of $y_b:$ this is the region where the symbolic flow equation has varying coefficients. The ${\bf B} \equiv 0$ region is to the right of $y_b.$ Given $(t',y')$ with $\t \leq t' \leq T \sqrt \e |\ln \e|$ and $|y'| \geq y_b^+,$ the characteristics passing by $(t',y')$ is issued from $(\t,y_0),$ with $|y_0| > y_b,$ since $|\d_\eta \mu_j| < 1.$ In particular, the characteristics is entirely contained in the ${\bf B} \equiv 0$ region. This implies that for $|y| \geq y_b^+,$ the representation \eqref{explicit:S} holds. The characteristics pictured here is incoming, corresponding to $\d_\eta \mu_j > 0;$ for outgoing characteristics the buffer zone $[y_b, y_b^+]$ is not necessary: in the case $\d_\eta \mu_j < 0,$ any $(t',y')$ with $|y'| > y_b$ sees its backward characteristics drawn entirely in the ${\bf B} \equiv 0$ region.} 
\label{fig:incoming}
\end{figure}

 Figure \ref{fig:incoming} illustrates the proof below.

\begin{proof}
 For notational simplicity, we let $(j_1,j_2) = (1,2)$ in this proof, and denote ${\bf B}$ the extra-diagonal terms in \eqref{2blocks}. The analysis applies to any resonant pair. For $(x,y,\xi,\eta)$ to satisfy $\chi_\Omega(x,y,\xi,\eta) = 0 ,$ we have either (a) $\tilde \chi_{long}(x)  = 0 $ or $\chi_{12}^\sharp(\xi,\eta) = 0 ,$ or (b) ${\bm \chi}_{trans}(y) =0.$ 
 
 In case (a), $y$ is arbitrary in $\R^2.$ Consider the partial differential equation satisfied by $S_{f,g}$ in this case. It's an equation in $(t,y) \in [\t, T \sqrt \e |\ln \e|] \times \R^2,$ where $(x,\xi,\eta)$ are parameters. By assumption on $(x,\xi,\eta),$ we have ${\bf B} = 0,$ for all $y.$ The system \eqref{2blocks} splits into two uncoupled transport equations  
 \be \label{case:a}
 \left( \d_t + \frac{i \mu_j}{\e} + \frac{\d_\eta \mu_j}{\sqrt \e} \d_y\right) S^j_{f,g} = f_j, \qquad S_{f,g}^j(\t;\t) = g_j, \quad y \in \R^2, \quad \t \leq t \leq T \sqrt \e |\ln \e|,
 \ee
  with $j = 1$ or $j = 2.$ Following the characteristics backward in time, we find that the unique solution to \eqref{case:a} is given by \eqref{explicit:S}.
  
 In case (b), we also have ${\bf B} = 0$ since $|y| \geq 2 R_{trans} |\ln \e|^{1/\kappa_a}.$ (Indeed, for such $y$ we have $\chi_{trans}(y) = 0,$ by definition of $\chi_{trans}$ in \eqref{def:chi:infty} hence $g_p = 0,$ by definition of $g_p$ in \eqref{def:g}, hence $b_{12}^\pm = 0:$ see \eqref{def:bjj'}). The cancellation ${\bf B} = 0$ is however not uniform in $y:$ for smaller values of $|y|,$ we have ${\bf B} \neq 0$ a priori. The equation in $S_{f,g}$ being a transport equation in $(t,y),$ the behavior of $S_{f,g}$ in the ${\bf B} = 0$ region is potentially affected by the behavior of $S_{f,g}$ in the ${\bf B} \neq 0$ region. 
 
 But by ${\bm \chi}_{trans}(y) = 0,$ we have in fact   $|y| \geq 2( R_{trans} |\ln \e|^{1/\kappa_a} + T |\ln \e|),$ and, since $|\d_\eta \mu_j| \leq 1,$ uniformly in $(\xi,\eta)$ (see Section \ref{sec:char0}), for $0 \leq \t \leq t \leq T \sqrt \e |\ln \e|$, this implies
  $$|y - \d_\eta \mu_j (t - \t)/\sqrt \e| > 2 R_{trans} |\ln \e|^{1/\kappa_a}.$$ 
Thus the characteristics is entirely drawn within the ${\bf B} = 0$ region (see Figure \ref{fig:incoming}), whence the explicit formulation \eqref{explicit:S}.
 \end{proof}

Given a triplet $(j_1,j_2,j_3)$ (as in Section \ref{sec:second:look}: this means $(j_1,j_2,j_3) = (1,3,5)$ or $(j_1,j_2,j_3) = (2,3,4)$), the associated truncation is 
\be \label{def:chi:Omega:triplet}
  \chi_\Omega(x,y,\xi,\eta) = \tilde \chi_{long}(x) {\bm \chi}_{trans}(y) \big( \chi_{j_1j_2}^\sharp(\xi,\eta) + \chi_{j_2j_3,-1}^\sharp(\xi,\eta)),\ee 
  and then we define $\Omega$ as in \eqref{def:Omega}. This $\Omega$ is compact just like the ones associated with resonant pairs.
  
 \begin{lem}[Far-field behavior for \eqref{3blocks}] \label{lem:away:3} For the solution $S_{f,g}$ to \eqref{3blocks} with a source $f$ and a datum $g:$ for each triplet $(j_1,j_2,j_3),$ we have an associated compact domain $\Omega$ (see above) and an explicit description outside of $\Omega$ that is analogous to \eqref{explicit:S} and implies in particular that $S_{f,g}$ does not grow in time outside of $\Omega,$ for $0 \leq \t \leq t \leq T \sqrt \e |\ln \e|.$ 
 \end{lem}

\begin{proof} For $(x,y,\xi,\eta)$ away from $\Omega$ defined in terms of \eqref{def:chi:Omega:triplet}, the system \eqref{3blocks} is a family of three uncoupled, constant-coefficient partial differential equations in $(t,y),$ which we solve by following the characteristics backward in time. As in the proof of Lemma \ref{lem:away}, the $y$ truncation ${\bm \chi}_{trans}$ ensures that within the small time window that we consider, the characteristics are entirely drawn in the constant-coefficient region. 
\end{proof}

\subsection{Existence and uniqueness for the symbolic flow} \label{sec:existence:uniqueness:flow}

We consider now the near-field behavior of $S,$ that is inside $\Omega$ defined in \eqref{def:Omega}. This means in particular that we consider a particular subsystem, whether in the form \eqref{2blocks} or in the form \eqref{3blocks}. It is convenient to use the notation ${\bf A}$ and ${\bf B}$ from the original system \eqref{resolvent0-in-proof} to denote also the operators in that particular subsystem, and we do so here. We have
$$ \d_t (\chi_\Omega^\sharp S_{f,g}) + \frac{1}{\e} \big( i {\bf A} - \sqrt \e {\bf B} \big) (\chi_\Omega^\sharp S_{f,g}) + \frac{1}{\sqrt \e} \d_\eta {\bf A} \d_y (\chi_\Omega^\sharp S_{f,g}) = \frac{1}{\sqrt \e} \d_\eta {\bf A} (\d_y \chi_\Omega^\sharp) S_{f,g} + \chi_\Omega^\sharp f,$$
with $\chi_\Omega \prec \chi_\Omega^\sharp.$
If we let 
$$ S_\Omega := \chi_{\Omega}^\sharp S_{f,g}, \qquad S^c_{\Omega} := (1 - \chi_\Omega^\sharp) S_{f,g},$$
  we obtain
\be \label{eq:S:Omega} \begin{aligned} 
 \left( \d_t + \frac{1}{\sqrt \e} \d_\eta {\bf A} \d_y\right) S_\Omega & + \left( \frac{1}{\e} \big( i {\bf A} - \sqrt \e {\bf B} \big) - \frac{1}{\sqrt \e} \d_\eta {\bf A} (\d_y \chi_\Omega^\sharp) \right) S_\Omega \\ & = \frac{1}{\sqrt \e} \d_\eta {\bf A} (\d_y \chi_\Omega^\sharp) S^c_\Omega + \chi_\Omega f.
\end{aligned} \ee  
The equation \eqref{eq:S:Omega} is linear and symmetric hyperbolic, since $\d_\eta {\bf A}$ is real diagonal. There are two source terms. First, we find $(\d_y \chi_\Omega^\sharp) S_\Omega^c,$ where $S_\Omega^c$ is explicitly given by Lemma \ref{lem:away}, in terms of $f$ and $g.$ The function $\d_y \chi_\Omega^\sharp$ is compactly supported. Thus $(\d_y \chi_\Omega^\sharp) S_\Omega^c$ is compactly supported, and has the pointwise regularity of $f$ and $g.$ The other source term is $\chi_\Omega^\sharp f:$ also compactly supported, with the same pointwise regularity as $f.$ 

The classical theory of linear symmetric hyperbolic system applies: there exists a unique solution $S_\Omega$ to \eqref{eq:S:Omega}, with a large degree of pointwise regularity in all variables, since the coefficients ${\bf A}$ and ${\bf B}$ are regular, and so are the source $f$ and datum $g.$ 

In conjunction with Lemmas \ref{lem:away} and \ref{lem:away:3}, this proves the existence and uniqueness of a solution $S_{f,g}$ to \eqref{symbolic:flow:source} and $S$ to \eqref{resolvent0-in-proof}.

\subsection{Control of the $L^2$ norm by the $L^\infty$ norm } \label{sec:S:L2:infty}

Another consequence of Lemmas \ref{lem:away} and \ref{lem:away:3} describing the far-field behavior of $S$ is the following bound, which gives a control of the $L^2(\R^2_y)$ norm of $S,$ pointwise in $(x,\xi,\eta),$ by its supremum in $y:$ 

\begin{lem} \label{lem:S:L2:Linfty} Given $(x,\xi,\eta) \in \R \times \R^3,$ we have, for $S_{f,g}$ solution to \eqref{symbolic:flow:source}, for $0 \leq \t \leq t \leq T \sqrt \e |\ln \e|$ and $|\a| < s_a - d/2:$ 
$$\begin{aligned} \| \d_{x,y}^\a & S_{f,g}(\t;t,x,\cdot,\xi,\eta) \|_{L^2(\R^2_y)}  \lesssim |\ln \e|^\star \| \d_{x,y}^\a S_{f,g}(\t;t,\cdot,\xi,\eta) \|_{L^\infty(\R^2_y)} \\ & + |\ln \e|^\star \left( \| \d_{x,y}^\a g(t,x,\cdot,\xi,\eta) \|_{L^2(\R^2_y)} + \int_\t^t \| \d_{x,y}^\a f(\t;t',x,\cdot,\xi,\eta) \|_{L^2(\R^2_y)} \right).\end{aligned}$$  
\end{lem} 

\begin{proof} Given $(x,\xi,\eta) \in \R \times \R^3,$ if $(\xi,\eta) \notin \Omega_{freq},$ or $x$ is not in the support of $\tilde \chi_{long},$ then the explicit representation of Lemma \ref{lem:away} or \ref{lem:away:3} holds. Transverse derivatives bring out $(t - \t)/\sqrt \e$ factors, which are $O(|\ln \e|^\star).$ Hence the desired bound in that case.

Given $(\xi,\eta) \in \Omega_{freq},$ and $x$ in the support of $\tilde \chi_{long},$ we decompose $\R^2_y$ into the support of ${\bm \chi}_{trans}$ and its complement. In the complement, the representation \eqref{explicit:S} holds, hence a control of the $L^2((\mbox{supp}{\bm \chi}_{trans})^c)$ norm of $\d_{x,y}^\a S$ by the $L^2(\R^2_y)$ norms of $\d_{x,y}^\a g$ and $\d_{x,y}^\a f.$ On the other hand, the support of ${\bm \chi}_{trans}$ has size $O(|\ln \e|^\star),$ hence the bound 
 $$  \| \d_{x,y}^\a  S_{f,g}(\t;t,x,\cdot,\xi,\eta) \|_{L^2(\mbox{\footnotesize supp}\,{\bm \chi}_{trans})} \lesssim |\ln \e|^\star \| \d_{x,y}^\a S_{f,g}(\t;t,x,\cdot,\xi,\eta) \|_{L^\infty(\R^2_y)},$$
 and the result follows. 
 \end{proof}

\section{Upper bounds for the symbolic flow} \label{sec:flow}

In this section, we will derive the optimal growth rate in time for the symbolic flow. We do so step by step, starting from suboptimal bounds. Keep in mind that, as shown in the previous section, the symbolic flow $S_{f,g}$ does not grow in time outside of $\Omega$, which is the compact support of $\chi_{\Omega}$ defined in \eqref{def:Omega}. The cut-off $\chi_{\Omega}$ is defined in \eqref{def:chi:Omega:pair} or \eqref{def:chi:Omega:triplet}, depending on the context. The frequency set $\Omega_{freq}$ is the projection of $\Omega$ onto $\R^{3}_{\xi, \eta}$.

\subsection{A rough upper bound} \label{bd:flow}
We consider first the symbolic flow equation in the form \eqref{2blocks}. For definiteness we write $j_1 = 1$ and $j_2 = 2,$ but the analysis applies to all resonant pairs. The bound that we prove here is not optimal but it is a stepping stone in the analysis.

Given $f,g$ as in \eqref{symbolic:flow:source} (so that $f$ represents a ``source" in the symbolic flow equation, and $g$ a datum), given a frequency domain $U_{freq} \subset \R^3_{\xi,\eta},$ we denote for $\a \in \N^3:$  
\be \label{def:Cinfty}
 N_{\infty,\a}(\t;t,f,g;\R^3 \times U_{freq}) 
 :=  \max_{\b \leq \a} \left( \| \d_x^\b g \|_{L^\infty(\R^3 \times U_{freq})} + \int_\t^t  \| \d_{x,y}^\b f(\t;t')\|_{L^\infty(\R^3 \times U_{freq})}\right),
 \ee 
and
\be \label{def:C2} \begin{aligned}
N_{2,\a}(\t;t,f,g;\R^3 & \times U_{freq} ) \\ & :=  \max_{\b \leq \a}  \left( \| \d_x^\b g \|_{L^\infty(\R_x \times U_{freq}, L^2(\R^2_y))} + \int_\t^t  \| \d_{x,y}^\b f(\t;t')\|_{L^\infty(\R_x \times U_{freq}, L^2(\R^2_y))}\right)\,.
\end{aligned}\ee
\begin{defi}[Growth rate] \label{def:growth:rate} Given $U_{freq} \subset \R^3_{\xi,\eta},$ we let $U = \R^3_{x,y} \times U_{freq},$ and say that $\g(U) \geq 0$ {\it is a growth rate for the symbolic flow} over $U$ if we have the bounds, for $|\a| < s_a - d/2:$ 
 \be \label{bd:rough} \begin{aligned} \| & \d_{x,y}^\a  S_{f,g}(\t;t,\cdot)\|_{L^\infty(U)} \\ &  \lesssim |\ln\e|^\star  \Big( N_{\infty,\a}(\t;t,f,g;\R^3 \times U_{freq})  +  N_{2,\a}(\t;t,f,g;\R^3 \times U_{freq}) \Big)  \\ & + \frac{|\ln \e|^\star}{\sqrt \e}  \int_\t^t e^{(t - t') \g(U)/\sqrt \e}  \Big( N_{\infty,\a}(\t;t',f,g;\R^3 \times U_{freq})  +  N_{2,\a}(\t;t',f,g;\R^3 \times U_{freq}) \Big)  \, dt' ,\end{aligned}\ee 
  where the implicit constant depends only on the $L^\infty(U)$ norms of $b_{12},$ $b_{21},$ but not on $\e,\t,t,x,y,\xi,\eta.$  
\end{defi}

In the above definition, we used notation $|\ln \e|^\star:= |\ln \e|^{N_\star}$ for some $N_\star > 0$ which depends on all parameters in the problem, in particular on $\a,$ but not on $\e,\t,t,x,y,\xi,\eta.$

\begin{lem}[Rough bound for the symbolic flow with a source] \label{lem:rough:S} For the symbolic flow equation in the form \eqref{2blocks} with a source $f$ and datum $g$ as in \eqref{symbolic:flow:source}, given any $U$ as in Definition {\rm \ref{def:growth:rate}},  
\be \label{def:gamma2}
\gamma_{j_1j_2}^+(U) := \Big( \| b_{j_1j_2}^+ \|_{L^\infty(U)} \cdot \| b_{j_2j_1}^- \|_{L^\infty( U)} \Big)^{1/2},
\ee
is a rate of growth. 
\end{lem}

\begin{proof} For notational simplicity, we let $(j_1,j_2) = (1,2)$ in this proof.
We coordinatize $S_{f,g} = (S^1,S^2)$ (keeping in mind Note \ref{note:dimensions}) and factor out the fast oscillations:
\be\label{S-new-1-2}
\tilde S^1(\tau;t): = e^{ i \mu_1 (t-\tau)/\e } S^1(\tau;t),\quad \tilde S^2(\tau;t): = e^{ i\mu_2(t-\tau)/\e} S^2(\tau;t).
\ee
The system in $(\tilde S^1, \tilde S^2)$ is 
\be\label{S-new-1}
\left\{\begin{aligned}
\d_t \tilde S^1 -  \frac{1}{\sqrt \e} \chi b_{12}^+ e^{i (\mu_1 - \mu_2) (t - \tau)/\e}  \tilde S^2  + \frac{1}{\sqrt\e} \d_\eta \mu_1  \cdot \d_y \tilde S^1 & = e^{i \mu_1 (t - \tau)/\e} f_1 =: \tilde f_1,\\
\d_t \tilde S^2 - \frac{1}{\sqrt\e} \chi  b_{21}^{-}  e^{ - i (\mu_1 - \mu_2) (t - \tau)/\sqrt \e}  \tilde S^1  + \frac{1}{\sqrt\e} \d_\eta \mu_2 \cdot \d_y \tilde  S^2& = e^{i \mu_2 (t - \tau)/\e} f_2 =: \tilde f_2,\\
\end{aligned}\right.
\ee
where $\chi$ is short for $\chi_{long}^\sharp,$ the longitudinal truncation in $x.$ We will keep this convention throughout the present Section on $S$ bounds. 
We let
\be\label{S-new-1-2-2}
\breve S^j(\tau;t,x,y,\xi,\eta)  : = \tilde S^j\left(\tau;t,x, y_j,\xi,\eta\right),  \quad y_j = y + \frac{\d_\eta \mu_j}{\sqrt \e} (t-\tau), \qquad j = 1,2,
 \ee
and
\be \label{velocity}
 \nu : = \nu_{12}, \quad \nu_{j_1j_2} := \d_\eta (\mu_{j_1} - \mu_{j_2})
\ee
wich yields
\be \label{sys:breve} \left\{ 
\begin{aligned} 
\d_t \breve  S^1(\tau;t,x,y,\xi,\eta)  - & \frac{e^{ i (\mu_1 - \mu_2) (t - \tau)/\e}}{\sqrt \e} (\chi b_{12}^+)\left( x, y_1,\xi,\eta \right)  \breve  S^2 (\tau; t, x, y + \frac{\nu t}{\sqrt \e},\xi,\eta)  \\ &  = \breve f_1\left(\tau;t,x, y,\xi,\eta\right), \\
\d_t \breve  S^2(\tau;t,x,y,\xi,\eta)  - & \frac{e^{- i (\mu_1 - \mu_2) (t - \tau)/\e}}{\sqrt \e} (\chi  b_{21}^-) \left( x, y_2,\xi,\eta \right)   \breve  S^1 (\tau; t, x, y - \frac{\nu t}{\sqrt \e},\xi,\eta) \\ & =  \breve f_{2} \left(\tau;t,x, y,\xi,\eta\right),
\end{aligned} \right.
\ee 
where $\breve f_j$ are defined in terms of $\tilde f_j$ just like $\breve S^j$ in terms of $\tilde S^j.$ 
With
\be \label{notation:rough} b_1 := \| b_{12}^+ \|_{L^\infty(U)}, \qquad b_2  := \| b_{21}^- \|_{L^\infty(U)}, \qquad s_j := \| S^j(\t;t) \|_{L^\infty(U)},\ee 
this gives
\be \label{S-new-4} \left\{ \begin{aligned}  
s_{1}(\tau;t) & \leq \| g_1 \|_{L^\infty} +    \int_{\tau}^{t} \, ( \frac{b_1}{\sqrt \e} s_2 (\tau; t') + \| f_1(\tau;t')\|_{L^\infty}) \,dt', \\
s_{2}(\tau;t) & \leq \| g_2 \|_{L^\infty} +    \int_{\tau}^{t} \, ( \frac{b_{2}}{\sqrt \e} s_1 (\tau; t') + \| f_2(\tau;t')\|_{L^\infty})  \,dt'.
\end{aligned}\right.\ee
Above the sup norms in $f$ and $g$ are naturally with respect to $U.$ At this point we may assume $b_1 \neq 0,$ $b_2 \neq 0$ (otherwise the bound in the Lemma obviously holds) and we deduce from \eqref{S-new-4} that
$\tilde s := \sqrt{b_2} s_1 + \sqrt{b_1} s_2$
satisfies 
$$ \begin{aligned}
\tilde s(\tau;t) & \leq (\sqrt{b_{1}}  + \sqrt{b_{2}} ) \| g \|_{L^\infty}   +  (\sqrt b_1 +  \sqrt b_2) \int_\t^t \| f(\t;t') \|_{L^\infty} ,dt \\ & \quad + \frac{\sqrt{ b_1 b_2} }{\sqrt \e} \int_{\tau}^{t} \, \tilde s (\tau; t') \, dt' . \end{aligned} $$
Gronwall's inequality concludes the proof in the case $\a = 0.$

The case $|\a| > 0$ follows by a simple induction argument, based on the equation satisfied by $\d_{x,y}^\a S_{f,g}.$ Indeed, we observe that $S_{\a,f,g} := \d_{x,y}^\a S_{f,g}$ solves
$$ \left\{ \begin{aligned} \d_t S_{\a,f,g} + \frac{1}{\e} (i {\bf A} - \sqrt \e {\bf B}) S_{\a,f,g}  + \frac{1}{\sqrt\e} [\d_{x,y}^\a, {\bf B}] S_{f,g} + \frac{1}{\sqrt\e} \d_\eta {\bf A} \d_y S_{\a,f,g}  = \d_{x,y}^\a f, \\ S_{\a,f,g}(\t;\t) = \d_{x,y}^\a g.
\end{aligned}\right.$$
That is (using notation $S_{f,g}$ to denote the solution to \eqref{symbolic:flow:source})
$$ \d_{x,y}^\a S_{f,g} = S_{f_\a, g_\a},$$
with notation
$$
 f_\a  := \d_{x,y}^\a f - \frac{1}{\sqrt \e} [\d_{x,y}^\a, {\bf B}] S_{f,g}, \ \  g_\a  := \d_{x,y}^\a g.
$$
Thus, we obtain, by the first step ($\a = 0$):
 \be \label{for:rough} \begin{aligned} \| \d_{x,y}^\a & S_{f,g}(\t;t,\cdot)\|_{L^\infty}  \lesssim |\ln\e|^\star \left( \| \d_{x,y}^\a g \|_{L^\infty} + \int_\t^t \| f_\a(\t;t'')\|_{L^\infty} \, dt''\right) \\ & \quad + \frac{|\ln \e|^\star}{\sqrt \e}  \int_\t^t \exp\left( \frac{(t - t') \g^+}{\sqrt \e} \,\right) \left(\| \d_{x,y}^\a g \|_{L^\infty} + \int_\t^{t'} \|  f_\a(\t;t'')\|_{L^\infty} \, dt'' \right) \, dt' ,\end{aligned}\ee
where $\g^+$ is short for $\g^+_{12}(U).$ The commutator $[\d_{x,y}^\a, {\bf B}] S_{f,g}$ in $f_\a$ is handled by the induction argument, since it involves strictly less than $\a$ spatial derivatives of $S_{f,g}.$ That is, by the induction hypothesis, 
\be \label{for:rough:2} \| f_\a(\t;t') \|_{L^\infty} \lesssim \| \d_{x,y}^\a f \|_{L^\infty} + \frac{|\ln \e|^\star}{\sqrt \e} e^{(t' - \t) \g^+/\sqrt \e} \sum_{\b \leq \a} \left( \| \d_{x,y}^\b g\|_{L^\infty} + \int_\t^{t'} \| \d_{x,y}^\b f \|_{L^\infty}\right).\ee
We then plug \eqref{for:rough:2} into \eqref{for:rough} and take advantage of the smallness of the time interval, to arrive at the result. 
\end{proof} 

\begin{rem} \label{rem:for:Duh} The existence of a growth rate for which the bounds of lemma {\rm \ref{lem:rough:S}} hold true does not rely on the structure of ${\bf B},$ only on the real diagonal structure of ${\bf A}.$ The sharp rate of growth found in Corollary {\rm \ref{cor:optimal:bound:S}} does however strongly depend on the structure of ${\bf B}.$ 
\end{rem}

Over frequency sets for which the equation takes the form \eqref{3blocks}, the same analysis applies:

\begin{lem} \label{cor:rough:2}
For the symbolic flow equation in the form \eqref{3blocks} with a source $f$ and datum $g$ as in \eqref{symbolic:flow:source}, given any $U$ as in Definition {\rm \ref{def:growth:rate}}, a rate of growth is given by 
\be \label{def:gamma3}
\gamma^+_{j_1j_2j_3}(U) :=
\Big( \| b_{j_1j_2}^+  \|_{L^\infty} \cdot \| b_{j_2j_1}^-  \|_{L^\infty} + \| b_{j_2j_3}^+ \|_{L^\infty} \cdot \| b_{j_3j_2}^- \|_{L^\infty} \Big)^{1/2},
\ee
with sup norms taken over $U.$
\end{lem}

\begin{proof} Essentially it suffices to follow the proof of Lemma \ref{lem:rough:S}, and account for the slightly different matrix structure of ${\bf B}.$ We give however the corresponding system in $\breve S^j$ here since it will be useful in the proof of Lemma \ref{lem:away:space:time:3}. We coordinatize $S_{f,g} = (S^1,S^2,S^3),$ and define as in the proof of Lemma \ref{lem:rough:S}:
 $$ \breve S^q = e^{i \mu_{j_q} (t - \t)/\e} S^q\left(\t;t,x, y_{q}, \xi, \eta \right), \quad y_{q} := y + \frac{\d_\eta \mu_{q}}{\sqrt \e} (t - \t), \quad q \in \{1, 2, 3\}.$$
 This leads to 
 \be \label{breve:3} \left\{
\begin{aligned}
\d_t \breve  S^1 - \frac{e^{ i (\mu_{j_1} - \mu_{j_2}) (t - \tau)/\e}}{\sqrt \e} \chi b_{j_1j_2}^+(y_1) \breve  S^2 (y_{12})  &  = \breve f_1, \\
\d_t \breve  S^2 - \frac{e^{i (\mu_{j_2} - \mu_{j_1}) (t - \tau)/\e}}{\sqrt \e} \chi  b_{21}^-(y_2)  \breve  S^1(y_{21}) -  \frac{e^{ i (\mu_2 - \mu_3) (t - \tau)/\e}}{\sqrt \e} \chi b_{j_2j_3}^+(y_2)  \breve  S^3(y_{23}) & =  \breve f_2,\\
\d_t \breve S^3 -  \frac{e^{ i (\mu_3 - \mu_1) (t - \tau)/\e}}{\sqrt \e} \chi  b_{j_3j_2}^-(y_3)  \breve  S^2 (y_{32}) & =  \breve f_3,
\end{aligned} \right.
\ee
where $\chi$ is short for $\chi_{long}^\sharp,$ and $y_{j_1j_2} := y + \nu_{j_1j_2} (t - \t)/\sqrt \e,$
and $\breve f_j$ are defined in terms of $f_j$ just like $\breve S^j$ in terms of $S^j.$ We can then argue as in the proof of Lemma \ref{lem:rough:S}.  
\end{proof}

\begin{cor}[Rough bound for the symbolic flow] \label{cor:rough}
The solution $S$ to the symbolic flow system \eqref{resolvent0-in-proof} satisfies the bounds, for $0 \leq \t \leq t \leq T \sqrt \e |\ln \e|$ and $|\a| < s_a - d/2:$ 
 $$ \|\d^{\a}_{x,y} S(\t;t) \|_{L^\infty(\R^6_{x,y,\xi,\eta})} \lesssim |\ln \e|^\star e^{(t  - t')\g^+/\sqrt \e},$$
 where 
 \be \label{def:gamma+}
 \g^+ := \max_{(j_1,j_2),(j_1,j_2,j_3)} \big( \g^+_{j_1j_2}(\R^6), \g^+_{j_1j_2j_3}(\R^6) \big),
 \ee 
 the maximum being taken over all resonant pairs $(j_1,j_2)$ (listed in Proposition {\rm \ref{prop:res}}) and all non-separated triplets $(j_1,j_2,j_3):$ that is, according to Proposition {\rm \ref{prop:separation}(3)}: $(j_1,j_2,j_3) = (1,3,5)$ and $(j_1,j_2,j_3) = (2,3,4).$ 
\end{cor}

\begin{proof} It suffices to use the block decomposition of Section \ref{sec:flow:blocks}, and apply Lemma \ref{lem:rough:S} with $f = 0,$ $g = {\rm Id},$ and $U = \R^6$ for blocks of size 2, and Lemma \ref{cor:rough:2} with $f = 0,$ $g = {\rm Id},$ and $U = \R^6$ for blocks of size 3.
\end{proof}

\subsection{Away from space-time resonances} \label{sec:away:space-time}

We consider first subsystems of the form \eqref{2blocks}. Given a resonant pair $(j,j'),$ recall that the domain $\Omega = \Omega_{sp} \times \Omega_{freq} \subset \R^3 \times \R^3$ was introduced in \eqref{def:Omega}, and that Lemma \ref{lem:away} gave an explicit description of $S_{f,g}$ outside of $\Omega.$   
 Consider now the domain
 \be \label{def:domain:away:space-time}
 U^{away}_{\delta} = \big\{ (\xi,\eta) \in \Omega_{freq}, \quad |\nu_{jj'}(\xi,\eta)| \geq \delta \big\},
 \ee
 for $\delta > 0,$ where notation $\nu_{jj'} = \d_\eta (\mu_{j} - \mu_{j'})$ is introduced in \eqref{velocity}. Just like we didn't use subscripts for $\Omega,$ we do not for $U^{away}_{\delta}.$ This is convenient, but we need to keep in mind that $U^{away}_{\delta}$ depends on the pair (for now) or triplet (later) that is being considered.

\begin{lem} \label{lem:away:space:time} For any resonant pair $(j,j'),$ any $\delta > 0,$ we have for the symbolic flow system in the form \eqref{2blocks} a rate of growth $\g_{jj'}(\R^3 \times U^{away}_\delta),$ in the sense of Definition {\rm \ref{def:growth:rate}}, which satisfies
$$ \lim_{\e \to 0} \g_{jj'}(\R^3 \times U^{away}_\delta)  = 0.$$
\end{lem}

\begin{proof} Here again, we let $(j,j') = (1,2),$ but the analysis applies to all resonant pairs. We shorten $\nu_{12} = \nu.$ Consider $(\xi,\eta)$ such that $|\nu(\xi,\eta) | \geq \delta.$ We go back to the proof of Lemma \ref{lem:rough:S}, specifically to the system in $\breve S^j,$ and use notation $y_j$ introduced just above \eqref{velocity}.

In the system in $\breve S^j,$ we find $b_{12}^+$ evaluated at $y_1,$ and $b_{21}^-$ evaluated at $y_2.$ We have
$$ |y_1 - y_2| \geq \frac{\delta (t - \t)}{\sqrt \e}.$$
This implies 
$$ |y_j| \geq \frac{\delta(t - \t)}{2 \sqrt \e}, \quad \mbox{for $j = 1$ or $j = 2.$}$$
Assume without loss of generality that $j = 1$ above. By condition (a3) bearing on the WKB amplitude (see page \pageref{page:a3}) and linearity of $b_{ij}^\pm$ in $g_\pm,$ this implies 
$$ |b_{12}^+(y_1)| \lesssim \exp\left( - r_a \left| \frac{\delta(t - \t)}{2 \sqrt \e}\right|^{\kappa_a}\right).$$
Given the contraint $t - \t \leq T \sqrt \e |\ln \e|,$ this implies
$$ |b_{12}^+(y_1)| \lesssim \exp\left( - C |\ln \e|^{\kappa_a}\right),$$
with $C = C(r_a, \delta, \kappa_a) > 0.$ In particular, the associated rate $\g_{12}^+(\R^3 \times \{ \xi,\eta\})$ \eqref{def:gamma2} tends to 0 as $\e \to 0,$ uniformly in $(\xi,\eta)$ with $|\nu(\xi,\eta)| \geq \delta.$  
 \end{proof} 

We now perform a similar ``away from space-time resonances" analysis for the subsystems \eqref{3blocks}. In this view, given $(j_1,j_2,j_3)$ equal either to $(1,3,5)$ or to $(2,3,4),$ consider the associated $\Omega$ defined in \eqref{def:chi:Omega:triplet} and \eqref{def:Omega}, its projection $\Omega_{freq}$ over the frequency domain, and 
\be \label{def:away:triplets}
 \begin{aligned} U^{away}_{\delta} :=  \big\{ (\xi,\eta) \in \Omega_{freq}, \quad |\nu_{j_1j_2}(\xi,\eta)| \geq \delta \quad \mbox{or} \quad |\nu_{j_2j_3}(\xi,\eta)| \geq \delta \big\}. \end{aligned}
 \ee

 \begin{lem} \label{lem:away:space:time:3} Given $(j_1,j_2,j_3)$ equal to $(1,3,5)$ or to $(2,3,4),$ we have for the symbolic flow system in the form \eqref{3blocks} a rate of growth $\g_{j_1j_2j_3}(\R^3 \times U^{away}_{\delta}),$ in the sense of Definition {\rm \ref{def:growth:rate}}, which satisfies
$$  \g_{j_1j_2j_3}(\R^3 \times U^{away}_{\delta}) \leq \g_\e + \max\big(\g_{j_1j_2}(\R^6), \g_{j_2j_3}(\R^6)\big),$$
with $|\ln \e|^\star \g_\e \to 0$ as $\e \to 0,$ for \emph{any} rate of growth $\g_{j_1j_2}$ for $(j_1,j_2)$ and $\g_{j_2j_3}$ for $(j_2,j_3).$  
\end{lem}

The key is that the above bound for $\g_{j_1j_2j_3}$ is valid for any admissible rate for $(j_1,j_2)$ and $(j_2,j_3).$ 

\begin{proof} We follow the proof of Lemma \ref{lem:away:space:time}, based on the $\breve S^j$ formulation \eqref{breve:3} of the system \eqref{3blocks}.

 Given $(\xi,\eta)$ such that $|\nu_{j_1j_2}(\xi,\eta)| \geq \delta,$ we argue as in the first part of the proof of Lemma \ref{lem:away:space:time} and find that either $|b_{j_1j_2}^+(y_1)| \lesssim \exp\left( - C |\ln \e|^{\kappa_a}\right),$ or $|b_{j_2j_1}^-(y_2)| \lesssim  \exp\left( - C |\ln \e|^{\kappa_a}\right),$ for any $y.$ Assume for now that it's $|b_{j_1j_2}^+|$ that is small. 
 
 We integrate in time the equation in $\breve S^1$ from \eqref{breve:3}, and find 
 \be \label{breveS1} \begin{aligned}
 & \breve  S^1(\tau;t,x,y, \xi,\eta) - \breve g_1 \\ \,\, & = \int_\t^t \frac{e^{ i (\mu_{j_1} - \mu_{j_2}) (t' - \t)/\e}}{\sqrt \e} (\chi b_{j_1j_2}^+)\left( x, y_1(\t;t'),\xi,\eta \right)  \breve  S^2 (\tau; t', x, y + \frac{\nu_{j_1j_2} (t' - \t)}{\sqrt \e},\xi,\eta) \, dt' \\ \,\, & + \int_\t^t \breve f_1\left(\tau;t',x, y_1(\t;t'),\xi,\eta\right) \, dt', 
 \end{aligned}\ee
 where as before $y_k(\t;t') := y + (\d_\eta \mu_{j_k}) (t' - \t)/\sqrt \e.$ Denoting $L^\infty$ the $L^\infty(\R^3 \times \{ (\xi,\eta) \})$ norm, and omitting the initial time $\t,$ this gives
 \be \label{bd:breveS1:3} \| S^1(t) \|_{L^\infty} \leq \| g_1 \|_{L^\infty} + \int_\t^t \| f_1 \|_{L^\infty} + \g_\e \max_{\t \leq t' \leq t} \| S^2(t') \|_{L^\infty},\ee 
 for some $\g_\e \geq 0$ which tends to $0$ as $\e \to 0,$ and such that $|\ln \e|^\star \g_\e \to 0$ as $\e \to 0.$ 
 Above, we used the bound $t \leq T |\ln \e|,$ and the fact that $|\ln \e| e^{- C |\ln \e|^{\kappa_a}} \to 0$ as $\e \to 0.$  
  
 And now we see the system in $(\breve S^j)_{j = 1,2,3}$ \eqref{breve:3} as a system in $(\breve S^j)_{j = 2,3},$ with a source term which includes the contribution of $\breve S^1.$ A rate of growth $\g_{j_2j_3} = \g_{j_2j_3}(\R^3 \times \{ (\xi,\eta) \})$ for a fixed $(\xi,\eta)$ with $|\nu_{j_1j_2}(\xi,\eta)| \geq \delta$ is also a rate of growth for the $(\breve S^2, \breve S^3)$ system, for the same $(\xi,\eta).$ (Note that it makes sense for consider a domain $\R^3 \times \{ (\xi,\eta)\},$ where the frequency domain is a singleton, for a growth rate as in Definition \ref{def:growth:rate}, since the frequency $(\xi,\eta)$ is only a parameter in the symbolic flow system.) Thus, denoting $\|S^{23}\|_{L^\infty} =  \| S^2 \|_{L^\infty} + \| S^3 \|_{L^\infty},$ with the same convention as above for $L^\infty,$ we obtain the bound 
 $$ \begin{aligned} \| S^{23} & \|_{L^\infty} \lesssim \|  g \|_{L^\infty} + \int_\t^t \| f \|_{L^\infty} + \frac{1}{\sqrt \e }\int_\t^t e^{(t - t') \g_{j_2j_3}/\sqrt \e} \Big( \|  g \|_{L^\infty} + \int_\t^t \| f \|_{L^\infty} \Big)\, dt' \\ & + \frac{1}{\sqrt \e} \int_\t^t e^{(t - t') \g_{j_2j_3}/\sqrt \e} \int_\t^{t'} \Big( \| g_1 \|_{L^\infty} + \int_\t^{t''} \| f_1 \|_{L^\infty} + \g_\e \max_{\t \leq t'' \leq t'} \| S^2(t'') \|_{L^\infty} \Big) \, dt'' \, dt'.\end{aligned}$$ 
  Above, we omitted the $|\ln \e|^\star$ from Definition \ref{def:growth:rate}, in order to gain space. From the above bound, it is easy to deduce, for a given $(\xi,\eta)$ with $|\nu_{j_1j_2}(\xi,\eta)| \geq \delta,$ the inequality 
  \be \label{for:3} \g_{j_1j_2j_3}(\R^3 \times \{ (\xi,\eta)\}) \leq |\ln \e| \g_\e + \g_{j_2j_3}(\R^3 \times \{ (\xi,\eta) \}),\ee
  where $|\ln \e|^\star \g_\e \to 0$ as $\e \to 0,$ as noted above.  
  Thus we covered the subcase $|b_{j_1j_2}^+|$ small in the case $|\nu_{j_1j_2}| \geq \delta.$ 
  
  Assuming $|\nu_{j_1j2}(\xi,\eta)| \geq \delta$ still, consider now the situation that $|b_{j_2j_1}^-|$ is small. The same arguments apply: we consider the subsystem in $(\breve S^2,\breve S^3),$ where $\breve S^1$ is a source term, which we can bound in terms of $\breve S^2.$ In the end it's the product $|b_{12}^+ b_{21}^-|$ that matters, and we find again the bound \eqref{for:3}.

\end{proof}

\subsection{Further separation analysis for the $(2,3,4)$ system} \label{sec:234} 

Here we focus on system \eqref{3blocks} with $(j_1,j_2,j_3) = (2,3,4)$ and use Proposition \ref{prop:separation}(4), to show that the rate of growth associated with the $(2,3)-(3,4)$ resonances is no greater than the maximum of the rate of growth for the $(2,3)$ resonance and the one for the $(3,4)$ resonance.

\begin{lem} \label{lem:234} For $\delta > 0$ small enough, with the phases $\Phi_{jj'}$ defined in Proposition \ref{prop:res}, the frequency set
$$ \{ |\Phi_{23}| \leq \delta \} \cap \{ |\Phi_{34,+k}| \leq \delta \} \cap \{|\nu_{23} \leq \delta \} \cap \{ |\nu_{34,+k}| \leq \delta \}$$
is empty.
\end{lem}

\begin{proof} Assume, by contradiction, that we can find $\delta$ arbitrarily small such that the frequency set in the statement contains $\zeta_{\delta}.$  By Proposition \ref{prop:res}, the phase functions $\Phi_{23}$ and $\Phi_{34}$ are proper. Thus from the fact that $\Phi_{23}(\zeta_{\delta})$ is bounded, we deduce that $\zeta_{\delta}$ is bounded. Up to a subsequence, it converges to $\zeta_\infty,$ which belongs to ${\mathcal S}_{23} \cap {\mathcal S}_{34,+k},$ by continuity of $\Phi_{jj'}$ and $\nu_{jj'}.$ But this is in contradiction with \eqref{st:sep:4} in Proposition \ref{prop:space-time:sep}.
\end{proof}

\begin{cor} \label{cor:234} If the constant $\delta_{res}$ that intervenes in the definition of the frequency cut-off around resonances (see Definition {\rm \ref{def:cut-offs}}) is small enough, depending on $\delta > 0$ from Lemma {\rm \ref{lem:234},} then for system \eqref{3blocks} with $(j_1,j_2,j_3) = (2,3,4),$ we have a rate of growth $\g_{234}(\R^6)$ such that 
$$ \g_{234}(\R^6) \leq \g_\e + \max(\g_{23}(\R^6), \g_{34}(\R^6)),$$
with $|\ln \e|^\star \g_\e \to 0$ as $\e \to 0,$ for \emph{any} rate of growth $\g_{23}$ for $(2,3)$ and $\g_{34}$ for $(3,4).$ 
\end{cor}

\begin{proof} Given $(\xi,\eta) \in \Omega_{freq},$ if $(\xi,\eta)$ does not belong to the support of $\chi_{23}$ or to the support of $\chi_{34,-1},$ then the $(2,3)-(3,4)$ system \eqref{3blocks} reduces to a $(2,3)$ system or a $(3,4)$ system \eqref{2blocks}, and the result follows immediately. 

Thus we may restrict to $(\xi,\eta)$ in the intersection of the support of $\chi_{23}$ and the support of $\chi_{34,-1}.$ Then, by \eqref{transp}, the phase functions $|\Phi_{23}|$ and $|\Phi_{34,-1}|$ are bounded from above by $c_0 \delta_{res},$ where $c_0 > 0$ is a fixed constant, and $\delta_{res} > 0$ was introduced in Definition \ref{def:cut-offs}. We choose $\delta_{res}$ such that $c_0 \delta_{res} < \delta,$ where $\delta$ is as in Lemma \ref{lem:234}.  

 With such a choice of $\delta$ and $\delta_{res},$ we observe that $\Omega_{freq} \setminus U^{away}_{\delta}$ is empty, by Lemma \ref{lem:234}. Hence we have $\g_{234}(\R^6) \leq \g_{234}(\R^3 \times U^{away}_{\delta}),$ and it suffices to appeal to Lemma \ref{lem:away:space:time:3}. 
\end{proof}

\subsection{Where we stand so far on the pursuit of an optimal bound for $S$} 

In this short Section we pause for a moment and consider where we stand in the derivation of an upper bound for $S:$
\begin{itemize}
\item we split (in Sections \ref{sec:flow:blocks} and \ref{sec:second:look})  the large original system in $S$ \eqref{resolvent0-in-proof} (equivalently, \eqref{symbolic:flow:source}) into subsystems for resonant pairs: those are \eqref{2blocks}, and subsystems for pairs of resonant pairs which cannot be separated in the sense of Proposition \ref{prop:separation}: those pairs are $(2,3)-(3,4)$ and $(1,3)-(3,5)$ and the corresponding subsystems are \eqref{3blocks}.  
\item With Lemmas \ref{lem:away:3} and Lemma \ref{lem:234} we effectively reduced the question of finding an optimal rate of growth for system \eqref{3blocks} in the case $(j_1,j_2,j_3) = (2,3,4)$ into finding an optimal rate of growth for systems \eqref{2blocks} for resonances $(2,3)$ and $(3,4)$ separately. This throws the triplet $(2,3,4)$ out of the equation.
\item For the remaining triplet $(1,3,5),$ we will eventually (in Section \ref{sec:interaction:coefficients}) see that the rough rate from Lemma \ref{cor:rough:2} is small enough. 
\item Thus our analysis of systems \eqref{3blocks} is complete, and we focus in the following on systems \eqref{2blocks}. For these we proved in Lemma \ref{lem:away:space:time} that not much growth occurs far from space-time resonances. The analysis near space-time resonances is the focus of the next two Sections. 
\end{itemize}

\subsection{In the small-trace region: little growth} \label{sec:small:trace}

We derive here bounds for $S$ for frequency parameters away from the frequency domain $U^{away}_{\delta}$ that was defined above in Section \ref{sec:away:space-time}.

The trace $\mbox{tr}\, b_{j_1j_2}^+ b_{j_2j_1}^-$ of the product of the interaction coefficients plays a crucial role in the argument. From \eqref{def:bjj'}, we have
 $$  \big| \mbox{tr}\, b_{j_1j_2}^+ b_{j_2j_1}^-(x,y,\xi,\eta) \big| = |g_{\pm}(x,y)|^2 |\mbox{tr}\, B_{1j_1j_2} B_{-1j_2j_1}(\xi,\eta)|.$$
 Given a resonant pair $(j_1,j_2),$ and $\delta > 0,$ we let 
 \be \label{def:Omega:delta} U^{small}_{\delta} = \big\{ (\xi,\eta) \in \Omega_{freq} \setminus U^{away}_{\delta}, \quad  \big| \mbox{tr}\, B_{1j_1j_2} B_{-1j_2j_1,+1}(\xi,\eta) \big| \leq \delta \big\}.\ee 
For some $\delta',$ we do not record much growth for $S$ in $\R^3 \times U^{small}_{\delta'},$ and this is due to the fact that near space-time resonances the trace $\mbox{tr}\, B_{1j_1j_2} B_{-1j_2j_1,+1}$ is small only when one of $B_{1j_1j_2}$ or $B_{-1j_2j_1,+1}$ is small:

\begin{lem}[Bound in the small-trace region] \label{lem:small-trace} For any resonant pair $(j_1,j_2),$ some $C > 0,$ any $\delta > 0,$ for some $\delta' \leq \delta,$ if the constant $\delta_{res}$ that appears in the definition of the frequency cut-off $\chi_{j_1j_2}$ is small enough depending on $\delta',$ the rate of growth $\g^+_{j_1j_2}(\R^3 \times U^{small}_{\delta'})$ as defined in \eqref{def:gamma2} satisfies
 $$ \g^+_{j_1j_2}(\R^3 \times U^{small}_{\delta'}) \leq C \delta^{1/2}.$$
\end{lem}

\begin{proof} As we did previously, we let $(j_1,j_2) = (1,2)$ in this proof. 

The computations of Section \ref{sec:interaction:coefficients} show that {\it at space-time resonances} $(\xi,\eta) \in {\mathcal S}_{12}$ the trace $\mbox{tr}\, B_{112} B_{-121,+1}$ is equal to zero only if one of $B_{112}$ or $B_{-121,+1}$ is equal to zero.

If we denote $\Phi = \Phi_{12|\Omega_{freq}}$ the phase function (see \eqref{def:phase:function})  in restriction to $\Omega_{freq},$ and similarly $\nu = (\d_\eta (\mu_1 - \mu_2))_{|\Omega_{freq}},$ and $\mbox{tr} := \big(\mbox{tr}\, B_{112} B_{-121,+1})_{|\Omega_{freq}},$ this means 
 $$ (\nu,\Phi, \mbox{tr})^{-1}( \{ 0 \}) \subset  \left( B_{112|\Omega_{freq}}\right)^{-1}( \{ 0 \}) \cup  \left( B_{-121,+1|\Omega_{freq}}\right)^{-1}( \{ 0 \}).$$
 By continuity of $B_{pjj'},$ this implies that for any $\delta > 0,$ the set
 \be \label{open} \left(B_{112|\Omega_{freq}}\right)^{-1} (B(0,\delta)) \cup  \left( B_{-121,+1|\Omega_{freq}}\right)^{-1}(B(0,\delta)),\ee
 where $B(0,\delta)$ is the ball of center 0 and radius $\delta$ in the target space of $B_{pjj'},$ is an open neigborhood of $(\nu,\Phi)^{-1}(\{ 0 \}).$ By compactness of $\Omega_{freq}$ and continuity of $(\nu,\Phi, \mbox{tr}),$ any open neighborhood of $(\nu,\Phi, \mbox{tr})^{-1}(\{ 0 \}),$ in particular \eqref{open}, contains an open set of the form $(\nu, \Phi, \mbox{tr})^{-1}(B(0,\delta')).$ Now recall that on the support of $\chi_{jj'}^\sharp$ we have the upper bound $|\Phi_{jj'}| \leq c_0 \delta_{res}.$ We choose $\delta_{res}$ such that $c_0 \delta_{res} < \delta'.$ Then, $U^{small}_{\delta'} \subset (\nu, \Phi, \mbox{tr})^{-1}(B(0,\delta')),$ which implies, by the above, that $U^{small}_{\delta'}$ is included in \eqref{open}: given a frequency $(\xi,\eta)$ in $U^{small}_{\delta'},$ one of $|B_{112}|$ or $|B_{-121,+1}|$ is smaller than $\delta.$ Thus the rough bound of Lemma \ref{lem:rough:S} gives an upper rate of growth 
 $$ \g^+_{12}(\R^3 \times U^{small}_{\delta'}) \leq C \delta^{1/2},$$
 where $C = (|g_\pm|_{\infty} \max(|B_{112}|_{L^\infty},|B_{-121}|_{L^\infty}))^{1/2}.$
 
 If it happens that $\delta' > \delta,$ then $U^{small}_{\delta} \subset U^{small}_{1,2,\delta'},$ and, $U^{small}_{1,2,\delta'}$ being included in the open set \eqref{open}, we deduce that $U^{small}_{1,2,\delta}$ itself is included in that open set, which implies that $\delta' = \delta$ is an admissible choice.
 \end{proof}

\subsection{In the large-trace region: an optimal rate of growth}\label{sec:large-trace} Given a resonant pair $(j_1,j_2),$ given $\delta > 0$ and an associated $\delta'$ as in Lemma \ref{lem:small-trace}, the frequency region that is left to consider is 
 \be\label{def:U-delta-large} \begin{aligned} U^{large}_{\delta'} & := \Omega_{freq} \setminus \big( U^{away}_{\delta'} \cup U^{small}_{\delta'} \big)  \\ & = \{ (\xi,\eta) \in \Omega_{freq}, \quad |\nu_{j_1j_2}| < \delta', \quad \big|\mbox{tr}\, B_{1j_1j_2} B_{-1j_2j_1,+1}\big| > \delta'\}.\end{aligned}\ee
In particular, $U^{large}_{\delta'}$ is open and relatively compact in $\R^3_{\xi,\eta}.$ 
We will use the following lemma from \cite{em4}:
 \begin{lem}\label{scalarize} Consider two smooth families $C_{12}$ and $C_{21}$ of square matrices with equal dimensions, defined over an open and relatively compact set ${\mathcal D}$ and such that
\be \label{cond:sym}
 {\rm rank}\, C_{12} \equiv {\rm rank}\,C_{21} \equiv 1 \quad \mbox{and} \quad {\rm tr}\, C_{12} C_{21} \neq 0, \qquad \mbox{over $\bar {\mathcal D}.$} 
 \ee
 Then, there exist smooth scalar maps $c_{12},$ $c_{21}$ and a smooth block-diagonal change of basis $P = \bp P_{11} &0\\0&P_{22}\ep$ such that ${\rm tr}\, C_{12} C_{21} = c_{12} c_{21},$ and, given $\t_{12}, \t_{21} \in \C:$
\be \label{eq:P}
 P^{-1}\left(\begin{array}{cc} 0 & \t_{12} C_{12} \\ \t_{21} C_{21} & 0 \end{array}\right) P =\left(\begin{matrix} 0 & \tilde C_{12} \\ \tilde C_{21} & 0
\end{matrix}\right), \qquad \tilde C_{ij} = \left(\begin{array}{cc} \t_{ij} c_{ij} & 0 \\ 0 & 0_{\C^{(N-1)\times(N-1)}} \end{array}\right),
\ee
and
\be \label{bd:P}
\sup_{z \in {\mathcal D}} |\d_z^\a P(z)| + |\d_{z}^\a P^{-1}(z)| < \infty.
\ee
\end{lem}

\begin{proof} See Lemma 3.19 in \cite{em4}. In Corollary \ref{lem:low} below, a little more information is required about matrix $P.$ The proof of Lemma 3.19 in \cite{em4} shows that the change of basis $P$ is defined by columns as
$$ P = \mbox{col}\,\big(e^\sharp, a_1^\sharp,\cdots,a_{n-1}^\sharp,f_\sharp,b_{1\sharp},\cdots,b_{n-1\sharp}\big) ,$$
where, for a vector $x \in \C^{N},$ we denote $x^\sharp = (x,0) \in \C^{2N}$ and $x_\sharp = (0,x) \in \C^{2N}.$ The vector $e$ generates the range of $C_{12}.$ The vector $f$ generates the range of $C_{21}.$
\end{proof}

\begin{cor} \label{cor:scalarize} For any resonant pair $(j_1,j_2),$ there exists a family of smooth Fourier multipliers $\tilde P$ defined over the closure of $U^{large}_{\delta'}$ such that
\be \label{bd:P-tilde}
\sup_{\zeta \in U^{large}_{\delta'}} |\d_\zeta^\a \tilde P(\zeta)| + |\d_\zeta^\a \tilde P^{-1}(\zeta)| < \infty, \qquad \a \in \N^d,
\ee
and, denoting $\tilde S := \tilde P^{-1} S$,
\be\label{eq-tildeS}
 \d_t  \tilde S + \frac{1}{\e} M \tilde S +  \frac{1}{\sqrt \e} \left(\begin{array}{cc} \d_\eta \mu_{j_1} & 0 \\ 0 & \d_\eta \mu_{j_2} \end{array}\right) \cdot  \d_y \tilde S(\t;t)  = 0, \quad (x,y,\xi,\eta) \in \R^3 \times U^{large}_{\delta'},
\ee
with $\mu_{j_1}, \mu_{j_2}$ defined in Section {\rm \ref{sec:second:look},} where
$$
 M : = \left(\begin{array}{cc} i\mu_{j_1} & \sqrt \e \tilde b_{j_1j_2}^+ \\ \sqrt \e \tilde b_{j_2j_1}^- & i\mu_{j_2} \end{array}\right),
 \qquad \tilde b_{jj'}^\pm := \chi g_\pm \left(\begin{array}{cc} s_{jj'} \sqrt{  |{\rm tr}\, B_{1j_1j_2} B_{-1j_2j_1,+1}|} & 0 \\ 0 & 0\end{array}\right),
$$
where the top left blocks of $\tilde b_{ij}^\pm$ are scalar, and $s_{jj'} \in \{-1,1\}$ are such that the product $s_{j_1j_2} s_{j_2j_1}$ is equal to the sign of ${\rm tr} \,b_{j_1j_2}^+ b_{j_2j_1}^-.$
\end{cor}

We will see in Section \ref{sec:interaction:coefficients} that for some resonant pairs the sign of ${\rm tr} \,b_{j_1j_2}^+ b_{j_2j_1}^-$ is positive, corresponding to a spectral instability for $M.$  

\begin{proof} We let $(j_1,j_2) = (1,2)$ in this proof. 
We apply Lemma \ref{scalarize} with ${\mathcal D} = U^{large}_{\delta'}$ defined in \eqref{def:U-delta-large} and $C_{12} = B_{112}$ and $C_{21} = B_{-121,+1}.$ We will use $\t_{12} = \chi g_1$ and $\t_{21} = \chi g_{-1}.$   

Indeed, by definition of $U^{large}_{\delta'}$ at the beginning of this Section, we have $\mbox{tr}\, C_{12} C_{21} \neq 0$ in the compact set $\bar {\mathcal D}.$ Thus Lemma \ref{scalarize} applies and yields a change of basis $P,$ which depends only on the frequency variables, and scalar symbols $c_{12},$ $c_{21}.$ These scalar symbols are bounded away from zero in ${\mathcal D}$ and we may consider the further change of basis
$\dsp{P_1 :=  \diag \left\{\frac{\sqrt{|c_{12}|}}{\sqrt{|c_{21}|}},1,\cdots, 1\right\}.}$
Then, $\tilde P := P P_1$ is an adequate change of basis such that
\be \label{eq:P-tilde}
\tilde P^{-1}  \left(\begin{array}{cc} 0 & \chi b_{12}^+ \\ \chi b_{21}^- & 0 \end{array}\right) \tilde P =  \left(\begin{matrix} 0 & \breve C_{12} \\ \breve C_{21} & 0 \end{matrix}\right), \quad \mbox{with} \quad \breve C_{ij} := \left(\begin{array}{cc} \t_{ij} \sqrt{|c_{12}c_{21}|}\sgn c_{ij} & 0 \\ 0 & 0 \end{array}\right).
\ee
By the block diagonal structure of $\tilde P$, and the fact that $\tilde P$ is independent of $(t,x,y)$, we obtain that $\tilde S = \tilde P^{-1} S$ solves 
\be\label{eq-tildeS-2}
\tilde P \d_t  \tilde S +  \frac{1}{\e} \tilde P \left(\begin{array}{cc} i \mu_1 &  \sqrt \e \chi b_{12}^+ \\  \sqrt\e \chi b_{21}^- & i \mu_2 \end{array}\right) \tilde S + \frac{1}{\sqrt \e } \tilde P \left(\begin{array}{cc} \d_\eta \mu_1 & 0 \\ 0 & \d_\eta \mu_2 \end{array}\right) \cdot  \d_y  \tilde S = 0.
\ee
We then apply $\tilde P^{-1}$ to \eqref{eq-tildeS-2} and use \eqref {eq:P-tilde}. This gives the announced equation \eqref{eq-tildeS}. 
\end{proof}

We now further decompose $U^{large}_{\delta'}$ into positive and negative regions for the trace:
\be \label{def:U:pm}
 U^{+}_{\delta'} := \big\{ (\xi,\eta) \in U^{large}_{\delta'}, \quad \mbox{tr}\, B_{1j_1j_2}, B_{-1j_2j_1,+} > \delta \big\}, \quad U^-_{\delta'} := U^{large}_{\delta'} \setminus U^+_{\delta'}. 
 \ee
The maximum growth is recorded in $U^+_{\delta'}:$ 

\begin{cor}[Refined bound for the symbolic flow] \label{cor:optimal:bound:S}
For $\e$ small enough, for any resonant pair $(j,j'),$ for the symbolic flow equation in the form \eqref{2blocks} with a source $f$ and datum $g$ as in \eqref{symbolic:flow:source}, a rate of growth given by 
\be \label{def:gamma:jj'}
 \g_{jj'}(\R^6) :=  \max_{(x,y) \in \R^3} |g_{\pm}(x,y)| \cdot \max_{(\xi,\eta) \in U^{+}_{\delta'}} \Re e \, \sqrt{{\rm tr}\, B_{1j_1j_2} B_{-1j_2j_1,+1} (\xi,\eta)}.
\ee
As a consequence, for the full symbolic flow $S_{f,g}$ solution of \eqref{symbolic:flow:source},  a rate of growth is given by 
 \be \label{def:gamma}
  \g(\R^6) := \max\left( \g^+_{135}(\R^3 \times (\R^3 \setminus U^{away}_{\delta})), \,\, \max_{(j,j')} \g_{jj'}(\R^6) \right),
 \ee
 for $\e$ and $\delta$ small enough. 
 We often write $\g$ for $\g(\R^6).$ In the definition of $\g,$ the maximum is taken over all resonant pairs $(j,j').$ 
\end{cor}

 The computations of Section \ref{sec:interaction:coefficients} reveal that in the definition of $\g,$ the maximum is strictly positive and attained at $(j_1,j_2) = (1,4)$ and $(j_1,j_2) = (2,5).$ The lower bound analysis of Section \ref{sec:lower} will show that this $\g$ is optimal, which is crucial in our analysis. 

 The proof below uses Corollary \ref{cor:234} in the treatment of the $(2,3)-(3,4)$ subsystem, so that technically there is an extra $\g_\e$ term in the growth rate. But since $|\ln \e| \g_\e \to 0$ as $\e \to 0,$ and the maximal time is $T \sqrt \e |\ln \e|,$ the exponential $e^{\g_\e t/\sqrt \e}$ is bounded uniformly in $\e$ and $t,$ and we may disregard that term in $\g.$ 
    
\begin{proof} By Corollary \ref{cor:234}, a growth rate for $(2,3)-(3,4)$ is bounded by the maximum of a growth rate for $(2,3)$ and one for $(3,4).$ By Lemma \ref{lem:away:space:time:3}, the same is true for $(1,3)-(3,5),$ away from space-time resonances. Thus both triplets are accounted for in \eqref{def:gamma}, and we have the global growth rate \eqref{def:gamma} if we are able to validate the growth rate \eqref{def:gamma:jj'} for any resonant pair $(j,j').$ As before, we let $(j,j') = (1,2)$ for notational simplicity. 

Given $\delta > 0,$ which we will eventually choose to be small enough, let $\delta' > 0$ as in Lemma \ref{lem:small-trace}. 

First we consider the $\R^3 \times U^{large}_{\delta'}$ region. The starting point is the equation in $\tilde S_{f,g} := \tilde P^{-1} S_{f,g}$ given in Corollary \ref{cor:scalarize}. We coordinatize (see Note \ref{note:dimensions}) $\tilde S_{f,g} =: (\tilde S^1, \tilde S^2, \tilde S^3, \tilde S^4),$ so that $\tilde S^1$ and $\tilde S^3$ are scalar. Thus $\tilde S^2$ and $\tilde S^4$ solve transport equations with source terms that involve only $f,$ and no amplification occurs for these components (i.e, $\g = 0$ is a rate of growth). We focus on $\tilde S^1$ and $\tilde S^3,$ and denote $S = (\tilde S^1, \tilde S^3).$  The corresponding system is
\be\label{S0-new}
\left\{\begin{aligned}
\big( \d_t  + \frac{i \mu_1 }{\e} + \frac{\d_\eta \mu_1 }{\sqrt\e} \d_y\big) \tilde S^1 + \frac{1}{\sqrt \e} \chi g_1 s_{12} \sqrt{|\mbox{tr}\,b_{12}^+ b_{21}^-|} \tilde S^3  & = (\tilde P^{-1} f)_1,\\
\big( \d_t  + \frac{i \mu_2 }{\e} + \frac{\d_\eta \mu_2 }{\sqrt\e} \d_y\big) \tilde S^3 + \frac{1}{\sqrt \e} \chi g_{-1} s_{21} \sqrt{|\mbox{tr}\,b_{12}^+ b_{21}^-|} \tilde S^1 & = (\tilde P^{-1} f)_3.
\end{aligned}\right.
\ee 
The key is that the coupling terms above are {\it scalar.} Thus for the above system the rough rate of growth deduced from Lemma \ref{lem:rough:S} is actually optimal.

Indeed, consider first $U^+_{\delta'} \subset U^{large}_{\delta'},$ where $U^+_{\delta'}$ is defined in \eqref{def:U:pm}. Starting from system \eqref{S0-new}, we can go through the arguments of the proof of Lemma \ref{lem:rough:S}, and find 
$$ \|  \d_{x,y}^\a S(\t;t) \|_{L^\infty(\R^3 \times U^{+}_{\delta'})} \lesssim |\ln \e|^\star e^{\tilde \g(t - \t)/\sqrt \e},$$
with $\tilde \g$ such that $\tilde \g \geq 0$ and $(\tilde \g)^2$ is equal to 
 $$ \max_{\begin{smallmatrix} (x,y) \in \R^3 \\ (\xi,\eta) \in U^{+}_{\delta'} \end{smallmatrix} } \big| \chi g_1(x,y) s_{12} \sqrt{\mbox{tr}\,b_{12}^+ b_{21}^-(\xi,\eta)} \,\big| \cdot \max_{\begin{smallmatrix} (x,y) \in \R^3 \\ (\xi,\eta) \in U^{+}_{\delta'} \end{smallmatrix} } \big| \chi g_{-1}(x,y) s_{21} \sqrt{\mbox{tr}\,b_{12}^+ b_{21}^-(\xi,\eta)}\, \big|.$$
 Since $g_{-1} = g_1^*$ \eqref{gpm1}, we find $\tilde \g$ to be equal to the announced growth rate $\g_{12}:$
 $$ \g_{12} := \max_{(x,y) \in \R^3} |g_{\pm}(x,y)| \cdot \max_{(\xi,\eta) \in U^{+}_{\delta'}} \sqrt{{\rm tr}\, B_{112} B_{-121,+1} (\xi,\eta)}.$$
 At this point we proved that $\g_{12}$ defined above is a growth rate in $\R^3 \times U^{+}_{\delta'}.$ It remains to put the pieces together and verify that $\g_{12}$ is a growth rate in the entirety of $\R^6.$ 
 
 Consider now $(\xi,\eta) \in U^-_{\delta'}.$ By Corollary \ref{cor:scalarize}, we may consider system \eqref{S0-new} again. Here $s_{12} = - s_{21}.$ Without loss of generality, assume $s_{12} = 1.$ We observe the cancellation 
\be \label{cancellation:negative:trace}
\Re e \, \frac{1}{\e} (M_0 S, S)_{L^2(\R^2_y)} =  0, \quad M_0 := \left(\begin{array}{cc} i \mu_1 & \mbox{tr}\,b^+_{j_1j_2} \\   - \mbox{tr}\,b^+_{j_1j_2} & i \mu_2 \end{array}\right),
\ee 
with $S$ here identified with $(\tilde S^1, \tilde S^3),$ the coordinates of $S$ that intervene in \eqref{S0-new}. The $L^2$ scalar product in \eqref{cancellation:negative:trace} is defined in terms of a Hermitian scalar product in $\C^{2m},$ where $m$ is the dimension of $\tilde S^1$ and $\tilde S^3$ (we refer to Note \eqref{note:dimensions} here; the exact value of $m$ is unimportant). Thus from system \eqref{S0-new} we deduce
$$ \Re e \, (\d_t S, S)_{L^2(\R^2_y)} = \Re e \, (\tilde f, S)_{L^2(\R^2_y)}, \qquad (\xi ,\eta) \in U^-_{\delta'},$$
implying no growth in time for $\| S \|_{L^2(\R^2_y)},$ for fixed $(x,\xi,\eta),$ with $(\xi,\eta) \in U^-_{\delta'}.$  

 Keeping $(\xi,\eta)$ in $U^-_{\delta'},$ we now verify that spatial derivatives $\d_{x,y}^\a S$ do not grow in $L^2(\R^2_y)$ norm either. That is, we now prove, by induction, the bound 
\be \label{induction:negative:trace}
 \| \d_{x,y}^\a S \|_{L^2(\R^2_y)} \lesssim |\ln \e|^\star N_{2,\a}(\t;t,f,g; \R^3 \times U^-_{\delta'}),
\ee 
with notation $N_{2,\a}$ introduced in \eqref{def:C2} just above Definition \ref{def:growth:rate}. 
We verified above \eqref{induction:negative:trace} for $\a = 0.$ 
By \eqref{cancellation:negative:trace}, applying first $\d_{x,y}^\a$ to \eqref{S0-new} and then taking the $L^2$ scalar product with $\d_{x,y}^\a S,$ we observe that 
 \be \label{07:10:20}  \d_t \| \d_{x,y}^\a S \|_{L^2(\R^2_y)} \leq  \frac{1}{\e} \| [\d_{x,y}^\a, M_0] S \|_{L^2(\R^2_y)} + \| \d_{x,y}^\a f \|_{L^2(\R^2_y)}.\ee
 The commutator $[\d_{x,y}^\a, M_0]$ involves only terms in $M_0$ that are $O(\sqrt \e).$ Thus we have by the induction hypothesis
 $$ \begin{aligned} \frac{1}{\e} \| [\d_{x,y}^\a, M_0] S \|_{L^2(\R^2_y)} & \lesssim \frac{1}{\sqrt \e} \max_{\b < \a} \| \d_{x,y}^\b S \|_{L^2(\R^2_y)} \lesssim \max_{\b < \a} \frac{1}{\sqrt \e} N_{2, \b}(\t;t,f,g; \R^3 \times U^-_{\delta'}).
\end{aligned}$$ 
 With \eqref{07:10:20}, this implies
$$ \d_t \| \d_{x,y}^\a S \|_{L^2(\R^2_y)} \lesssim  \| \d_{x,y}^\a f \|_{L^2(\R^2_y)} + \max_{\b < \a} \frac{1}{\sqrt \e} N_{2, \b}(\t;t,f,g; \R^3 \times U^-_{\delta'}),$$
and, integrating in time over $[\t;t]$ with $t = O(\sqrt \e |\ln \e|),$ we find \eqref{induction:negative:trace} for any $\a.$   

The bound \eqref{induction:negative:trace} is pointwise in $(x,\xi,\eta).$ Given $\a_a$ an index of maximal size for which \eqref{induction:negative:trace} holdss true, by Sobolev embedding we deduce the pointwise bound, for $\a < \a_a - [d/2] - 1:$ 
\be \label{bound:U:minus} 
 \| \d_{x,y}^\a S\|_{L^\infty(\R^3 \times U^-_{\delta'})} \lesssim |\ln \e|^\star N_{2,\a}(\t;t,f,g;\R^3 \times U^-_{\delta'}).
\ee
That is, we have a growth rate in $\R^3 \times U^-_{\delta'},$ in the sense of Definition \ref{def:growth:rate}, which is equal to 0.

  We now choose $\delta > 0$ to be small enough so that $C \delta^{1/2} < \g,$ where $C$ is the constant that appears in Lemma \ref{lem:small-trace}. Lemma \ref{lem:away:space:time} asserts that for $\e$ small enough, $\g_{12}$ is also a rate of growth in $\R^3 \times U^{away}_{\delta'}.$ 
 
 Thus it appears that for $\e$ small enough $\g_{12}$ is a rate of growth in 
 $$\R^3 \times (U^{away}_{\delta'} \cup U^{small}_{\delta'} \cup U^{large}_{\delta'}) = \R^3 \times \Omega_{freq}.$$ Finally, Lemma \ref{lem:away} shows that there is no growth outside of $\Omega,$ so that $\g_{12}$ is also a rate of growth in $\R^6 \setminus \Omega.$ Since 
  $$ (\R^3 \times \Omega_{freq}) \cup (\R^6 \setminus \Omega) = \R^6,$$
 we verified $\g_{12}$ is indeed a growth rate (over $\R^6$) for the $(1,2)$ subsystem for $S.$ The same goes for every other resonant pair, and we arrive at \eqref{def:gamma:jj'}. 
 \end{proof}

\subsection{Bounds for frequency derivatives} 

 Frequency derivatives of $S$ enjoy the same growth rate but each frequency derivative brings out an $\e^{-1/2}$ factor. Crucially, the prefactor is smaller for $\d_\eta$ derivatives near space-time resonances:  
 
\begin{lem} \label{lem:derS}
The symbolic flow $S$ from \eqref{resolvent0-in-proof} enjoys the bounds, with the rate $\gamma$ defined in \eqref{def:gamma}:
 \be \label{der:S:bd:0} \e^{|\b|/2} \Big( \| \d_{\xi,\eta}^\b \d_{x,y}^\a S(\t;t) \|_{L^\infty}  + \| \d_{\xi,\eta}^\b \d_{x,y}^\a S(\t;t) \|_{L^\infty(\R^4_{x,\xi,\eta}, L^2(\R^2_y)}  \Big) \lesssim |\ln \e|^\star e^{\g (t - \t) /\sqrt \e}.\ee 
 Near $\eta = 0,$ we have the improved bound for transverse frequency derivatives: 
\be \label{der:S:bd:1} \begin{aligned} \e^{|\b|/4} \Big( |\d_\eta^\b \d_{x,y}^\a S(\t;t,x,y,\xi,\eta)| & + \|\d_\eta^\b \d_{x,y}^\a S(\t;t,x,\cdot,\xi,\eta)\|_{L^2(\R^2_y)} \Big) \\ & \lesssim |\ln \e|^\star e^{\g (t - \t) /\sqrt \e}, \qquad |\eta| \leq \e^{1/4}, \end{aligned}\ee
uniformly in $(\e,\t,t,x,y,\xi,\eta),$ for $0 \leq \t \leq t \leq T \sqrt \e |\ln \e|$ and $|\eta| \leq \e^{1/4}.$ 
 \end{lem}

\begin{proof} For $\b = 0,$ both bounds hold true for all $\a$ by Corollary \ref{cor:optimal:bound:S} and Lemma \ref{lem:S:L2:Linfty}.

For a proof of the first bound \eqref{der:S:bd:0}, we work by induction on $|\b|,$ as we assume now that \eqref{der:S:bd:0} holds true for up to $k$ frequency derivatives of $S,$ for some positive integer $k,$ for any number of spatial derivatives of $S.$ Let $\b \in \N^3$ with $|\b| = k+1.$ We now work by induction on $|\a|.$ By applying $\d_{\xi,\eta}^\b \d_{x,y}^\a$ to the symbolic flow equation \eqref{resolvent0-in-proof}, we obtain 
 \be \label{for:detaS} \d_{\xi,\eta}^\beta \d_{x,y}^\a S = S_{f_{\a,\b},0}, \qquad \dsp{f_{\a,\b} := - \frac{1}{\e} [\d_{\xi,\eta}^\b \d_{x,y}^\a,i {\bf A} - \sqrt \e {\bf B}] S - \frac{1}{\sqrt \e} [\d_{\xi,\eta}^\b, \d_\eta {\bf A}] \d_{x,y}^\a \d_y S},\ee 
 with notation $S_{f,g}$ introduced in \eqref{symbolic:flow:source}. 
 
 We initiate the induction on $|\a|:$ the source $f_{0,\b}$ contains strictly less than $|\b|$ frequency derivatives of the symbol, and we can use the induction hypothesis. This gives
 $$ \| f_{0,\b}(\t;t)\|_{L^\infty} + \| f_{0,\b}(\t;t)\|_{L^\infty(L^2)} \lesssim \e^{-1 - (|\b| - 1)/2} |\ln \e|^\star e^{\g (t - \t)/\sqrt \e},$$
 where $L^\infty(L^2)$ is short for $L^\infty(\R^4_{x,\xi,\eta}, L^2(\R^2_y)).$ 
 Then we use Corollary \ref{cor:optimal:bound:S}, which gives 
 $$ \|  \d_{\xi,\eta}^\beta S(\t;t)\|_{L^\infty} = \| S_{f_{0,\b},0}(\t;t) \|_{L^\infty} \lesssim \e^{1/2} \e^{-1 - (|\b| - 1)/2} |\ln \e|^\star e^{\g (t - \t)/\sqrt \e}.$$
 Lemma \ref{lem:S:L2:Linfty} then gives the same bound for the $L^\infty(L^2)$ norm of $S_{f_0,\b},$ 
 and the case $\a = 0$ is proved. 
 
 Now we assume that the bound \eqref{der:S:bd:0} is known for $\d_{\xi,\eta}^\b \d_{x,y}^\a S,$ for all $|\a| \leq k'-1,$ and use \eqref{for:detaS} with $|\a| = k'.$ 
The first term in $f_{\a,\b}$ is a sum of terms of the form 
\be \label{first:alpha} \e^{-1} (\d_{\xi,\eta}^{\b_1} \d_{x,y}^{\a_1} ( i {\bf A} - \sqrt \e {\bf B}))(\d_{\xi,\eta}^{\b_2} \d_{x,y}^{\a_2} S), \quad \mbox{with $|\a_1| + |\b_1| > 0.$}\ee 
 If $|\b_1| > 0,$ then this term is amenable to the first induction hypothesis (where the induction is on $|\b|$). If $|\a_1| > 0,$ then this term is amenable to the second induction hypothesis (induction on $|\a|$). The second term in $f_{\a,\b}$ contains at most $|\b| - 1$ frequency derivatives of the source, and we use the first induction hypothesis. Thus $f_{\a,\b}$ satisfies the same bound as $f_{0,\b}$ above, and another application of Corollary \ref{cor:optimal:bound:S} and Lemma \ref{lem:S:L2:Linfty} allows us to complete the induction, just like in the case $\a = 0.$

 We turn to \eqref{der:S:bd:1},  the bound near $\eta = 0.$ We use the same scheme of proof, that is an induction on $\a$ within an induction of $\b,$ but this time we pay attention to the size of the source $f_{\a,\b}.$ Here frequency derivatives are in the transverse frequency variables: in the definition of $f_{\a,\b},$ we replace $\d_{\xi,\eta}^\b$ with $\d_\eta^\b.$ 
 
 As noted above, if $\b = 0$ then \eqref{der:S:bd:1} holds true for all $\a.$ Now we assume that it does hold true for all $|\b| \leq k,$ for some $k \geq 0,$ and for all $\a.$ Consider $\b \in \N^3$ with $|\b| = k+1.$ First we consider the case $\a = 0,$ and the first term in $f_{0,\b}$ \eqref{for:detaS}. It contains terms of the form $\e^{-1} (\d_\eta^{\b_1} (i {\bf A} - \sqrt \e {\bf B})) (\d_\eta^{\b_2} S).$ If $|\b_1| = 1,$ then we observe the following:
 \be \label{cancel:small:eta}
 |\d_\eta^{\b_1} (i {\bf A} - \sqrt \e {\bf B})| \lesssim \e^{1/4}, \qquad |\b_1| = 1, \qquad |\eta| \lesssim \e^{1/4}.
 \ee
 Indeed, we have $|\d_\eta \l_j| = O(|\eta|)$ for $j \neq 3,$ and $\l_3$ and its derivatives are $ O(\sqrt \e |(\xi,\eta)|),$ uniformly in frequency (see Section \ref{sec:spectral:dec}, in particular Section \ref{sec:parallel:modes}). This takes care of the ${\bf A}$ terms in \eqref{cancel:small:eta}. The ${\bf B}$ terms are $O(\sqrt \e)$ and so are their frequency derivatives. This verifies \eqref{cancel:small:eta}.

 As a consequence, by the induction hypothesis, we have
 \be \label{for:der:S:bd:2}
  \begin{aligned} \e^{-1} \big| (\d_\eta^{\b_1} (i {\bf A} - \sqrt \e {\bf B})) (\d_\eta^{\b_2} S ) \big| & \lesssim \e^{-1} \e^{1/4} \e^{-(|\b| - 1)/4} |\ln \e|^\star e^{(t - \t) \g/\sqrt\e}, \\ &  \mbox{if $|\b_1| = 1$ and $|\b_2| = |\b| - 1,$ for $|\eta| \leq \e^{1/4}.$} \end{aligned}
  \ee 
 In the case $|\b_1| \geq 2,$ we cannot use \eqref{cancel:small:eta}, but then we simply use the induction hypothesis, which implies 
 \be \label{for:der:S:bd:3} \begin{aligned} \e^{-1} \big| (\d_\eta^{\b_1} (i {\bf A} - \sqrt \e {\bf B})) (\d_\eta^{\b_2} S ) \big| & \lesssim \e^{-1} \e^{-(|\b| - 2)/4} |\ln \e|^\star e^{(t - \t) \g/\sqrt\e}, \\ &  \mbox{if $|\b_1| \geq 2$ and $|\b_2| \leq |\b| - 2,$ for $|\eta| \leq \e^{1/4}.$} \end{aligned}
  \ee 
  
  The second term in $f_{0,\b}$ contains terms of the form $(\d_\eta^{\b_1} {\bf A}) (\d_\eta^{\b_2} \d_y S),$ with $|\b_1| \geq 2$ and $|\b_1| + |\b_2| = |\b| + 1.$ In particular, $|\b_2| \leq |\b| - 1,$ so that, by the induction hypothesis, we have
  \be \label{for:der:S:bd:2:2}
  \begin{aligned} \e^{-1/2} \big| (\d_\eta^{\b_1} {\bf A}) (\d_\eta^{\b_2} \d_y S) \big| & \lesssim \e^{-1/2} \e^{-(|\b| - 1)/4} |\ln \e|^\star e^{(t - \t)\g/\sqrt \e}, \\ & \mbox{with $|\b_1| \geq 2,$ $|\b_2| \leq |\b| - 1,$ and $|\eta| \leq \e^{1/4}.$}  
  \end{aligned}
  \ee
 
 With \eqref{for:der:S:bd:2}, \eqref{for:der:S:bd:3} and \eqref{for:der:S:bd:2:2}, we obtain
 \be \label{fob} | f_{0,\b} | \lesssim \e^{-1/2} \e^{-|\b|/4} |\ln \e|^\star e^{(t - \t)\g/\sqrt \e}, \qquad \mbox{for $|\eta| \leq \e^{1/4}$}.\ee  
 
 We now apply Corollary \ref{cor:optimal:bound:S} and Lemma \ref{lem:S:L2:Linfty}, as we did for the proof of \eqref{der:S:bd:0}, and this gives, thanks to the time integration in the bound of Lemma \ref{lem:rough:S}, the bound \eqref{der:S:bd:1} in the case $\a = 0$ and $|\b| = k + 1.$ 
 
 Finally we assume that for any $\b$ with $|\b| = k+1,$ for some $k' \geq 1,$ the bound \eqref{der:S:bd:1} holds for all $|\a| \leq k'-1,$ and, for some $\a \in \N^3$ with $|\a| = k',$ we bound $f_{\a,\b}.$ As in the proof of \eqref{der:S:bd:0}, the first term is \eqref{first:alpha}. If in \eqref{first:alpha}, we have $|\b_1| > 0,$ then we use the induction on $\beta,$ and bound the corresponding term exactly as in \eqref{for:der:S:bd:2} and \eqref{for:der:S:bd:3}. If we have $|\a_1| > 0,$ then we can use the induction on $\a.$ Thus we find the term in \eqref{first:alpha} to be controlled just like $f_{0,\b}$ \eqref{fob} above. The second term in $f_{\a,\b}$ has the same upper bound, by the induction on $\b.$ We conclude again with Corollary \ref{cor:optimal:bound:S} and Lemma \ref{lem:S:L2:Linfty}.
 \end{proof}

\subsection{Bounds on the symbolic flow operator} \label{sec:bd:symbolic:flow:operator}

The above pointwise bounds on the symbolic flow translate into bounds, in $L^2$ and $\fl1,$ for the action of the symbolic flow {\it operator}:

\begin{prop} \label{prop:opS} We have the bound for the symbolic flow $S$ from \eqref{resolvent0-in-proof}, for $0 \leq \t \leq t \leq T \sqrt \e |\ln \e|:$
 $$ \| \op_\e(S(\t;t)) \|_{L^2 \to L^2} +  \| \op_\e(S(\t;t)) \|_{\fl1 \to \fl1} \lesssim |\ln \e|^\star e^{(t - \t) \g /\sqrt \e},$$
 where $\g$ is the growth rate defined in \eqref{def:gamma}. 
\end{prop}

\begin{proof} By Proposition \ref{prop:H}, the $L^2 \to L^2$ norm of $\op_\e(S(\t;t))$ is bounded by a constant (depending only on dimensions) times
$$ \sup_{\begin{smallmatrix} (\xi,\eta) \in \R^3 \\ 0 \leq |\a| \leq 4 \end{smallmatrix} } \| \d_{x,y}^\a S(\t;t) \|_{L^1(\R^3_{x,y})}.$$
By Proposition \ref{prop:op:fl1}, the $\fl1\to \fl1$ norm of $\op_\e(S(\t;t))$ is bounded by a constant (depending only on dimensions) times 
$$ \sup_{\begin{smallmatrix} (\xi,\eta) \in \R^3 \\ 0 \leq |\a| \leq 2 \end{smallmatrix} } \| \d_{x,y}^\a S(\t;t) \|_{L^2(\R^3_{x,y})}.$$  
We decompose %
$$ S = S \chi_{\Omega}^\sharp + S (1 - \chi_{\Omega}^\sharp),$$
where $\chi_\Omega \prec \chi_\Omega^\sharp,$ with $\chi_\Omega$ defined in \eqref{def:chi:Omega:pair} and \eqref{def:chi:Omega:triplet}.  By property of $\Omega$ (see Section \ref{sec:far-field}), we may choose the size of the support of $\chi_\Omega^\sharp$ to be $O(|\ln \e|^\star).$ In particular,
$$ \| \d_{x,y}^\a (\chi_{\Omega}^\sharp S(\t;t) ) \|_{L^p} \lesssim |\ln \e|^\star \| \d_{x,y}^\a S(\t;t) \|_{L^\infty}, \quad p \in  \{1,2\}.$$
and by Corollary \ref{cor:optimal:bound:S}, the above is controlled by $|\ln \e|^\star e^{ \g (t - \t) /\sqrt \e}.$ 

The symbol $(1 - \chi_{\Omega}^\sharp) S$ is supported away from $\Omega,$ and by Lemma \ref{lem:away}, we have an explicit description of $S$ in that region. Indeed, applying Lemma \ref{lem:away} with $f = 0$ and $g = {\rm Id},$ we find (keeping in mind Note \ref{note:dimensions})
$$ \begin{aligned} (1 - \chi_\Omega^\sharp) S  = (1 - \chi_\Omega^\sharp) e^{ - i (t - \t) \mu_j/\e} = e^{ - i (t - \t) \mu_j/\e} - \chi_\Omega^\sharp(x,y,\xi,\eta) e^{ - i (t - \t) \mu_j/\e}.\end{aligned}$$
Since $\chi_\Omega^\sharp$ can be chosen to be a tensor product, Lemma \ref{lem:tensor} applies to the operators associated with both symbols in the above right-hand side, and we obtain
$$ \| \op_\e\big( (1 - \chi_\Omega^\sharp) S \big) \|_{L^2 \to L^2} \lesssim 1, \quad\mbox{and} \quad \| \op_\e\big( (1 - \chi_\Omega^\sharp) S \big) \|_{\fl1 \to \fl1} \lesssim |\ln \e|^\star, $$
where the $|\ln \e|^\star$ factor comes from the $\fl1$ norm of ${\bm \chi}_{trans}^\sharp.$ 
\end{proof}

\subsection{Action of the symbolic flow on an oscillating datum} \label{sec:flow:action:datum}

The action of the symbolic flow operator on the product of a fast oscillation times a fast decaying function has a simple form: 

\begin{lem} \label{lem:for:datum} For the symbolic flow $S$ from \eqref{resolvent0-in-proof}, for $0 \leq t \leq T \sqrt \e |\ln \e|,$ for any $f$ in the Schwartz class, any $\xi_0 \in \R:$
$$ \big\| \op_\e(S(0;t) \big( e^{i x \xi_0/\e} f \big) - e^{i x \xi_0/\e} S(0;t,x,y,\xi_0,0) f \big\|_{X} \lesssim \e^{1/4} |\ln\e|^\star e^{t \g/\sqrt \e},$$
with $X = L^2$ and $X = \fl1,$ the implicit constants depending on various $W^{k,p}$ norms of $f.$  
\end{lem}

Lemma \ref{lem:for:datum} will be used in Lemma \ref{lem:look:datum} and Corollary \ref{lem:low}, in both cases with $f = \phi,$ where $\phi$ is the smooth and compactly supported amplitude which carries the initial perturbation \eqref{initial:datum}-\eqref{varphi:initial:perturbation}. 

\begin{proof} We use the fact that derivatives of $S$ are not too large, as described in Lemma \ref{lem:derS}, and that $\hat f$ is very small for large frequencies. Precisely, Taylor expanding the symbol, we find
 $$ e^{- i x \xi_0/\e}\op_\e(S(0;t) \big( e^{i x \xi_0/\e} f \big) - S(0;t,x,y,\xi_0,0) f = f_x + f_y,$$
 with notation
 $$ \begin{aligned} f_x & := \e  \int_0^1 \op_1(\d_\xi S(0;t,x,y,\xi_0 + \e \t \xi, \sqrt \e \t \eta)) \d_x f \, d\t \\ f_y &:= \sqrt \e \int_0^1 \op_1(\d_\eta S(0;t,x,y,\xi_0 + \e \t \xi, \sqrt \e \t \eta)) \cdot \d_y f \, d\t.
 \end{aligned}$$
For the $f_x$ term we first decompose, with $\chi_\Omega$ as in the proof of Proposition \ref{prop:opS} (throughout this proof, in a harmless abuse of notation we identify $\chi_\Omega$ with its extension $\chi_\Omega^\sharp$): 
\be \label{dec:dxiS} \begin{aligned} \d_\xi S =  \chi_\Omega \d_\xi S + (1 - \chi_\Omega) \d_\xi S = \chi_\Omega \d_\xi S +  (1 - \chi_\Omega) \left( \frac{- i t  \d_\xi \mu_j}{\e} e^{-i t  \mu_j/\e} \right)_j,\end{aligned}\ee
the second equality by definition of $\Omega$ by Lemma \ref{lem:away}. 

As in the proof of Proposition \ref{prop:opS}, for the $\chi_\Omega \d_\xi S$ part we can use Corollary \ref{cor:optimal:bound:S} in conjunction with Propositions \ref{prop:H} and \ref{prop:op:fl1}: 
$$ \begin{aligned} \| \op_1 \Big( \, (\chi_{\Omega} \d_\xi S(0;t)) & (x,y,\xi_0 + \e \t \xi, \sqrt \e \t \eta) \, \Big) \|_{L^2 \to L^2} \\ & \lesssim \sup_{ \begin{smallmatrix} (\xi,\eta) \in \R^3  \\ 0 \leq |\a| \leq 4 \end{smallmatrix}} \| \d_\xi \d_{x,y}^\a (\chi_\Omega S(0;t))(\cdot,\xi_0 + \e \t \xi, \sqrt \e \t \eta)) \|_{L^1(\R^3_{x,y})} \\ & \lesssim |\ln \e|^\star \sup_{ \begin{smallmatrix} (\xi,\eta) \in \R^3  \\ 0 \leq |\a| \leq 4 \end{smallmatrix}} \| \d_\xi \d_{x,y}^\a S(0;t)(\cdot,\xi, \eta)) \|_{L^\infty(\R^3_{x,y})} \\ & \lesssim \e^{-1/2} |\ln \e|^\star e^{t \g/\sqrt \e}. \end{aligned}$$
For the $\fl1 \to \fl1$ bound, an $L^2$ norm of the symbol replaces the above $L^1$ norm, as per Proposition \ref{prop:op:fl1}. This gives 
$$ \begin{aligned} \| \op_1 \Big( \, (\chi_{\Omega} \d_\xi S(0;t)) & (x,y,\xi_0 + \e \t \xi, \sqrt \e \t \eta) \, \Big) \|_{\fl1 \to \fl1} \\ & \lesssim \sup_{ \begin{smallmatrix} (\xi,\eta) \in \R^3  \\ 0 \leq |\a| \leq 4 \end{smallmatrix}} \| \d_\xi \d_{x,y}^\a (\chi_\Omega S(0;t))(\cdot,\xi_0 + \e \t \xi, \sqrt \e \t \eta)) \|_{L^2(\R^3_{x,y})} \\ & \lesssim \sup_{ \begin{smallmatrix} (\xi,\eta) \in \R^3  \\ 0 \leq |\a| \leq 4 \end{smallmatrix}} \| \d_\xi \d_{x,y}^\a S(0;t)(\cdot,\xi, \eta)) \|_{L^\infty(\R^3_{x,y})} \\ & \lesssim \e^{-1/2} |\ln \e|^\star e^{(t  - \t) \g/\sqrt \e}. \end{aligned}$$
The far-field contribution $(1 - \chi_\Omega) \d_\xi S$ remains. It is given explicity on the second line of \eqref{dec:dxiS}. We observe that 
$$ \left| \frac{ i t \d_\xi \mu_j}{\e} e^{i t \mu_j/\e} \right| \lesssim \e^{-1/2} |\ln \e|,$$
in the time interval under consideration. By the fact that $\chi_\Omega$ is a tensor product and Lemma \ref{lem:tensor}, this implies
$$ \big\| \op_1\big( ((1 - \chi_\Omega) \d_\xi S)(0;t)(x,y,\xi_0 + \e \t \xi, \sqrt \e \t \eta) \big) \big\|_{L^2 \to L^2} \lesssim \e^{-1/2} |\ln \e|,$$
and 
$$ \big\| \op_1\big( ((1 - \chi_\Omega) \d_\xi S)(0;t)(x,y,\xi_0 + \e \t \xi, \sqrt \e \t \eta) \big) \big\|_{\fl1 \to \fl1} \lesssim \e^{-1/2} |\ln \e|^\star.$$
At this point we proved
 $$ \| f_x\|_{L^2} + \| f_x \|_{\fl1} \lesssim \e^{1/2}  |\ln \e|^\star e^{t \g/\sqrt \e},$$
 since $f$ is rapidly decaying.

 For the $f_y$ term we use the bound \eqref{der:S:bd:1} from Lemma \ref{lem:derS}.  We need another argument, though, since \eqref{der:S:bd:1} gives us a good control of $\d_\eta \d_{x,y}^\a S$ only for $|\eta|$ small. For large $\eta,$ we use the fast decay of $\widehat{\d_y f},$ which derives from the assumption that $f$ be in the Schwartz space. 

 In this view, we introduce a decomposition of $\d_y f$ into 
 $$ \d_y f = \chi_{LF}(\e^{1/4} D_y) \d_y f + (1 - \chi_{LF}(\e^{1/4} D_y)) \d_y f,$$
 where $\chi_{LF}$ is a smooth low-frequency cut-off, identically equal to $1$ in a small ball and identically equal to 0 outside of the unit ball. By fast decay of $f,$ we have
 \be \label{decay:schw} \| (1 - \chi_{LF}(\e^{1/4} D_y) \d_y f \|_{L^2} \leq C_m \e^m, \qquad \mbox{for any $m > 0,$ for some $C_m > 0.$}\ee
 We use the decomposition
 $$ \d_\eta S = \chi_\Omega \d_\eta S + (1 - \chi_\Omega) \left( \frac{-i t \d_\eta \mu_j}{\e} e^{-i t \mu_j /\e} \right)_j,$$
 analogous to \eqref{dec:dxiS}. With this decomposition of the symbol in $f_y,$ and the above decomposition of $f,$ there are four terms to bound. 
 
 First we observe that
 \be \label{chiLF} \op_\e(M) \chi_{LF}(\e^{1/4} D_y) = \op_\e( \chi_{LF}(\e^{-1/4} \cdot) M), \quad \mbox{for any symbol $M,$}\ee
 by definition of the anisotropic semiclassical quantization \eqref{aniso}. By definition of $\chi_{LF},$ this means that the composition to the right by $\chi_{LF}(\e^{1/4} D_y)$ implies a localization of the symbol to the domain $|\eta| \leq \e^{1/4}.$ Thus
 \be \label{near:low} \begin{aligned} \big\| \op_\e\big( (\chi_\Omega \d_\eta S(0;t))& (x,y,\xi_0 + \cdot, \cdot) \big) \chi_{LF}(\e^{1/4} D_y) \d_y f \big\|_{L^2} \\ & \lesssim  \sup_{\begin{smallmatrix} (\xi,\eta) \in \R^3 \\ 0 \leq |\a| \leq 4 \end{smallmatrix} } \| \d_{x,y}^\a\big(  \chi_\Omega \d_\eta S\big) \chi_{LF}(\e^{-1/4} \eta) \|_{L^1(\R^3_{x,y})} \\ & \lesssim  |\ln \e|^\star \sup_{\begin{smallmatrix} \xi \in \R \\ |\eta| \leq \e^{1/4} \\ 0 \leq |\a| \leq 4 \end{smallmatrix} } \| \d_{x,y}^\a \d_\eta S \|_{L^\infty(\R^3_{x,y})} \\ & \lesssim \e^{-1/4} |\ln \e|^\star e^{t \g/\sqrt \e}, \end{aligned} \ee 
 the second-to-last inequality by size of the support of $\chi_\Omega,$ and the last by \eqref{der:S:bd:1} in Lemma \ref{lem:derS}. The $\fl1$ bound is handled by the same arguments (see the $f_x$ term above). If we take into account the $\sqrt \e$ prefactor in $f_y,$ the above takes care of the first (near-field/low-transversal frequency) term in $f_y.$

 In view of \eqref{chiLF}, the far-field/low-transversal frequency term in $f_y$ is 
 $$ \op_\e\Big( (1 - \chi_\Omega) \frac{i t \d_\eta \mu_j}{\e} e^{i t \mu_j /\e} \chi_{LF}(\e^{-1/4} \cdot) \Big) \d_y f,$$
 where the shift $+\xi_0$ in $\xi$ is implicit.  
 As in the proof of Lemma \ref{lem:derS}, we observe that $|\d_\eta \mu_j| \lesssim \max(|\eta|, \sqrt \e).$ Thus, in the time interval under consideration
 $$ \left| \chi_{LF}(\e^{-1/4} \eta) \frac{i t \d_\eta \mu_j}{\e} e^{i t \mu_j /\e} \right| \lesssim \e^{-1/4} |\ln \e|:$$ 
 we gained a factor $\sqrt \e |\ln \e|$ from $t$ and a factor $\e^{1/4}$ from $\eta \chi_{LF}(\e^{-1/4} \eta).$ 
 By Lemma \ref{lem:tensor}, this implies
$$ \big\| \op_\e\big( ((1 - \chi_\Omega) \d_\eta S)(0;t)(x,y,\xi_0 + \cdot) \big) \chi_{LF}(\e^{1/4} D_y) \d_y f \big\|_{L^2} \lesssim \e^{1/4} |\ln \e|,$$
and 
$$ \big\| \op_\e\big( ((1 - \chi_\Omega) \d_\xi S)(0;t)(x,y,\xi_0 + \cdot) \big) \chi_{LF}(\e^{1/4} D_y) \d_y f \big\|_{\fl1} \lesssim \e^{1/4} |\ln \e|^\star.$$

The high-transversal frequency terms remain. We handle the near-field/high-transversal frequency term very much like we did in \eqref{near:low}: this gives
\be \label{near:high} \begin{aligned} \big\| \op_\e\big( (\chi_\Omega \d_\eta S(0;t))& (x,y,\xi_0 + \cdot, \cdot) \big) (1 - \chi_{LF}(\e^{1/4} D_y)) \d_y f \big\|_{L^2} \\ & \lesssim  \sup_{\begin{smallmatrix} (\xi,\eta) \in \R^3 \\ 0 \leq |\a| \leq 4 \end{smallmatrix} } \| \d_{x,y}^\a\big(  \chi_\Omega \d_\eta S\big) \|_{L^1(\R^3_{x,y})}  \| (1 - \chi_{LF}(\e^{1/4} D_y)) \d_y f \|_{L^2} \\ & \lesssim  |\ln \e| \sup_{\begin{smallmatrix} (\xi,\eta) \in \R^3 \\ 0 \leq |\a| \leq 4 \end{smallmatrix} } \| \d_{x,y}^\a \d_\eta S \|_{L^\infty(\R^3_{x,y})} \| (1 - \chi_{LF}(\e^{1/4} D_y)) \d_y f \|_{L^2} \\ & \lesssim \e^{-1/2} |\ln \e|^\star e^{t \g/\sqrt \e} \| (1 - \chi_{LF}(\e^{1/4} D_y)) \d_y f \|_{L^2}. \end{aligned} \ee 
Note above the $\e^{-1/2}$ prefactor instead of $\e^{-1/4}$ in \eqref{near:low}. But now we can use \eqref{decay:schw} in \eqref{near:high}, so that this term is appropriately small. The $\fl1$ bound is similar, with an $L^2$ norm of the symbol which is promptly bounded by an $L^\infty$ norm. 

For the fourth, far-field/high-transversal frequency term in $f_y,$ we use Lemma \ref{lem:tensor}  and \eqref{decay:schw} again.
\end{proof}

 \section{Estimates in the {\it in}} \label{sec:Duh}

We go back to equation \eqref{eq:Uflat}, which we reproduce here:
$$ \d_t v_{in} + \frac{i}{\e} \op_\e({\bf A}) v_{in} = \frac{1}{\sqrt \e} \op_\e({\bf B}) v_{in} +  {\bf F},
$$
 where ${\bf F}$ is equivalent to $\dot F,$ in the sense of Definition \ref{def:equiv:remainder}.
The datum $v_{in}(0)$ is partially described in \eqref{datum:vin}. A complete description of $v_{in}(0)$ will be given in Section \ref{sec:lower}. 
With $L$ defined by 
\be \label{def:L}
 L := i {\bf A} - \sqrt \e {\bf B},
\ee
Theorem \ref{th:duh} (stated and proved in Section \ref{app:duh}) applies to the initial-value problem for the system in $v_{in},$ as we verify now:

\begin{prop} \label{prop:duh} The solution $v_{in}$ to \eqref{eq:Uflat} issued from $v_{in}(0)$ is
\be \label{rep:W}
 v_{in}(t) = \op_\e({\bf S}(0;t)) v_{in}(0) + \int_0^t \op_\e({\bf S}(t';t))(\Id + \e^\sigma R_1) \Big( {\bf F} + \e^{\sigma} R_2(\cdot) v_{in}(0)\Big)(t')\, dt',
  \end{equation}
 where $\sigma > 0$ and the remainders $R_1$ and $R_2$ satisfy
$$ \begin{aligned}
 & \sup_{0 \leq t \leq T \sqrt \e |\ln \e|} 
  \| R_1 z_1 (t)\|_{\e,s_1}  \lesssim |\ln \e|^\star  
  \sup_{0 \leq t \leq T \sqrt \e |\ln \e|} \|  z_1 (t)\|_{\e,s_1}, \\
& \sup_{0 \leq t \leq T \sqrt \e |\ln \e|}   \| R_2(t) z_{2} \|_{\e,s_1} \lesssim |\ln \e|^{\star} \|  z_{2} \|_{\e,s_1},
\end{aligned}$$ 
for all $z_{1} \in L^{\infty}(0, T\sqrt \e |\ln \e|, H^{s_1}(\R^3)),$ $z_2 \in H^{s_1}(\R^{3}).$
 The symbol ${\bf S}$ above is defined in \eqref{def:S} in Section {\rm \ref{app:duh},} with $S_0 = S$ and $L$ defined in \eqref{def:L}. The symbol $S$ is the symbolic flow studied in Section {\rm \ref{sec:flow}.} The exponent $\sigma$ defined in \eqref{def:sigma} can be made arbitrarily large if the WKB solution is smooth. 
\end{prop}

\begin{proof} We invoke Theorem \ref{th:duh}, which relies on Assumptions \ref{ass:B}, \ref{ass:BS}, and \ref{ass:new:for:S}.

 The symbol $L$ defined in \eqref{def:L} in terms of ${\bf A}$ and {\bf B} (both defined in Lemma \ref{lem:step3bis}) satisfies Assumption \ref{ass:B}, by property of ${\bf A}$ (order zero, independent of $(x,y)$), and ${\bf B}$ (order 0, uniformly bounded in in $(\e,t)$). 
 
  Assumption \ref{ass:BS} holds too, by Corollary \ref{cor:optimal:bound:S} and Lemma \ref{lem:S:L2:Linfty}.
  
  Assumption \ref{ass:new:for:S} holds by Lemma \ref{lem:away}.
  
  Thus Theorem \ref{th:duh} applies, and gives the representation \eqref{rep:W}.
\end{proof}

We note that $u_{in}$ and $v_{in}$ have equivalent weighted Sobolev and $\fl1$ norms, in the following sense:
\be \label{uin:vin} \begin{aligned}
 c_1 \| u_{in} \|_{\e,s} & \leq \| v_{in} \|_{\e,s} \leq c_2 \| u_{in} \|_{\e,s}, \\ c_1 \| \d_\e^m u_{in} \|_{\fl1} & \leq \| \d_\e^m v_{in} \|_{\fl1} \leq c_2 \| \d_\e^m u_{in} \|_{\fl1},   \end{aligned}
 \ee
 for some $c_1 > 0,$ $c_2 > 0,$ for all $m < s - 3/2.$ (The notation $\d_\e^m$ is introduced in \eqref{fl1:m}.)
 Indeed, $v_{in}$ was defined in terms of $u_{in}$ in the course of Section \ref{sec:6}, by changes of variables which involve a spectral decomposition, multiplications by $e^{ \pm i \theta},$ and projections and localizations. The corresponding operators are linear bounded in weighted Sobolev and $\fl1$ norms (see Lemma \ref{lem:tensor} and the paragraph just above Definition \ref{def:equiv:remainder}).

\begin{cor} \label{lem:final:upper}
 For $0 \leq t \leq t_\star(\e) \leq T \sqrt \e |\ln \e|,$ for $u_{in}$ defined in \eqref{def:U:in:out:high}, we have  
 $$ \begin{aligned} \| u_{in}(t) \|_{\e,s_1}  & \lesssim \e^K |\ln \e|^* e^{t \g/\sqrt \e} \\ &  + \e^{-1/2 + K'} \int_0^t e^{\g (t - t')/\sqrt \e} (\| u_{out} \|_{\e,s_1} + \| u_{high} \|_{\e,s_1 + 1})(t') \, dt'.\end{aligned}$$
 where the rate of growth $\gamma$ is defined in \eqref{def:gamma}. 
 \end{cor}
 
 A choice for $s_1$ is made at the end of the proof of Proposition \ref{prop:fl1:final}. In particular, $s_1$ is chosen much smaller than $s.$ 

\begin{proof} The datum for $u_{in}$ is $O(\e^K)$ in Sobolev weighted norm. Proposition \ref{prop:opS} gives an $L^2 \to L^2$ bound for $\op_\e(S).$ Since 
$$ \d^\a_\e \op_\e(S) f = \sum_{\a_1 + \a_2 = \a} \op_\e(\d_\e^{\a_1} S) \d_\e^{\a_2} f,$$
and since $S$ and $\d^\a_\e S$ enjoy the same bounds (see Corollary \ref{cor:optimal:bound:S}), we have
$$ \| \op_\e(S) f \|_{\e,s_1} \lesssim |\ln \e|^\star e^{t \g/\sqrt \e} \| f \|_{\e,s_1}, \quad \mbox{for all $f \in H^{s_1},$ all $s_1 \leq s.$}$$
 By Lemma \ref{lem:bd-S}, the correctors $S_q$ satisfy the same spatial derivative bounds. This implies 
 \be \label{correctors:up} \| \op_\e({\bf S}) f  \|_{\e,s_1} \lesssim |\ln \e|^\star e^{\g (t - \t)/\sqrt \e} \| f \|_{\e,s_1}, \quad \mbox{for all $f \in H^{s_1},$ all $s_1 \leq s.$}
 \ee 
The remainder ${\bf F}$ is equivalent to $F_{in}.$ 
 Thus by the Sobolev estimate \eqref{bd:F:sob} and Lemmas \ref{lem:dot-checkv} (equivalence of Sobolev norms for $\dot u$ and $\check u$) and \ref{lem:reconstruct} (control of $\check u$ by $u_{in},$ $u_{out}$ and $u_{high}$), we have
 \be \label{bd:bfF:Sob} \| {\bf F}\|_{\e,s_1} \lesssim \e^{-1/2 + K'} (\| u_{in} \|_{\e,s_1} + \| u_{out} \|_{\e,s_1} + \| u_{high} \|_{\e,s_1 + 1}) + \e^{K_a}.\ee 
 With \eqref{uin:vin} and the Sobolev bound on $R_1$ and $R_2$ (see Proposition \ref{prop:duh}), we then obtain 
$$ \begin{aligned} \| u_{in}(t) \|_{\e,s_1}  & \lesssim \e^K |\ln \e|^* e^{t \g/\sqrt \e} \\ &  + \e^{-1/2 + K'} \int_0^t e^{ (t - t') \g/\sqrt \e} (\| u_{in} \|_{\e,s_1} + \| u_{out} \|_{\e,s_1} + \| u_{high} \|_{\e,s_1 + 1})(t') \, dt',\end{aligned}$$
which implies the result, by Gronwall's lemma.
\end{proof}

\begin{lem} \label{lem:fl1:in} 
For $t \leq t_\star(\e):$
 \be \label{bd:fl1:uin} \begin{aligned} \| & u_{in}(t) \|_{\fl1}  \lesssim \e^K e^{t \g/\sqrt\e} \\ & + \e^{-1/2 + K'} \int_0^t e^{(t - t')\g/\sqrt \e} ( \max_{0 \leq t'' \leq t'} \|  u_{out}(t'') \|_{\fl1} + \max_{0 \leq t'' \leq t'} \|  u_{high}(t'') \|_{\fl1}) \, dt' \\ & + \e^{\sigma - 1 + K'} |\ln \e|^\star \int_0^t e^{(t - t') \g/\sqrt \e} \| \dot u(t') \|_{\e,s_1 + 1} \, dt',\end{aligned}\ee 
 where $\sigma > 0$ is introduced in Proposition {\rm \ref{prop:duh}.} 
Besides, using notation \eqref{fl1:m}, for $1 \leq m \leq s_1-3:$ 
$$ \begin{aligned}  \| \d_\e^{m} & u_{in}(t) \|_{\fl1} \lesssim \e^K e^{t \g/\sqrt \e} \\ & + |\ln \e|^\star \e^{-1/2 + K'} \int_0^t e^{(t - t') \g/\sqrt \e} (\max_{\begin{smallmatrix} 0 \leq t' \leq t \end{smallmatrix} }  \| \d_\e^{m} u_{out}(t') \|_{\fl1} + \max_{\begin{smallmatrix} 0 \leq t' \leq t \end{smallmatrix} }  \| \d_\e^{m + 1} u_{high}(t') \|_{\fl1}) \, dt'
\\ & +  \e^{\sigma - 1 + K'} |\ln \e|^\star \int_0^t e^{(t - t') \g/\sqrt \e} \| \dot u(t') \|_{\e, m + 3} \, dt'. 
\end{aligned}$$  
\end{lem}

The exponent $\sigma$ is defined in \eqref{def:sigma}, in terms of $q_0,$ which corresponds to the order of the Taylor expansion in the proof of Corollary \ref{lem:duh-remainder}. The only limitation on $q_0$ is the regularity of the WKB solution.

\begin{proof} By Proposition \ref{prop:opS} and the bound \eqref{op:Sq} in Lemma \ref{lem:bd-S}, we have
$$ \| \op_\e({\bf S}) \|_{\fl1 \to \fl1} \lesssim |\ln \e|^\star e^{t \g /\sqrt \e}.$$ 
Thus from the integral representation \eqref{rep:W} and \eqref{uin:vin}, we deduce the bound
  $$ \begin{aligned} \| u_{in}(t) \|_{\fl1} & \lesssim \e^K |\ln \e|^\star e^{t \g/\sqrt\e} \\ & \quad + |\ln \e|^\star \int_0^t e^{(t - t')\g/\sqrt \e} \left\| (\Id + \e^\sigma R_1) \Big( {\bf F} + \e^{\sigma} R_2(\cdot) v_{in}(0)\Big)(t') \right\|_{\fl1} \, dt'.\end{aligned}$$
  Consider first the $\| {\bf F} \|_{\fl1}$ norm above. We use \eqref{bd:dotF:fl1:apriori}, and the fact that ${\bf F}$ is equivalent to $F_{in}.$ This gives
  $$ \| {\bf F }\|_{\fl1} \lesssim \e^{-1/2 + K'} \| \dot u(t) \|_{\fl1} + \e^{K_a}, \qquad t \leq t_\star(\e).$$
  Next consider the $R_1 {\bf F}$ term. Here $R_1$ is a remainder coming from Proposition \ref{prop:duh} (itself a consequence of Theorem \ref{th:duh} in the appendix). We use the Sobolev embedding \eqref{embed:loss} and the bound on $R_1$ from Proposition \ref{prop:duh}:
  $$ \e^\sigma \| R_1 {\bf F}\|_{\fl1} \lesssim \e^{\sigma - 1/2} \| R_1 {\bf F} \|_{\e,s_1} \lesssim \e^{\sigma - 1/2} |\ln \e|^\star \| {\bf F} \|_{\e,s_1}, \quad \mbox{since $s_1 > 3/2.$}$$
  Thus by \eqref{bd:bfF:Sob}, 
 $$ \e^{\sigma} \| R_1 {\bf F}\|_{\fl1} \lesssim \e^{K_a} + \e^{\sigma - 1 + K'} |\ln \e|^\star \| \dot u \|_{\e,s_1 + 1}. $$ 
  For the $R_2$ term the same reasoning gives
 $$ \| ({\rm Id} + \e^{\sigma} R_1) \e^{\sigma} R_2 v_{in}(0) \|_{\fl1} \lesssim \e^{K}.$$
   Summing up, we obtained
 $$ \begin{aligned} \| u_{in} \|_{\fl1} \lesssim \e^K |\ln \e|^\star e^{t \g/\sqrt\e} & + \e^{-1/2 + K'} |\ln \e|^\star \int_0^t e^{(t - t')\g/\sqrt \e} \| \dot u \|_{\fl1} \, dt' \\ & + \e^{\sigma - 1 + K'} |\ln \e|^\star \int_0^t e^{(t - t')\g/\sqrt \e} \| \dot u \|_{\e,s_1 + 1} \, dt' 
 \end{aligned}$$ 
and the bound \eqref{bd:fl1:uin} follows.

Next we apply $\d_\e^{m}$ to \eqref{rep:W}. Looking at the proof of Proposition \ref{prop:opS}, we see that $\| \op_\e(\d_\e^\a S) \|_{\fl1 \to \fl1} \lesssim |\ln\e|^\star e^{t \g/\sqrt \e}.$ Thus we find  
$$  \| \d_\e^{m} u_{in} \|_{\fl1} \lesssim \e^K e^{t \g/\sqrt \e} + |\ln \e|^\star \sum_{0 \leq m' \leq m} \int_0^t e^{(t - t') \g/\sqrt \e} \big( \| \d_\e^{m'} {\bf F} \|_{\fl1} + \e^{\sigma - 1/2} \| {\bf F} \|_{\e,m' + 2} \big) \, dt'.$$  
 Then we use for $\d_\e^m {\bf F}$ the bound \eqref{final:F:fl1}, which gives
$$ %
  \| \d_\e^m {\bf F} \|_{\fl1}  \lesssim \e^{K_a} + \e^{-1/2 + K'} \big( \| \d_\e^m u_{in} \| + \| \d_\e^m u_{out} \|_{\fl1} + \| \d_\e^{m+1} u_{high} \|_{\fl1} \big), \qquad t \leq t_\star(\e).%
 $$ %
Finally, for the Sobolev norm of ${\bf F}$ that appears above, we use \eqref{bd:bfF:Sob}.
\end{proof} 

\section{A refined Sobolev bound} \label{sec:refined:bound}

We now put together the upper bounds for $u_{high}$ (Section \ref{sec:high}), $u_{out}$ (Section \ref{sec:out}) and $u_{in}$ (Section \ref{sec:Duh}), and derive an upper bound for $\dot u$ that is much improved compared to the rough bound of Section \ref{sec:weak:Sobolev:bound}.

\begin{prop}[Refined Sobolev bound] \label{cor:the:upper:bound} For an appropriate choice of the spatial truncation $\chi_{long},$ depending on $T$ and $K',$ if $s - s_1$ is large enough (meaning $s_a$ large enough), depending on $T,$ $K'$ and $\Gamma - \gamma,$ for $t \leq t_\star(\e)$ the solution $\dot u$ to the initial-value problem \eqref{ivp:dotu} satisfies the bound: 
 $$ \| \dot u(t) \|_{\e,s_1} \lesssim \e^K |\ln \e|^\star e^{t \g/\sqrt \e}.$$
 \end{prop}

The rate $\g$ given in Proposition \ref{cor:the:upper:bound} is much smaller than the rate $\G$ given in Proposition \ref{prop:weak:bound}. The key is that $\g$ is optimal, as seen in Section \ref{sec:lower} below. The condition bearing on $s - s_1$ is inequality \eqref{s:s1} in the proof below.

\begin{proof} Consider first the upper bound for $u_{out}$ from Proposition \ref{prop:out}: for all $s_1 \leq s,$ 
$$ \begin{aligned} \| u_{out}(t) \|_{\e,s_1}^2 & \lesssim  \e^{2 K_a} e^{t \g_{out}/\sqrt \e} +  \e^{K'} |\ln \e|^\star e^{t \g_{out}/\sqrt \e} \max_{0 \leq t' \leq t} \big( \| u_{in}(t') \|^2_{\e,s_1} + \| u_{high}(t') \|^2_{\e,s_1} \, \big).
\end{aligned}
$$%
For $t \leq T \sqrt \e |\ln \e|,$ we have $e^{t \g_{out}/\sqrt \e} \leq \e^{-T \g_{out}}.$ We may choose the spatial truncation $\chi_{long}$ to have such a large support that $\g_{out}$ is so small that $K' - T \g_{out} > 0.$ In particular, the above {\it out} bound then implies 
\be \label{out:in:proof:2}
\| u_{out}(t) \|_{\e,s_1}^2 \lesssim \e^{2 K_a - T \g_{out}} + \e^{K' - T \g_{out}} |\ln \e|^\star \max_{0 \leq t' \leq t} \big( \| u_{in}(t') \|^2_{\e,s_1} + \| u_{high}(t') \|^2_{\e,s_1} \big).\ee 
Next consider the upper bound for $u_{high}$ from Proposition \ref{prop:high}: for all $s_1, s_2$ with $s_1 + s_2 \leq s,$ 
\be \label{high:in:proof} \begin{aligned} \| u_{high}(t)\|^2_{\e,s_1} &  \lesssim \e^{2 K_a} + \e^{2 K + s_2} |\ln \e|^\star e^{t (2 \Gamma + \g_{high})/\sqrt \e}  \\ & +  \e^{K'} |\ln \e|^\star  e^{t \g_{high}/\sqrt \e} \max_{ 0 \leq t' \leq t} (\| u_{in}(t') \|^2_{\e,s_1} + \| u_{out}(t') \|^2_{\e,s_1}).\end{aligned}
\ee 
We now {\it choose} $\g_{high}$ (this parameter is indeed arbitrarily small, see the proof of Proposition \ref{prop:high}) so that
$$ K' - T(\g_{out} + \g_{high}) > 0.$$
 Then, when we plug \eqref{out:in:proof:2} into \eqref{high:in:proof}, we find 
 \be \label{high:in:proof:1.6} \begin{aligned}
 \max_{0 \leq t' \leq t} \| u_{high}(t') \|_{\e,s_1}^2   \lesssim \e^{2 K_a}&  + \e^{2 K + s_2} |\ln \e|^\star e^{t (2 \Gamma + \g_{high})/\sqrt \e} \\ & +  \e^{K' - T \g_{high}} |\ln \e|^\star \max_{0 \leq t' \leq t} \| u_{in}(t') \|^2_{\e,s_1}. 
\end{aligned}\ee 
We now plug \eqref{high:in:proof:1.6} into \eqref{out:in:proof:2}, and find 
\be \label{out:in:proof:3} 
\begin{aligned}
 \| u_{out}(t) \|_{\e,s_1}^2  & \lesssim \e^{2 K_a - T \g_{out}} + \e^{2 K + s_2 + K' - T \g_{out}} |\ln \e|^\star e^{t (2 \Gamma + \g_{high})/\sqrt \e} \\ &  +  \e^{K' - T \g_{out}} |\ln \e|^\star \max_{0 \leq t' \leq t} \| u_{in}(t') \|^2_{\e,s_1}. 
\end{aligned}\ee
At this point we took out the {\it out} term in the {\it high} bound, and conversely the {\it high} term in the {\it out} bound. 

We want to use both \eqref{high:in:proof:1.6} and \eqref{out:in:proof:3} into the {\it in} bound from Corollary \ref{lem:final:upper}, which we reproduce here:
\be \label{in:in:proof} \begin{aligned} \| u_{in}(t) \|_{\e,s_1}  & \lesssim \e^K |\ln \e|^* e^{t \g /\sqrt \e} \\ &  + \e^{-1/2 + K'} \int_0^t e^{(t - t') \g/\sqrt \e} (\| u_{out} \|_{\e,s_1} + \| u_{high} \|_{\e,s_1 + 1})(t') \, dt'.\end{aligned}\ee
  First, we have to modify a little the above {\it in} bound, so that it features squares of norms. 
We let 
 $$ U_{\star, s'}(t) := e^{-  t \g/\sqrt \e} \max_{0 \leq t' \leq t} \| u_{\star}(t') \|_{\e,s'}, \qquad s' \leq s, \quad \star \in \{ in, out, high\}.$$
In \eqref{in:in:proof}, we use Cauchy-Schwarz on the small time interval $[0, t] \subset [0, T \sqrt \e |\ln \e|],$ and find
\be \label{in:in:proof:1.5} \begin{aligned}
 U_{in,s_1}^2(t) & \lesssim \e^{2 K} |\ln \e|^\star + \e^{-1/2 + 2K'} \int_0^t (U_{out,s_1}^2(t') + U_{high,s_1 + 1}^2(t')) \, dt'.
\end{aligned}  \ee
Using $\max |\cdot|^2 = (\max |\cdot|)^2,$ the {\it in} bound \eqref{in:in:proof} takes the form
$$ %
\begin{aligned} \max_{0 \leq t' \leq t} & \| u_{in}(t') \|_{\e,s_1}^2 \lesssim \e^{2 K} |\ln \e|^* e^{2 t \g/\sqrt \e} \\ & +  \e^{-1/2 + 2 K'} |\ln \e|^\star \int_0^t e^{2 (t - t') \g/\sqrt \e} (\max_{0 \leq t'' \leq t'} \| u_{out}(t'') \|_{\e,s_1}^2 + \max_{0 \leq t'' \leq t'} \| u_{high}(t'')\|_{\e,s_1 + 1}^2 ) dt'.
\end{aligned}
$$ %
We use \eqref{out:in:proof:3} (with $s_2 = s - s_1$) in the above, and find 
$$ %
 \begin{aligned} 
\max_{0 \leq t' \leq t} & \| u_{in}(t') \|_{\e,s_1}^2 \lesssim \e^{2 K} |\ln \e|^* e^{2 t \g/\sqrt \e}  + \e^{2 (K + K') + s - s_1} |\ln \e|^\star e^{2 t \Gamma/\sqrt\e} \\ & +  \e^{-1/2 + 3 K' - T \g_{out}} |\ln \e|^\star \int_0^t e^{2 (t - t') \g/\sqrt \e} \max_{0 \leq t'' \leq t'} \| u_{in}(t'') \|_{\e,s_1}^2 \\ & +  \e^{-1/2 + 2 K'} |\ln \e|^\star \int_0^t e^{2 (t - t') \g/\sqrt \e} \max_{0 \leq t'' \leq t'} \| u_{high}(t'')\|_{\e,s_1 + 1}^2  dt'.
\end{aligned}
$$ %
We now use \eqref{high:in:proof:1.6} in the above. For the {\it high} bound \eqref{high:in:proof:1.6}, the index is $s_1 + 1,$ so that $s_2 = s - s_1 - 1$ is admissible. This gives  
$$ %
 \begin{aligned} 
\max_{0 \leq t' \leq t} & \| u_{in}(t') \|_{\e,s_1}^2 \lesssim \e^{2 K} |\ln \e|^* e^{2 t \g/\sqrt \e}  + \e^{2 (K + K') + s - s_1 -1 - T \g_{high}} |\ln \e|^\star e^{2 t \Gamma/\sqrt\e} \\ & +  \e^{-1/2 + 3 K' - T \g_{out}} |\ln \e|^\star \int_0^t e^{2 (t - t') \g/\sqrt \e} \max_{0 \leq t'' \leq t'} \| u_{in}(t'') \|_{\e,s_1}^2 \\ & +  \e^{-1/2 + 3 K' - T \g_{high}} |\ln \e|^\star \int_0^t e^{2 (t - t') \g/\sqrt \e} \max_{0 \leq t'' \leq t'} \| u_{in}(t'')\|_{\e,s_1 + 1}^2  dt'.
\end{aligned}
$$
Going back to the $U_{in,s'}$ notation and using the smallness of $\g_{out},$ we find 
$$ %
 \begin{aligned} 
 U_{in,s_1}^2(t) \lesssim \e^{2 K} |\ln \e|^* & + \e^{2 (K + K') + s - s_1 -1 - T \g_{high}} |\ln \e|^\star e^{2 t (\Gamma - \g)/\sqrt\e} \\ &  +  \e^{3 K' - T \g_{high}} |\ln \e|^\star U_{in,s_1 + 1}^2(t).\end{aligned}
 $$ %
We now choose $s$ and $s_1$ such that $s > 1 + s_1 + 2 T (\Gamma - \gamma).$ (A stronger constraint on $s - s_1$ awaits.) Then, for $t \leq T \sqrt \e |\ln \e|,$  
$$ \e^{2 (K + K') + s - s_1 -1 - T \g_{high}} |\ln \e|^\star e^{2 t (\Gamma - \g)/\sqrt\e} \lesssim \e^{2 K} |\ln \e|^\star,$$
implying 
$$ %
U_{in,s_1}^2(t) \lesssim \e^{2 K} |\ln \e|^*  +  \e^{3 K' - T \g_{high}} |\ln \e|^\star U_{in,s_1 + 1}^2(t),
$$ %
and, by a straightforward induction,
$$ %
U_{in,s_1}^2(t) \lesssim \e^{2 K} |\ln \e|^\star + \e^{(3 K' - T \g_{high})(s-s_1)} U_{in,s}^2(t).
$$ %
At this point we use the rough bound of Proposition \ref{prop:weak:bound}, which implies
$$ U_{in,s}(t) \lesssim \e^K |\ln \e|^\star e^{t (\Gamma - \g)/\sqrt \e}.$$
This gives
\be \label{in:final}
U_{in,s_1}^2(t) \lesssim \e^{2 K} |\ln \e|^\star + \e^{2 K + (3 K' - T \g_{high})(s-s_1)} e^{2 t (\Gamma - \gamma)/\sqrt \e}.
\ee 
Under the condition
\be \label{s:s1}
 s > s_1 + \frac{2 T (\Gamma - \gamma)}{3 K'},
\ee
we can choose $\g_{high}$ small enough so that \eqref{in:final} implies the bound 
$$ U_{in,s_1}(t) \lesssim \e^{K} |\ln \e|^\star,$$
which translates into the expected bound for $\| u_{in} \|_{\e,s_1}.$ By \eqref{high:in:proof:1.6} and \eqref{s:s1}, this implies the same bound for $\| u_{high} \|_{\e,s_1},$ and by \eqref{out:in:proof:3} and \eqref{s:s1} the same for the {\it out} component.
\end{proof}

\section{An upper bound in the $\fl1$ norm} \label{sec:upper:fl1}

We now put together the upper bounds in $\fl1$ norm from Sections \ref{sec:high}, \ref{sec:out} and \ref{sec:Duh}.

\begin{prop}[Upper bound in the $\fl1$ norm] \label{prop:fl1:final} For $t \leq t_\star(\e),$ if $s_1$ satisfies the condition given in \eqref{cond:m} below, the solution $\dot u$ to the initial-value problem \eqref{ivp:dotu} satisfies the bound: 
 $$ \| \dot u(t) \|_{\fl1} + \| \d_\e \dot u(t)\|_{\fl1} \lesssim \e^K |\ln \e|^\star e^{t \g/\sqrt \e}.$$
  \end{prop}

The rate $\g$ is the same as in Proposition \ref{cor:the:upper:bound}; it is defined in \eqref{def:gamma}.

\begin{proof}
We let 
$$ V_{\star,m}(t) := e^{-t \g/\sqrt \e} \max_{\begin{smallmatrix} 0 \leq t' \leq t \\ 0 \leq |\a| \leq m \end{smallmatrix} } \| \d_\e^{\a} u_{\star}(t') \|_{\fl1}.$$
Then the high-frequency bounds of Proposition \ref{prop:high:fl1} take the form
 \be \label{fl1:2.1}
V_{high,m}(t)  \lesssim \e^{K_a} + \e^{K'} |\ln \e|\big( V_{high,m+1}(t) + V_{in,m}(t) + V_{out,m}(t) \big).
\ee
The {\it out} bounds of Proposition \ref{prop:out:fl1} take the form
\be \label{fl1:2.2} 
 V_{out,m}(t) \lesssim \e^{K_a - T \g_{out}} + \e^{K' - T \g_{out}} (V_{in,m}(t) + V_{high,m + 1}(t)).
 \ee
With the refined Sobolev bound of Proposition \ref{cor:the:upper:bound}, the {\it in} bounds of Lemma \ref{lem:fl1:in} imply  
\be \label{fin:fl1:in:2} \begin{aligned}  V_{in,m}(t) \lesssim \e^K + \e^{K'} |\ln \e|^\star  (V_{out,m}(t) + V_{high,m+1}(t)).
\end{aligned}
\ee 
We denote $V_m = V_{in,m} + V_{out,m} + V_{high,m},$ and obtain, by a straightforward induction,
$$ V_1 \lesssim \e^K + \e^{(m-1)(K' - T \g_{out})} |\ln \e|^\star V_m.$$
For $m \leq s_1 -2,$ we can use the Sobolev embedding \eqref{embed:loss} then the refined Sobolev bound of Proposition \ref{cor:the:upper:bound}. This gives
$$ V_1 \lesssim \e^K + \e^{(m-1)(K' - T \g_{out})} \e^{-1/2 + K} |\ln \e|^\star.$$
Thus (with $m = s_1 -3$) the condition bearing on $s_1$ appears to be 
\be \label{cond:m}
 (s_1-4) (K' - T \g_{out}) \geq \frac{1}{2}.
\ee  
Under this condition, we find $V_m \lesssim \e^K |\ln \e|^\star,$ which translates into the result. We choose $s_1$ to be the smallest integer such that \eqref{cond:m} holds true. 
\end{proof}

\begin{rem} \label{rem:Sobolev:indices} We find in Proposition {\rm \ref{cor:the:upper:bound}} the constraint \eqref{s:s1} bearing on $s.$ In Proposition {\rm \ref{prop:fl1:final}}, we have a constraint bearing on $s_1.$ Taking into account the definitions of $T$ \eqref{def:Tstar}, the upper bound for $s$ in terms of $s_a$ \eqref{def:s}, and the constraint $K_a \geq K + 1/2,$ this means for $s_a$ the lower bound 
\be \label{lower:sa}
  s_a > 6 + 2 K + \frac{1}{2K'} + \frac{2(\Gamma - \gamma)}{3 \gamma} \cdot \frac{K - K'}{K'}.
  \ee
The proof of the symbolic flow Theorem {\rm \ref{th:duh}} introduces another constraint on $s_a,$ in the form $s_a > q_0 + s_1 + 4  + d/2,$ with $q_0$ satisfying \eqref{def:sigma}. With the constraint on $s_1$ from Proposition {\rm \ref{prop:fl1:final}}, this means
\be \label{lower:sa:2}
s_a > 7 + \frac{3}{2} + 2 C_0 + \frac{1}{2K'} + 2 (K - K'),
\ee
where $C_0$ is the constant depending only on the spatial dimension which appears in Corollary {\rm \ref{lem:duh-remainder}}. In the small $K'$ and large $K$ limit, the first lower bound \eqref{lower:sa} implies the second \eqref{lower:sa:2}. 
\end{rem}

\section{Existence up to the optimal amplification time} \label{sec:existence}

Much of the analysis so far has been devoted to showing that {\it so long as the solution (in weighted $W^{1,\infty}$ norm) is controlled by $\e^{K'},$ it is actually controlled by $\e^K e^{t \g/\sqrt \e}.$}

Precisely, recall that $t_\star(\e)$ was defined in the last paragraph of Section \ref{sec:tstar} as the largest time $t \in [0, T \sqrt \e |\ln \e|] \cap [0, T_\star(\e))$ for which 
$$ \max_{0 \leq t' \leq t} ( \| \dot u(t) \|_{\fl1} + \| \d_\e \dot u(t) \|_{\fl1}) \leq \e^{K'}.$$
Here $T_\star(\e)$ is the maximal existence time. %
Now we know from the above Section \ref{sec:upper:fl1} that for $t \leq t_\star(\e),$ we have the improved bound
\be \label{the:fl1:bound} \| \dot u(t) \|_{\fl1} + \| \d_\e \dot u(t) \|_{\fl1} \leq C |\ln \e|^{C'} \e^{K} e^{t \g/\sqrt \e},\ee
for some constants $C > 0$ and $C' > 0$ which are independent of $\e$ and $t.$ 
We observe that
$$ \e^K e^{t \g/\sqrt \e} = \e^{K'} \quad \mbox{for} \quad t = T \sqrt \e |\ln \e|,$$
with $T$ defined in \eqref{def:Tstar}. Let $T_\e$ be defined by 
\be \label{def:T0}
 T_\e := T - \frac{C' \ln |\ln \e| + \ln C}{\g} \sqrt \e = \frac{K - K'}{\g} \sqrt \e |\ln \e| -  \frac{C' \ln |\ln \e| + \ln C}{\g} \sqrt \e,
\ee
where $C$ and $C'$ are the constants that appear in \eqref{the:fl1:bound}.
Then, 
$$ C |\ln \e|^{C'} \e^K e^{t \g/\sqrt \e} = \e^{K'}, \quad \mbox{for} \quad t = T_\e \sqrt \e |\ln \e|.$$
This means in particular that 
\be \label{t0} T_\e \sqrt \e |\ln \e| \leq t_{\star}(\e).\ee 
Indeed, if \eqref{t0} were not true, then the weighted $W^{1,\infty}$ norm over $[0, t_\star(\e)]$ would be strictly smaller than $\e^{K'}.$ By the continuation criterion\footnote{Equivalently, by the strict inequality $t_\star(\e) < T_\star(\e),$ with notation from Section \ref{sec:tstar}.} for solutions to first-order, quasilinear symmetric hyperbolic systems, we could extend the solution a little beyond $t_\star(\e),$ and by continuity of the maximum in time of the weighted $W^{1,\infty}$ norm, this norm would not reach $\e^{K'}$ even a little after $t_\star(\e),$ contradicting the definition of $t_\star(\e).$ 

 In conclusion, the perturbative solution $\dot u$ is defined up to time $T_\e \sqrt \e |\ln \e|,$ and satisfies the bound \eqref{the:fl1:bound} over $[0, T_\e \sqrt \e |\ln \e|],$ and also the weighted Sobolev bound of Proposition \ref{cor:the:upper:bound} within that interval.

\section{Lower bound} \label{sec:lower}

In order to conclude the proof of Theorem \ref{th:main} it now suffices to show that the datum $v_{in}(0)$ is maximally amplified by the solution operator $\op_\e(S(0;t)).$

In the following Lemma, we denote $v_{in}(0)_j \in \C^{14},$ with $j \in \{1,2_{out}, 2_{in}, 3, 4_{out}, 4_{in},5\}$ the components of $v_{in}(0),$ with $v_{in}$ defined in \eqref{b5}.

We define vectors $\vec e_0$ and $\vec f_0$ by their components in the decomposition of Section \ref{sec:coord}, so that $B_{|\vec e_0}$ denotes the magnetic field ($B$) component of $\vec e_0,$ etc., as follows:
\be \label{def:e1}
B_{|\vec e_0} = 0, \quad E_{|\vec e_0} = 0, \quad v_{e|\vec e_0} = (0,1,0), \quad n_{e|\vec e_0} = 0, \quad v_{i|\vec e_0} = 0, \quad n_{i|\vec e_0} = 0,
\ee
and
\be \label{def:e2}
B_{|\vec f_0} = 0, \quad E_{|\vec f_0} = 0, \quad v_{e|\vec f_0} = (1,0,0), \quad n_{e|\vec f_0} = 0, \quad v_{i|\vec f_0} = 0, \quad n_{i|\vec f_0} = 0.
\ee
Note that, with these definitions of $\vec e_0$ and $\vec f_0,$ the divergence equations \eqref{div:eq} are trivially satisfied by the datum \eqref{initial:datum}.

In accordance with the computations of Section \ref{sec:interaction:coefficients}, we choose $\xi_0 \in \R$ (describing the oscillation in the initial perturbation) to be 
\be \label{def:xi0}
\xi_0 := \left\{\begin{aligned} \xi_{14}^-, & \quad \mbox{if $k > \sqrt 3,$} \\ \xi_{14}^+, & \quad \mbox{if $k < - \sqrt 3.$}
\end{aligned}\right.
\ee

\begin{lem} \label{lem:look:datum} With the above choice \eqref{def:e1} and \eqref{def:e2} for the vectors $\vec e_0$ and $\vec f_0$ that appear in the initial datum \eqref{def:phie}, and \eqref{def:xi0} for the frequency $\xi_0$ that appear in the initial datum, the datum $v_{in}(0)$ satisfies
$$ \begin{aligned}
 \left(\begin{array}{c} v_{in}(0)_{1} \\ v_{in}(0)_{4_{out}} \end{array}\right) = & \frac{\e^K}{\sqrt 2} e^{i x \xi_0/\e} \phi(x,y)\left(\begin{array}{c} \l_1(\xi_0 + k)^{-1} e_1(\xi_0 + k,0) \\ e_4(\xi_0,0) \end{array}\right) \\
 & + \mbox{other oscillating terms} + O(\e^{K + 1/2}).
 \end{aligned}$$
Above, the vector $e_1$ is a transverse eigenvector of the hyperbolic symbol, and $e_4$ is a longitudinal eigenvector of the hyperbolic symbol, as defined in \eqref{e1:e5} and \eqref{e2:e4}. 
 The other components are given by 
$$ v_{in}(0)_j = \mbox{other oscillating terms} + O(\e^{K + 1/2}), \qquad j \notin \{1, 4_{out}\},$$
where
\begin{itemize}
\item ``other oscillating terms" refer to terms of the form $\e^K e^{i x \xi_1/\e} \phi_1(x,y),$ for some $\xi_1 \in \R$ and some compactly supported and vector-valued $\phi_1,$ which will be seen not to be maximally amplified by the solution operator, in the sense that
$$ \| \op_\e(S(0;t))(\e^K e^{ix \xi_1/\e} \phi_1) \|_{L^2} \lesssim \e^K |\ln \e|^\star e^{t \g_1/\sqrt \e} + \e^{K + 1/4} |\ln \e|^\star \e^{t \g/\sqrt \e},$$ for $t \leq T \sqrt \e |\ln \e|,$
with $\g_1 < \g,$ where $\g$ is the optimal growth rate given in Corollary {\rm \ref{cor:optimal:bound:S}};
\item remainders $O(\e^{K + 1/2})$ are understood in $L^2$ norm.
\end{itemize}
\end{lem}

\begin{proof} By definition of $\phi^\e$ in \eqref{def:phie} and description of $u_{in}(0)$ in \eqref{datum:uin}, we have
$$ \begin{aligned} u_{in}(0) & = \e^K e^{i x (\xi_0 + k)/\e} \op_\e(\chi_{low}(\xi_0 + k + \cdot)) (\phi \vec e_0) \\ & + \e^K e^{i x \xi_0/\e} \op_\e(\chi_{low}(\xi_0  + \cdot)) (\phi \vec f_0) \\ & + \e^{K} e^{- i x (\xi_0 + k)/\e} \op_\e(\chi_{low}(- \xi_0 - k + \cdot)) (\phi \vec e_0) \\ & + \mbox{o.o.t.} + O(\e^{K + 1/2}),\end{aligned}$$
where $\mbox{o.o.t.}$ stands for ``other oscillating terms" and $O(\cdot)$ is understood in $\| \cdot \|_{\e,s}$ norm. The term in the third line in the above right-hand side belongs to the $\mbox{o.o.t.}$ category. We will describe its fate in detail; the other $\mbox{o.o.t.}$ terms are similar.

By description of $v_{in}(0)$ in \eqref{datum:vin}, the above implies
$$  \begin{aligned}
 v_{in}(0)_{1}  &  =  \e^K e^{i x \xi_0/\e} \op_\e( (\chi_{low} \Pi_1) (\xi_0 + k + \cdot)) (\phi \vec e_0) \\  & + \e^K e^{i x (\xi_0 - k)/\e} \op_\e( (\chi_{low} \Pi_1) (\xi_0 + \cdot)) (\phi \vec f_0) \\ &  + \e^K e^{- i  x (- \xi_0 - 2 k)/\e} \op_\e( (\chi_{low} \Pi_1)(- \xi_0 - k + \cdot)(\phi \vec e_0) \\ & + \mbox{o.o.t.} + O(\e^{K  + 1/2}).\end{aligned}$$

Since $\phi$ is smooth and compactly supported,
$$ \op_\e((\chi_{low} \Pi_1)(\xi_0 + k + \cdot)) (\phi \vec e_0) = ((\chi_{low} \Pi_1)(\xi_0 + k) (\phi \vec e_0) + O(\e^{1/2}),$$
with $O(\cdot)$ understood in $\| \cdot \|_{\e,s}$ norm still. The support of $\chi_{low}$ being large enough around the resonant set ${\mathcal R}$ (see Figure \ref{fig:r:chi}), we have
$$\chi_{low}(\xi_0) = 1, \qquad \chi_{low}(\xi_0 + k) = 1.$$
By definition of $\vec e_0$ in \eqref{def:e1} and the spectral computations of Section \ref{sec:spectral:dec} and \ref{sec:interaction:coefficients} in which eigenvectors of the hyperbolic symbol are given in detail, we have:
$$ \Pi_1(\xi_0 + k) \vec e_0 = \frac{1}{\sqrt 2 \l_1(\xi_0 + k)} e_1(\xi_0 + k), \quad \mbox{and} \quad \Pi_1(\xi_0) \vec f_0 = 0.$$
Thus
$$ \begin{aligned} v_{in}(0)_1 & = \frac{\e^K e^{i x \xi_0/\e}}{\sqrt 2 \l_1(\xi_0 + k)} e_1(\xi_0 + k) + \e^K e^{- i  x (- \xi_0 - 2 k)/\e} (\chi_{low} \Pi_1)(- \xi_0 - k)(\phi \vec e_0) \\ & + \mbox{o.o.t.}  + O(\e^{K + 1/2}).\end{aligned}$$
We focus now on the second term in the above right-hand side. The goal is to show that this term is not maximally amplified, as in the statement of the present Lemma. This term has the form $\e^K  e^{- i  x (- \xi_0 - 2 k)/\e} \phi_1(x,y),$ for some smooth and compactly supported $\phi_1.$ We use Lemma \ref{lem:for:datum}:
$$ \begin{aligned} \| \op_\e(S(0;t))   \e^K  e^{- i  x (- \xi_0 - 2 k)/\e} \phi_1(x,y) & - \e^K e^{- i  x (- \xi_0 - 2 k)/\e} S(0;t,x,y,-\xi_0 - 2k,0) \phi_1 \|_{L^2} \\ & \lesssim \e^{K + 1/4} |\ln \e|^\star e^{t \g/\sqrt \e}.\end{aligned}$$

Thus it only remains to show that $S(0;t,x,y,-\xi_0 - 2k,0)$ grows at an exponential rate smaller than the optimal rate $\g.$ 

If $k > \sqrt 3,$ then $\xi_0 = \xi_{14}^-.$ We see on Figure \ref{figure:trace} and \eqref{1425} that the only other frequency at which the maximal growth rate is attained is $\xi_{25}^+ = (\o^2 - 2 \o)^{1/2}.$ We have 
$$ - \xi_0 - 2 k = (\o^2 - 2\o)^{1/2} - k \notin \{ \xi_{25}^+, \xi_{14}^-\},$$  
 so that the flow at $(-\xi_0 - 2 k,0)$ is not maximally amplified. (If $k < - \sqrt 3,$ then it's the same story.) 

This completes the description of $v_{in}(0)_1.$ The computations for $v_{in}(0)_{4_{out}}$ are similar; we use $\Pi_4(\xi_0) \vec f_0 = \frac{1}{\sqrt 2} e_4(\xi_0),$ and $\Pi_4(\xi_0) \vec e_0 = 0.$

Now onto the other components of $v_{in}(0).$ In these components, the above arguments show that only the terms with oscillations $e^{i \xi_0 x/\e}$ are susceptible to be maximally amplified. Looking at the components of $v_{in}$ as defined in \eqref{b5}, we see that those terms are
\begin{itemize}
\item the term with oscillations in $e^{i (\xi_0 + k) x/\e}$ in $U_2^{in},$ since $U_2^{in}$ comes in with a $e^{- i \theta}$ prefactor; this term is denoted $w_2;$
\item and the term with oscillations in $e^{i \xi_0 x/\e}$ in $U_3,$ since $U_3$ is not shifted; this term is denoted $w_3.$
\end{itemize}
We have
$$ \begin{aligned} w_{2} & = \e^K e^{- i k x/\e} \op_\e(\chi_{23,-1} \chi_{low} \Pi_2) (e^{i (\xi_0 + k) x/\e} \phi \vec e_0) \\ &  = \e^K e^{i x \xi_0/\e} \chi_{23}(\xi_0) \Pi_2(\xi_0 + k) \vec e_0 \phi + O(\e^{K + 1/2}).\end{aligned}$$
By definition of $\vec e_0$ in \eqref{def:e1} and property of longitudinal eigenvectors (see Section \ref{sec:spectral:dec}: their electronic velocities at frequency $\zeta \in \R^3$ are parallel to $\zeta$), we have $\Pi_2(\xi_0 + k) \vec e_0 = 0.$ Thus $w_{2} = O(\e^{K + 1/2}).$

The other term is
$$ \begin{aligned} w_3 = \e^K \op_\e(\Pi_3)(e^{i x \xi_0/\e} \vec f_0 \phi)  = \e^K e^{i x \xi_0/\e} \Pi_3(\xi_0) \vec f_0 \phi + O(\e^{K + 1/2}).\end{aligned}$$
Recall that $\Pi_3$ (defined in Section \ref{sec:char0}) comprises eigenmodes that belong to the kernel of the operator (described in Section \ref{sec:kernel}) and longitudinal acoustic modes that are $O(\sqrt \e |\zeta|)$ (those are described in Section \ref{sec:parallel:modes}). The eigenmodes in the kernel have electronic velocities which are transverse to the frequency. Thus at frequency $(\xi_0,0),$ those electronic velocities are perpendicular to the electronic velocity in $\vec f_0$ defined in \eqref{def:e2}. Besides, as aboserved in Section \ref{sec:(1,3):coeffs}, the longitudinal acoustic eigenmodes have electronic velocities $O(\sqrt \e).$ Thus $w_3 = O(\e^{K + 1/2}),$ which concludes the proof.
\end{proof}

\begin{cor} \label{lem:low} Given a perturbation datum as in Lemma {\rm \ref{lem:look:datum}} with $\phi$ such that $\phi(x_0,y_0) \neq 0,$ we have the lower bound, for $0 \leq t \leq T\sqrt \e |\ln \e|:$
$$ \| \op_\e(S(0;t)) v_{in}(0) \|_{L^2(B((x_0,y_0), \e^\b)} \geq C \e^{K + 3 \b/2} e^{t \g/\sqrt \e}.$$
Above, $(x_0,y_0) \in \R^3$ is the argmax of $|g_\pm|$ over $\R^3,$ and $\b < 1/6.$
\end{cor}

\begin{proof} The computations of Section \ref{sec:interaction:coefficients} show that the maximal rate $\g$ defined in \eqref{def:gamma} is positive and associated with the $(1,4)$ and $(2,5)$ space-time resonances. We choose the initial datum to be polarized along $(1,4).$ (The other choice, polarizing along $(2,5),$ would of course have been equally acceptable.) 

Thus we may focus on the solution $S$ to system \eqref{2blocks} with $(j,j') = (1,4),$ and the component $\left(\begin{array}{c} v_{in}(0)_1 \\ v_{in}(0)_{4_{out}} \end{array}\right) =: v_{in}(0)_{1,4_{out}}$ of the datum. 

By Lemmas \ref{lem:for:datum} and \ref{lem:look:datum},  
$$ \begin{aligned} e^{- i x \xi_0/\e} & \op_\e(S(0;t))  v_{in}(0)_{1,4_{out}} = \frac{\e^K}{\sqrt 2} S(0;t,x,y,\xi_0,0) \phi \left(\begin{array}{c} \l_1(\xi_0 + k)^{-1} e_1(\xi_0 + k,0) \\ e_4(\xi_0,0) \end{array}\right) \\ & + \e^K e^{t \g_1/\sqrt \e} f_1^\e + \e^{K + 1/4} e^{t \g/\sqrt \e} f^\e,\end{aligned}$$
where $\g_1$ is a suboptimal rate: $\g_1 < \g$ and $f_1^\e$ and $f^\e$ are bounded in $L^2,$ uniformly in $\e$ and $t$ within the observation interval.

The focus is on the space-time resonance $(\xi_0,0),$ with $\xi_0$ defined in \eqref{def:xi0}. The trace $\mbox{tr}\, B_{114} B_{-141,+1}$ is maximal at $(\xi_0,0).$ We use Corollary \ref{cor:scalarize} (and its proof): there exists a change of basis $\tilde P$ such that $\tilde S := \tilde P^{-1} S$ solves, at $(x,\xi,\eta) = (x,\xi_0,0),$ the system
$$ \big( \d_t + \frac{i \l_4(\xi_0,0)}{\e} \big) \tilde S + \frac{1}{\sqrt \e} \left(\begin{array}{cc} 0 & \tilde b_{14} \\ \tilde b_{41} & 0 \end{array}\right) \tilde S = 0, \qquad \tilde S(\tau;\tau) = \tilde P^{-1},$$
with
$$ \tilde b_{14} := g_1(x,y) \left(\begin{array}{cc} \g_0 & 0 \\ 0 & 0\end{array}\right), \qquad \tilde b_{41} := g_{-1}(x,y) \left(\begin{array}{cc} \g_0  & 0\\ 0 & 0 \end{array}\right) = \tilde b_{12}^\star,$$
where 
$$ \g_0 = \sqrt{\mbox{tr}\, B_{114} B_{-141,+1}(\xi_0,0)}, \quad\mbox{so that} \quad \g = |g_\pm(x_0,y_0)| \g_0,$$
where $(x_0,y_0)$ is the argmax of the modulus of $g_1$ and $g_{-1}$ over $\R^3$ (recall, the complex conjugate of $g_{-1}$ is equal to $g_1,$ by reality of the initial datum for the WKB approximate solution). 
  The top left blocks of $\tilde b_{ij}$ are scalar. Above, we used $\d_\eta \l_4(\xi_{14}^\pm,0) = 0.$  
This gives 
$$ e^{i (t - \t) \l_4(\xi_0,0)/\e} \tilde S(\t;t,x,y,\xi_0,0) = P_0 \tilde P^{-1} + \sum_\pm e^{\pm (t - \t) |g_\pm(x,y)| \g_0/\sqrt \e} P_\pm \tilde P^{-1},$$
where $P_0$ and $P_\pm$ are the eigenprojectors of matrix
$$\left(\begin{array}{cc} 0 & J_1 \\ J_1 & 0 \end{array}\right) \in \C^{28 \times 28}, \quad \mbox{with} \quad J_1 = \left(\begin{array}{cc} 1 & 0_{\C^{13}} \\ 0_{\C^{13}} & 0_{\C^{13 \times 13}} \end{array}\right).$$
Thus
$$ \begin{aligned} e^{-i x \xi_0/\e} & \op_\e(S(0;t)) v_{in}(0)_{1,4_{out}} \\ & = \frac{\e^K}{\sqrt 2} e^{t |g_\pm(x,y)| \g_0/\sqrt \e} \phi(x,y) \tilde P P_+ \tilde P^{-1} \left(\begin{array}{c} \l_1(\xi_0 + k)^{-1} e_1(\xi_0 + k,0) \\  e_4(\xi_0,0) \end{array}\right) \\ & + \e^{K+1/4} e^{t \g/\sqrt \e} f^\e + \e^{K} e^{t \g_1/\sqrt \e} f_1^\e, \end{aligned}$$
with $\g_1,$ and $f^\e, f^\e_1$ as above. Let $\b > 0.$ Since $g_{\pm} \in C^2,$ for some $C > 0,$ given $(x,y) \in B((x_0,y_0),\e^\b),$  we have 
$$|g_{\pm}(x,y)| \geq |g_{\pm}(x_0,y_0)| - C \e^{2 \b}.$$ This implies
$$ \begin{aligned} \|\op_\e(S(0;t) v_{in}(0)_{1,4_{out}}\|_{L^2(B((x_0,y_0),\e^\b))} & \geq C_0 \e^K e^{t (|g(x_0,y_0) - C \e^{2\b}) \g_0} |\phi|_{L^2(B((x_0,y_0),\e^\b))} \\ & - C_1 (\e^{K + 1/4} e^{t \g /\sqrt \e} + \e^K e^{t \g_1/\sqrt \e}),\end{aligned}$$
where $C_1 > 0$ and 
$$C_0 := \frac{1}{\sqrt 2} \left| \tilde P P_+ \tilde P^{-1} \left(\begin{array}{c} \l_1(\xi_0 + k)^{-1} e_1(\xi_0 + k,0) \\ e_4(\xi_0,0) \end{array}\right) \right|.$$
We need to verify that $C_0 > 0.$ Since $\tilde P$ is invertible, it suffices to check that 
\be \label{verif:lower} z := \tilde P^{-1} \left(\begin{array}{c} \l_1(\xi_0 + k)^{-1} e_1(\xi_0 + k,0) \\   e_4(\xi_0,0) \end{array}\right) \notin \mbox{ker}\, P_+.\ee Given $z' \in \C^{28},$ the norm of $P_+ z'$ is a non-zero constant times $|z'_1 + z'_{15}|.$ At this point we go back to the proof of Lemma \ref{scalarize}: the vector $z$ in the left-hand side of \eqref{verif:lower} is equal to the coordinates of $\tilde P z$ in the basis $\big(e^\sharp, a_1^\sharp,\cdots,a_{n-1}^\sharp,f_\sharp,b_{1\sharp},\cdots,b_{n-1\sharp}\big),$ where the vector $e$ generates the range of $\Pi_1(\xi_0 + k,0) B(\vec e_1) \Pi_4(\xi_0,0),$ and the vector $f$ generates the range of $\Pi_4(\xi_0,0) B(\vec e_{-1}) \Pi_1(\xi_0 + k,0),$ and $e^\sharp = (e,0),$ $f_\sharp = (0,f).$ Thus
 $$  |P_+ z | = |z_1 + z_{15}| = \big|\l_1(\xi_0 + k) \langle e_{1}(\xi_0 + k,0), e \rangle + \langle e_{\parallel}(\xi_0,0), f \rangle\big|,$$
 with $\langle \cdot, \cdot \rangle$ the Hermitian scalar product in $\C^{14}.$
The Hermitian products $\langle e_1(\xi_0 + k,0), e\rangle$ and $\langle e_4(\xi_0,0), f\rangle$ are equal to one. Hence $0 < \l_1(\xi_0 + k,0) + 1$ implies that $P_+ z \neq 0,$ which in turns implies that $C_0$ is indeed positive. 

Finally, we choose $\phi$ so that $\phi(x_0,y_0) \neq 0.$ Then $\| \phi \|_{L^2(B((x_0,y_0), \e^\b)} \geq c_0 \e^{d \b/2},$ for some $c_0 > 0$ and $\e$ small enough. Thus we obtained
$$ \|\op_\e(S(0;t) v_{in}(0)_{1,2}\|_{L^2(B((x_0,y_0),\e^\b))} \geq  \e^K e^{t \g/\sqrt \e} \big(c_0 C_0  \e^{d \b/2} - C_1 \e^{1/4} \big),$$
and it suffices to choose $\b < 1/6$ (recall, the spatial dimension is $d=3$) in order to conclude, since the other components of $v_{in}(0)$ are not maximally amplified by $\op_\e(S(0;t)),$ as proved in Lemma \ref{lem:look:datum}.
\end{proof}

\section{Endgame} \label{sec:endgame}

From the representation \eqref{rep:W} via the symbolic flow, we find, with $(x_0,y_0)$ and $\b$ as in Lemma \ref{lem:low}:  
$$ \begin{aligned} \| v_{in}(t) \|_{L^2(B((x_0,y_0), \e^{\b})} & \geq \| \op_\e({\bf S}(0;t)) v_{in}(0) \|_{L^2(B((x_0,y_0),\e^{\b}))} \\ & -  \int_0^t \Big\| \op_\e({\bf S}(t';t))(\Id + \e^\sigma R_1) \Big( {\bf F} + \e^{\sigma} R_2(\cdot) v_{in}\Big)(t') \Big\|_{L^2(\R^3)} \, dt'.\end{aligned}$$
By Lemma \ref{lem:low} and \eqref{correctors:up}:
$$  \| \op_\e({\bf S}(0;t)) v_{in}(0) \|_{L^2(B((x_0,y_0),\e^{\b}))} \geq C \e^{K + \b d/2} e^{t \g/\sqrt \e}.$$
A bound for ${\bf F}$ is given by \eqref{bd:bfF:Sob}, in terms of $\|\dot u\|_{L^2}.$ With the refined upper bound of Corollary \ref{cor:the:upper:bound}, this gives
$$ \| {\bf F} \|_{L^2} \lesssim 	\e^{-1/2 + K' + K} |\ln \e| e^{t \g/\sqrt \e} + \e^{K_a}.$$
We may assume $K_a > K + K' - 1/2.$ By Theorem \ref{th:duh}, the operator norms of the remainders $R_1$ and $R_2$ are controlled by $|\ln \e|^\star.$ An upper bound for the operator norm of $\op_\e({\bf S})$ is given in \eqref{correctors:up}. Thus for $t \leq \sqrt\e T |\ln \e|:$
$$ \begin{aligned} \int_0^t \Big\| \op_\e({\bf S}(t';t))(\Id & + \e^\delta R_1) \Big( {\bf F} + \e^{\delta} R_2(\cdot) v_{in}\Big)(t') \Big\|_{L^2(\R^3)} \, dt' \\ & \lesssim |\ln \e|^\star \int_0^t e^{\g (t - t')/\sqrt \e} \e^{-1/2 + K + K'} e^{t' \g/\sqrt \e} \, dt'
\\ & \lesssim \e^{K + K'} |\ln \e|^\star \e^{t \g /\sqrt \e}.\end{aligned}$$
We obtained finally
$$ \| v_{in}(t) \|_{L^2(B((x_0,y_0), \e^{\b}))} \geq \e^K e^{t \g/\sqrt \e} \big(C \e^{\b d/2} - C' \e^{K'} |\ln \e|^\star\big),$$
for some $C > 0$ and some $C' > 0$ which depend neither on $\e$ nor on $t.$ Let $\b = K'/d$ (without loss of generality, we may assume $K' < 1/2,$ so that the condition of Lemma \ref{lem:low} is satisfied). Then, for $\e$ small enough,
$$ C \e^{\b d/2} - C' \e^{K'} |\ln \e|^\star \geq C \e^{K'/2},$$
so that 
$$ \left \| v_{in}(T_\e \sqrt \e |\ln \e|) \right \|_{L^2(B(x_0,y_0,\e^{K'/3}))}\geq C \e^{K' + K'/2},$$
with $T_\e$ defined in \eqref{def:T0}. Since
$$  \left \| v_{in}(T_\e \sqrt \e |\ln \e|) \right \|_{L^2(B(x_0,y_0,\e^{K'/3}))} \lesssim \e^{K'/2} \max_{|x - x_0| + |y - y_0| \leq \e^{K'/3}} |  v_{in}(T_\e \sqrt \e |\ln \e|,x,y)|,$$
we have
\be \label{final:for:vin}  \|  v_{in}(T_\e \sqrt \e |\ln \e|,\cdot) \|_{L^\infty} \geq C \e^{K'},\ee
for some constant $C > 0$ which is independent of $\e$ and $t.$

We finally go back to the original unknown $\dot u.$ To this effect, we observe that $v_{in}$ has the form
 $$ v_{in} = \sum_{q \in {\mathcal H}} e^{i q \theta} \op_\e(H_q) \dot u,$$
 where $H_q \in S^0$ does not depend singularly on $\e,$ and ${\mathcal H}$ is finite. By Proposition \ref{prop:pointwise}, this implies
 \be \label{final-1} \|v_{in} \|_{L^\infty} \lesssim |\dot u|_{L^\infty} \big(1 + |\ln \e| + | \ln \| \dot u\|_{\e,s}|\big).\ee
By the refined bound of Proposition \ref{cor:the:upper:bound} and \eqref{t0},
 $$ \sup_{0 \leq t \leq T_\e \sqrt \e |\ln \e|} \big| \ln \| \dot u(t) \|_{\e,s} \big| \lesssim |\ln \e|.$$
 Putting \eqref{final:for:vin} and \eqref{final-1} together, we arrive at
 $$ C \e^{K'}(1 + |\ln \e|)^{-1} \lesssim \|\dot u(T_\e \sqrt \e |\ln \e|)\|_{L^\infty},$$
 which implies the result, since $K'$ is arbitrarily small. 

\medskip

This concludes the main part of the proof of Theorem \ref{th:main}. In Section \ref{sec:wkb}, we sketch the WKB computations. In Section \ref{sec:spectral:dec}, we describe the eigenvalues and eigenvectors of the leading hyperbolic operator. Section \ref{sec:interaction:coefficients} contains a detailed computation of the interaction coefficients at space-time resonances. In Section \ref{app:symb}, classical results on pseudo-differential operators are gathered. Section \ref{app:duh} gives an anisotropic formulation of the symbolic flow method. The final Sections comprise a notation index, an index, and lists of parameters.

\section{WKB computations} \label{sec:wkb}

The fundamental phase is $(\o,(k,0,0)) \in \R \times \R^3.$ We often consider $k$ to be a vector in $\R^3,$ that is we identify when convenient $k \in \R$ with the vector $(k,0,0) \in \R^3.$

The wavenumber $k$ is given in the initial WKB datum \eqref{data}. The associated characteristic frequency is encoded in the polarization condition for the WKB profile $a:$ we assume that for all $(x,y) \in \R^3,$ we have
\be \label{pola:a} \Pi_1(k) a(x,y) = a(x,y),\ee
where $\Pi_1$ is the spectral projector introduced in Section \ref{sec:resonances}. This implies that the initial oscillations at frequency $k$ are propagated by the hyperbolic system at the time frequency $\o := \l_1(k) = (1 + k^2)^{1/2}.$

In the rescaled spatial frame \eqref{def:tildeu}, the Zakharov ansatz is
\be \label{ansatz} u_a(\e,t,x,y) = \sum_{0 \leq k \leq 2 K_a} \e^{k/2} \sum_{q \in {\mathcal H}_k} e^{i q \theta} u_{a,k,q}(t,x,y),\ee
where ${\mathcal H}_k \subset \Z$ is finite, $\theta := (k x - \o t)/\e,$ and the amplitudes $u_{a,k,q}$ are independent of $\e.$ We denote $(B_{k,q}, E_{k,q}, v_{e,k,q}, n_{e,k,q}, v_{i,k,q}, n_{i,k,q})$ the components of $u_{a,k,q}.$ We plug \eqref{ansatz} into the Euler-Maxwell system and find a cascade of equations. The first are the equations at order $O(\e^{-1}):$
$$\begin{aligned}  - i q \o B_{0,q} + i q k \times E_{0,q} & = 0,\nonumber \\
 - i q \o E_{0,q} - i q k  \times B_{0,q} -  {v}_{e0,q} &  = 0,\nonumber \\
 - i q \o v_{e 0,q} + i q k \theta_e n_{e0q} +  {E}_{0,q} & = 0,\nonumber \\
- i q \o {n}_{e 0,q} + i q \theta_e k \cdot {v}_{e 0,q} &= 0, \nonumber \\
- i q \o {v}_{i 0,q}  &= 0,\nonumber \\
- i q \o {n}_{i 0, q} &= 0.\nonumber
 \end{aligned}$$

Given any $k \neq 0,$ and $(k, \o) = (k, \l_1(k))$ on the variety, we have ${\mathcal H}_0 = \{ -1 ,0 ,1\},$ by symmetry of the characteristic variety (see Figure \ref{fig1} on page \pageref{fig1}), and since $\theta_e < 1$ and $|\nabla \l_j| < 1,$ uniformly in $(\xi,\eta),$ for all $j \in \{1,\dots,5\}.$

It will be convenient to express the components of the leading WKB profile in terms of the leading electronic velocity $v_{e0},$ which is orthogonal to $(1,0,0)$, and which we choose to be proportional to $(0,1,0).$ We have, for $p \in \{ - 1, 1 \}:$
\be \label{for:interaction:coefficients}
\begin{aligned} v_{e0,p} & = g_{(p)}(t,x,y) \left(\begin{array}{c} 0 \\ 1 \\ 0 \end{array}\right), \\ E_{0,p} & =  i p \o v_{e0,p}, \\ B_{0,p} & = g_{(p)}(x,y) \left(\begin{array}{c} i p k \\ 0 \\ 0 \end{array}\right) \times \left( \begin{array}{c} 0 \\ 1 \\ 0 \end{array}\right).
\end{aligned}
\ee
Besides, the leading mean modes are all equal to zero: 
$$ B_{0,0} = 0, \quad E_{0,0} = 0, \quad v_{e0,0} = 0,$$
and
$$ n_{e0} = 0, \quad v_{i0} = 0, \quad n_{i0}  = 0.
$$ 
 The above defines vectors $\vec e_p,$ with $p \in \{ - 1 , 1 \},$ such that $\tilde u_{0p} = g_{(p)} \vec e_p$ \eqref{def:gp}. We denote 
 \be \label{def:vv}
 \vec v := \left(\begin{array}{c} 0 \\ 1 \\ 0 \end{array}\right), \quad \vec v{\,'} := \left(\begin{array}{c} 0 \\ 0 \\ 1 \end{array}\right).
 \ee
 Then, the vector $\vec e_p$ takes the form
 \be \label{explicit:ep}
 \vec e_p = \Big( i k {\vec v}{\,'},\,\,  i p \o \vec v, \,\, \vec v, \,\, 0, \,\,  0, \,\, 0 \Big).
 \ee

The next equations are the equations at order $O(\e^{-1/2}).$ These involve the large semilinear source terms. After lengthy computations, detailed in \cite{em1}, and which involve equations up to order $O(\e^{1/2}),$ we arrive at the Zakharov system (Z): 
$$ \mbox{(Z)} \left\{\begin{aligned} \big(i(\d_t + \frac{k}{\o} \d_x) + \frac{1}{2 p \o} \Delta_Y\big) E_{0,p} - \frac{1}{2 p \o \theta_e^2} E_{0,p} & = \frac{1}{2 p \o} n_{1,0} E_{0,p}, \\ (\d_t^2 - (\a_{ie} + 1)^2 \Delta_Y) n_{1,0} & = - \frac{2}{\o^2} \Delta_Y |E_{0,p}|^2,\end{aligned}\right.$$
where $n_{1,0} = n_{e1,0} = n_{i1,0}$ is the (electronic or ionic) mean mode of the fluctuation of density in the first corrector $u_1$ of the WKB ansatz.

 The equations for the corrector terms $(E_{2,p}, n_{2,0})$ are the {\it linearized} (Z) equations at $(E_{1,p}, n_{1,0})$ applied to $(E_{2,p}, n_{2,0}),$ with extra source terms depending on $(E_{1,p}, n_{1,0}).$ The same goes for higher-order terms: we find the (Z) system linearized at $(E_{1,p}, n_{1,0})$ at every step with extra source terms depending on the lower-order terms of the ansatz. In particular, an existence time for (Z) is an existence time for the whole WKB approximate solution. 

\section{Spectral decomposition} \label{sec:spectral:dec}

We compute here the eigenvalues and eigenprojectors of the symbol of the linear hyperbolic operator $A(i \xi , i \eta)$ introduced in Section \ref{sec:hypop}. Denoting $\zeta = (\xi,\eta),$ we find that eigenmodes $\l$ and associated eigenvectors $(B, E, v_e, n_e, v_i, n_i)$ of $A(i \zeta)$ satisfy
 $$ \left\{\begin{aligned}
  i \zeta \times E = \l B, & \qquad - i \zeta \times B - v_e + \frac{\sqrt \e}{\theta_e} v_i  = \l E, \\ E + i \theta_e \zeta n_e = \l v_e, & \qquad i \theta_e \zeta \cdot v_e  = \l n_e, \\ -\frac{\sqrt \e}{\theta_e} E + i \a^2_{ie} \sqrt \e \zeta n_i = \l v_i, & \qquad i \sqrt \e \zeta \cdot v_i = \l n_i.  \end{aligned}\right.$$

\subsection{Away from the kernel} \label{sec:notkernel} If $\l \neq 0,$ the above is equivalent to
$$ \left\{\begin{aligned}
 B = \frac{i}{\l} \zeta \times E,  \quad n_e = \frac{i \theta_e}{\l} \zeta \cdot v_e, \quad n_i =  \frac{i \sqrt \e}{\l} \zeta \cdot v_i, \\  \frac{1}{\l} \zeta \times (\zeta \times E) - v_e + \frac{\sqrt \e}{\theta_e} v_i = \l E, \\ E - \frac{\theta^2_e}{\l} (\zeta \cdot v_e) \zeta = \l v_e, \\ -\frac{\sqrt \e}{\theta_e} E - \frac{\a^2_{ie} \e}{\l} (\zeta \cdot v_i) \zeta = \l v_i.
 \end{aligned}\right.$$

\subsubsection{Transverse modes} \label{sec:transverse:modes}
If $v_e$ and $v_i$ are orthogonal to $\zeta,$ then the above system becomes
\be \label{spectral:orthogonal} \left\{\begin{aligned}
 B = i \zeta \times v_e, \quad v_e = \l^{-1} E, \quad n_e = 0, \quad v_i = - \sqrt\e \theta_e^{-1} \l ^{-1} E, \quad n_i = 0, \\
 \Big(1 + |\zeta|^2 + \e \theta_e^{-2} + \l^2 \Big) E = 0.
 \end{aligned}\right.\ee
We find two transverse modes $i \l_1^\e = - i \l^\e_5$ with 
$$ %
 \l^\e_1 = -\l^\e_5 = \big(1 + \e \theta_e^{-2} + |\zeta|^2 \big)^{1/2}.
$$ %
Clearly
\be \label{orthogonal:modes}
 \l^\e_1 = -\l^\e_5 = \big(1 + |\zeta|^2 \big)^{1/2} + \frac{O(\e)}{\big(1 + |\zeta|^2 \big)^{1/2} } = \l_1(\zeta)  + \frac{O(\e)}{\big(1 + |\zeta|^2 \big)^{1/2} },
\ee
where the remainder $O(\e)$ is real and uniform in $\zeta\in \R^{3}$.

The first equations in \eqref{spectral:orthogonal} are {\it polarization} conditions which specify the eigenspace. The last equation in \eqref{spectral:orthogonal} is an {\it eikonal} equation which describes the eigenmodes.

 These transverse modes are ``electromagnetic" modes, since they correspond to $n_e = 0,$ $n_i = 0.$ Each of these modes has multiplicity two, since $\zeta^\perp$ defines a plane in $\R^3.$
 
 We denote $e_\perp$ the transverse eigenvectors, so that
\be \label{e1:e5}
 e_1(\zeta) = e_\perp(\l_1(\zeta), \zeta), \quad e_5(\zeta) = e_\perp(\l_5(\zeta), \zeta).
\ee 

\subsubsection{Longitudinal modes} \label{sec:parallel:modes} If $v_e$ and $v_i$ are parallel to $\zeta$ then the system at the start of Section \ref{sec:notkernel} becomes
\be \label{spectral:parallel}
 \left\{\begin{aligned}
 B =0,  \quad  v_e = (\l + \frac{\theta_e^2}{\l} |\zeta|^2)^{-1} E, \quad n_e = \frac{i \theta_e}{\l} \zeta \cdot v_e, \\ v_i = - \frac{\sqrt \e}{\theta_e} \cdot \frac{\l}{\l^2 +\a^2_{ie} \e |\zeta|^2} E,  \quad n_i = \frac{i \sqrt \e}{\l}  \zeta \cdot v_i, \\ \Big( 1 + \l^2 + \theta_e^2 |\zeta|^2 + \frac{\e}{\theta_e^2} \cdot \frac{\l^2 + \theta_e^2 |\zeta|^2}{\l^2 + \e \a^2_{ie} |\zeta|^2}  \,   \Big) E = 0. &
 \end{aligned}\right.\ee
 The longitudinal eikonal equation is
 \ba
 \l^{4} + \l^{2} \Big(1+\th_{e}^{2} |\zeta|^{2} + \e \a_{ie}^{2} |\zeta|^{2} + \frac{\e}{\th_{e}^{2}} \Big) + \e |\zeta|^{2} (\a_{ie}^{2} + \a_{ie}^{2} \th_{e}^{2} |\zeta|^{2} + 1) = 0. 
 \ea
We find two solutions $ i\l_2^\e$ and $i \l_4^\e = - i \l_2^\e$ such that 
\be \label{parallel:modes:24}
\l_{2}^\e = \left( 1 + \theta_e^2  |\zeta|^2  + \frac{O(\e)}{1 + \theta_e^2  |\zeta|^2}\right)^{1/2}  = \big( 1 + \theta_e^2  |\zeta|^2 \big)^{1/2} +  \frac{O(\e)}{ \big( 1 + \theta_e^2  |\zeta|^2 \big)^{3/2} },
\ee
where the remainder $O(\e)$ is real and {\it uniform in $\zeta\in \R^{3}.$} The other two solutions, which we denote $\pm i \l^\e_{3}$ are such that, for $\e$ small enough with respect to $\a_{ie}$ and $\theta_e,$ 
\be \label{parallel:modes:3}
 \l^\e_{3} = i \left( \sqrt \e \Big(  \a_{ie}^{2} |\zeta|^{2} + \frac{1}{\th_{e}^{2}} -  \frac{1}{\th_{e}^{2}  (1+\th_{e}^{2} | \zeta|^{2})}\Big)^{\frac 12} - \frac{O(\e)}{(1 + |\zeta|^2)^{3}} \right).
 \ee
with an $O(\e)$ remainder which is real (in fact, non-negative) and uniform in $\zeta$. The low-frequency behavior of these ``acoustic" modes $\pm \l_3^\e$ is given by 
 \be \label{parallel:3:small:xi}
 (\l^\e_{3})^{2} = - \e (\a_{ie}^{2}  + \frac{1}{1+\th_{e}^{2} |\zeta|^{2}}) |\zeta|^{2}  -  \frac{O(\e^{2})}{(1+\th_{e}^{2} |\zeta|^{2})}, \qquad |\zeta| \leq \frac{1}{\a_{ie} \th_{e}},
 \ee
 where the remainder $O(\e^{2})$ is positive and uniform in $|\zeta| \leq \frac{1}{\a_{ie} \th_{e}}.$
 The four longitudinal eigenmodes are simple hence analytical with analytical eigenvectors given by the first two lines of \eqref{spectral:parallel}. Those modes are ``electronic" modes since their $n_e$ component is non-zero (and their $n_i$ component is equal to zero) and their $B$ component is equal to zero. They correspond to the response of the plasma to the incoming (for us, initial) electromagnetic wave.

We denote $e_\parallel$ the longitudinal eigenvectors associated with the slow Klein-Gordon eigenvalues, so that
\be \label{e2:e4}
 e_2(\zeta) = e_\parallel(\l_2(\zeta), \zeta), \qquad e_4(\zeta) = e_\parallel(\l_4(\zeta), \zeta).
\ee
We denote $e_{3\parallel}$ the longitudinal acoustic eigenvectors.

\subsection{The kernel} \label{sec:kernel}

Since components $v_e$ and $v_i$ of an eigenmode associated with $\l \not\equiv 0$ are either both parallel to $\zeta$ or both orthogonal to $\zeta,$ it remains only to describe the kernel of $A(i \zeta).$ The kernel has dimension 6. A family of independent eigenvectors spanning the kernel is given by (for $\zeta \neq 0$): 
$$ \begin{aligned} e_{3,1}:  & \quad B \parallel \zeta, \quad \mbox{all other components equal to 0;} \\
 e_{3,j}, 2 \leq j \leq 5: & \quad \mbox{$\dsp{ - i \zeta \times B - v_e + \frac{\sqrt \e}{\theta_e} v_i = 0,}$ \,\, $E = 0,$ \,\, $n_e = n_i = 0,$ \quad any $v_e, v_i \perp \zeta;$}\\
 e_{3,6}: & \quad \mbox{$n_e = \a^2_{ie} n_i, \quad E=-i\theta_e\zeta n_e,$ \quad and all other components equal to 0.}
 \end{aligned}$$

 \section{Interaction coefficients at space-time resonances} \label{sec:interaction:coefficients}

The goal here is to show that the rate $\g$ from \eqref{def:gamma} is positive, and identify the resonant pair at which the maximum in \eqref{def:gamma} is attained.

\subsection{Reduction: it suffices to compute the growth rate at space-time resonances}

 The growth rate $\g$ \eqref{def:gamma} involves a maximum over $U^{large}_{\delta'},$ not over the set ${\mathcal S}_{jj'}$ of $(j,j')$-space-time resonances. It suffices however to compute the rate at space-time resonances, by virtue of
 \be \label{rate:space-time}
 \lim_{\delta' \to 0} d(U^{large}_{\delta'}, {\mathcal S}_{jj'}) = 0,
 \ee
where $d$ is the distance function. %

By \eqref{rate:space-time} and continuity of $\mbox{tr}(B_{1jj} B_{-1j'j}),$ we find that the rate function computed by taking the maximum over ${\mathcal S}_{jj'}$ converges to $\g_{jj'}.$ Since in our analysis $\delta$ can be chosen to be arbitrarily small, and $\delta' \leq \delta$ (see Lemma \ref{lem:small-trace}), this implies that we may replace $U^{large}_{\delta'}$ with ${\mathcal S}_{jj'}$ in the computation of $\g$ \eqref{def:gamma}.

{\it Verification of \eqref{rate:space-time}.} Otherwise, for some $\e' > 0,$ for all $n,$ we could find $\zeta_n \in U^{large}_{1/n}$ such that $d(\zeta_n, {\mathcal S}_{jj'}) > \e'.$ Since $\zeta_n \in U^{large}_{1/n},$ we have for $n$ large enough $|\Phi_{jj'}(\zeta_n)| \leq 1/n.$ Since $\Phi_{jj'}$ is proper, this implies that $\zeta_n$ is bounded, and, passing to the limit in a subsequence, we find convergence to an element of ${\mathcal S}_{jj'},$ by continuity of $\Phi_{jj'}$ and $\nu_{jj'},$ which contradicts the lower bound $d(\zeta_n, {\mathcal S}) > \e'.$

\subsection{The leading source term at space-time resonances} \label{B:at:space-time} The goal is to compute $B_{pjj'},$ then the trace $\mbox{tr} \, B_{1jj'} B_{-1j'j}.$ The interaction coefficients $B_{pjj'}$ are defined in \eqref{def:interaction:coefficient}.

 First, we go back to the definition of the ``source" term $B\footnote{We are aware that our notation is potentially confusing, with the letter $B$ being used both for the magnetic field and the linearized bilinear source term. But we believe that in practive the context makes it clear which $B$ we are dealing with.},$ defined in Section \ref{sec:Bua}, in terms of the current density, Lorentz force, and convective terms. Here we may overlook the ionic components of the vectors altogether, since the ionic terms do not contribute to $B.$ So we are looking at vectors $u = (B, E, v, n) \in \C^{10}.$ 

 For the contribution of ${\mathcal B}$ to $B,$ we use the definition of ${\mathcal B}$ in \eqref{def:unB} and the explicit expression of $\vec e_p$ in \eqref{explicit:ep}, in terms of the vectors $\vec v$ and $\vec v^{\, '}$ defined in \eqref{def:vv}:
$$ {\mathcal B}(u, \vec e_p) + {\mathcal B}(\vec e_p, u)  = \left(\begin{array}{c} 0 \\ n \vec v \\ - \theta_e \big( \vec v \times B +  i  p k  v \times \vec v^{\,'} \big) \\ 0 \end{array} \right), \qquad u = (B, E, v, n) \in \C^{10}.$$

 Next we compute the contribution of the electronic convective term ${\mathcal J}_e$ (defined in Section \ref{sec:conv}) to $B.$ We have, given $u = (B, E, v, n) \in \C^{10}:$ 
$$ {\mathcal J}_e(u, i p k) \vec e_p + (1 - \chi_{high,+p}) {\mathcal J}_e(\vec e_p,i \zeta) u  = \left(\begin{array}{c} 0 \\ 0 \\ \theta_e \big( (v \cdot i pk ) \vec v + (1 - \chi_{high,+p}(\zeta)) (\vec v \cdot i \zeta) v \big) \\ 0 \end{array}\right).$$ 
Space-time resonances occur only for frequencies $\zeta \in \R \times \{0\}_{\R^2}$ (see Proposition \ref{prop:separation}(4)). In particular, for those frequencies we have, according to the definition of $\vec v$ in \eqref{def:vv}: 
$$ \vec v \cdot \zeta = 0, \qquad \zeta \in \R \times \{0\}_{\R^2}.$$
Thus the leading source term $B$ defined in Section \ref{sec:Bua} is 
$$ B  := \sum_{p \in \{-1,1\}} e^{i p \theta} g_p(x,y) B_p,$$
with
\be \label{Bpu} B_p u = \left(\begin{array}{c} 0 \\ n \vec v \\ - \theta_e \big(  \vec v \times B +  ip  k v \times \vec v^{\,'} + (v \cdot i pk ) \vec v \big) \\ 0 \end{array}\right), \quad u = (B, E, v, n) \in \C^{10},\ee 
where $\vec v$ and $\vec v^{\,'}$ are defined in \eqref{def:vv}.

\subsection{Action of the leading source term over longitudinal and transverse eigenvectors}  \label{sec:B-perp-parallel}
Given $\zeta = (\xi,0,0)$ (recall, we focus on space-time resonances, for which $\eta = 0$), the plane $\zeta^\perp$ is generated by $\vec v$ and $\vec v{\,'}$ \eqref{def:vv}. We denote $e_\perp$ and $e'_{\perp}$ the transverse eigenvectors of $A(i \zeta)_{|\e = 0 },$ described in Section \ref{sec:transverse:modes}. The associated eigenvalues are $\l_\perp(\zeta) = \pm i (1 + |\zeta|^2)^{1/2},$ up to $O(\sqrt \e),$ uniformly in frequency (so that $\l_\perp = i \l_1 + O(\sqrt \e)$ or $\l_\perp = - i \l_1 + O(\sqrt \e) = i \l_5 + O(\sqrt \e)$). According to Section \ref{sec:transverse:modes}, we have
\be \label{eperp} e_\perp = \Big( i \zeta \times \vec v , \,\, \l_\perp(\zeta) \vec v, \,\, \vec v, 0 \Big), \qquad e'_\perp = \Big( i \zeta \times \vec v^{\,'} , \,\, \l_\perp(\zeta) \vec v^{\,'}, \,\, \vec v^{\,'}, 0 \Big),\ee
with norms
$$ |e_\perp| = |e'_\perp| = \sqrt 2 |\l_\perp(\zeta)|.$$
We have $\vec v \cdot i p k = \vec v^{\,'} \cdot i p k = 0 ,$ since $k  = k(1,0,0).$ Thus the cancellation
$$ B_p e'_\perp = 0.$$
Besides, in $B_p e_\perp$ the third term in the electronic velocity cancels, and 
 \be \label{for:0711:1} B_p e_\perp = - i \theta_e \left(\begin{array}{c} 0 \\ 0 \\  (p k + \zeta) \vec v \times \vec v^{\,'} \\ 0 \end{array}\right).\ee
In particular,
\be \label{vBperp}
\mbox{the electronic velocity in $B_p e_\perp$ (third coordinate above) is parallel to $ \vec v^{\,''} := (1,0,0).$}
\ee
Denote $e_\parallel = e_\parallel(\zeta)$ a longitudinal mode of $A(i \zeta),$ as described in Section \ref{sec:parallel:modes}, with associated eigenvalues $\l_\parallel(\zeta) = \pm i (1 +
 \theta_e^2 |\zeta|^2)^{1/2},$ up to $O(\sqrt \e),$ uniformly in frequency. According to Section \ref{sec:parallel:modes}, we have, at $\zeta = (\xi,0):$ 
 \be \label{eparallel} e_\parallel = \Big(0, \,\, \big(\l_\parallel + \frac{\theta_e^2 |\xi|^2}{\l_\parallel}\big)  \vec v^{\,''}, \,\,  \vec v^{\,''}, \,\, \frac{i \theta_e \xi}{\l_\parallel}\Big), \qquad \vec v^{\,''} = \left(\begin{array}{c} 1 \\ 0 \\ 0 \end{array}\right),\ee
 with Euclidian norm
 $|e_\parallel|  =
 \sqrt 2.$
We compute, at $\zeta = (\xi,0):$ 
\be \label{for:0711:2} B_p e_\parallel = \left(\begin{array}{c} 0 \\ \frac{i \theta_e \xi}{\l_\parallel(\xi)} \vec v \\ 0  \\ 0 \end{array}\right).\ee

\subsection{Interaction coefficients: the $(1,2)$ and $(4,5)$ resonances} \label{sec:(1,2):coeffs}

With notation from the previous Section, the orthogonal eigenprojectors $\Pi_1$ and $\Pi_2$ of $A_{|\e = 0 }$ are 
$$ \Pi_1 u = |e_\perp|^{-2} \langle u, e_\perp \rangle e_\perp + |e'_\perp|^{-2} \langle u, e'_\perp \rangle e'_\perp, \qquad \Pi_2 u = |e_\parallel|^{-2} \langle u, e_\parallel \rangle e_\parallel,
$$ 
where $\langle \cdot, \cdot \rangle$ denotes the Hermitian scalar product in $\C^{10}$ (by convention, linear to the left and anti-linear to the right). Taking into account the description of Section \ref{sec:B-perp-parallel}, this gives
$$ \begin{aligned} B_{-121,+1}  u & = \frac{\langle B_{-1} e_{\perp,+1}, e_\parallel\rangle}{|e_\parallel|^2 |e_{\perp,+1}|^2} \langle u, e_{\perp,+1}\rangle e_\parallel, \end{aligned}$$ 
and
\be \label{trace:10:09:2020}
 \mbox{tr}\, B_{112} B_{-121,+1} = \frac{\langle B_1 e_\parallel, e_{\perp,+1} \rangle \langle B_{-1} e_{\perp,+1}, e_\parallel \rangle}{|e_\parallel|^2 |e_{\perp,+1}|^2},
\ee 
where $e_\perp$ is the transverse eigenvector \eqref{eperp} with $\l_\perp = i \l_1,$ and $e_\parallel$ is the longitudinal eigenvector \eqref{eparallel} with $\l_\parallel = i \l_2.$ %
According to Section \ref{sec:B-perp-parallel}, we have, at $\zeta = (\xi,0) \in \R^3:$ 
$$ \langle B_1 e_\parallel, e_{\perp,+1} \rangle = \frac{- i \theta_e \xi \l_1(\xi + k)}{\l_2(\xi)}, \quad \mbox{and} \quad  \langle B_{-1} e_{\perp,+1}, e_\parallel \rangle = - i \theta_e \xi.$$
so that
\be \label{trace:10:09:2020:2}
 \mbox{tr}\, B_{112} B_{-121,+1} = \frac{ - \theta_e^2 \xi^2  \l_1(\xi + k)}{4 \l_2(\xi)^3 \l_1(\xi + k)^2} = \frac{ - \theta_e^2 \xi^2}{4 \l_1(\xi + k) \l_2(\xi)^3}.
 \ee
The trace is negative, implying stability (no growth for the symbolic flow; see the proof of Corollary \ref{cor:the:upper:bound}). We focus on a space-time resonance $\zeta = (\xi,0) \in {\mathcal R}_{12},$ so that
\be \label{12} \l_1(\xi + k) = (1 + |\xi + k|^2)^{1/2} = \o + \l_2(\xi) = \o + (1 + \theta_e^2 |\xi|^2)^{1/2}.\ee
In the small $\theta_e$ limit this gives
$$
 |\xi + k| = (\o^2 + 2 \o)^{1/2} + O(\theta_e^2),
$$
with solutions
\be \label{12resonances}
 \xi_{1,2}^+ = -k + (\o^2 + 2 \o)^{1/2} + O(\theta_e^2), \qquad \xi_{1,2}^- = - k - (\o^2 + 2 \o)^{1/2} + O(\theta_e^2),
\ee 
as we saw in the proof of Proposition \ref{prop:space-time:sep}. Thus
\be \label{for:rate:12}
 \mbox{tr}\, B_{112} B_{-121,+1}(\xi_{1,2}^+) = \frac{ - \theta_e^2 (- k \pm (\o^2 +2 \o )^{1/2})^2}{4 (\o + 1)} + O(\theta_e^4).
\ee
Similarly, the computations for the $(4,5)$ resonance lead to 
\be \label{trace:11:09:2020}
 \mbox{tr}\, B_{145} B_{-154,+1} =  \frac{ - \theta_e^2 (\xi + k)^2}{4 \l_4(\xi + k)^3 \l_5(\xi)}.
 \ee
 In the small $\theta_e$ limit, we find the $(4,5)$ space-time resonances to be $(\xi^{\pm}_{4,5},0)$ with
$$ \xi^{\pm}_{4,5} = \pm (\o^2 + 2 \o)^{1/2} + O(\theta_e^2).$$
The trace in \eqref{trace:11:09:2020} is negative: the $(4,5)$ resonance does not generate any instability. We observe that
$$ \mbox{tr}\, B_{145} B_{-154,+1}(\xi_{4,5}^+) = \mbox{tr}\, B_{112} B_{-121,+1}(\xi_{1,2}^-) + O(\theta_e)^2,$$
and
$$\mbox{tr}\, B_{145} B_{-154,+1}(\xi_{4,5}^-) = \mbox{tr}\, B_{112} B_{-121,+1}(\xi_{1,2}^+) + O(\theta_e^2).$$

\subsection{Interaction coefficients: the $(1,4)$ and $(2,5)$ space-time resonances} \label{sec:insta:14:25}

 For the $(1,4)$ resonance, we have, similarly to \eqref{trace:10:09:2020}: 
\be \label{rate:14}
 \mbox{tr}\, B_{114} B_{-141,+1} = \frac{\langle B_1 e_\parallel, e_{\perp,+1} \rangle \langle B_{-1} e_{\perp,+1}, e_\parallel \rangle}{|e_\parallel|^2 |e_{\perp,+1}|^2}.
\ee 
The difference with \eqref{trace:10:09:2020} is that in \eqref{rate:14}, the vector $e_\perp$ is the transverse eigenvector \eqref{eperp} with $\l_\perp = i \l_1,$ and $e_\parallel$ is the longitudinal eigenvector \eqref{eparallel} with $\l_\parallel = i \l_4 = - i \l_2.$ According to Section \ref{sec:B-perp-parallel}, we have, at $\zeta = (\xi,0) \in \R^3:$ 
$$ \langle B_1 e_\parallel, e_{\perp,+1} \rangle = \frac{i \theta_e \xi \l_1(\xi + k)}{\l_2(\xi)}, \quad \mbox{and} \quad  \langle B_{-1} e_{\perp,+1}, e_\parallel \rangle = - i \theta_e \xi.$$
Note the sign difference in $\langle B_1 e_\parallel, e_{\perp,+1} \rangle$ compared to the computation for the $(1,2)$ resonance in Section \ref{sec:(1,2):coeffs} just above. 
This gives 
\be \label{trace:27:09:2020}
 \mbox{tr}\, B_{114} B_{-141,+1} = \frac{\theta_e^2 \xi^2}{4 \l_2(\xi)^3 \l_1(\xi + k)} > 0.
 \ee
As seen in the proof of Proposition \ref{prop:space-time:sep}, the $(1,4)$ space-time resonances occur at $(\xi_{1,4}^\pm,0)$ with 
\be \label{xi14pm} \xi_{1,4}^\pm = - k \pm (\o^2 -2 \o)^{1/2} + O (\theta_e^2).\ee 
Thus 
\be \label{rate:14:2}
 \mbox{tr}\, B_{114} B_{-141,+1}(\xi_{1,4}^\pm) = \frac{\theta_e^2 (-k \pm (\o^2 - 2 \o)^{1/2})^2}{4 (\o - 1)} + O(\theta_e^4).
\ee 
For the $(2,5)$ resonance, we find symmetrically
\be \label{trace:27:09:2020:2}
\mbox{tr}\, B_{125} B_{-152,+1} = \frac{\theta_e^2 (\xi + k)^2}{4 \l_2(\xi + k)^3 \l_1(\xi)} > 0.
 \ee
 The $(2,5)$ space-time resonances occur at $\xi_{2,5}^\pm$ such that
 $$ \xi_{2,5}^\pm = \pm (\o^2 - 2 \o)^{1/2} + O(\theta_e^2).$$
There is a symmetry between $(1,4)$ and $(2,5),$ just like between $(1,2)$ and $(4,5):$ 
\be \label{1425} \mbox{tr}\, B_{125} B_{-152,+1}(\xi_{2,5}^+) = \mbox{tr}\, B_{114} B_{-141,+1}(\xi_{1,4}^-) + O(\theta_e)^2,\ee 
and
$$\mbox{tr}\, B_{125} B_{-152,+1}(\xi_{2,5}^-) = \mbox{tr}\, B_{114} B_{-141,+1}(\xi_{1,4}^+) + O(\theta_e^2).$$

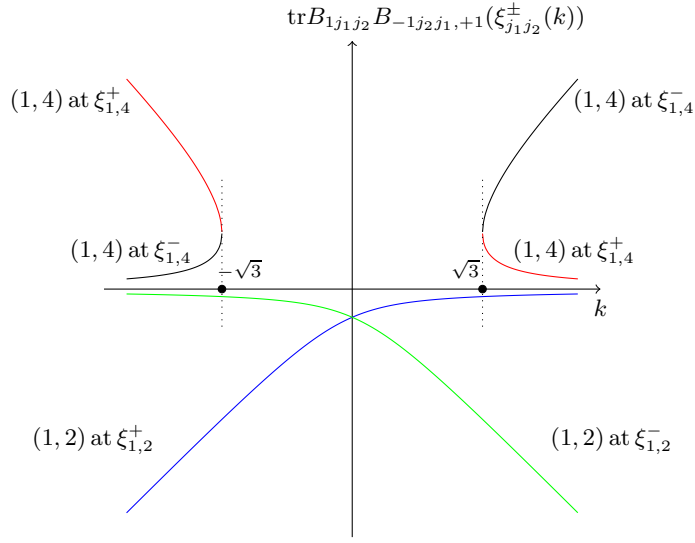
\begin{figure} 
\begin{tikzpicture}

\draw[->] (-3.3,0) -- (3.3,0) ; 
\draw (3.3,0) node[anchor=north] {$\footnotesize{\mbox{$k$}}$};  
\draw[->] (0,-3.3) -- (0,3.3) ;
\draw (-1,3.6) node[anchor=west] {\footnotesize $\mbox{tr} B_{1j_1j_2} B_{-1j_2j_1,+1}(\xi^\pm_{j_1j_2}(k))$} ; 
\draw[domain=1.73206:3, samples=100,red]
plot( \x, { ( -\x + ( 1 + \x^2 - 2*(1 + \x^2)^(1/2) )^(1/2) )^2/(4*((1+\x^2)^(1/2) - 1) ) } ) ; 
\draw (2.8,2.5) node[anchor=west] {\footnotesize $(1,4) \,\mbox{at} \,\xi_{1,4}^-$} ;  
\draw (-2.8,2.5) node[anchor=east] {\footnotesize $(1,4) \,\mbox{at} \,\xi_{1,4}^+$} ;

\draw[domain=1.73206:3, samples=100] 
plot( \x, { ( -\x - ( 1 + \x^2 - 2*(1 + \x^2)^(1/2) )^(1/2) )^2/(4*((1+\x^2)^(1/2) - 1)) } ) ;
\draw (-2,.5) node[anchor=east] {\footnotesize $(1,4) \,\mbox{at} \,\xi_{1,4}^-$} ; 
\draw (2,.5) node[anchor=west] {\footnotesize $(1,4) \,\mbox{at} \,\xi_{1,4}^+$} ;  
\draw[domain=-1.732051:-3, samples=100,red]
plot( \x, { (-\x + ( 1 + (-\x)^2 - 2*(1 + (-\x)^2)^(1/2) )^(1/2) )^2/(4*((1+(-\x)^2)^(1/2) - 1) ) } ) ; 
\draw[domain=-1.732051:-3, samples=100] 
plot( \x, { (-\x - ( 1 + (-\x)^2 - 2*(1 + (-\x)^2)^(1/2) )^(1/2) )^2/(4*((1+(-\x)^2)^(1/2) - 1)} ) ;
\draw[dotted] (1.732051,-.5) -- (1.732051,1.5) ; 
\draw[dotted] (-1.732051, -.5) -- (-1.732051,1.5) ; 
\filldraw (1.732051,0) circle (0.05cm) ; 
\filldraw (-1.732051,0) circle (0.05cm) ; 
\draw (1.5,0) node[anchor=south] {\tiny $\sqrt 3$} ; 
\draw (-1.5,0) node[anchor=south] {\tiny $-\sqrt 3$} ; 
\draw[domain=-3:3, samples=100,blue]
plot( \x, { - (  ( 1 + (\x)^2 + 2*(1 + (\x)^2)^(1/2)  )^(1/2) - \x )^(2)/(4*(1 + (1 + (\x)^2)^(1/2) ) } ) ; 
\draw (-2.5,-2) node[anchor=east] {\footnotesize $(1,2) \, \mbox{at} \, \xi_{1,2}^+$} ; 
\draw (2.5,-2) node[anchor=west] {\footnotesize $(1,2) \, \mbox{at} \, \xi_{1,2}^-$} ; 
\draw[domain=-3:3, samples=100,green]
plot( \x, { - ( ( 1 + (\x)^2 + 2*(1 + (\x)^2)^(1/2))^(1/2) + \x )^(2)/(4*(1 + (1 + (\x)^2)^(1/2)) ) } ) ; 

\end{tikzpicture}

\caption{The graphs of $k \to \mbox{tr}\, B_{1j_1j_2} B_{-1j_2j_1,+1}(\xi^\pm_{j_1,j_2}(k)),$ where $(\xi^\pm_{j_1j_2},0)$ are the $(j_1,j_2)$-space-time resonances, are shown on this figure, for $(j_1,j_2) = (1,4)$ and $(j_1,j_2) = (1,2),$ and, by symmetry, for $(j_1,j_2) = (2,5)$ and $(j_1,j_2) = (4,5).$ We kept only the leading term in $\theta_e,$ and let $\theta_e = 1.$ The trace for $(1,4)$ is given by \eqref{rate:14:2}. The curves in red correspond to the trace evaluated at $\xi_{1,4}^+$ \eqref{xi14pm}. The curves in black correspond to the trace evaluated at $\xi_{1,4}^-.$ We see that the trace is positive for $|k| > \sqrt 3,$ indicating instability. The greatest instability is recorded at $\xi_{1,4}^-$ for $k > \sqrt 3$ and at $\xi_{1,4}^+$ for $k < - \sqrt 3.$ As noted at the end of Section \ref{sec:insta:14:25}, the (leading term of the) trace for $(2,5)$ at $\xi_{2,5}^+$ is equal to the trace for $(1,4)$ at $\xi_{1,4}^-,$ and conversely the trace for $(2,5)$ at $\xi_{2,5}^-$ is equal to the trace for $(1,4)$ at $\xi_{1,4}^+.$ The trace for $(1,2)$ at $\xi_{1,2}^+$ (equivalently, for $(4,5)$ at $\xi_{4,5}^-$) is shown in blue. The trace for $(1,2)$ at $\xi_{1,2}^-$ (or for $(4,5)$ at $\xi_{4,5}^+$) is shown in green. The resonances $(1,2)$ and $(4,5)$ are stable.} \label{figure:trace} 
\end{figure}

\subsection{Interaction coefficients: the $(2,4)$ resonance} \label{sec:(2,4):coeffs}

In view of \eqref{eparallel} and \eqref{for:0711:2}, we see that
$$ \langle B_1 e_4, e_{2,+1} \rangle = 0, \qquad \langle B_{-1} e_{2,+1}, e_4 \rangle = 0,$$
denoting $e_2$ the $e_\parallel$ eigenvector associated with $\l_2,$ and $e_4$ the $e_\parallel$ eigenvector associated with $\l_4.$ That is, the $(2,4)$ resonance is transparent in the sense of Joly, M\'etivier and Rauch \cite{JMR-TMB}. According to the $S$ bounds of Section \ref{sec:flow}, this implies in particular that it does not generate any instability.

\subsection{Interaction coefficients: the $(1,3)$ and $(3,5)$ resonances} \label{sec:(1,3):coeffs}

 The $(1,3)$ space-time resonant frequencies $(\xi_{1,3},0)$ satisfy
\be \label{13:eq} (1 + |\xi + k|^2)^{1/2} = \o.\ee
Recall indeed that we approximate the acoustic modes $\pm \l_{3}^\e$ described in \eqref{parallel:modes:3} by $\l_3 \equiv 0.$ The error is significant only for large $|\zeta|.$ Since the resonant set is bounded, this approximation is good enough for our analysis.

The solutions to \eqref{13:eq} are 
$$\xi_{1,3} = - 2 k \quad \mbox{and} \quad \xi_{1,3} = 0.$$
We observe the cancellation
\be \label{cancellation:13} \langle B_p e_\perp, e_{3,j} \rangle = 0, \qquad 1 \leq j \leq 6, \quad p \in \{-1,1\},\ee 
which is due to \eqref{vBperp} and the description of the vectors in the kernel in Section \ref{sec:kernel}.

It remains to consider the eigenvectors $e_{3\parallel}$ associated with the acoustic longitudinal eigenvalues $\pm i \l_3^\e,$ described in \eqref{parallel:modes:3}. The associated eigenvectors are longitudinal modes, described in Section \ref{sec:parallel:modes}. 

We observe that the $v_e$ component of these eigenvectors is $O(\sqrt\e).$ Indeed, for $\zeta$ far from zero, this is a direct consequence of \eqref{spectral:parallel} and \eqref{parallel:modes:3}. In the small-frequency limit, this is found by taking the limit as $\zeta \to 0$ of the eigenvector described in \eqref{spectral:parallel}. Using \eqref{parallel:3:small:xi}, we find that the eigenvectors associated with $\pm \l^\e_{3}$ do converge in the $\zeta \to 0$ limit, which is not obvious in the first place since these acoustic modes coalesce at $\zeta = 0.$ The limiting eigenvectors are %
 $$ e_{3\parallel}(0,0) = \left( 0, 0, \frac{\pm \sqrt{\e (\a_{ie}^2 + 1)}}{\theta_e^2}, \, \frac{i}{\theta_e} \right).$$
 (The other, ionic components of the full $e_{3\parallel}$ eigenvector at $(\xi,\eta) = (0,0)$ are equal to zero.)  

Thus the $v_e$ component of $e_{3\parallel}$ is $O(\sqrt \e)$ everywhere, and as a consequence of the description of $B_p e_\perp$ in \eqref{for:0711:1}, this implies 
\be \label{int:coeff:13} \big| \langle B_{-1} e_{\perp,+1}, e_{3\parallel} \rangle \big| \lesssim \sqrt \e.\ee

By symmetry, the same holds for the $(3,5)$ space-time resonances $\xi^\pm_{3,5}.$ Those satisfy $0 = \o - (1 + |\xi|^2)^{1/2},$ so that $\xi_{3,5}^\pm = \pm k.$ Indeed, the cancellation \eqref{cancellation:13} holds, since $\l_5$ is a transverse mode just like $\l_1.$ Then we use again the fact that $v_e = O(\sqrt \e)$ for the $e_{3\parallel}$ modes. 

These observations imply that the rough rate of growth \eqref{def:gamma3} associated with the triplet $(1,3,5)$ is $O(\sqrt \e).$ 

\subsection{Interaction coefficients: the $(2,3)$ and $(3,4)$ resonances} \label{sec:(2,3):coeffs}

The $(2,3)$ space-time resonances occur for frequencies $\zeta = (\xi,0)$ with
$$ (1 + \theta_e^2 |\xi + k|^2)^{1/2} = \o ,$$
so that $\xi$ is somewhat large (since $\theta_e$ is small): 
$$ \xi_{2,3}^\pm  = -k \pm k/\theta_e.
$$ 
We note that $\l_2(\xi_{2,3}^\pm + k) = \o.$
From the description of $B_p e_\parallel$ in \eqref{for:0711:2}, we note that in $\langle B_{-1} e_{2,+1}, e_3 \rangle,$ only the $E$ component of $e_3$ intervenes. Actually, only the component along $\vec v$ of the $E$ component in $e_3$ comes into play, where $e_3$ is an eigenvector in the kernel (described in Section \ref{sec:kernel}) or an acoustic eigenvector (described in Section \ref{sec:parallel:modes}). But amongst the eigenvectors in the kernel, only $e_{3,6}$ has a non-zero $E$ component, and that component is orthogonal to $\vec v.$ The acoustic eigenvectors also have $E$ components orthogonal to $\vec v.$ As a consequence,
$$ \langle B_{-1} e_{2,+1}, e_3 \rangle = 0,$$
for any eigenvector $e_3$ which is either in the kernel or acoustic. By symmetry, the same applies to $(3,4).$ By \eqref{def:gamma3}, this implies a zero growth rate for $(2,3).$

\subsection{Comparison of growth rates and Raman backscattered waves} \label{sec:rates:comparison}

From \eqref{cancellation:13}, \eqref{int:coeff:13}, the definition of $\g^+_{135}$ in \eqref{def:gamma3}, and \eqref{rate:space-time}, we deduce that the growth rates  
 $ \gamma_{135}^+(\R^3 \times (\R^3 \setminus U^{away}_{\delta})),$ and $\g_{13}(\R^6)$ and $\g_{35}(\R^6)$ can all be made arbitrarily small if $\e$ and $\delta$ are small enough.

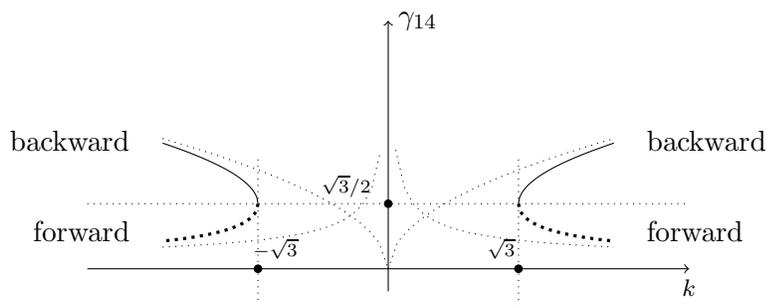
\begin{figure} 
\begin{tikzpicture}

\draw[dotted] (-4,.866) -- (4,.866)  ; 

\draw[->] (-4,0) -- (4,0) ; 
\draw (0,3.3) node[anchor=west] {$\g_{14}$} ; 
\draw (4,0) node[anchor=north] {$\footnotesize{\mbox{$k$}}$};  
\draw[->] (0,-.3) -- (0,3.3) ;

\draw[domain=1.73206:3, samples=100,very thick,dotted]
plot( \x, { sqrt( (1/4)*( -\x + sqrt( 1 + \x^2 - 2*sqrt(1 + \x^2) ) )^2/( sqrt( 1+\x^2) - 1) )) } )  ; 
\draw[domain=1.73206:3, samples=100]
plot( \x, { sqrt( (1/4)*( -\x - sqrt( 1 + \x^2 - 2*sqrt(1 + \x^2) ) )^2/( sqrt( 1+\x^2) - 1) )) } )  ;

\draw[domain=-1.732051:-3, samples=100]
plot( \x, { sqrt( (1/4)*( -\x + sqrt( 1 + (-\x)^2 - 2*sqrt(1 + (-\x)^2) ) )^2/( sqrt( 1+(-\x)^2) - 1) ) ) } ) ;
\draw[domain=-1.732051:-3, samples=100,very thick, dotted]
plot( \x, { sqrt( (1/4)*( -\x - sqrt( 1 + (-\x)^2 - 2*sqrt(1 + (-\x)^2) ) )^2/( sqrt( 1+(-\x)^2) - 1) ) ) } ) ;
\draw (3.3,1.7) node[anchor=west] {backward} ; 

\draw (3.3,.5) node[anchor=west] {forward} ; 

\draw (-3.3,1.7) node[anchor=east] {backward} ; 

\draw (-3.3,.5) node[anchor=east] {forward} ; 

\filldraw (0,.866) circle (0.05cm) ; 

\draw (-.1,1.1) node[anchor=east] {\tiny{$\sqrt 3/2$}} ; 
\draw[dotted] (1.732051,-.5) -- (1.732051,1.5) ; 
\draw[dotted] (-1.732051, -.5) -- (-1.732051,1.5) ; 
\filldraw (1.732051,0) circle (0.05cm) ; 
\filldraw (-1.732051,0) circle (0.05cm) ; 
\draw (1.5,0) node[anchor=south] {\tiny $\sqrt 3$} ; 
\draw (-1.5,0) node[anchor=south] {\tiny $-\sqrt 3$} ; 

\draw[domain=0:3,dotted] plot( \x, {sqrt(\x)}); 
\draw[domain=-3:0,dotted] plot(\x, {sqrt(-\x)});
\draw[domain=0.1:3,dotted] plot (\x, {(1/2)*(\x)^(-1/2)}); 
\draw[domain=-3:-0.1,dotted] plot (\x, {(1/2)*(-\x)^(-1/2)});
\end{tikzpicture}

\caption{The growth rate $\g_{14}$ (or $\g_{25}$) is pictured as a function of the WKB wavenumber $|k| > \sqrt 3.$ This rate is associated with backscattered Raman waves. The smaller rate of growth associated with forward Raman waves is shown in dotted lines. In this figure, $\theta_e = 1,$ and $\max_{\R^3} |a| = 1,$ so that the graphs here are simply the square roots of the graphs of Figure \ref{figure:trace} that are above the $k$ axis. In the large $|k|$ limit, the backscattered Raman rate $\g_{14}$ grows like $|a|_{L^\infty} \sqrt{|k| \theta_e} ,$ and the forward scattered Raman rate decays like $|a|_{L^\infty} \sqrt{\frac{\theta_e}{2 |k|}}.$} 
 \label{fig:growth-rate} 
\end{figure}

From the previous Sections we deduce that
$$ \g_{jj'}(\R^6) = 0, \qquad (j,j') \in \{ (1,2), (4,5), (2,3), (3,4), (2,4) \},$$
and
$$ \g_{135}(\R^6) \lesssim \sqrt \e.$$ 
The $(1,4)$ and $(2,5)$ resonances remain, and the maximal growth rate $\g$ defined in \eqref{def:gamma+} is equal to $\g_{14}$ and $\g_{25}:$ 
 $$ \g = \g_{14}(\R^6) = \g_{25}(\R^6), \quad \mbox{for any $|k| > \sqrt 3.$}$$ 
The graph of this growth rate as a function of $k,$ with $|k| > \sqrt 3,$ is pictured on Figure \ref{fig:growth-rate}, where for the sake of representation we chose a maximum of the norm of the WKB amplitude equal to 1, and $\theta_e = 1.$

 If $k > 0$ and $\o > 0,$ then the WKB wave propagates to the right; if $k < 0$ and $\o > 0,$ the WKB wave propagates to the left. 
 
 Consider for instance the case $k > \sqrt 3.$ Then, we see on Figure \ref{figure:trace} that  the maximal rate is attained at $\xi_{1,4}^-$ \eqref{xi14pm}. By \eqref{1425}, it is also attained at $\xi_{2,5}^+,$ but we may focus on $(1,4).$
 
 We choose the initial perturbation in \eqref{initial:datum} to involve the branches 1 and 4 on the characteristic variety, and $\xi_0 = \xi_{14}^-.$ The corresponding phases in the initial perturbation are $(\xi_{14}^- + k, 0, \l_1(\xi_{14}^- + k,0))$ and $(\xi_{14}^-, 0, \l_4(\xi_{14}^-,0)).$ We have $\xi_{14}^- < 0$ and $\xi_{14}^- + k < 0.$ In particular, the plane wave carried by the phase $(\xi_{1,4}^- + k, 0, \l_1(\xi_{14}^- + k,0)),$ and which is exponentially amplified in time, travels to the left: it is ``backscatterred", that is, travels in a direction opposite to the WKB solution. The electronic plasma wave carried by the phase $(\xi_{14}^-,0, \l_4(\xi_{14}^-,0))$ travels to the right. 

 In the case $k < - \sqrt 3,$ the situation is symmetrical, with the maximal rate attained at $\xi_{14}^+.$ The electromagnetic wave associated with $\l_1$ that is present in the initial perturbation is exponentially amplified and propagates to the right, while the electronic plasma wave associated with $\l_4$ in the initial perturbation (also exponentially amplified) and the WKB solution travel to the left.

\subsection{Forward Raman waves} \label{sec:forward:Raman} 

 Forward Raman waves are associated with $\xi_{14}^+$ for $k > \sqrt 3$ and $\xi_{14}^-$ for $k < - \sqrt 3.$ Indeed, if for instance $k > 0,$ then $\xi_{14}^+ < 0$ and $\xi_{14}^+ + k > 0,$ so that the wave carried by $(\xi_{14}^+ + k, 0,  \l_1(\xi_{14}^+ + k,0))$ travels to the right, just like the WKB solution, and so does the wave carried by $(\xi_{14}^+,0,\l_4(\xi_{14}^+,0)).$ We see on Figure \ref{fig:growth-rate} that forward Raman waves have associated rates of growth which are positive, although smaller than the backward Raman rates. 
 
 If we chose an initial perturbation that will ride those forward Raman waves, we would expect the proof to break down before the optimal observation time, since those amplified waves grow at a rate that is not maximal among all possible rates. This would imply a lower bound on $K',$ in the spirit of Theorem 2.11 in \cite{em4} (``all non-transparent resonances are amplified") describing the instability generated by any incompatible (that is, non-transparent) resonance, even those not associated with the maximal growth rate.

\section{Symbols and operators} \label{app:symb}

We denote $\hat u(\xi,\eta)$ the Fourier transform of $u(x,y).$ The dual Fourier variable of $(x,y) \in \R \times \R^2$ is often denoted $\zeta = (\xi,\eta) \in \R \times \R^2.$ We consider waves with typical frequencies $1/\e$ in $x$ and $1/\sqrt \e$ in the transverse direction $y.$ Accordingly, the relevant quantization of operators is
\be \label{aniso}
(\op_\e(a) u)(x,y) = \int_{\R^3} e^{i (x \xi + y \cdot \eta)} a(x,y,\e \xi,\e^{1/2}\eta) \hat u(\xi,\eta) \, d\xi \, d\eta.
\ee
A symbol $a,$ possibly matrix-valued, is said to belong to the classical class $S^m = S^m_{1,0}$ if it satisfies the bounds
\be \label{sm10}
 \sup_{\begin{smallmatrix} x,y \in \R^3 \\ \xi,\eta \in \R^3 \end{smallmatrix}} \langle \xi,\eta\rangle^{(m - |\b|)/2} |\d_{x,y}^\a \d_{\xi,\eta}^\b a(x,y,\xi,\eta)| < \infty, \qquad \mbox{for all $\a,\b \in \N^3,$}
\ee
with $\langle \xi,\eta\rangle := (1 + |(\xi,\eta)|^2)^{1/2}.$ 
Pseudo-differential operators with symbols in $S^0$ are linear bounded from $L^2$ to itself. The space of operators associated with symbols in $S^m$ also enjoys a form of stability by composition. Those are well-known results; in our case, since we're dealing with singular symbols in $\e,$ the specific symbolic norms that come in those results matter very much. They are made precise in the next three statements.

\begin{prop} \label{prop:H} Given $a \in S^0,$ we have the bound 
$$  \| \op_\e(a) \|_{L^2 \to L^2} \lesssim \sum_{|\a| \leq 4} \sup_{(\xi,\eta)\in\R^3} \| \d_{x,y}^\a a(\cdot,\xi,\eta)\|_{L^1(\R^3_{x,y})},$$
 the implicit constant depending only on dimensions.
\end{prop}

\begin{proof} This is Theorem 18.8.1' in H\"ormander's treatise \cite{H}. 
\end{proof}

Stability by composition is expressed by the equality
 \be \label{compo:e} \begin{aligned}
 \op_\e(a_1) \op_\e(a_2) & = \op_\e(a_1 a_2) + \sum_{1 \leq |\a| \leq q } \e^{|\a|/2} \op_\e\left( \frac{(-i)^{|\a|}}{\a!} \d_\eta^\a a_1 \d_y^\a a_2 \right) \\ & + \sum_{1 \leq |\a| \leq [q/2] } \e^{|\a|} \op_\e\left( \frac{(-i)^{|\a|}}{\a!} \d_\xi^\a a_1 \d_x^\a a_2 \right) + \e^{[q/2]} \op_\e(R_{q/2}(a_1,a_2)),
 \end{aligned} \ee
and the following result about the remainder $R:$

\begin{prop} \label{prop:composition} For all $q,$ all $a_1 \in S^{m_1},$ all $a_2 \in S^{m_2},$ with $m_1, m_2 \in \R$ and $[q/2] - m_1 - m_2 \leq 0,$ the remainder $R_{q/2}(a_1,a_2)$ belongs to $S^{m_1 + m_2 - [q/2]},$ and satisfies the bound, for some $C(d) > 0$ which depends only the dimension $d:$
$$ \begin{aligned} \| \op_\e(R_{q/2}(&a_1,a_2))\|_{L^2 \to L^2} \\ & \lesssim \sup_{\begin{smallmatrix} |\a_1| \leq C(d) \\ |\b_1| \leq q + C(d) \\ (x,y,\xi,\eta) \in \R^6 \end{smallmatrix}} \langle \xi, \eta\rangle^{|\b_1| - m_1}   \big|\d_{x,y}^{\a_1} \d_{\xi,\eta}^{\b_1} a_1 \big| \sup_{\begin{smallmatrix} |\a_2| \leq q + C(d) \\ |\b_2| \leq C(d) \\ (x,y,\xi,\eta) \in \R^6 \end{smallmatrix}} \langle \xi, \eta\rangle^{|\b_2| - m_2}  \big|\d_{x,y}^{\a_2} \d_{\xi,\eta}^{\b_2} a_2 \big|.\end{aligned}$$
\end{prop}

In our context $d=3;$ for clarity we wrote $d$ in the statement of Proposition \ref{prop:composition}.

\begin{proof}[Proof of Proposition {\rm \ref{prop:composition}}] This is a classical result. Details of the remainder bound are found for instance by putting together the results of Theorem 1.1.5, Theorem 1.1.20, Remark 4.1.2 and Remark 4.1.4 in Lerner's book \cite{Le}. The extension to (anisotropic) semi-classical quantization follows easily by introduction of dilations and weights, as described for instance in the appendix of \cite{em4}.
\end{proof}

In \eqref{compo:e} we observe that $(a_1,a_2) \to \op_\e(R_{q/2}(a_1,a_2))$ is bilinear. In particular, for $\a \in \N^3,$
$$ \d_{x,y}^\a \op_\e(R_{q/2}(a_1,a_2)) = \sum_{\a_1 + \a_2 + \a_3 =\a} C_{\a_1,\a_2,\a_3} \op_\e(R_{q/2}(\d_{x,y}^{\a_1} a_1, \d_{x,y}^{\a_2} a_2)) \d_{x,y}^{\a_3} ,$$
for some $C_{\a_1,\a_2,\a_3} \in \N,$ implying, with Proposition \ref{prop:composition}, the Sobolev bound 
\be \label{composition:sobolev} \begin{aligned}
\| & \op_\e(R_{q/2}(a_1,a_2) u \|_{\e,s} \\ & \lesssim \sup_{\begin{smallmatrix} |\a_1| \leq s + C(d) \\ |\b_1| \leq q + C(d) \\ (x,y,\xi,\eta) \in \R^6 \end{smallmatrix}} \langle \xi, \eta\rangle^{|\b_1| - m_1}  \big|\d_{x,y}^{\a_1} \d_{\xi,\eta}^{\b_1} a_1 \big| \cdot \sup_{\begin{smallmatrix} |\a_2| \leq s + q + C(d) \\ |\b_2| \leq C(d) \\ (x,y,\xi,\eta) \in \R^6 \end{smallmatrix}} \langle \xi, \eta\rangle^{|\b_2| - m_2}  \big|\d_{x,y}^{\a_2} \d_{\xi,\eta}^{\b_2} a_2 \big| \cdot \| u \|_{\e,s}, \end{aligned}
\ee
for $s \in \N,$ $u \in H^s,$ $a_j \in S^{m_j}$ and $q$ large enough so that $q/2 - m_1 - m_2 \geq 0.$

The next result gives an elementary bound for the $\fl1 \to \fl1$ norm of pseudo-differential operators:

\begin{prop} \label{prop:op:fl1} Given a symbol $a,$ if $\sup_{(\xi,\eta) \in \R^3} \| \d_{x,y}^\a a(\cdot,\xi,\eta)\|_{L^2(\R^3_{x,y})} < \infty$ for $|\a| \leq 2,$ then $\op_\e(a)$ maps $\fl1$ to itself, and  
$$ \| \op_\e(a) \|_{\fl1 \to \fl1} \lesssim \sup_{\begin{smallmatrix} (\xi,\eta) \in \R^3 \\ 0 \leq |\a| \leq 2 \end{smallmatrix} } \| \d_{x,y}^\a a \|_{L^2(\R^3_{x,y})},$$
the implicit constant depending only on dimensions.
\end{prop}

\begin{proof} Given $f \in \fl1,$ we compute
$$ \begin{aligned} {\mathcal F} \big( \op_\e(a) f \big)(\xi,\eta) & = \int e^{- i (x,y) \cdot (\xi,\eta)} \int e^{i (x,y) \cdot (\xi',\eta')} a(x,y,\e \xi', \sqrt \e \eta') \hat f(\xi',\eta') \, d\xi d\eta \, dx\,dy \\ & = \int \left( \int e^{ -i (x,y) \cdot (\xi - \xi', \eta - \eta')} a(x,y,\e \xi',\sqrt \e \eta') \, dxdy \right) \hat f(\xi',\eta') \, d\xi' d\eta'\\ & = \int \hat a(\xi - \xi', \eta - \eta', \e \xi', \sqrt \e \eta') \hat f(\xi',\eta') d\xi' d\eta', \end{aligned}$$  
where $\hat a$ refers to the Fourier transform of $a$ in its spatial variables. Thus
$$ \begin{aligned} \| \op_\e(a) f \|_{\fl1} & \lesssim \int_{\R^6} |\hat a(\xi - \xi',\eta- \eta', \e \xi', \sqrt \e \eta')| |\hat f(\xi',\eta')| d\xi' d\eta' d\xi d\eta \\ & \lesssim \int_{\R^3} \left( \int_{\R^3} |\hat a(\xi,\eta,\e \xi', \sqrt \e \eta')| d\xi d\eta \right) |\hat f(\xi',\eta')| d\xi' d\eta' \\ & \lesssim \sup_{\zeta' \in \R^3} \|\hat a(\cdot, \zeta')\|_{L^1(\R^3)} \| f \|_{\fl1}.\end{aligned}$$
Now for fixed $\zeta \in \R^3,$ we have
$$ \| \hat a(\cdot,\zeta)\|_{L^1} \lesssim \| a(\cdot, \zeta) \|_{H^2(\R^3)}, \quad \mbox{since $2 > 3/2,$}$$
and the result follows. 
   \end{proof}

We use a pointwise bound at the very end of the proof:
\begin{prop} \label{prop:pointwise} Given $a \in S^0,$ given $u \in H^s(\R^3)$ with $s > d/2,$ we have
\be \label{sup:w}
  \|\op_\e(a) u \|_{L^\infty} \lesssim C(a) \big(1 + |\ln \e| +| \ln \| u \|_{\e,s}|\big),
 \ee
 for some constant $C(a) > 0.$
\end{prop}

\begin{proof} By introduction of dilations and weighted norms (see for instance the appendix of \cite{em4}), we derive from estimate (B.1.1) in Appendix B of \cite{Tay} the bound
$$   |\op_\e(a) u|_{L^\infty} \lesssim \| a \|_{C(d)} |u|_{L^\infty} \Big(1 + |\ln \e| + \Big| \ln\Big(\frac{N_{\e,s}(u)}{|u|_{L^\infty}}\Big)\, \Big| \, \Big),$$
where $N_{\e,s}(\cdot)$ is the weighted Sobolev norm associated with the anisotropic quantization \eqref{aniso}:
$$ N_{\e,s}(u) := \big\| (1 + |\e \xi|^2 + |\sqrt \e \eta|^2)^{s/2} \hat u(t,\xi,\eta) \big\|_{L^2(\R^3)}.$$
In particular,
$$ N_{\e,s}(u) \leq \| u \|_{\e,s},$$
where the (weighted in $x$ only) Sobolev norm $\| \cdot \|_{\e,s}$ is introduced in \eqref{weighted:norm}. 
The Sobolev embedding takes the form
$$ \| u \|_{L^\infty} \lesssim \e^{-1/2} \| u \|_{\e,s}.$$
Thus we have
$$ \Big| \, \ln\Big(\frac{\|u\|_{\e,s}}{\|u\|_{L^\infty}}\Big) \, \Big| \leq |\ln \|u\|_{\e,s}| + |\ln \|u\|_{L^\infty}| \lesssim |\ln \e| + |\ln \| u \|_{\e,s}|,$$
and the result follows.
\end{proof}

\begin{rem}[Action of a pseudo-differential operator on highly-oscillating functions] \label{rem:shift}
The action of a pseudo-differential operator $\op_\e(a)$ on a function with space-time oscillations in $\theta = (k x - \o t)/\e$ (where $k$ is fixed, equal to the fundamental wavenumber in the WKB oscillations, and $\o = (1 + k^2)^{1/2}$) induces a shift in the argument of the symbol:
$$
 \op_\e(a) (e^{i p\theta} f) = e^{ip\theta} \op_\e(a_{+p}) f, \qquad \theta = (k x - \o t)/\e,
$$
with notation
$$
a_{+p}(x,y,\xi,\eta) := a(x,y,\xi + p k, \eta), \qquad p \in \Z.
$$
\end{rem}

Many symbols that we meet in the proof are tensor products: $b(x,y,\xi,\eta) = b_1(x,y) b_2(\xi,\eta),$ where $b_1$ is a term in the WKB approximate solution, and $b_2$ is bounded. For the associated operators, it is straightforward to derive bounds in weighted Sobolev and $\fl1$ norms:

\begin{lem} \label{lem:tensor} Given $q \in \N,$ given a fast oscillation $e^{i q \theta},$ with $\theta = (k x - \o t)/\e,$ given an $\e$-independent Sobolev map $b(x,y) \in W^{s,\infty},$ with $s \in \N,$ given a bounded map $b_2 \in L^\infty(\R^3_{\xi,\eta}),$ we have for all $u \in H^s$ the bound 
$$ \| \op_\e( e^{i q \theta} b_1 b_2) u \|_{\e,s} \lesssim \| b_2 \|_{L^\infty} \cdot \sum_{s_1 + s_2 = s} \| b_1 \|_{W^{s,\infty}} \| u \|_{\e,s_2}.$$
If $\d_\e^{\a'} b_1 \in \fl1$ and $\d_\e^{\a'} u \in \fl1$ for all $\a' \leq \a,$ then 
$$ \| \d_\e^\a \op_\e(e^{i q \theta} b_1 b_2) u \|_{\fl1} \lesssim \| b_2 \|_{L^\infty} \cdot \sum_{\a_1 + \a_2 = \a} \| \d_\e^{\a_1} b_1 \|_{\fl1} \| \d_\e^{\a_2} u \|_{\fl1}  $$
\end{lem}

\begin{proof} Both bounds are elementary. For the first one, we observe
$$ \d_\e^\a (\op_\e(e^{i q \theta} b_1 b_2 u) = \sum_{\a_1 + \a_2 + \a_3} C_{\a_1,\a_2,\a_3} \d_\e^{\a_1} (e^{i q \theta}) \d_\e^{\a_2} b_1 \op_\e(b_2) \d_\e^{\a_3} u,$$
for some $C_{\a_1,\a_2,\a_3} \in \N,$ and then we use
$$\| \d_\e^{\a_1} (e^{i q \theta}) \d_\e^{\a_2} b_1 \op_\e(b_2) \d_\e^{\a_3} u \|_{L^2} \lesssim \| \d_\e^{\a_2} b_1 \|_{L^\infty} \| \op_\e(b_2) \d_\e^{\a_3} u \|_{L^2},$$
and finally the fact that the $L^2 \to L^2$ norm of the Fourier multiplier $\op_\e(b_2)$ is controlled by the sup norm of $b_2.$ 

For the second bound, we use first the fact that the multiplication by $e^{i q \theta}$ corresponds to a translation in Fourier, so that 
$$\| \big( \d_\e^{\a_1} e^{i q \theta} \big) f \|_{\fl1} \lesssim \| f \|_{\fl1}.$$
Then we use the convolution bound 
$$ \| \d_\e^{\a_2} b_1 \op_\e(b_2) \d_\e^{\a_3} u \|_{\fl1} \leq \| \d_\e^{\a_2} b_1 \|_{\fl1}  \| \op_\e(b_2) \d_\e^{\a_3} u \|_{\fl1},$$
and finally the fact that the $L^1 \to L^1$ norm of the Fourier multiplier $\op_\e(b_2)$ is controlled by the sup norm of $b_2.$ 
\end{proof} 
   
\section{The symbolic flow theorem in anisotropic quantization} \label{app:duh}

  Consider the initial value problem
 \begin{equation} \label{buff01}
  \begin{aligned} \d_t u + \frac{1}{\e} \op_\e(L) u  = f_0, \qquad u(0) \in H^{s_1}(\R^3),\end{aligned}
  \end{equation}
 where $f_0 \in L^\infty([0,  T \sqrt \e |\ln \e|], H^{s_1}(\R^3)).$ Here as before $x \in \R,$ $y \in \R^2,$ and $\R^3 = \R_x \times \R^2_{y}.$

 \begin{assump}  \label{ass:B} The family of symbols $L = L(\e,t,x,y,\xi,\eta) \in S^0$ is bounded in $S^0,$ in the sense that it satisfies, for all $\a \in \N^3$ with $|\a| < s_a - d/2,$ for all $ \b \in \N^3,$ 
$$ \langle \xi,\eta \rangle^{|\b|} \big| \d_{x,y}^\a \d_{\xi,\eta}^\b L(\e,t,x,y,\xi,\eta) \big| < \infty,$$
uniformly in $\e \in (0,\e_0)$ for some $\e_0 > 0,$ in $t \in[0, T \sqrt \e |\ln \e|],$ and $(x,y,\xi,\eta) \in \R^{2d}.$
Moreover, we have the decomposition $L = L_0 + \sqrt\e L_1,$ with $L_0$ and $L_1$ each separately bounded in $S^0,$ and $L_0$ independent of $(x,y).$ 
\end{assump}

Associated with \eqref{buff01}, we consider the linear partial differential equation
\be \label{resolvent0fg} \d_t S_0(\t;t) + \frac{1}{\e} L(t) S_0(\t;t) + \frac{(-i)}{\sqrt\e} \d_\eta L_0(t)  \d_y S_0(\t;t)  = f_1, \qquad S_0(\t;\t) = f_2,\ee
where $f_1$ and $f_2$ depend on $(\t;t,x,y,\xi,\eta),$ and are differentiable up to order $s_a - d/2,$ with bounded derivatives, and finite $L^2(\R^2_y)$ norms, uniformly in $(x,\xi,\eta).$ (That is: $f_1$ and $f_2$ are such that the norms in the upcoming Assumption \ref{ass:BS} are finite.) 
For the solution $S_{0,f_1,f_2}$ to \eqref{resolvent0fg}, we assume an exponential growth in time:

 \begin{assump} \label{ass:BS} There exists a unique solution $S_{0,f_1,f_2}$ to \eqref{resolvent0fg}, which, for some $\g > 0,$ for all $\a \in \N^3$ with $|\a| < s_a - d/2,$ satisfies the bound 
$$ \begin{aligned}\|\d_{x,y}^\a & S_{0,f_1,f_2}(\t;t)\|_{L^\infty(\R^4_{x,\xi,\eta}, L^2(\R^2_y))}  + \|\d_{x,y}^\a S_{0,f_1,f_2}(\t;t)\|_{L^\infty(\R^6)} \\ & \lesssim |\ln\e|^\star  \Big( N_{\infty,\a}(\t;t,f_1,f_2;\R^6)  +  N_{2,\a}(\t;t,f_1,f_2;\R^6) \Big)  \\ & + \frac{|\ln \e|^\star}{\sqrt \e}  \int_\t^t e^{(t - t') \g(U)/\sqrt \e}  \Big( N_{\infty,\a}(\t;t',f_1,f_2;\R^6)  +  N_{2,\a}(\t;t',f_1,f_2;\R^6) \Big)  \, dt' ,\end{aligned}$$ 
  where the implicit constant depends on norms of $L,$ but not on $\e,$ and is uniform in time for  $0 \leq \t \leq t \leq T \sqrt \e |\ln \e|.$ 
 Above,  $|\ln \e|^*$ means $|\ln \e|^{N_\star}$ for some $N_\star > 0$ possibly depending on all parameters, but independent of $(\e,\t,t,x,y,\xi,\eta).$ 
 
 The norms $N_{\infty,\a}$ and $N_{2,\a}$ are defined in \eqref{def:Cinfty} and \eqref{def:C2}.
\end{assump}

We also assume a specific spatial far-field behavior for $S_{0,f_1,f_2}:$ 

\begin{assump} \label{ass:new:for:S} If $f_1$ and $f_2$ are independent of $(x,y)$ outside of a possibly $\e$-dependent bounded domain in $\R^3 \times \R^3,$ of area $O(|\ln \e|^\star),$ then outside of some bounded domain of the same area, $S_{0,f_1,f_2}$ is also independent of $(x,y)$  and uniformly bounded in $(\xi,\eta).$ 
\end{assump}

\begin{rem} \label{rem:ass:L} Assumptions {\rm \ref{ass:BS}} and {\rm \ref{ass:new:for:S}} are implicit in $L.$ As argued in Remark {\rm \ref{rem:for:Duh}}, the proof of Lemma {\rm \ref{lem:rough:S}} shows that if $L$ has the form 
  $L = i \underline A + \sqrt \e \underline B,$
  where $\underline A$ is real, diagonal and independent of $(x,y),$ then Assumption {\rm \ref{ass:B}} implies Assumption {\rm \ref{ass:BS}}. If moreover $\underline B$ vanishes outside of a bounded domain of size $O(|\ln \e|^\star),$ then Lemma {\rm \ref{lem:away}} shows that Assumption {\rm \ref{ass:new:for:S}} holds as well. 
\end{rem}

We denote
\be \label{def:S0}
S_0 := S_{0,0,{\rm Id}},
\ee
the solution to \eqref{resolvent0fg} with $f_1 =0$ and $f_2 = {\rm Id}.$ We see in the following computation (where we use the composition result \eqref{compo:e}) the role played by the $\d_\eta L_0 \d_y S_0$ term in \eqref{resolvent0fg}, and also how $\op_\e(S_0)$ fails to be an approximate solution operator for \eqref{buff01}:
\be \label{the:cancellation} \begin{aligned} - \d_t \op_\e(S_0) & = \frac{1}{\e} \op_\e(L S_0) + \frac{(-i)}{\sqrt \e} \op_\e(\d_\eta L_0 \d_y S_0) \\ & = \frac{1}{\e} \Big( \op_\e(L) \op_\e(S_0) - (-i) \sqrt \e \op_\e(\d_\eta L_0 \d_y S_0)  + \e(\dots) \Big) + \frac{(-i)}{\sqrt \e} \op_\e(\d_\eta L_0 \d_y S_0) \\ & = \frac{1}{\e} \op_\e(L) \op_\e(S_0) + O(1),
\end{aligned}\ee
and now the goal is to change the above $O(1)$ into something small. In this view, we introduce correctors $\{ S_q\}_{1 \leq q \leq q_0},$ for some $q_0 \in \N$ to be determined, defined
  as the solutions of the triangular system of linear partial differential equations
 \begin{equation} \label{resolventk}
 \left\{\begin{aligned} \d_t S_q  + \frac{1}{\e} L S_q + \frac{(-i)}{\sqrt \e} \d_\eta L_0  \d_y S_q + \frac{1}{\sqrt\e} D_q  & = 0,
 \\  S_{q}(\t;\t) & = 0. \end{aligned}\right.
  \end{equation}
where the ``source" term $D_q = D_q(L, (S_{q'})_{0 \leq q' \leq q-1})$ is chosen so that cancellations like the one observed in \eqref{the:cancellation} occur at every order.

For $q = 1,$ this means that $D_1$ must contain all the terms of order $O(1)$ in the expansion 
\be \label{expansion:duh} \frac{1}{\e} \op_\e(L S_0) = \frac{1}{\e} \op_\e(L) \op_\e(S_0) + \frac{1}{\sqrt \e} (\cdots).\ee
Thus we let
\be \label{def:D1}
 D_1 := \sum_{|\a| = 2} \frac{(-i)^{|\a|}}{\a!} \d_\eta^\a L_0 \d_y^\a S_0 + (-i) \d_\eta L_1 \d_y S_0 + (-i) \d_\xi L_0 \d_x S_0.
\ee
Then we obtain, formally, 
$$ - \d_t \op_\e(S_0 + \e^{1/2} S_1) = - \frac{1}{\e} \op_\e(L) \op_\e(S_0 + \e^{1/2} S_1) + O(\sqrt \e).$$ 
We define similarly $D_q$ to include all terms of order $\e^{(q-1)/2}$ in the expansion \eqref{expansion:duh}:
\ba \label{def:D}
 D_q & := \sum_{2 \leq |\a| \leq q + 1} \frac{(- i)^{|\a|}}{\a!} \d_\eta^\a L_0 \d_y^\a S_{q + 1 - |\a|} + \sum_{1 \leq |\a| \leq q } \frac{(- i)^{|\a|}}{\a!} \d_\eta^\a L_1 \d_y^\a S_{q - |\a|} \\ & + \sum_{1 \leq \a \leq [(q+1)/2]}   \frac{(- i)^{\a}}{\a!} \d_\xi^\a L_0 \d_x^\a S_{q + 1 - 2 |\a|}  + \sum_{1 \leq \a \leq [q/2]} \frac{(- i)^{|\a|}}{\a!} \d_\xi^\a L_1 \d_x^\a S_{q - 2 |\a|},  
 \ea
for $1 \leq q \leq q_0,$ where a sum over an empty is by convention equal to zero. With $S_q$ and $D_q$ defined in \eqref{resolventk} and \eqref{def:D}, and ${\bf S}$ defined as
\be \label{def:S} {\bf S} := \sum_{0 \leq q \leq q_0} \e^{q/2} S_q
\ee
 we will see that $\op_\e({\bf S})$ is an approximate solution operator for \eqref{buff01}. But first, we prove bounds for $S_q$ and the associated operators: 

 \begin{lem} \label{lem:bd-S} Under Assumptions {\rm \ref{ass:B}} and {\rm \ref{ass:BS}}, the leading term $S_0$ defined in \eqref{def:S0} and correctors $S_q$ defined in \eqref{resolventk} satisfy for all
 $q \in [0,q_0],$ all $\a \in \N^3$ with $|\a| < s_a - d/2 - q,$ and all $\b \in \N^3,$ the bounds:
  \be \label{bd:Sq} \langle \xi,\eta \rangle^{|\b|} |\d_{x,y}^\a \d_{\xi,\eta}^\b S_q(\t;t) | \lesssim \e^{-|\b|/2} |\ln \e|^* \exp( \g (t - \t) /\sqrt \e), \qquad 0 \leq \t \leq t \leq T \sqrt \e |\ln \e|,\ee
  with the same convention for $|\ln \e|^*$ as in Assumption {\rm \ref{ass:BS}}. With Assumption {\rm \ref{ass:new:for:S}}, this implies the bounds, for $0 \leq \t \leq t \leq T \sqrt \e |\ln \e|,$ if $s_a > q_0 + 4 + d/2:$ 
  \be \label{op:Sq} \| \op_\e(S_q(\t;t))\|_{L^2 \to L^2} + \| \op_\e(S_q(\t;t))\|_{\fl1 \to \fl1} \lesssim |\ln \e|^* \exp( \g (t - \t) / \sqrt \e),\ee
  and, if $s_a > q_0 + s_1 + 4 + d/2:$
  \be \label{op:Sq:s}
   \| \op_\e(S_q(\t;t)) z \|_{\e,s_1}  \lesssim |\ln \e|^* \exp( \g (t - \t) / \sqrt \e) \| z\|_{\e,s_1}, \qquad \mbox{for all $z \in H^{s_1}.$}
  \ee
 \end{lem}

 \begin{proof} This is similar to the proof of Lemma \ref{lem:derS}. We work by induction on $q,$ and, within the induction on $q,$ perform a double induction on $|\a|$ and $|\b|$ (on $|\b|$ only in the case $q = 0$). For $q = 0$ and $\b = 0,$ the bound \eqref{bd:Sq} is given in Assumption \ref{ass:BS}. For $|\b| > 0,$ we observe that
 $$ \d_{x,y}^\a \d_{\xi,\eta}^\b S_0 = S_{0,\Gamma_{\a,\b}(S_0),0}, \qquad \Gamma_{\a,\b} := - [ \d_{x,y}^\a \d_{\xi,\eta}^\b , \e^{-1} L + (-i) \e^{-1/2} \d_\eta L_0 \d_y].$$
 By the induction hypothesis (bearing on $|\b|$), and Assumption \ref{ass:BS}, we have
 $$ |\Gamma_{\a,\b}(S_0)| \lesssim \e^{-(|\b| + 1)/2} |\ln \e|^\star e^{\g (t - \t)/\sqrt \e}.$$
 Hence, by Assumption \ref{ass:BS} again, used here with $f = \Gamma_{\a,\b}(S_0)$ and $g = 0,$ we obtain
 $$ | \d_{x,y}^\a \d_{\xi,\eta}^\b S_0| \lesssim \e^{-|\b|/2} |\ln \e|^\star e^{\g(t - \t)/\sqrt \e}.$$
 At this point we proved the bound \eqref{bd:Sq} for $q = 0$ and all $\b$ (and all $\a$). We assume now that the bound \eqref{bd:Sq} holds for all $q' \leq q-1,$ for some $q \geq 1.$ Then, by definition of the correctors \eqref{resolventk}, for $q \geq 1$ we have 
 $$ \d_{x,y}^\a \d_{\xi,\eta}^\b S_q = S_{0,f_{\a,\b,q},0}, \qquad f_{\a,\b,q} := \Gamma_{\a,\b}(S_q) - \e^{-1/2} \d_{x,y}^\a \d_{\xi,\eta}^\b D_q.$$ 
 If $\a = \b = 0,$ then $\Gamma_{\a,\b} = 0,$ and $D_q$ involves only $S_{q'},$ with $q' \leq q-1.$ Hence, by the induction hypothesis and Assumption \ref{ass:BS}, we find the expected bound for $S_q.$ 
 
 For any $\a,$ and $\b,$ since $D_q$ involves only lower-order correctors $S_{q'},$ with $q' \leq q-1,$ the $D_q$ term in $f_{\a,\b,q}$ contributes an appropriate upper bound, via the induction hypothesis on $q,$ Assumption \ref{ass:BS} and the smallness of the time interval under consideration.
 
 Thus in the induction on $|\a|$ and $|\b|,$ we may focus on the commutator term in $f_{\a,\b,q}.$ 
We assume that the bound \eqref{bd:Sq} for $S_q,$ for all $\b'$ with $|\b'| \leq |\b| - 1,$ for some $|\b| \geq 1.$ Most terms in $\Gamma_{\a,\b}(S_q)$ contains strictly less than $|\b|$ frequency derivatives of $S_q,$ hence the induction hypothesis on $|\b|$ applies to those. The other terms have the form 
 $$ \e^{-1/2} \d_{x,y}^{\a_1} L_1 \d_{\xi,\eta}^\b \d_{x,y}^{\a_2} S_q, \qquad \a_1 + \a_2 = \a, \quad |\a_1| > 0.$$
 Here an induction on $|\a|$ is made possible by the fact that $\a_1$ is strictly positive. We obtain
 $$ |f_{\a,\b,q}| \lesssim \e^{-(|\b| + 1)/2} e^{(t - \t) \g/\sqrt \e},$$
 and conclude with Assumption \ref{ass:BS}.

Thus \eqref{bd:Sq} is proved and we move on to a proof of \eqref{op:Sq}. In the case $q = 0,$ with bound \eqref{bd:Sq} and Assumption \ref{ass:new:for:S}, the bounds \eqref{op:Sq} are proved exactly as in the proof of Proposition \ref{prop:opS}. 

Then it suffices to observe that Assumption \ref{ass:new:for:S} carries over to the correctors $S_q,$ with $q \geq 1,$ since the source $D_q$ involves only $S_{q'}$ with $q' \leq q-1.$ Thus we may repeat again the proof of Proposition \ref{prop:opS} in the case $q \geq 1,$ and arrive at \eqref{op:Sq}.

The bound \eqref{op:Sq:s} is then derived by differentiation with respect to the spatial directions. The use of Proposition \ref{prop:H} for $\d_\e^{\a} S_q,$ with $q \leq q_0$ and $|\a| \leq s_1 + 4$ requires $s_a > q_0 + s_1 + 4 + d/2.$
\end{proof} 

 \begin{cor} \label{lem:duh-remainder} Under Assumptions {\rm \ref{ass:B}}, {\rm \ref{ass:BS}} and {\rm \ref{ass:new:for:S}}, we have, for ${\bf S}$ defined in \eqref{def:S}, if $s_a > q_0 + s_1 + 4 + d/2,$ 
  \begin{equation} \label{RS} \op_\e(\d_t {\bf S}) = - \frac{1}{\e} \op_\e(L) \op_\e({\bf S}) + \rho,
  \end{equation}
  where for $0 \leq \t \leq t \leq T \sqrt \e |\ln \e|,$ for all $z \in H^{s_1}:$ 
  \begin{equation} \label{est:RS} \| \rho(\t;t) z \|_{\e,s_1} \lesssim \e^{q_0/2 - C_0} |\ln \e|^* \exp\big( \g (t - \t) / \sqrt \e \big) \| z \|_{\e,s_1},\end{equation}
  for some constant $C_0 > 0$ which depends only on the spatial dimension.
  \end{cor}

  \begin{proof}
  By definition of ${\bf S}$ and \eqref{resolventk},
  \begin{equation} \label{c1}
   - \d_t \op_\e({\bf S}) = \sum_{0 \leq q \leq q_0} \e^{(q-2)/2} \op_\e(L S_q) + \sum_{0 \leq q \leq q_0} \e^{(q-1)/2}  \op_\e ( (-i) \d_\eta L_0 \d_y S_q +  D_q).
   \end{equation}
By composition of operators (Proposition \ref{prop:composition}):
$$ \sum_{0 \leq q \leq q_0} \e^{(q-2)/2} \op_\e(L S_q) = \e^{-1} \op_\e(L) \op_\e({\bf S}) - {\rm I}_y - {\rm I}_x - \rho,$$
with notations
$$
\begin{aligned}
{\rm I}_y  & := \sum_{\begin{smallmatrix} 0 \leq q \leq q_0 \\ 1 \leq |\b| \leq q_0  +1 - q \end{smallmatrix}} \e^{(q - 2 + |\b|)/2} \frac{(-i)^{|\b|}}{\b!} \d_\eta^\b L \d_y^\b S_q, \\
 {\rm I}_x  & := \sum_{\begin{smallmatrix} 0 \leq q \leq q_0 \\ 1 \leq \b \leq [(q_0  + 1 - q )/2] \end{smallmatrix}} \e^{(q-2)/2 + |\b|} \frac{(-i)^{\b}}{\b!} \d_\xi^\b L \d_x^\b S_{q}, \\
 \rho & := \sum_{0 \leq q \leq q_0} \e^{(q-2)/2} \e^{[(q_0 + 1 - q)/2]}  R_{(q_0 - q  )/2}(L, S_q),
\end{aligned}
$$
where notation $R_j(a,b)$ for the remainder in the composition of two operators is introduced in \eqref{compo:e}.
Reindexing, we obtain
$$ \begin{aligned} {\rm I}_y & = \sum_{\begin{smallmatrix} 0 \leq q \leq q_0 \\ 1 \leq |\b| \leq q  +1 \end{smallmatrix}} \e^{(q-1)/2} \frac{(-i)^{|\b|}}{\b!} \d_\eta^\b L \d_y^\b S_{q + 1  - |\b|} \\ & = \sum_{0 \leq q \leq q_0}  \e^{(q-1)/2} (-i) \d_\eta L \d_y S_{q} +  \sum_{\begin{smallmatrix} 1 \leq q \leq q_0 \\ 2 \leq |\b| \leq q  +1 \end{smallmatrix}} \e^{(q-1)/2} \frac{(-i)^{|\b|}}{\b!} \d_\eta^\b L \d_y^\b S_{q + 1  - |\b|}.
\end{aligned}$$
We can decompose $L = L_0 + \sqrt \e L_1$ and reindex once more, to find
$$ \begin{aligned} {\rm I}_y & = \sum_{0 \leq q \leq q_0} \e^{(q-1)/2} (-i) \d_\eta L_0 \d_y S_q + \sum_{0 \leq q \leq q_0} \big(\mbox{terms in $\d_\eta/\d_y$ in $D_q$}\big). \end{aligned}$$
The same holds for $I_x$ in relation to the $\d_\xi/\d_x$ terms in $D_q.$ Thus we obtain the cancellation 
\be \label{cancel:flow} \sum_{0 \leq q \leq q_0} \e^{(q-1)/2} \big( (-i) \d_\eta L_{0} \d_y S_q + D_q\big) - {\rm I}_x - {\rm I}_y = 0,\ee
and the equality \eqref{RS} follows, with $\rho$ defined above. 
 According to \eqref{composition:sobolev} and bound \eqref{bd:Sq} in Lemma \ref{lem:bd-S},
$$ \| \op_\e(R_{(q_0 - q )/2}(L, S_q)) z \|_{\e,s_1} \lesssim \e^{-C(d)/2} | \ln \e|^* e^{\g (t - \t)/\sqrt \e} \| z\|_{\e,s_1},,$$
for all $z \in H^{s_1},$ where $C(d)$ is the constant, depending only the spatial dimension $d = 3,$ that appears in Proposition \ref{prop:composition}. Summing up in $q$ as in the above definition of $\rho,$ we arrive at \eqref{est:RS}.
\end{proof}

Denote
\be \label{X}
 X := C^0(0, T\sqrt \e |\ln \e|, H^{s_1}(\R^3)), \qquad \| f \|_X :=  \sup_{0 \leq t \leq T \sqrt \e |\ln \e|} \| f(t) \|_{\e,s_1}.
\ee

  \begin{theo} \label{th:duh} Under Assumptions {\rm \ref{ass:B}}, {\rm \ref{ass:BS}} and {\rm \ref{ass:new:for:S}}, if $s_a > q_0 + s_1 + 4 + d/2,$ where $q_0$ is large enough so that 
  \be \label{def:sigma}
 \sigma := \frac{q_0}{2} - C_0 - T \g > 0, \end{equation}
  then the initial value problem \eqref{buff01} has a unique solution $u \in X,$ which satisfies the representation
  \begin{equation} \label{buff0.01} u = \op_\e({\bf S}(0;t)) u(0) + \int_0^t \op_\e({\bf S}(t';t))(\Id + \e^{\sigma} R_1) \Big( f_0 + \e^{\sigma} R_2(\cdot) u(0)\Big)(t')\, dt',
  \end{equation}
 with ${\bf S}$ defined in \eqref{def:S}. The constant $C_0$ depends only on the spatial dimension and was introduced in Corollary {\rm \ref{lem:duh-remainder}.}
  The remainders $R_1$ and $R_2$ satisfy the bound
  \begin{equation} \label{bd-R121}
  \| R_1 \|_{X \to X} \lesssim |\ln \e|^{\star}, \qquad \| R_2(t) z \|_{\e,s_1} \lesssim |\ln \e|^\star \| z \|_{\e,s_1},
  \end{equation}
for all $z \in H^{s_1}$ and all $t \leq \sqrt \e T |\ln \e|.$ %
 \end{theo}


In \eqref{def:sigma}, the integer $q_0$ enters in the definition of $S$ \eqref{def:S}, the constant $C_0$ comes from Corollary \ref{lem:duh-remainder} and the rate $\g$ is introduced in Assumption \ref{ass:BS}. The notation $|\ln \e|^\star$ is introduced in Assumption \ref{ass:BS}.

\begin{proof} Let $g \in X.$ By Corollary \ref{lem:duh-remainder}, the map 
 $$
 u := \op_\e( S(0;t)) u(0) + \int_0^t \op_\e( S(t';t)) g(t') \, dt'
 $$
   solves \eqref{buff01} if and only if for all $0 \leq t \leq T \sqrt \e |\ln \e|,$ we have the identity
\begin{equation} \label{cond-g1}
  (\Id + r) g(t)  = f_0(t) - \rho(0;t) u(0),
   \end{equation}
 where $r$ is the linear integral operator
 $$r: \ v \in X \to \Big(t \to \int_0^t  \rho(\t;t) v(\t) \, d\t\Big) \in X.$$
 Above, $\rho$ is the remainder in Corollary \ref{lem:duh-remainder}. We now choose the index $q_0$ that enters in the definition of $S$ such that $\sigma$ defined in \eqref{def:sigma} is positive.
 Then, by estimate \eqref{est:RS}, for all $0 \leq t \leq T \sqrt \e |\ln \e|,$ the operator $r$ maps $X$ to itself, with the bound
  \be \label{bd-rf1} \| r \|_{X \to X} \lesssim \e^{\sigma} |\ln \e|^*.
  \ee
  In particular, for $\e$ small enough the operator $\Id + r$ is invertible,
  with inverse $(\Id + r)^{-1}$ bounded as an operator from $X$ to itself, uniformly in $\e.$ As a consequence,
  we can solve \eqref{cond-g1} in $X$
 and obtain the representation formula \eqref{buff0.01}, with
 $$ \e^\sigma R_1 :=  (\Id + r)^{-1} - \Id, \qquad \e^\sigma R_2(t) := - \rho(0;t).$$
Bounds \eqref{bd-R121} then directly follow from \eqref{bd-rf1} and \eqref{est:RS}. Since $L \in S^0,$ $\op_\e(L)$ is linear bounded $L^2 \to L^2$ (Proposition \ref{prop:H}), hence, by spatial regularity of $L,$ also linear bounded $H^{s_1} \to H^{s_1}.$ Thus \eqref{buff01} is a differential equation in $H^{s_1},$ and uniqueness is a consequence of the Cauchy-Lipschitz theorem.
\end{proof}

\section{Notation index} \label{sec:notation:index}

 \setlength{\columnseprule}{0.001cm}
 \begin{multicols}{2}

$\prec:$ \eqref{sharp:tilde} 

\smallskip

$|\ln \e|^\star:$ unspecified power of $|\ln \e|,$ page \pageref{defi:uniform:remainder}

\smallskip

$\nabla_\e:$ weighted spatial gradient, \eqref{nablae}

\smallskip

$\d_\e:$ another weighted spatial gradient, \eqref{nablae}

\smallskip

$A(\nabla_\e):$ the main differential operator in (EM), Section \ref{sec:hypop}

\smallskip

${\bf A}:$ the diagonal order-zero symbol resulting from the changes of variables of Section \ref{sec:6}, see Lemma \ref{lem:step3bis}

\smallskip

$\vec e_p$ for $p \in \{-1,1\}:$ fixed vectors such that $\tilde u_{0p} = g_{(p)} \vec e_p,$ see \eqref{def:gp}, \eqref{for:interaction:coefficients} and \eqref{explicit:ep}

\smallskip

$g_{(p)}:$ untruncated WKB amplitude \eqref{def:gp}

\smallskip

$g_p:$ truncated initial WKB amplitude \eqref{def:g}

\smallskip

$\g:$ optimal rate of growth, \eqref{def:gamma}

\smallskip

$\g^+_{j_1j_2}(U):$ suboptimal rate of growth, \eqref{def:gamma2}

\smallskip

$\g^+_{j_1j_2j_3}(U):$ suboptimal rate of growth, \eqref{def:gamma3}

\smallskip

$\g^+:$ suboptimal rate of growth, \eqref{def:gamma+}

\smallskip

$\G:$ rough (suboptimal) rate of growth, \eqref{def:Gamma} 

\smallskip

${\mathcal J}:$ the convective terms in (EM), Section \ref{sec:conv} 

\smallskip

${\mathcal J}_e:$ the electronic convective terms, Section \ref{sec:conv} 

\smallskip

${\mathcal J}_i:$ the ionic convective terms, Section \ref{sec:conv}

\smallskip

$J:$ the high-frequency component of the electronic convective terms \eqref{def:J}

\smallskip

$\check J:$ the order-one symbol created by the first normal form reduction \eqref{def:checkJ}

\smallskip

$k:$ the initial WKB wavenumber defined in \eqref{data}; often identified with $k (1,0,0)$

\smallskip

$K:$ arbitrarily large parameter, page \pageref{th:main}

\smallskip

$K':$  arbitrarily small parameter, page \pageref{th:main}

\smallskip

$\l_j:$ approximate eigenvalues of the symbol of the main differential operator \eqref{approximate:spectral:decomposition}

\smallskip

$L:$ order-zero symbol in the abstract symbolic flow theorem of Section \ref{app:duh}; see Assumption \ref{ass:B}

\smallskip

$L_0:$ the leading term in $L$, see Assumption \ref{ass:B} 

\smallskip

$\tilde \mu_j:$ defined in terms of the diagonal entries of ${\bf A};$ the precise definition depends on the context, see Section \ref{sec:flow:blocks}

\smallskip

$\mu_j:$ like the $\tilde \mu_j$ but without the truncation; the precise definition is context-dependent, \eqref{def:mu}

$\nu_{jj'} = \d_\eta (\mu_j - \mu_{j'})$ \eqref{velocity}

\smallskip

$N_{\infty,\a}:$ sup norm of the datum and source in the resolvent flow equation \eqref{def:Cinfty}

\smallskip

$N_{2,\a}:$ transverse spatial $L^2$ norm of the datum and source in the resolvent flow equation \eqref{def:C2}

\smallskip

$\Pi_j:$ approximate spectral projectors \eqref{approximate:spectral:decomposition}

\smallskip
 $q_0:$ the number of correctors required for the symbolic flow, page \pageref{def:S0} 
 
 \smallskip

${\mathcal R}_{jj'}:$ the frequency set of $(j,j')$ resonances, Definition \ref{def:resonances}

\smallskip

${\mathcal R}:$ the total resonant set, Definition \ref{def:resonances}

\smallskip 

$R_A:$ order-one symbol that appears as a remainder in the approximation of $A$ by $A|_{\e = 0};$ \eqref{spectral:dec:symbols}

\smallskip

$S:$ the symbolic flow \eqref{resolvent0-in-proof}

\smallskip

$S_{f,g}:$ symbolic flow with source \eqref{symbolic:flow:source}

\smallskip

$S_q:$ correctors to the symbolic flow, page \pageref{resolventk}

 \smallskip
 
${\bf S}:$ the complete symbolic flow \eqref{def:S0}

\smallskip

${\mathcal S}_{jj'}:$ the set of $(j,j')$ space-time resonances \eqref{def:space-time}

 \smallskip
 
 $\s:$ exponent associated with the symbolic flow remainders, \eqref{def:sigma}

 \smallskip

 $\theta = (k x - \o t)/\e,$ page \pageref{def:theta}

 \smallskip
 
 $U^{away}_{\delta}:$ frequency domain relative to a given resonant pair or triplet \eqref{def:domain:away:space-time} or \eqref{def:away:triplets}

 \smallskip
 
 $U^{small}_\delta:$ frequency domain relative to a given resonant pair or triplet \eqref{def:Omega:delta} 
 
 \smallskip
 
 $U^{large}_{\delta'}:$ frequency domain for a given resonant pair, Section \ref{sec:large-trace} 

\smallskip

 $U^{\pm}_{\delta'}:$ frequency domain for a given resonant pair,\eqref{def:U:pm}

 \smallskip
 
$\vec v, \vec v^{\,'}:$ fixed vectors defined in \eqref{def:vv}

\smallskip

$\vec v^{\,''} = \vec v \times \vec v^{\,'}$ \eqref{vBperp}

\smallskip

$\Phi_{jj'}:$ the phase function associated with the $(j,j')$ resonance \eqref{def:phase:function}

\smallskip 

$\xi_{jj'}^\pm:$ longitudinal frequency at which the optimal growth rate associated with resonance $(j,j')$ is attained; Section \ref{sec:interaction:coefficients}

\smallskip
 
 $\o = (1 + k^2)^{1/2}:$ the characteristic frequency associated with $k,$ page \pageref{wkb:intro} and Figure \ref{fig:fund-phase}

 \smallskip
 
 $\Omega:$ space-frequency domain defined in \eqref{def:Omega}, based on \eqref{def:chi:Omega:pair} for resonant pairs, and \eqref{def:chi:Omega:triplet} for triplets

 \smallskip
 
 $\Omega_{freq}:$ projection of $\Omega$ over $\R^3_{\xi,\eta}$ \eqref{def:Omega} 
 
 \end{multicols}

\section{Index} \label{sec:index}

 \setlength{\columnseprule}{0.001cm}
 \begin{multicols}{2}

 {\it anisotropic quantization:} \eqref{aniso}

 \smallskip

  {\it equivalent source:} Definition \ref{def:equiv:remainder} page \pageref{def:equiv:remainder}

\smallskip

{\it $\fl1$ space:} Section \ref{sec:fl1}

\smallskip

 {\it growth rate}: Definition \ref{def:growth:rate} page \pageref{def:growth:rate}
 
  \smallskip

 {\it phase function:} Proposition \ref{prop:res} page \pageref{prop:res}
 
  \smallskip

 {\it resonance:} Definition \ref{def:resonances} page \pageref{def:resonances}

 \smallskip

 {\it resonant pair:} Definition \ref{def:resonances} page \pageref{def:resonances}

 \smallskip

 {\it space-time resonance:} \eqref{def:space-time}
 
   \smallskip

 {\it symbolic flow:} leading term: \eqref{resolvent0-in-proof} and \eqref{resolvent0fg}-\eqref{def:S0}; correctors: \eqref{resolventk}; full symbolic flow:  \eqref{def:S}
    
     \smallskip

 {\it triplet:} \eqref{triplets}

  \smallskip

{\it uniform remainder:} Definition \ref{defi:uniform:remainder} page \pageref{defi:uniform:remainder}

 \smallskip

 {\it weighted Sobolev norm:} \eqref{weighted:norm}

 \smallskip

\end{multicols}

\section{Parameter lists and further notation} \label{sec:lists}

\subsection{Time parameters}

We work with $t \geq 0,$ $x \in \R$ and $y \in \R^2.$ 

\begin{itemize}

\item $T_\star(\e)$ is the maximal existence time for our main initial-value problem (EM)-\eqref{initial:datum}, given by the classical quasi-linear symmetric hyperbolic theory. It is introduced in Section \ref{sec:tstar}.

\item $t_\star(\e)$ is the maximal time for which the perturbative unknown is controlled by $\e^{K'}$ in weighted $\fl1$ norm. A precise definition is given in \eqref{a:priori:new}.

\item $T = (K - K')/\g$ is such that the amplification takes place around time $T \sqrt \e |\ln\e|.$ It is introduced in \eqref{def:Tstar}.

\item $T_\e$ is the precise amplification time. It is just a bit smaller than $T,$ and introduced near the end of the analysis in \eqref{def:T0}.

\end{itemize}

\subsection{The main unknown $u$ and its avatars}

\begin{itemize}

\item $u = u(\e,t,x,y)$ is the local-in-time solution to system (EM) given in Section \ref{sec:EM} with initial datum \eqref{initial:datum}.

\item $\tilde u$ is $u$ in a rescaled time frame in which we change $y$ into $\sqrt \e y.$ See \eqref{def:tildeu}.

\item $\dot u$ is defined in \eqref{def:dotu} as $\dot u := \tilde u - \tilde u_a.$

\item $\check u$ is defined in terms of $\dot u$ in \eqref{def:checku:good}, as a result of the first normal form reduction.

\item $\check v$ is defined in terms of $\dot u$ in \eqref{def:checkv}, as a result of the second normal form reduction.

\item $u_{in}, u_{out},$ and $u_{high}$ are defined in Section \ref{sec:in:out:high}. 

\item $v_{in}$ is defined in terms of $u_{in}$ in the course of Section \ref{sec:6}. 

\end{itemize}

\subsection{Sobolev indices} 

\begin{itemize}
\item $s_a$ is the Sobolev index of regularity of the WKB amplitude $a;$ see below in Section \ref{sec:index:wkb}.
\item $s_1$ is the Sobolev index such that the refined bound of Proposition \ref{cor:the:upper:bound} holds in $\|\cdot\|_{\e,s_1}$ norm. A lower bound on $s_1$ is given in \eqref{cond:m}. See also Remark \ref{rem:Sobolev:indices}. 
\item $s = s_a - 2K_a - 1$ is the Sobolev index of regularity of the exact solution. It is such that the rough bound of Proposition \ref{prop:weak:bound} holds in $\|\cdot\|_{\e,s}$ norm. In particular, $s$ is much larger than $s_1$ if $K$ is large and $K'$ is small. See Remark \ref{rem:Sobolev:indices}.  
\end{itemize}

\subsection{The main ``source'' term: convection, current density and Lorentz force} 

\begin{itemize}

\item ${\mathcal B}:$ the current density and Lorentz force, defined in \eqref{def:unB}.

\item $\dot {\mathcal B}:$ the linearized current density and Lorentz force and the low-frequency contribution of the linearized convective terms, defined in \eqref{def:B0}. The linearization occurs at the WKB solution $u_a.$ 

\item $B:$ the leading term in $\dot {\mathcal B},$ defined just below \eqref{def:g}. Not to be confused with the magnetic field in the (EM) system, which we also denote $B.$ 
 
\item $B_{pjj'}:$ these {\it interaction coefficients} are defined in \eqref{def:interaction:coefficient}.

\item $\check B:$ the avatar of $B$ that results from the first normal form reduction. It is defined in \eqref{def:checkB}. 

\item $\check B_{(-1)}:$ the avatar of $B$ that results from the second normal form reduction, defined in \eqref{def:checkB-1:informal} and \eqref{def:checkB-1}.

\item $b_{jj'}^\pm:$ these interaction coefficients are defined in \eqref{def:bjj'}.

\item ${\bf B}:$ the final avatar of $B$ in the {\it in} equation, defined just below \eqref{eq:Uflat}. 

\end{itemize} 

\subsection{The WKB solution}  \label{sec:index:wkb}

Most parameters associated with the WKB approximate solution bear the subscript $a,$ as in ``approximate". 

\begin{itemize}

\item $u_a$ is the ($\e$-dependent family of) WKB solution(s).

\item $a(x,y) \in \R^{14}$ is the leading initial WKB amplitude, introduced in Section \ref{sec:ZEM}. 

\item $s_a$ is the Sobolev index of regularity of the initial WKB amplitude $a.$ We have $s_a \to \infty$ as $K \to \infty$ and $K' \to 0.$ See Remark \ref{rem:Sobolev:indices}.

\item $K_a$ is the order of precision of the WKB solution $u_a.$ It satisfies $K_a \geq K + 1/2.$ 

\item $p \in \{-1,1\}$ is used to denote the harmonics of the leading term in the WKB solution.

\item $R_a$ is the remainder associated with the WKB solution. See Section \ref{sec:new:wkb}.

\end{itemize}

\subsection{Cut-offs}

\begin{itemize} 

\item Notation $\prec$ is introduced in \eqref{sharp:tilde}: we have $\chi_1 \prec \chi_2$ if $\chi_2 \equiv 1$ on the support of $\chi_1.$

\item If $\chi$ is a cut-off, then we use, sometimes without introduction, associated cut-offs $\chi^\flat,$ $\chi^\sharp$ and $\tilde \chi$ such that $\chi^\flat \prec \chi \prec \chi^\sharp \prec \tilde \chi.$ See Section \ref{sec:cutoffs}.

\item $\chi_0$ is a smooth $\R^d \to [0,1]$ ``plateau" cut-off. With $d=2,$ it is used in the definition of $\chi_{trans}$ on page \pageref{chi0}. With $d = 1,$ it is used in the definition of the frequency cut-offs $\chi_{jj'}$ in Definition \ref{def:cut-offs}. 
 
 \item $\chi_{long}$ is a cut-off in the longitudinal direction $x,$ introduced in Section \ref{sec:in:out:high}. 

\item $\chi$ is a shorthand for $\chi_{long}^\sharp,$ which we use exclusively in Section \ref{sec:flow}. 

\item  $\chi_{trans}$ is a cut-off in the transverse directions $y,$ defined in \eqref{def:chi:infty}.

\item
${\bm \chi}_{trans}$ is another cut-off in $y,$ which takes into account a buffer zone. It is introduced in \eqref{def:bm:chi:trans}.

\item $\chi_{jj'}$ is a frequency cut-off around the $(j,j')$ resonance, introduced in Definition \ref{def:cut-offs}.

\item $\chi_{low}$ is a low-frequency cut-off. It is introduced in \eqref{def:chilow}.

\item $\chi_{high}$ is a high-frequency cut-off. It is introduced in \eqref{def:chiHF} The link between $\chi_{low}$ and $\chi_{high}$ is pictured on Figure \ref{fig:r:chi}.

\item 
$\chi_\Omega$ is a space-frequency cut-off that is introduced in \eqref{def:chi:Omega:pair} or \eqref{def:chi:Omega:triplet}, depending on the context, in view of the description of the far-field behavior of $S$ in Lemma \ref{lem:away}.

\end{itemize}

 \end{document}